\documentclass[12pt,a4paper,leqno]{report}
\usepackage{pdftricks}
\begin{psinputs}
\usepackage{pstricks}
\usepackage{multido}
\end{psinputs}
\usepackage{graphicx}
\usepackage{amssymb}
\usepackage{amsthm}
\usepackage{epstopdf}
\usepackage{enumerate}
\usepackage{pst-all}
\newtheorem{theorem}{Theorem}[chapter]
\newtheorem{lemma}[theorem]{Lemma} 
\newtheorem{remark}[theorem]{Remark} 
\newtheorem{definition}[theorem]{Definition} 
\newtheorem{corollary}[theorem]{Corollary} 
\newtheorem{claim}[theorem]{Claim} 
\newtheorem{example}[theorem]{Example}
\def\da{\downarrow}
\def\nda{\not\downarrow}
\def\A{\mathcal{A}}
\def\B{\mathcal{B}}
\def\C{\mathcal{C}}
\def\D{\mathcal{D}}
\def\K{\mathcal{K}}
\def\M{\mathbb{M}}
\def\a{\alpha}
\def\g{\gamma}
\def\o{\omega}
\def\raj{\upharpoonright}

\medskip
\medskip

\begin{document}

\begin{titlepage}
\thispagestyle{empty}

\vspace*{4cm}

\begin{center}
\huge{Finding groups in Zariski-like structures}

\vspace{1cm}

\Large{Kaisa Kangas}

\vspace{11cm}

\end{center}
\begin{flushright}
\today\\
Licentiate's thesis\\
University of Helsinki\\
Department of Mathematics and Statistics\\
Advisor: Tapani Hyttinen  \\
\end{flushright}
\end{titlepage}

\begin{abstract}
We study quasiminimal classes, i.e. abstract elementary classes (AECs) that arise from a quasiminimal pregeometry structure.
For these classes, we develop an independence notion, and in particular, a theory of independence in $\M^{eq}$.
We then generalize Hrushovski's Group Configuration Theorem to our setting.
In an attempt to generalize Zariski geometries to the context of quasiminimal classes, we give the axiomatization for Zariski-like structures, and as an application of our group configuration theorem, show that groups can be found in them assuming that the pregeometry obtained from the bounded closure operator is non-trivial.
Finally, we study the cover of the multiplicative group of an algebraically closed field and show that it provides an example of a Zariski-like structure.
\end{abstract}

\newpage
 
 \vspace{5cm}
 ACKNOWLEDGEMENTS
 
 \vspace{0.5cm}
 
I express my gratitude to my advisor Tapani Hyttinen for his ideas and guidance. 
Jonathan Kirby deserves thanks for discussions on the subject, in particular for suggesting transferring torus equations into a canonical form using a coordinate change. 
I would also like to thank Alexander Engstr\"om for acting as one of the referees for the thesis.
I am grateful to Jenny and Antti Wihuri Foundation and to the Finnish National Doctoral Programme in Mathematics and its Applications for financially supporting my work.

 \newpage
 
 \tableofcontents

 
\chapter{Introduction}

In \cite{HrZi93} and \cite{HrZi}, E. Hrushovski and B. Zilber introduced the concept of \emph{Zariski geometry}, a structure that generalizes the Zariski topology of an algebraically closed field. 
One of the results in \cite{HrZi} is that in a non locally modular, strongly minimal set in a Zariski geometry, an algebraically closed field can be interpreted. 
This result plays an important role in Hrushovski's proof of the geometric Mordell-Lang Conjecture (\cite{Hr}, see also e.g. \cite{bous1}), where model-theoretic ideas were applied to solve a problem from arithmetic geometry.
The field is acquired by first finding an Abelian group and then using it to construct the field.
At both steps, the Group Configuration Theorem originally presented by Hrushovski in his Ph.D. thesis (see e.g. \cite{pi}) is utilized. 
This theorem roughly states that whenever a certain kind of configuration of elements can be found, there exists a group.
 
The origin of this thesis was the question whether Zariski geometries, and the theorem from \cite{HrZi} stating the existence of a group, could be generalized from the context of first-order logic to that of quasiminimal classes, i.e. abstract elementary classes (AECs) that arise from a quasiminimal pregeometry structure (see \cite{monet}).
The results presented here will be included in joint papers with T. Hyttinen.
From the beginning, we had the idea that covers of the multiplicative group of an algebraically closed field together with the PQF-topology (see \cite{Lucy}) should serve as an example of the generalized Zariski geometries.
This eventually led to the axiomatization of \emph{Zariski-like} structures, presented in Chapter 4.
The road was not completely straightforward, as we first had to generalize Hrushovski's Group Configuration Theorem to the context of quasiminimal classes.
For this, we developed an independence calculus that has all the usual properties of non-forking and works in our context.
 
Quasiminimal classes are uncountably categorical.
They have both the amalgamation property (AP) and the joint embedding property (JEP), and thus also have a model homogeneous universal monster model, which we will denote by $\M$.
These classes are also excellent in the sense of B. Zilber (this is different from the original notion of excellence due to S. Shelah).
In the second chapter, we develop the independence notion for them.
We first isolate some properties of AECs (axioms AI-AVI presented in Chapter 2) and prove that under them the class has a perfect theory of independence (ideas used here originate from \cite{HyLess} and \cite{Hy2}).
This somewhat resembles the elementary case of strongly minimal structures, where the independence notion and Morley ranks can be obtained from the pregeometry associated to the model theoretic algebraic closure operator (acl).
In the quasiminimal case, we replace the algebraic closure operator by the bounded closure operator (bcl).

In our context, we cannot construct $\M^{eq}$ so that it would be both $\omega$-stable (in the sense of AECs) and have elimination of imaginaries.
Since $\omega$-stability is vital, we build the theory so that we can always move from $\M$ to $\M^{eq}$ and then, if needed, to $(\M^{eq})^{eq}$ and so on.
We then show that the properties expressed by axioms AI-AVI are preserved when moving from $\M$ to $\M^{eq}$, and finally that the axioms are satisfied by quasiminimal classes.

In Chapter 3, we show, generalizing Hrushovski, that from a group configuration a Galois definable rank $1$ group can be constructed.
Since $\M^{eq}$ does not necessarily have elimination of imaginaries in our setting, this group is found in $(\M^{eq})^{eq}$ rather than in $\M^{eq}$.
Essentially the first trick used in Hrushovski's original proof does not work in our context  (we would need to take rather arbitrary countable sets as elements of $\M^{eq}$, which is not possible), but otherwise the proof generalizes nicely to our context.
To overcome the problem, we move from the pregeometry to the canonical geometry associated to it and work there.
This is possible since for all (singletons) $a \in \M$, $\textrm{bcl}(a) \setminus \textrm{bcl}(\emptyset)$ is indeed in our $\M^{eq}$ (note that in the elementary case, $\textrm{acl}(a) \setminus \textrm{acl}(\emptyset)$ need not be in $\M^{eq}$).

In Chapter 4, we look at possibilities of generalizing Zariski geometries to our context. 
We give the axioms (ZL1)-(ZL9) for a \emph{Zariski-like} structure, and then apply our group configuration theorem to show that a group can be found there. 
We also point out that Zariski geometries satisfy our axioms, so we indeed have a generalization.
We work with quasiminimal classes and formulate the axioms within this context.
In the original context of Zariski geometries, a single structure is used as a starting point.
It is assumed that a collection of topologies arises from the structure, and the axiomatization is given for the closed sets in these topologies.
Then, a saturated elementary extension of the original structure is taken and the work is carried out there.
Unlike in the elementary case, we do not start from a single structure, but formulate our axioms to generalize the setting obtained after moving into the elementary extension.
Thus, we are able to use properties of quasiminimal classes to our advantage.

Instead of arbitrary closed sets, we have decided to look at irreducible closed sets (which, for simplicity, we call just irreducible sets) and state our axioms for them.
In the case of Zariski geometries, the irreducible $\emptyset$-closed sets satisfy the axioms.
The notion of a closed set could also be useful, as can be seen in the example of covers of the multiplicative group of an algebraically closed field, treated in Chapter 5, where there is a natural notion of a closed set.
However, we don't feel our insight is strong enough to formulate the axioms for arbitrary closed sets.
In \cite{Zi}, B. Zilber has given one axiomatization for closed sets in a non-elementary case, which he calls analytic Zariski structures, but we have chosen a somewhat different route.
Partially because of not using the more general concept of a closed set, some of our axioms come from Assumptions 6.6. in \cite{HrZi} rather than from the axiomatization (Z0)-(Z3) for Zariski geometries.
In our axiomatization, axioms (ZL1)-(ZL6) give meaning to the key axioms (ZL7)-(ZL9).
If, in (ZL9), we take $\kappa$ to be finite and choose $S=\{\kappa\}$, then we get just the axiom (Z3) of Zariski geometries (the dimension theorem).
In the elementary case, (ZL9) is the immediate consequence of (Z3) and Compactness.
Axioms (ZL7) and (ZL8) come from Assumptions 6.6 in \cite{HrZi}.
 
In Chapter 5, we study the cover of the multiplicative group of an algebraically closed field, a class originally introduced by Zilber.
It can be obtained from complex exponentiation $exp: (\mathbb{C},+) \to (\mathbb{C}^*, \times)$, or more precisely, from the exact sequence $0 \to \mathbb{Z} \to (V,+) \to (F^*, \times) \to 1$, where $F$ is an algebraically closed field of characteristic $0$, and $V$ is a vector space over $\mathbb{Q}$.
In particular, we show that the irreducible $\emptyset$-closed sets in the PQF-topology (see \cite{Lucy}) satisfy our axioms for Zariski-like structures, and thus the cover provides an example of such a structure.
This class is quasiminimal by \cite{Zilber}.
Prior to \cite{monet}, the uncountable categoricity of the class was known by \cite{bays}.

The main result of \cite{HrZi} is that every very ample Zariski geometry arises from the Zariski topology of a smooth curve over an algebraically closed field.
In addition to improving our axiomatization, the final goal in our study of Zariski-like structures might be to prove an analogue to this theorem, i.e. that all non-trivial Zariski-like structures resemble in some sense the cover presented in Chapter 5.
This would mean that on the level of the canonical geometry we would be back in the elementary case (pregeometries can be very complicated). 
A result like this would be in line with the existing studies of geometries in non-elementary cases.
However, since the existence of a non-classical group (see \cite{HLS} and \cite{Hy3} for locally modular cases) is still open, to prove something like this seems very difficult, and if it turns out that there are non-classical groups, the playground is completely open.

Since Zariski geometries serve as the starting point of our work, and since some results on them are needed in Chapter 5, we now provide a brief introduction to them.
 
\section{Zariski geometries}

Zariski geometries were introduced by Hrushovski and Zilber in \cite{HrZi93} and \cite{HrZi}.
In this section we present the definiton of a Zariski geometry and some basic properties of Zariski geometries.  
All results on Zariski geometries that are presented in this section can be found in \cite{HrZi}. 
More information on Zariski geometries can also be found in  \cite{Marker}  or \cite{Zi}.
The former reference contains some illustrative and relatively easily approachable material.
 
Zariski geometries are structures that generalize the idea of the Zariski topology on an algebraically closed field.
Let $F$ be an algebraically closed field.
Then, we can define a topology on $F^n$ for each $n$ as follows.
Let $S \subset F[x_1, \ldots, x_n]$.
We say that the set
$$\{x \in F^n \, | \, f(x)=0 \textrm{ for all } f \in S \}$$
is the \emph{vanishing set} of the polynomial set $S$.
We say that a set $V \subset F^n$ is \emph{Zariski closed} if it is the vanishing set of some set of polynomials.
The Zariski closed sets form a topology on $F^n$ called the \emph{Zariski topology}.  
The Zariski topology is \emph{Noetherian}, i.e. there are no infinite descending sequences of closed sets. (see e.g. \cite{Hu} for details.)

\begin{definition}
Let $X$ be a topological space, and let $C \subseteq X$ be a closed set.
We say $C$ is \emph{irreducible} if there are no closed sets $C_1, C_2 \subsetneq C$ such that $C=C_1 \cup C_2$.
\end{definition}  
  
The proof of the following lemma can be found from e.g. \cite{Hu}.

\begin{lemma}
Let $X$ be a Noetherian topological space, and let $C \subset X$ be closed.
Then, there are finitely many irreducible closed sets $C_1, \ldots, C_n$ such that $C=C_1 \cup \ldots \cup C_n$.
Moreover, if we choose $C_1, \ldots, C_n$ so that $C_i \not\subseteq C_j$ for $i \neq j$, then $C_1, \ldots, C_n$ are unique up to permutation.
\end{lemma}  

\begin{definition}
The sets $C_1, \ldots, C_n$ from the lemma are called the \emph{irreducible components} of $C$.
\end{definition}
  
For a Noetherian topology, we define the dimension of a set as follows.
\begin{definition}
If $X$ is a Noetherian space and $C \subseteq X$ is irreducible, closed and nonempty, then we define the \emph{dimension} of $C$ inductively as follows:
\begin{itemize}
\item $\textrm{dim}(C) \ge 0$,
\item $\textrm{dim}(C)=\textrm{sup }\{\textrm{dim}(F)+1 \, | \, F \subsetneq C, \textrm{$F$ closed, irreducible and nonempty } \}$.
\end{itemize}

If $C\subseteq X$ is an arbitrary closed set, then the dimension of $C$ is the maximum dimension of its irreducible components.

If $A \subseteq X$ is an arbitrary set, then the dimension of $A$ is the dimension of its closure.
\end{definition}

In the following, we use the concept of dimension in the sense of the definition.
 
\begin{definition}
A \emph{Zariski geometry} is an infinite set $D$ together with a family of Noetherian topologies on $D, D^2, D^3, \ldots$ such that the following axioms hold:

(Z0) Coherence and separation:
\begin{enumerate}[(i)]
\item If $f: D^n \to D^m$ is defined by $f(x)=(f_1(x), \ldots, f_m(x))$, where $f_i: D^n \to D$ is either constant or a coordinate projection for each $i=1, \ldots, m$, then $f$ is continuous.
\item Each diagonal $\Delta_{i,j}^n=\{(x_1, \ldots, x_n)\in D^n \, | \, x_i=x_j\}$ is closed.
\end{enumerate}

(Z1) Weak quantifier elimination: If $C \subseteq D^n$ is closed and irreducible, and $\pi: D^n \to D^m$ is a projection, then there is a closed $F \subsetneq \overline{\pi(C)}$ such that $\overline{\pi(C)} \setminus F \subseteq \pi(C)$.

(Z2) Uniform one-dimensionality:
\begin{enumerate}[(i)]
\item $D$ is irreducible.
\item Let $C \subseteq D^n \times D$ be closed and irreducible.
For $a \in D^n$, let $C(a)=\{x \in D \, | \, (a,x) \in C \}$.
There is a number $N$ such that for all $a \in D^n$, either $|C(a)| \le N$ or $C(a)=D$.
In particular, any proper closed subset of $D$ is finite.
\end{enumerate}

(Z3) Dimension theorem: Let $C \subseteq D^n$ be closed and irreducible.
Let $W$ be a non-empty irreducible component of $C \cap \Delta_{i,j}^n$.
Then, $\textrm{dim } C \le \textrm{dim } W+1$.
\end{definition}

The Dimension theorem (Z3) is the key structural condition that allows us to interpret an algebraically closed field in a non locally modular Zariski geometry.

\begin{remark}
\begin{enumerate}[(i)]
\item It follows from (Z0) that if $C_1, C_2$ are closed, then $C_1 \times C_2$ is closed.
Indeed, $C_1 \times C_2=\pi_1^{-1}(C_1) \cap \pi_2^{-1}(C_2)$ where $\pi_1, \pi_2$ are the suitable projections.
\item If $C \subset D^n \times D^m$ is closed, and $a \in D^n$, then $C(a)=f^{-1}(C)$, where $f(x)=(a,x)$ for $x \in D^m$.
Thus, $C(a)$ is closed by (Z0).
Also, if $a \in D$, then $g: D \to D^2$, $g(x)=(a,x)$ is a continuous function.
Since the diagonal of $D^2$ is closed, $g^{-1}(\Delta_{1,2}^2)=\{a\}$ is closed.
Thus singletons are closed.
\item It can be shown that $\textrm{dim } C_1 \times C_2=\textrm{dim } C_1 +\textrm{dim }C_2$, so in particular $\textrm{dim }D^n=n$ (see \cite{HrZi}, Chapter 2).
Thus, every set has finite dimension.
\end{enumerate}
\end{remark}

An algebraically closed field $F$ together with the Zariski topology for each $F^n$ satisfies the axioms, and even a more general result can be proved:
If $D$ is a smooth quasi-projective algebraic curve, then $D$, equipped with the Zariski topologies on $D^n$, is a Zariski geometry. (see \cite{Marker} for details).

The following lemma is proved completely similarly as Lemma 2.2. in \cite{HrZi}.
In Chapter 5, present the same result for the so-called PQF-topology on a cover of the multiplicative group of an algebraically closed field (Lemma \ref{PQFcartesianirred}).
The proof is essentially similar also in this case.
 
 \begin{lemma}\label{closirr}
Let $C_1, C_2$ be closed and irreducible.
Then, $C_1 \times C_2$ is irreducible.
In particular, $D^n$ is irreducible.
\end{lemma}

Now we can look at the irreducible components of cartesian products. 
 
\begin{lemma}\label{components}
Let $C$ and $F$ be two closed sets, and let $C_1, \ldots, C_n$ be the irreducible components of $C$, and $F_1, \ldots, F_m$ the irreducible components of $F$.
Then,  the irreducible components of $C \times F$ are $C_i \times F_j$ ($1 \le i \le n$, $1 \le j \le m$). 
\end{lemma}
\begin{proof}
By Lemma \ref{closirr}, $C_i \times F_j$ is closed and irreducible for $1 \le i \le n$, $1 \le j \le m$.
Clearly, 
\begin{displaymath}
C \times F=\bigcup_{1 \le i \le n, 1\le j\le m}  C_i \times F_j,
\end{displaymath}
and $C_i \times F_j \neq C_{i'} \times F_{j'}$ for $(i,j) \neq (i',j')$.
\end{proof}
  
We also have a stronger version of (Z3):

\begin{theorem}\label{dimthmzariski}
Let $C_1, C_2$ be closed, irreducible subsets of $D^n$.
Then, every irreducible component of $C_1 \cap C_2$ has dimension at least $\textrm{dim }C_1 +\textrm{dim }C_2 -n$. (Lemma 2.5 in \cite{HrZi})
\end{theorem}  

By Lemma \ref{closirr}, $D^k$ is irreducible for every $k$, and by (Z0) (i),
the set $\Delta_{i,j}^n$ is isomorphic with $D^{n-1}$.
Thus, (Z3) follows  from Theorem \ref{dimthmzariski}: If $C$ is a closed set and $W$ an irreducible component of $C \cap \Delta_{i,j}^n$, then 
\begin{displaymath}
\textrm{dim } W \ge \textrm{dim }C +\textrm{dim }\Delta_{i,j}^n -n= \textrm{dim }C+(n-1)-n=\textrm{dim }C-1.
\end{displaymath}
  
Suppose now $D$ is a countable Zariski geometry.
Let $\mathcal{L}_D$ be the language where we have an $n$-ary predicate for each closed subset of $D^n$.
Let $T_D$ be the $\mathcal{L}_D$-theory of $D$.
We note that since singletons are closed sets, each element of $D$ has its own predicate.
 
\begin{theorem}
\begin{enumerate}[(i)]
\item $T_D$ admits elimination of quantifiers.
\item $T_D$ is $\omega$-stable, and the Morley rank of a definable set $X$ equals the dimension of its closure.
In particular, $D$ is strongly minimal. (\cite{HrZi}, section 2)
\end{enumerate}
\end{theorem}

Let $M$ be an elementary extension of $D$.
Define a topology on $M$ so that the basic closed sets are those sets $X$ for which there is a closed $C \subseteq D^m \times D^n$ for some $m,n$, and $a \in M^m$ such that $X=C(a)$, i.e.
\begin{displaymath}
X=\{b \in M^n \, | \, M \models C(a,b)\}.
\end{displaymath}
It turns out that with respect to this topology, $M$ is a Zariski geometry (\cite{HrZi}, Proposition 4.1).

From now on, we will replace $D$ by a saturated elementary extension.
Thus, we assume that there is a Zariski geometry $D_0$ such that $D$ is a saturated elementary extension of $D_0$ in the language $\mathcal{L}_{D_0}$ and that the topology on $D$ is obtained from the topology on $D_0$ as described above. 
It is this situation that we generalize when presenting our axioms for Zariski-like structures in Chapter 4.
There, we give a more general framework with axioms that are satisfied by the irreducible closed sets of $D_0$ after moving into the saturated elementary extension.
 
\begin{definition}
Let $A \subset D$.
We say that a set $X$ is \emph{$A$-closed} if $X=C(a)$ for some $C \in \mathcal{L}_{D_0}$ and some $a \in A^n$ for some $n$.
For $x \in D^n$, we define the \emph{locus} of $x$ over $A$ to be the smallest $A$-closed set containing $x$. 

If $C$ is an irreducible closed set, we say that an element $a \in C$ is \emph{generic} (over $A$) if $C$ is the locus of $a$ (over $A$).
\end{definition}
 
\begin{remark}
We note that the Morley rank of a tuple $a$ over a set $A$ coincides with the dimension of the locus of $a$ over $A$:
$$\textrm{MR}(a/A)=\textrm{dim}(C),$$
where $C$ is the locus of $a$ over $A$.
\end{remark}

Suppose $A \subset B \subset D$, and let $a \in D^n$.
We say that $a$ is \emph{independent} from $B$ over $A$ if $\textrm{MR}(a/A)=\textrm{MR}(a/B)$.
In Chapter 2, we will present a notion of independence that can be applied in a more general setting.   
 
\subsection{Regular points}

In the Zariski geometry context, we often need to work inside some closed set $C \subset D^n$ rather than inside $D^n$ itself.  
When doing so, we use a generalized version of Theorem \ref{dimthmzariski} that states that if $C_1$ and $C_2$ are closed, irreducible subsets of $C$, then all ``nice enough" irreducible components of $C_1 \cap C_2$ have dimension at least $\textrm{dim }C_1 +\textrm{dim }C_2 -\textrm{dim } C$.  
Unfortunately, this does not hold for all irreducible components of $C_1 \cap C_2$, but it holds for components that pass through a regular point of $C$ (Lemma 5.4 in \cite{HrZi}).
In a sense, regular points are the Zariski geometry analogue of smooth points on a variety. 

\begin{definition}
Let $C \subset D^n$ be a closed set.
We define the \emph{codimension} of $C$ in $D^n$, denoted $\textrm{codim}_{D^n} C$, to be the number
$$\textrm{codim}_{D^n} C=\textrm{dim } D^n-\textrm{dim C}=n-\textrm{dim }C.$$
\end{definition}

\begin{definition}  
Let $C \subseteq D^n$ be an irreducible closed set and let $p \in C$.
Denote $\Delta_C=\{(x,y) \in C \times C : x=y \}$.
We say that $p$ is a \emph{regular} point of $C$ if there is a closed irreducible set $G \subseteq D^n \times D^n$ such that
\begin{enumerate}[(i)]
\item $\textrm{codim}_{D^n \times D^n} G = \textrm{dim } C$
\item $\Delta_C$ is the unique irreducible component of $G \cap C \times C$ passing through $(p,p)$.
\end{enumerate}
\end{definition}  

\begin{lemma}\label{regonD}
Any $a \in D$ is regular on $D$.
\end{lemma}
\begin{proof}
Now $\textrm{codim}_{D \times D} (\Delta_D)=1=\textrm{dim } D$, so we may choose $G=\Delta_D$.
\end{proof} 
 
\begin{lemma}\label{regonloc}
Any point is regular on its own locus.
\end{lemma}
\begin{proof}
Let $a \in D^n$, and let $C$ be the locus of $a$.
We prove that $a$ is regular on $C$.
Let Let $k=\textrm{dim}(C)$, and suppose for the sake of convenience that the first $k$ coordinates of a generic point of $C$ are independent.
Denote $\Delta_C=\{(x,y)\in C \times C \, | \, x=y\}$.
Let 
\begin{displaymath}
G=\{(x_1, \ldots, x_n, y_1, \ldots, y_n) \in D^n \times D^n \, | \, x_1=y_1, \ldots, x_k=y_k\}.
\end{displaymath}
Now $\textrm{dim}(G)=2n-k$ and $\textrm{codim}_{D^n \times D^n}=k$.
We have to prove that $\Delta_C$ is the unique irreducible component of $G \cap C\times C$ passing through $(a,a)$.
Clearly  $\textrm{dim}(G \cap C\times C)=k=\textrm{dim}(\Delta_C)$, so $\Delta_C$ is indeed an irreducible component of $G \cap C \times C$.

Suppose now there is some other irreducible component $F$ of $G \cap C \times C$ such that $(a,a) \in F$.
Then, $F \cap \Delta_C \neq \emptyset$.
Moreover, as $F \neq \Delta_C$, there is some $b \in C$ such that $(b,b) \in (G \cap C \times C) \setminus F$.
Denote
\begin{displaymath}
C'=\{x \in C \, | \, (x,x) \in F \cap \Delta_C\}.
\end{displaymath}
Then, $b \in C \setminus C'$, and $C'$ is closed as $C'=f^{-1}(F \cap \Delta_C)$, where $f$ is such that $f(x)=(x,x)$ (continuous by (Z0)). 
But then $a \in C' \subsetneq C$ which contradicts the fact that $C$ is the locus of $a$.
\end{proof}

While the definition of regular points is non-intuitive, we will show that if $V$ is an irreducible variety, then every non-singular point is regular in the sense defined above (\cite{Marker}, section 2).
We first remind that a point $p$ on an irreducible variety $V$ is non-singular if the dimension of the tangent space at $p$ equals the dimension of the variety $V$.
The dimension of the tangent space at $p$ can be calculated as the dimension of the linear subspace defined by $J_p$, the Jacobian matrix of the partial derivatives at $p$ of any defining equations for $V$ chosen so that the corresponding polynomials generate the ideal of all polynomials vanishing on $V$ (see e.g. \cite{Cox} or \cite{Hu} for details). 

To illustrate the idea, we first consider the case where our variety is a plane curve $C$ defined by the equation $F(X,Y)=0$.
Let $(x,y)$ be a regular point of $C$.
Then, at least one of the partial derivatives of $F$ at $p$ is nonzero, so we may assume $\frac{\partial F}{\partial Y} \neq 0$.
Let $G=\{(X,Y,Z,W) \, | \, X=Z\}$.
Now $G \cap (C \times C)$ has dimension $1$, and thus $\Delta_C$ is an irreducible component.
We claim that it is the unique component containing $(x,y,x,y)$.
For this, it suffices to show that $(x,y,x,y)$ is non-singular on $G \cap (C \times C)$, as any point on two components is singular. 
The (possibly reducible) variety $G \cap (C \times C)$ is given by the equations
\begin{eqnarray*}
F(X,Y)&=&0, \\
F(Z,W)&=&0, \\
X-Z &=&0.
\end{eqnarray*}
The Jacobian matrix at $(x,y,x,y)$ is
\begin{displaymath}
J=\left( \begin{array}{cccc}
\frac{\partial F}{\partial X}(x,y) & \frac{\partial F}{\partial Y}(x,y) & 0 & 0\\
 0 & 0 & \frac{\partial F}{\partial X}(x,y) & \frac{\partial F}{\partial Y}(x,y) \\
 1 & 0 & -1 & 0
\end{array} \right).
\end{displaymath}
Since $\frac{\partial F}{\partial X}\neq 0$, the rows are linearly independent.
Thus, the tangent space at $(x,y,x,y)$ has dimension $4-3=1$ and $(x,y,x,y)$ is non-singular as desired.  
 
Let now $V \subseteq K^n$ be an irreducible variety of dimension $m$.
Suppose $V$ is defined by the equations 
\begin{equation}\label{defeqs}
F_1(\overline{X})= \ldots =F_l(\overline{X})=0,
\end{equation}
where the polynomials $F_1, \ldots, F_l$ generate the ideal of all polynomials vanishing on $V$.
If $p=\overline{x}$ is a smooth point of $V$, then the matrix
\begin{displaymath}
J=\left(\frac{\partial F_i}{\partial X_j}(\overline{x}) \right)
\end{displaymath}
has rank $n-m$.  

Renumbering equations and variables if necessary, we may assume that the minor
\begin{displaymath}
M=\left(\frac{\partial F_i}{\partial X_j} \right)(\overline{x}) \quad 1 \le i \le n-m, \quad m+1 \le j \le n
\end{displaymath}
is a nonsingular matrix.

Let $G=\{(\overline{x}, \overline{y}) \in K^{2n} \, x_i=y_i \textrm{ for $1 \le i \le m$} \}$.
Now $G \cap (V \times V)$ has dimension $m$ unless there are algebraic dependencies between the first $m$ coordinates.
If such dependencies exist, we may without loss of generality assume that the list (\ref{defeqs}) contains equations in the variables $x_1, \ldots, x_m$ only, giving these dependencies. 
If $F_i$ is one of the corresponding polynomials, then $\frac{\partial F_i}{\partial x_j}=0$ for $m+1 \le j \le n$.
Thus, $F_i$ gives a row in $J$ that has zeros at the indices $m+1 \le j \le n$.
Clearly we cannot have $1 \le i \le m-n$, as the nonsingular minor $M$ would then contain a zero row.
On the other hand, the nonsingularity of $M$ implies that the first $m-n$ rows of $J$ are linearly independent.
Thus, as $J$ has rank $m-n$, all the other rows are linear combinations of the first $m-n$ rows. 
But this means that if we would have $i>m-n$, then the rows of $M$ would be linearly dependent which is also impossible. 
Thus, there are no algebraic dependencies between the first $m$ coordinates, and $G \cap (V \times V)$ has dimension $m$.

As before, $\Delta_V$ is an irreducible component of $G \cap (V \times V)$.
We show that $(p,p)$ is a nonsingular point of $G \cap (V \times V)$ which again proves that $G \cap (V \times V)$ is the unique component containing it. 

To calculate the dimension of the tangent space at $(p,p)$, we must consider the $(2l+m) \times 2n$ matrix $J'$ where  
\begin{itemize}
\item For $1 \le i \le l$, the $i$:th row of $J'$ is
\begin{displaymath}
\left( \begin{array}{cccccc}
\frac{\partial F_i}{\partial X_1}(\bar{x}) & \ldots & \frac{\partial F_i}{\partial X_n}(\bar{x}) & 0 & \ldots & 0 
\end{array} \right)
\end{displaymath}
and the $(l+i)$:th row is
\begin{displaymath}
\left( \begin{array}{cccccc}
0 & \ldots & 0 & \frac{\partial F_i}{\partial X_1}(\bar{x}) & \ldots & \frac{\partial F_i}{\partial X_n}(\bar{x})  \end{array} \right),
\end{displaymath}
\item For $i \le m$, the $(2l+i)$:th row has $1$ in the $i$:th column and $-1$ in the $(n+i)$:th column.
\end{itemize}
The rows $1, \ldots, n-m$, $l+1, \ldots, l+n-m$, $2l+1, \ldots, 2l+m$ form a maximal linearly independent set, and thus $J$ has rank $2(n-m)+m=2n-m$.
Hence, the tangent space at $(p,p)$ has dimension $2n-(2n-m)=m$, as desired.

\subsection{Specializations}

The concept of a specialization plays an important role in finding an algebraically closed field from a non locally modular strongly minimal set in a Zariski geometry, and we will also be using it in our framework of Zariski-like structures.
 
 \begin{definition}\label{SpecHrZi}
 Let $D$ be a Zariski geometry.
 If $A \subset D$, we say that a function $f: A \to D$ is a \emph{specialization} if for any $a_1, \ldots, a_n \in A$ and for any $\emptyset$-closed set $C \subseteq D^n$, it holds that if $(a_1, \ldots, a_n) \in C$, then $(f(a_1), \ldots, f(a_n)) \in C$.
 
 If $A=(a_i : i \in I)$, $B=(b_i:i\in I)$ and the indexing is clear from the context, we write $A \to B$ if the map $a_i \mapsto b_i$, $i \in I$, is a specialization.
\end{definition}
 
 \begin{remark}
 It is easy to see that the following hold (tp denotes the first-order type):
 \begin{itemize}
 \item If $\textrm{tp}(a/\emptyset)=\textrm{tp}(a'/\emptyset)$, then $a \to a'$.
 \item If $a \to a'$ and $a' \to a''$, then $a \to a''$.
 \item Let $a=(a_i:i \in I)$, $\iota : I \to I$ a permutation of the index set, $a \iota = (a_{\iota (i)} : i \in I)$.
 If $a \to a'$, then $a \iota \to a' \iota$.
 \item If $a \in D$ is a generic singleton, then  $a \to a'$ holds for any singleton $a' \in D$.
 \item If $a \to a'$, then either $\textrm{tp}(a/\emptyset)=\textrm{tp}(a'/\emptyset)$ or $MR(a/\emptyset)>MR(a'/\emptyset)$.
 \end{itemize}
 \end{remark}
 
 \begin{definition}
 We define $\textrm{rk}(a \to a')=MR(a/\emptyset)-MR(a'/\emptyset)$.
 \end{definition}

The Dimension Theorem (Z3) can be reformulated in terms of specializations as follows (Lemma 4.13 in \cite{HrZi}).

\begin{lemma}\label{dimthmspecHrZi}
Let $a=(a_1, \ldots, a_n)$, $a''=(a_1'', \ldots, a_n'')$, 
$a \to a''$, and suppose $a_1 \neq a_2$, $a_1''=a_2''$.
Then there exists $a'=(a_1', \ldots, a_n')$ such that $a_1'=a_2'$, $a \to a' \to a''$, and $\textrm{rk}(a \to a')=1$.
\end{lemma}

\begin{proof}
Let $C$ be the locus of $a$.
Then, $a'' \in C \cap \Delta_{12}^n$.
Hence, $a''$ must lie on some irreducible component $W$ of $C \cap \Delta_{12}^n$.
By (Z3), $\textrm{dim}(W) \ge \textrm{dim}(C)-1$.
As $a_1 \neq a_2$, we have $C \cap \Delta_{12} \subsetneq C$, and thus $\textrm{dim}(W) < \textrm{dim}(C)$.
Thus,  $\textrm{dim}(W)=\textrm{dim}(C)-1$.
Choose $a''$ to be a generic point of $W$.
Then, $a''$ is as wanted. 
\end{proof} 

It is this version of the Dimension Theorem that is used (together with Compactness) when finding in a  non locally modular Zariski geometry the configuration that yields a group.
In the more general setting in which we will be working, we don't have Compactness.
There, the axiom (ZL9) captures Lemma \ref{dimthmspecHrZi} and the traces of compactness needed for the argument.
Also, (ZL9) implies Lemma \ref{dimthmspecHrZi}.
   
In the Zariski geometry setting, the concepts of regular and good specializations allow us to take regular points into account when working with specializations. 
In Chapter 4, we will present the concepts of strongly regular and strongly good specializations that generalize these notions.
We first recall the definition of the model theoretic algebraic closure.

\begin{definition}
Let $b \in D^n$.
We say $b$ is \emph{algebraic} over $A$ if there is some formula $\phi(x, a)$, where $a \in A^m$ for some $m$, such that the set $\{x \in D^n \, | \, \phi(x,a)\}$ is finite and $\phi(b,a)$ holds.

For $A \subseteq D$, the \emph{algebraic closure} of $A$, denoted $\textrm{acl}(A)$, is the set of all elements of $D$ algebraic over $A$.
\end{definition}

If $D$ is an algebraically closed field, then the model theoretic notions of an algebraic element and the algebraic closure of a set coincide with the field theoretic ones (see e.g. \cite{MarIntro}).
      
\begin{definition}\label{goodHrZi}
A specialization $a \to a'$ is called \emph{regular} if $a'$ is regular on the locus $a$.  

A \emph{good} specialization is defined recursively as follows.
Regular specializations are good.
Let $a=(a_1, a_2, a_3)$, $a'=(a_1', a_2', a_3')$, and $a \to a'$.
Suppose:
\begin{enumerate}[(i)]
\item $(a_1, a_2) \to (a_1', a_2')$ is good.
\item $a_1 \to a_1'$ is an isomorphism.
\item $a_3 \in \textrm{acl}(a_1)$.
\end{enumerate}
Then, $a \to a'$ is good.
\end{definition}
 
We now list some properties of regular specializations that will be utilized when forming the definition of a strongly regular specialization in Chapter 4.
 
 \begin{lemma}\label{viim}  
\begin{enumerate}[(i)]
\item If $aa' \to bb'$ is a specialization, and $a \to b$, $a' \to b'$ are regular specializations, and if $a$ is independent from $a'$ over $\emptyset$, then $aa' \to bb'$ is regular.
\item If $a$ is a generic element of $D$,  then $a \to a'$ is always regular.  
\item Isomorphisms are regular.
\end{enumerate}
\end{lemma}
 
 \begin{proof}
 For (i), we need to prove that $(b,b')$ is regular on the locus of $(a,a')$.
 Let $C_1$ be the locus of $a$ and $C_2$ be the locus of $a'$.
 Suppose $C_1 \subseteq D^n$, $C_2 \subseteq D^m$, $\textrm{dim}(C_1)=r_1$, and $\textrm{dim}(C_2)=r_2$.
 As $a$ is independent from $a'$ over $\emptyset$, it holds that the locus of $a$ over $a'$ is $C_1$.
 The independence relation is symmetric (see e.g. \cite{MarIntro}), so the locus of $a'$ over $a$ is $C_2$.
Thus, the locus of $(a,a')$ is $C_1 \times C_2$.
By our assumptions, there are closed, irreducible sets $G_1 \subseteq D^n \times D^n$ and $G_2 \subseteq D^m \times D^m$ such that $\textrm{codim}(G_1)=r_1$, $\textrm{codim}(G_2)=r_2$, $\Delta_{C_1}$ is the unique irreducible component of $G_1 \cap (C_1 \times C_1)$ passing through $(b,b)$, and $\Delta_{C_2}$ is the unique irreducible component of $G_2 \cap (C_2 \times C_2)$ passing through $(b',b')$.
Now
\begin{displaymath}
\textrm{codim}(G_1 \times G_2)=r_1+r_2=\textrm{dim}(C_1 \times C_2). 
\end{displaymath}
 As coordinate permutations are isomorphisms, it suffices to show that $G_1 \times G_2$ is the unique irreducible component  of $(C_1 \times C_1) \times (C_2 \times C_2)$ passing through $(b,b,b',b')$, but this follows from Lemma \ref{components}.

Parts (ii) and (iii) follow directly from Lemmas \ref{regonD} and \ref{regonloc}, respectively. 
\end{proof}

The concept of a good specialization is used in the following two lemmas that are utilized when proving that a group can be interpreted in a non locally modular Zariski geometry.
In our setting, the analogues of these lemmas will be the axioms (ZL7) and (ZL8). 
  
\begin{lemma}\label{amalgHrZi}
Let $a \to a'$ be a good specialization of rank $\le 1$.
Then any specializations $ab \to a'b'$, $ac \to a'c'$ can be amalgamated: there exists $b^{*}$, independent from $c$ over $a$ such that $\textrm{tp}(b^{*}/a)=\textrm{tp}(b/a)$, and $ab^{*}c \to a'b'c'$.
(Lemma 5.14 in \cite{HrZi})
\end{lemma}

\begin{lemma}\label{yhdistHrZi}  
Let $(a_i : i \in I)$ be independent over $b$ and indiscernible over $b$, where the set $I$ is infinite.
Suppose $(a_i' : i \in I)$ is indiscernible over $b'$, and $a_i b \to a_i' b'$ for each $i \in I$.
Further suppose $\textrm{rk}(b \to b') \le 1$ and $b \to b'$ is good.
Then, $(ba_i : i \in I) \to (b'a_i' : i \in I)$.
(Lemma 5.15 in \cite {HrZi})
\end{lemma}

\chapter{Independence in Abstract Elementary Classes}
 
In this chapter, we will develop an independence notion within the context of abstract elementary classes satisfying certain axioms. We will then show that it has all the usual properties of non-forking.
The ideas used originate from \cite{HyLess} and \cite{Hy2}.

First, we need to present some basic definitions.
 
\begin{definition}
Let $L$ be a countable language, let $\K$ be a class of $L$ structures and let $\preccurlyeq$ be a binary relation on $\K$.
We say $(\K, \preccurlyeq)$ is an \emph{abstract elementary class} (AEC for short) if the following hold.
\begin{enumerate}[(1)]
\item Both $\K$ and $\preccurlyeq$ are closed under isomorphisms.
\item If $\A, \B \in \K$ and $\A \preccurlyeq \B$, then $\A$ is a substructure of $\B$.
\item The relation $\preccurlyeq$ is a partial order on $\K$.
\item If $\delta$ is a cardinal and $\langle \A_i \, | \, i<\delta \rangle$ is an $\preccurlyeq$-increasing chain of structures, then
\begin{enumerate}[a)]
\item $\bigcup_{i<\delta} \A_i \in \K;$
\item for each $j<\delta$, $\A_j \preccurlyeq \bigcup_{i<\delta} \A_i$;
\item if $\B \in \kappa$ and for each $i<\delta$, $\A_i \preccurlyeq \B$, then $\bigcup_{i<\delta} \A_i \preccurlyeq \B$.
\end{enumerate}
\item If $\A, \B, \C \in \K$, $\A \preccurlyeq \C$, $\B \preccurlyeq \C$ and $\A \subseteq \B$, then $\A \preccurlyeq \B$.
\item There is a L\"owenheim-Skolem number $LS(\K)$ such that if $\A \in \K$ and $B \subseteq \A$, then there is some structure $\A' \in \K$ such that $B \subseteq \A' \preccurlyeq \A$ and $\vert \A' \vert = \vert B \vert +LS(\K)$.
\end{enumerate}

If $\A \preccurlyeq \B$, we say that $\A$ is an \emph{elementary substructure} of $\B$.
\end{definition} 

It is easy to see that the class $(\K, \preccurlyeq)$ of all models of some first-order theory $T$, where $\preccurlyeq$ is interpreted as the elementary submodel relation, is an AEC.
 
We also consider the following example, presented in \cite{kir}.

\begin{example}\label{equivalence}
Let $\K$ be the class of all models $M=(M,E)$ such that $E$ is an equivalence relation on $M$ with infinitely many classes, each of size $\aleph_0$.
For any set $X$, we define the closure of $X$ to be 
$$\textrm{cl}(X)=\bigcup\{x/E \, | \, x \in XÊ\}.$$
We define $\preccurlyeq$ so that $\A \preccurlyeq \B$ if and only if $\A \subseteq \B$ and $\A=\textrm{cl}(\A)$.
Then, it is easy to see that $(\K, \preccurlyeq)$ is an AEC.
\end{example} 

\begin{definition}
Let $\A, \B \in \K$.
We say a function $f: \A \to \B$ is an \emph{elementary embedding}, if there is some $\C \in \K$ such that $\C \preccurlyeq \B$ and $f$ is an isomorphism from $\A$ to $\C$.
\end{definition}
 
\begin{definition}
We say a class of structures $\K$ has the \emph{amalgamation property} (AP for short)
if for all $\A, \B \in \K$ and any map $f: \A \to \B$ such that $f: \A' \to \B$ is an elementary embedding for some
$\A' \preccurlyeq \A$, there exists some $\C \in \K$ such that $\B \subseteq \C$ and an elementary embedding $g: \A \to \C$ such that $f \subseteq g$.
\end{definition} 
 
\begin{definition}
We say a class of structures $\K$ has the \emph{joint embedding property} (JEP for short) if for all $\A, \B \in \K$,
there is some $\C \in \K$ such that $\B \preccurlyeq \C$ and an elementary embedding $f: \A \to \C$.
\end{definition} 

\begin{definition}
Let $\M \in \K$, and let $\delta$ be a cardinal.
We say $\M$ is \emph{$\delta$- model homogeneous} if whenever $\A, \B \preccurlyeq \M$ are such that $\vert \A \vert, \vert \B \vert < \delta$ and $f: \A \to \B$ is an isomorphism, there is some automorphism $g$ of $\M$ such that $f \subseteq g$.
\end{definition}

\begin{definition}
Let $\M \in \K$, and let $\delta$ be a cardinal.
We say $\M$ is \emph{$\delta$-universal} if for every $\A \in \K$ such that $\vert \A \vert < \delta$ there is an elementary embedding $f: \A \to \M$.
\end{definition}

We note that if $\M \in \K$ is both $\delta$- model homogeneous and $\delta$-universal, then for any $\A, \B \in \K$ such that $\A \preccurlyeq \B$ and $\vert \B \vert <\delta$, and any elementary embedding
 $f: \A \to \M$, there is an elementary embedding $g: \B \to \M$ such that $f \subseteq g$.
Indeed, by $\delta$-universality, there is some elementary embedding $g': \B \to \M$.
Then, $g'(\A)$ and $f(\A)$ are isomorphic, so let $h: g'(\A) \to f(\A)$ be an isomorphism.
By $\delta$ -model homogeneousness, $h$ extends to an automorphism $h'$ of $\M$.
Thus, $g=(h' \raj \B) \circ g'$ is as wanted.

From the above observation it follows that if all the structures we are considering are small compared to some cardinal $\delta$ and our class $\K$ contains a structure $\M$ of size $\delta$ that is both $\delta$- model homogeneous and $\delta$-universal, we can view all the other structures we are considering as elementary substructures of $\M$.
Let now $\delta$ be a cardinal bigger than any structure we will be considering, and let us call a $\delta$- model homogeneous and $\delta$-universal structure $\M \in \K$ of size $\delta$ a \emph{monster model} for $\K$.
We may now think we are always working inside the monster model $\M$.
This means that every structure we will be considering will be an elementary substructure of $\M$ of cardinality less than $\delta$, every set we will be considering will be a subset of $\M$ of cardinality less than $\delta$, and every tuple we will be considering will be a tuple of elements of $\M$.
 
From now on we suppose that $(\mathcal{K} ,\preccurlyeq )$ is an AEC with AP and JEP
and with arbitrarely large structures, $LS(\mathcal{K} )=\omega$
and $\mathcal{K}$ does not contain finite models.
Moreover, we suppose that $\mathcal{K}$ has a monster model which we will denote by $\mathbb{M}$.
Then, every time we use the term \emph{model}, we mean a structure $\A \in \K$ such that $\A \preccurlyeq \M$.
Also, whenever we write $\A \in \K$, we assume that actually $\A \preccurlyeq \M$.
If $A$ is a set, we usually write ``$a \in A$" as shorthand for ``$a \in A^n$ for some natural number $n$". 
If $a$ and $b$ are finite tuples, we will write $ab$ for the concatenation $a \frown b$.
Also, for a set $A$ and a tuple $a$, we will write $Aa$ for $A \cup a$.  

It is easy to see that in Example \ref{equivalence}, all closed models of the same cardinality are isomorphic.
Indeed, two models of the same cardinality have the same number of classes in the equivalence relation $E$.
When constructing the isomorphism, you just map equivalence classes onto equivalence classes.
It is then easy to see that $(\K, \preccurlyeq)$ satisfies the requirements listed above.
For a monster model, one can just choose any closed structure that is large enough.

We will list six axioms (AI-AVI) and show that if these axioms hold for $\mathcal{K}$, then Lascar non-splitting will satisfy the usual properties of an independence notion.

\begin{definition}
Suppose $A \subset \mathbb{M}.$
We denote by $\textrm{Aut}(\mathbb{M}/A)$ the subgroup of the automorphism group of $\mathbb{M}$ consisting of those automorphisms $f$ that satisfy $f(a)=a$ for each $a \in A$.

We say that $a$ and $b$ have the same \emph{Galois type} over $A$ if there is some $f \in \textrm{Aut}(\mathbb{M}/A)$ such that $f(a)=b$.
We write $t^{g}(a/A)=t^{g}(a/A;\mathbb{M} )$ for the Galois-type of $a$ over $A$.

We say that $a$ and $b$ have the same $\emph{weak type}$ over $A$ if for all finite subsets $B \subseteq A$, it holds that $t^{g}(a/B)=t^{g}(b/B)$.
We write $t(a/A)$ for the weak type of $a$ over $A$.

We often denote types by letters $p, q$, etc, and write e.g. $p=t(a/A)$.
Then, we say that the element $a$ \emph{realizes} the type $p$, or that $a$ is a \emph{realization} of $p$.
\end{definition}

\begin{definition}
Let $A$ and $B$ be sets such that $A \subseteq B$ and $A$ is finite.
We say that $t(a/B)$ \emph{splits} over $A$ if there are $b,c \in B$ such that $t(b/A)=t(c/A)$ but $t(ab/A) \neq t(ac/A)$.

We write $a\da^{ns}_{B}C$ (``$a$ is free from $C$ over $B$ in the sense of non-splitting") if there is
some finite $A\subseteq B$ such that $t(a/B \cup C )$ does not split
over $A$. By $A\da^{ns}_{B}C$ we mean that $a\da^{ns}_{B}C$
for each $a\in A$. 
\end{definition}

We note that if $A \subseteq B \subseteq C$ for some finite $B$, and $t(a/C)$ does not split over $A$, then $t(a/C)$ does not split over $B$ either.
Indeed, if $t(a/C)$ would split over $B$, then we could find $b,c \in C$ such that $t(b/B)=t(c/B)$ but $t(ab/B) \neq t(ac/B)$.
Since $B$ is finite, there is some tuple $d \in B$ such that $B=Ad$.
Now, $t(bd/A)=t(cd/A)$ but $t(abd/A) \neq t(acd/A)$, so the tuples $bd$ and $cd$ witness the splitting of $t(a/C)$ over $A$, a contradiction.

It is now easy to see that $\da^{ns}$ is monotone, i.e. that if $A \subseteq B \subseteq C \subseteq D$, then $a \da_{A}^{ns} D$ implies $a \da_{B}^{ns} C$.

In the context of Example \ref{equivalence}, $a \da^{ns}_B C$ means that if $a=(a_1, \ldots, a_n)$ and  if for some $1 \le i \le n$ there 
is some $c \in C$ such that $(a_i, c) \in E$, then there is also some $b \in B$ such that $(a_i, b) \in E$.
Moreover, if $a_i=c$ for some $1 \le i \le n$ and some $c \in C$, then $c \in B$.

\section{Our axioms} 
 
For the sake of readability, instead of first presenting all the definitions needed and then giving the axioms AI-AVI in the form of a simple list, we will now start listing the axioms and give the related definitions, lemmas and remarks in midst of them.

\textbf{AI:
Every countable model $\A \in \K$ is s-saturated, i.e. for any $b \in \M$ and any finite $A \subseteq \A$, there is $a \in \A$ such that $t(a/A)=t(b/A)$.}

We note that the AEC $(\K, \preccurlyeq)$ of Example \ref{equivalence} satisfies AI.
Indeed, for a tuple $b=(b_1, \ldots, b_n) \in \M$ and a finite set $A \subseteq \A$, we find a tuple $b'=(b_1', \ldots, b_n') \in \A$ such that $t(b'/A)=t(b/A)$ as follows.
Let $1 \le i \le n$.
If $(b_i, a) \in E$ for some $a \in A$, then $b_i \in \A$ since $\A$ contains $\textrm{cl}(A)=\bigcup \{a/E \, | \, a \in A\}$,
and we may choose $b_i'=b_i$.
If it holds for every $a \in A$ that $(b_i, a) \notin E$, then choose $b_i' \in \A$  so that $(b_i', a) \notin E$ holds for all $a \in A$ (such an element can be found since $A$ is finite and $E$ has infinitely many classes).
Moreover, one needs to take care that for $1 \le i<j \le n$, $b_i'=b_j'$ if and only if $b_i=b_j$.
 
\begin{lemma}\label{vaplaajykskas}
Let $\B$ be a model.
If $a \da^{ns}_\B A$, $b \da^{ns}_\B A$ and $t(a/\B)=t(b/\B)$, then  $t(a/A)=t(b/A)$.
\end{lemma}
 
\begin{proof}
Let $c \in A$ be arbitrary.
We need to show that $t(ac/\emptyset)=t(bc/\emptyset)$.
Let $B_a \subset \B$ be a finite set such that $t(a/ \B \cup A)$ does not split over $B_a$, and let $B_b \subset \B$ be a finite set such that $t(b/ \B \cup A)$ does not split over $B_b$.
Then, neither $t(a/ \B \cup A)$ nor $t(b/ \B \cup A)$ splits over $B=B_a \cup B_b$.  
By AI, there is some $d \in \B$ such that $t(d/B)=t(c/B)$. 
We have
$$t(ac/\emptyset)=t(ad/\emptyset)=t(bd/\emptyset)=t(bc/\emptyset),$$
where the first and the last equality follow from non-splitting.
The middle equality holds since $d \in \B$ and $t(a/\B)=t(b/\B)$.
\end{proof}

\begin{lemma}\label{vaplaajol}
Suppose  $\A$ and $\B$ are countable models, $t(a/\A )$ does not split over some finite $A\subseteq\A$,
and   $\A \subseteq \B$.
Then there is
some $b$ such that  $t(b/\A )=t(a/\A )$ and $b\da^{ns}_{A} \B$.
\end{lemma}

\begin{proof} 
As both $\A$ and $\B$ are countable and contain $A$, we can, using AI and back-and-forth methods, construct an automorphism $f\in Aut(\M /A)$ such that $f(\A )=\B$.
Choose $b=f(a).$
\end{proof}

\begin{definition}
We say that a model $\B =\A a\cup\bigcup_{i<\o}a_{i}$, where $a_i$ is a singleton for each $i$,
is \emph{$s$-primary} over $\A a$ if for all $n<\o$, there is
a finite $A_{n}\subset \A$
such that for all $(a', a_0', \ldots, a_n') \in \M$ such that $t(a'/\A)=t(a/\A)$, $t(a', a_0', \ldots, a_n'/A_{n})=t(a, a_0, \ldots, a_n/A_{n})$
implies  
$t(a', a_0', \ldots, a_n'/\A)=t(a, a_0, \ldots, a_n/\A)$
\end{definition}
 
\textbf{AII: For all $a$ and countable $\A$, there is
an $s$-primary model
$\B =\A a\cup\bigcup_{i<\o}a_{i}$ ($\le\M$) over $\A a$}.

We denote a countable $s$-primary model $\B= \A a\cup\bigcup_{i<\o}a_{i}$ over $\A a$ that is as above by $\A [a]$.
 
Also AII is satisfied in Example \ref{equivalence}.
Indeed, for a model $\A$ and a tuple $b=(b_1, \ldots, b_m)$, we choose 
$$\A [b]=\A \cup \bigcup \{b_i/E \, | \, 1 \le i \le m\}$$
with any enumeration.
Then, for $A_n$, we choose $\{b_1, \ldots, b_m, a_0, \ldots, a_n\} \cap \A$.
(Note that $\A$ is closed and that the type of an element is determined by identity and 
its $E$ -equivalence class.)

 \begin{lemma}\label{primary}
 Let $\A$ be a countable model, and let $t(b/\A )=t(a/\A )$. 
Then, there is an isomorphism $f:\A[a] \rightarrow\A[b]$ such that $f\raj\A =id$ and $f(a)=b$.
\end{lemma}

\begin{proof}
Let $\A[a]=\A a\cup\bigcup_{i<\o}a_{i}$ and $\A[b]=\A b\cup\bigcup_{i<\o}b_{i}$.
Now there is some finite $A_0 \subset \A$ such that it holds for any $a', a_0'$ that if $t(a/\A)=t(a'/\A)$ and $t(a',a_0'/A_0)=t(a,a_0/A_0)$, then $t(a,a_0/\A)=t(a',a_0'/\A)$.
As $t(b/\A )=t(a/\A )$, there is an automorphism $F \in \textrm{Aut}(\mathbb{M}/\A)$ such that $F(a)=b$.
Let $a_0'=F(a_0)$.
By AI, there is some $i$ such that $t(b_i/A_0 b)=t(a_0'/A_0 b)$, and in particular $t(b_i, b/A_0)=t(a_0, a/A_0)$.
Thus, $t(b_i, b / \A)=t(a_0, a/ \A)$.
Moreover, choose $i$ so that it is the least possible. 
Let $f_0: \A \cup \{a, a_0\} \to \A \cup \{b,b_i\}$ be such that $f_0 \raj \A= \textrm{id}$, $f_0(a)=b$ and $f_0(a_0)=b_i$.

Construct inductively functions $f_k$ for $k \in \o$
such that $\A a \subseteq \textrm{dom}(f_k) \subseteq \A[a]$, 
$\A b \subseteq \textrm{ran}(f_k) \subseteq \A[b]$, $f_k \raj \A=id$, $f_k(a)=b$, $\textrm{dom}(f_k) \setminus \A$ is finite and for all $c \in \textrm{dom}(f_k)$, it holds that  $t(c/\emptyset)=t(f_k(c)/\emptyset)$.
Moreover, take care that
$t(A/ \A)=t(B/ \A)$, where $A=\textrm{dom}(f_{k}) \setminus \A$ and $B=\textrm{ran}(f_{k}) \setminus \A$.

This is done as follows.
Suppose we have constructed $f_k$.
Let $i$ be least such that $a_i \notin \textrm{dom}(f_k)$.
If $b_j \in \textrm{ran}(f_k)$ for all $j<i$, we start looking for an image for $a_i$.
Otherwise, we will consider $b_j$ for the least $j$ such that $b_j \notin \textrm{ran}(f_k)$ and start looking for a pre-image.
Here we treat the former case, the latter is similar.
Let $n$ be greatest possible such that $a_n \in \textrm{dom}(f_k)$.
Let $A_n \subset \A$ be the finite subset such that if $t(a//A)=t(a'/\A)$ and $t(a,a_0, \ldots, a_n/A_n)=t(a',A_0', \ldots, a_n'/A_n)$, then  $t(a',a_0', \ldots, a_n'/\A)=t(a, a_0, \ldots, a_n/\A)$. 
Denote $A'=A \setminus \{a, a_0, \ldots, a_{i-1}\}$.
Choose now the least $j$ such that 
$$t(a, a_0, \ldots, a_{i-1}, a_i, A'/A_n)=t(f_k(a), f_k(a_0), \ldots, f_k(a_{i-1}), b_j, f(A')/A_n).$$
Then, also 
$$t(a, a_0, \ldots, a_{i-1}, a_i, A'/\A)=t(f_k(a), f_k(a_0), \ldots, f_k(a_{i-1}), b_j, f(A')/\A),$$
and we may set $f_{k+1}=f_k \cup \{(a_i,b_j)\}$.
 
Denote $f=\bigcup_{k<\o} f_k$.
Then, $f$ is the desired isomorphism.
\end{proof} 
 
In particular, it follows from the above lemma that for a countable model $\A$, $t(a/\A )$ determines $t^{g}(a/\A )$.

\begin{definition}
We say $a$ \emph{dominates} $B$ over $A$ if the following holds for all $C$: If there is a finite $A_0\subseteq A$ such that $t(a/AC)$ does not split over $A$, then $B \da^{ns}_{A}C$.
\end{definition}

\begin{lemma}\label{domination}
If $\A$ is a countable model, then the element $a$ dominates $\A[a]$ over $\A$. 
\end{lemma}

\begin{proof}
Let $A \subset \A$ be finite, and let $B$ be such that $t(a/\A B)$ does not split over $A$.
It suffices to show that for each $n$, it holds that
$$a, a_0, \ldots, a_n \da^{ns}_{\A} B.$$  
We make a counterassumption and suppose that $n$ is the least number such that
$$a, a_0, \ldots, a_n \nda^{ns}_{\A} B.$$
Let $C \subset \A$ be a finite set so that $A \subseteq C$, $A_\gamma  \subseteq C$ for each $\g \le n$, and $t(a, a_0, \ldots, a_{n-1}/ \A B)$ does not split over $C$.
By the counterassumption, there are  $c,d \in \A \cup B$ such that $t(c/C)=t(d/C)$ but $t(c, a, a_0, \ldots, a_n/C) \neq t(d, a, a_0, \ldots, a_n/C)$.
By AI, there is some $d' \in \A$ so that $t(d'/C)=t(d/C)$.
Then, either  $t(d', a, a_0, \ldots, a_n/C) \neq t(d, a, a_0, \ldots, a_n/C)$ or $t(d', a, a_0, \ldots, a_n/C) \neq t(c, a, a_0, \ldots, a_n/C)$. 
We may without loss suppose the latter.
Since $t(a, a_0, \ldots, a_{n-1}/\A B)$ does not split over $C$, we have that 
$$t(c/C, a, a_0, \ldots, a_{n-1})=t(d'/C, a, a_0, \ldots, a_{n-1})$$
(otherwise $c$ and $d'$ would witness the splitting of $t(a, a_0, \ldots, a_{n-1}/ \A B)$ over $C$).
Thus, there is some $f \in \textrm{Aut}(\M/C, a, a_0, \ldots, a_{n-1})$ such that $f(c)=d'$.
Denote $a_n'=f(a_n)$.
Then, $t(a_n', a_0, \ldots, a_{n-1}/A_n)=t(a_n, a_0, \ldots, a_{n-1}/A_n)$, and thus $t(a_n'/\A, a, a_0, \ldots, a_{n-1})=t(a_n/\A,a, a_0, \ldots, a_{n-1})$.
In particular, $$t(a_n,d' / C,a, a_0, \ldots, a_{n-1})=t(a_n', d'/ C,a, a_0, \ldots, a_{n-1}),$$
as $d' \in \A$.
But 
$$t(a_n, d'/C,a, a_0, \ldots, a_{n-1}) \neq t(a_n,c/C,a, a_0, \ldots, a_{n-1})=t(a_n', d'/C,a, a_0, \ldots, a_{n-1}),$$  
a contradiction.
\end{proof} 
 
\begin{definition}
Let $\alpha$ be a cardinal and $\A_i \preccurlyeq \M$ for $i<\alpha$, and let $A=\bigcup_{i<\alpha} \A_i$.
We say that $f:A\rightarrow\M$ is \emph{weakly elementary}
if for all $a\in A$, $t(a/\emptyset )=t(f(a)/\emptyset )$ and for all $i<\alpha$, $f(\A_i) \preccurlyeq \M$.
\end{definition}

\begin{definition}
We say a model $\A$ is \emph{$s$-prime} over $A=\bigcup_{i<\alpha} \A_i$, where $\alpha$ is a cardinal and $\A_i$ is a model for each $i$,
if for every model $\B$ and
every weakly elementary $f:A\rightarrow\B$,  
there is an elementary embedding $g:\A\rightarrow\B$ such that
$f\subseteq g$.
\end{definition}

\textbf{AIII: Let $\A, \B, \C$ be models.
If $\A\da^{ns}_{\B}\C$ and $\B=\A\cap\C$,
then there is a unique (not only up to isomorphism)
$s$-prime model $\D$ over $\A\cup\C$. Furthermore,
if $\C'$ is such that $\C \subseteq \C'$ and $\A\da_{\B}\C'$, then
$\D\da_{\C}\C'$.}

It follows that
if also $\A',\B', \C'$ and $\D'$ are as in AIII,
$f:\A\rightarrow\A'$ and $g:\C\rightarrow\C'$
are isomorphisms and $f\raj\B =g\raj\B$, then there is an
isomorphism $h:\D\rightarrow\D'$ such that
$f\cup g\subseteq h$.

\begin{remark}
 Note that if $\A$, $\B$ and $\C$ are models such that $\A\da^{ns}_{\B}\C$ and $\B\subseteq\A\cap\C$, then we must have $\B=\A \cap \C$:
Suppose not, and let  $a \in \A \cap \C \setminus \B$ and let $B \subset \B$ be a finite subset.
Then, there is some $a' \in \B$ such that $t(a'/B)=t(a/B)$.
However, $t(a/\A) \neq t(a'/\A)$ and so $t(a/\C)$ splits over $B$ with $a$ and $a'$ as the witnesses.
This contradicts the assumption that $\A\da^{ns}_{\B}\C$.
\end{remark}

The class $\K$ of Example \ref{equivalence} satisfies AIII.
Indeed, since both $\A$ and $\C$ are closed, also $\A \cup \C$ is closed, so $\A \cup \C \in \K$.
Then, $\A \cup \C$ is $s$-prime over $\A \cup \C$.
Suppose now $\C'$ is such that $\C \subseteq \C'$ and $\A \da_\B \C'$. 
We claim that $\A \cup \C \da_\C^{ns} \C'$.
Suppose not.
Then, there is some $d \in \A \cup \C$ such that $d \nda_\C^{ns} \C'$.
This means that $d=(d_1, \ldots, d_n)$ and there is some $c \in \C' \setminus \C$ such that for some $1 \le i \le n$, 
$d_i=c$ (note that since $\C'$ is closed, $(d_i, c) \in E$ implies $d_i \in \C'$).
Clearly we must then have $d \notin \C$, so $d \in \A$.
But since $\A \da_\B \C'$,  we have $d \in \B \subseteq \C$, a contradiction.

\begin{definition}
Let $\A$ be a model, $A\subseteq\A$ finite and
$a\in\M$. The game $GI(a,A,\A )$ is played as follows:
The game starts at the position $a_{0}=a$ and $A_{0}=A$.
At each move $n$, player I first chooses $a_{n+1}\in\M$ and
a finite subset $A'_{n+1}\subseteq\A$ such that
$t(a_{n+1}/A_{n})=t(a_{n}/A_{n})$, $A_{n}\subseteq A'_{n+1}$
and $t(a_{n+1}/A'_{n+1})\ne t(a_{n}/A'_{n+1})$.
Then player II chooses a finite subset $A_{n+1}\subseteq\A$
such that $A'_{n+1}\subseteq A_{n+1}$.
Player II wins if player I can no longer make a move.
 \end{definition}

\textbf{AIV: For each $a \in \M$, there is a number $n<\o$ such that for any countable model $\A$ and any finite subset $A \subset \A$, player II has a winning strategy in $GI(a,A,\A )$ in $n$ moves.}

Also AIV is satisfied in Example \ref{equivalence}.
Consider $GI(a,A,\A)$, and suppose $a=(a_{01}, \ldots, a_{0m})$. 
Assume player I has succeeded in his first move and played a tuple $a_1=(a_{11}, \ldots, a_{1m})$  such that $t(a_1/A)=t(a/A)$ and a set $A_1'$ such that $A \subseteq A_1' \subset \A$ and $t(a/A_1') \neq t(a_1/A_1')$.
The model $\A$ is closed, and thus, if for some $1 \le i \le m$ there is an element $a_{1i}' \in \A$ such that $(a_{1i}, a_{1i}') \in E$, then $a_{1i} \in \A$.
We may without loss suppose that there is a number $k \le m$ so that for $i \le k$, we can find some
$a_{1i}' \in \A$ such that $(a_{1i}, a_{1i}') \in E$, and that for $i>k$, there is no such element.
As her first move, player II plays the set $A_1=A_1' \cup (\A \cap a_1)$.
After this, player I must play some element $a_2=(a_{21}, \ldots, a_{2m})$  and some set $A_2' \subset \A$ so that $t(a_2/A_1)=t(a_1/A_1)$, $A_1 \subseteq A_2'$ and $t(a_2/A_2') \neq t(a_1/A_2')$.
For this to be possible, he must choose $a_{2i}=a_{1i}$ for $i \le k$.
The only way to ensure that $t(a_2/A_2') \neq t(a_1/A_2')$ is to for some $i>k$  choose $a_{2i}$ so that there is some $a_{2i}' \in A_2'$ such that $(a_{2i}, a_{2i}') \in E$.
Now, II plays using the same strategy as before.
Thus, after his first move, Player I can survive at most $m-k$ moves.
 
\begin{lemma}\label{ansA}
Let $a \in \M$ be arbitrary, and let $\A$ be a model.
Then, $a \da^{ns}_\A \A$.
\end{lemma}

\begin{proof}
It suffices to show that there is a finite $A\subseteq\A$ such that $t(a/\A )$ does not split over $A$.
Suppose not.
Assume first that $\A$ is countable.
We claim that then player I can survive $\o$ moves in $GI(a,A,\A )$ for any finite subset $A \subset \A$, which contradicts AIV.
Suppose we are at move $n$ and that $t(a_n/ \A)$ splits over every finite subset of $\A$ containing $A_n$.
In particular, it splits over $A_n$.
Let $b,c$ be tuples witnessing this splitting.
Let $f \in \textrm{Aut}(\M/A_n)$ be such that $f(b)=c$ and $f(\A)=\A$ (note that we may find such an automorphism as all countable models are $s$-saturated).
Now player I chooses $a_{n+1}=f(a_n)$ and $A_{n+1}=A_n \cup \{c\}$.
Then,  $t(a_n/A_n)=t(a_{n+1}/A_n)$ but $t(a_{n+1} c/ A_n)=t(a_n b/A_n) \neq t(a_n c/A_n)$ and thus $t(a_{n+1}/A_{n+1}) \neq t(a_{n}/A_{n+1})$.
As $t(a_n/\A)$ splits over every finite subset of $\A$ containing $A_n$, the same is true for $t(a_{n+1}/ \A)$.
 
Let now $\A$ be arbitrary and suppose that $t(a/\A)$ splits over every finite $A \subset \A$.
Let $\B$ be a countable submodel of $\A$.
Then, $\B$ contains only countably many finite subsets.
For each finite $B \subset \B$, we find some tuples $b, c \in \A$ witnessing the splitting of $t(a/\A)$ over $B$.
We now enlarge $\B$ into a countable submodel of $\A$ containing all these tuples.
After repeating the process $\omega$ many times we have obtained a countable counterexample.
\end{proof}

\begin{lemma}\label{notypes}
For all models $\A$, the number of weak types
$t(a/\A )$ for $a\in\M$, is $\vert\A\vert$.
\end{lemma}

\begin{proof}
We prove this first for countable models.
Suppose, for the sake of contradiction, that there is a countable model $\A$ and elements $a_i \in \M$, $i<\o_1$ so that $t(a_i/\A) \neq t(a_j/\A)$ if $i \neq j$.
As countable models are $s$-saturated, there are only countably many types over a finite set.
In particular, by the pigeonhole principle, we find an uncountable set $J \subseteq \o_1$ so that $t(a_i/\emptyset)$ is constant for $i \in J$.
After relabeling, we may set $J=\o_1$.
For each $i$, there is a number $n<\o$ such that player II wins  $GI(a_i, \emptyset, \A)$ in $n$ moves.
Using again the pigeonhole principle, we may assume that the number $n$ is constant for all $i<\o_1$.

Now we start playing $GI(a_i, \emptyset, \A)$ simultaneously for all $i< \o_1$.
Since the $a_i$ have different weak types over $\A$, for each $i$ of the form $i=2 \alpha$ for some $\alpha< \o_1$, we can find a finite set $A_\alpha \subset \A$ such that $t(a_{2 \alpha} /A_\alpha) \neq t(a_{2\alpha+1}/A_\alpha)$.
We write $A_0^i=A_\alpha$ for $i=2 \alpha$ and $i=2 \alpha+1$.
As there are only countably many finite subsets of $\A$, we find an uncountable $I \subseteq \o_1$ so that for all $i \in I$, $A_0^i=A$ for some fixed, finite $A \subset \A$.
In $GI(a_i, \emptyset, \A)$ for $i \in I$, on his first move player I plays $a_{2\alpha+1}$ and $A$ if $i=2 \alpha$ for some $\alpha<\o_1$, and $a_{2\alpha}$ and $A$ if $i=2 \alpha+1$ for some $\alpha <\o_1$.
All the rest of the games he gives up. 
Now, in each game $GI(a_i, \emptyset, \A)$ player II plays some finite $A_1^i \subset \A$ such that $A \subseteq A_1^i$.
Again, there is an uncountable $I_1' \subseteq I$ such that for $i \in I_1'$, we have $A_1^i=A_1$ for some fixed, finite $A_1$.
As there are only countably many types over $A_1$, we find an uncountable $I_1 \subset I_1'$ so that $t(a_i/A_1)=t(a_j/A_1)$ for all $i,j \in I_1$.
Again, player I gives up on all the games except for those indexed by elements of $I_1$.
Continuing like this, he can survive more than $n$ moves in uncountably many games.
This contradicts AIV.

Suppose now $\A$ is arbitrary.
Denote $X=\mathcal{P}_{<\o}(\A)$.
Then, $\vert X \vert= \vert \A \vert$.
For each $A \in X$, choose a countable model $\A_A \preccurlyeq \A$ such that $A \subset \A_A$.
By Lemma \ref{ansA}, for each weak type $p=t(a/\A)$, there is some $A_p \in X$ so that $a \da_{A_p}^{ns} \A$, and hence also $a \da_{\A_{A_p}}^{ns} \A$.
By Lemma \ref{vaplaajykskas}, $t(a/\A_{A_p})$ determines $t(a/\A)$ uniquely.
As there are only countably many types over countable models, the number of weak types over $\A$ is
$$\vert X \vert \cdot \omega=\vert \A \vert.$$
\end{proof}

\begin{lemma}\label{typegalois}
For any $a \in \M$ and any model $\A$, the weak type $t(a/\A )$ determines the Galois type $t^{g}(a/\A )$.
\end{lemma}

\begin{proof}
Suppose $t(a/\A)=t(b/\A)$.
By Lemma \ref{ansA}, we can find a countable submodel $\B$ of $\A$ so that $a \da^{ns}_\B \A$ and $b \da^{ns}_\B \A$.
By Lemma \ref{primary}, there is some $f \in \textrm{Aut}(\M/ \B)$ such that $f(\B[a])=\B[b]$ and $f(a)=b$.
Moreover, by Lemma \ref{domination}, $\B[a] \da^{ns}_\B \A$ and $\B[b] \da^{ns}_\B \A$.
We claim that the map $g=(f \raj \B[a]) \cup \textrm{id}_\A$ is weakly elementary.
For this, it suffices to show that $t(c/\A)=t(f(c)/\A)$ for every $c \in \B[a]$.
But $t(c/\B)=t(f(c)/\B)$, $c \da_{\B}^{ns} \A$, and $f(c) \da_{\B}^{ns} \A.$
Thus, by Lemma \ref{vaplaajykskas}, $t(c/\A)=t(f(c)/\A)$.

By AIII, there are unique $s$-prime models $\D_a$ and $\D_b$, over $\B[a] \cup \A$ and $\B[b] \cup \A$, respectively.
The map $g$ extends to an automorphism $h \in \textrm{Aut}(\M/\A)$ so that $h(\D_a) \subseteq \D_b$. 
The $s$-prime models are unique and preserved by automorphisms, thus we must have $h(\D_a)=\D_b$.
Since $h(a)=b$, we have $t^g(a/\A)=t^g(b/\A)$.
\end{proof}
  
\textbf{AV: If $\A$ and $\B$ are countable models,
$\A\subseteq\B$ and $a\in\M$,
and $\B\da^{ns}_{\A}a$, then $a\da^{ns}_{\A}\B$.}

We note that AV is satisfied in Example \ref{equivalence}.
Indeed, suppose $ \B \da^{ns}_\A a$ but $a \nda_\A^{ns} \B$ and write $a=(a_1, \ldots, a_n)$.
Then, there is some $b \in \B$ and some $1\le i \le n$ so that $(b,a_i) \in E$  
but $(a_i, c) \notin E$ for all $c \in \A$.
But now $b \nda_\A^{ns} a$, since otherwise there would exist some $c \in \A$ such that $(b,c) \in E$ and hence $(a_i,c) \in E$ as $E$ is an equivalence relation.
Thus, we must have $a \da_\A^{ns} \B$, as wanted.

\begin{lemma}\label{symlemma}
Let $A,C\subseteq\M$ and let $\B\subseteq A\cap C$ be a model.
If $A\da_{\B}^{ns}C$, then $C\da^{ns}_{\B}A$.
\end{lemma}

\begin{proof}  
We note first that for any finite tuples $a, c \in \M$,  and for any countable model $\B$ it holds that if $a \da_\B^{ns} c$, then $c \da_\B^{ns} a$. 
Indeed, then by dominance in $s$-primary models, it holds that $\B[a] \da_\B^{ns} c$, and thus by AV, $c \da_\B^{ns} \B[a]$, and in particular,  $c \da_\B^{ns} a$.

Let now $\B$ be arbitrary, and suppose $a \da_\B^{ns} c$ but $c \not\da_\B^{ns} a$.
Then, there is some finite $B \subset \B$ so that $t(a/\B c)$ does not split over $B$.
However, $t(c/ \B a)$ splits over $B$.
Let $b, d \in \B a$ be tuples witnessing this.
If $\B' \preccurlyeq \B$ is a countable model containing $B$, $b \cap \B$ and $d \cap \B$, then $a \da_{\B'}^{ns} c$ but 
$c \not\da_{\B'}^{ns} a$, which contradicts what we have just proved.

Suppose now $A\da_{\B}^{ns}C$ but $C\not\da^{ns}_{\B}A$.
Then, there is some $c \in C$ so that $c\not\da^{ns}_{\B}A$, and this is witnessed by some finite $a \in A$, i.e. $c\not\da^{ns}_{\B}a$.
But we have $a\da_{\B}^{ns}C$ and hence $a\da_{\B}^{ns}c$, a contradiction.
 \end{proof}

\begin{remark}\label{symremark}
Note that from Lemma \ref{symlemma} it  follows that for any $a, b \in \M$ and any model $\A$, it holds that $a \da_\A^{ns} b$ if and only if $b \da_\A^{ns} a$. 
\end{remark}

\textbf{AVI: For all models $\A, \B$ and $\D$ such that $ \A \subseteq \B \cap \D$,
there is a model $\C$
such that $t(\C/\A )=t(\B/\A )$ and $\C\da^{ns}_{\A}\D$.}

It follows that AVI holds also without the assumption
that $\B$ and $\D$ are models, as we can always find models extending these sets.

We also note that the class $\K$ in Example \ref{equivalence} satisfies AVI.
Let $\{b_i \, | \, i<\kappa\}$ contain exactly one representative for each $E$-class that intersects $\B \setminus \A$.
Let $B=\{b_i' \, | \, i<\kappa\}$ be such that $B \cap \D=\emptyset$ and $(b_i',b_j') \notin E$ for $i<j<\kappa$.
Then, we may choose $\C=\textrm{cl}(\A \cup B)$.

\begin{lemma}\label{transitivity}
If $\B$ is a model, $A \subseteq \B$ and
$\B \subseteq C$,
then $a\da^{ns}_{A}C$ if and only if $a\da^{ns}_{A}\B$ and
$a\da^{ns}_{\B}C$.
\end{lemma}
 
\begin{proof}
If $a\da^{ns}_{A}C$, then $a\da^{ns}_{A}\B$ and
$a\da^{ns}_{\B}C$ follow by monotonicity.

Suppose now $a\da^{ns}_{A}\B$ and
$a\da^{ns}_{\B}C$.
Let $A_0 \subset A$ and $B_0 \subset \B$ be finite sets so that $A_0 \subseteq B_0$, $t(a/\B)$ does not split over $A_0$ and $t(a/C)$ does not split over $B_0$.
Suppose $a\not\da^{ns}_{A}C$.
Then, $t(a/C)$ splits over $A_0$.
Let $b, c \in C$ witness the splitting, i.e. $t(b/A_0)=t(c/A_0)$ but $t(ab/A_0) \neq t(ac/A_0)$.
By AI, there are $b', c' \in \B$ so that $t(b'/B_0)=t(b/B_0)$ and $t(c'/B_0)=t(c/B_0)$.
Since $t(a/C)$ does not split over $B_0$, we have $t(ab'/B_0)=t(ab/B_0)$ and $t(ac'/B_0)=t(ac/B_0)$.
Thus,
$$t(ab'/A_0)=t(ab/A_0)\neq t(ac/A_0)=t(ac'/A_0),$$
a contradiction since $t(a/\B)$ does not split over $A_0$.
\end{proof} 

\begin{lemma}\label{vaplaajol2}
Suppose $\A$ is a model, $t(a/\A)$ does not split over some finite $A \subset \A$ and $B$ is such that $\A \subseteq B$.
Then, there is some $b$ such that $t(b/\A)=t(a/\A)$ and $b \da_A^{ns} B$.
\end{lemma}

\begin{proof}
Let $\B$ be a model such that $B \subseteq \B$.
Let $\C$ be a model containing $\A a$.
By AVI, there is a model $\C'$ such that $t(\C/\A)=t(\C'/\A)$ and $\C' \da^{ns}_\A \B$.
In particular, there is some $b \in \C'$ such that $t(b/\A)=t(a/\A)$ and $b \da^{ns}_ \A \B$.
Let $A' \subseteq \A$ be a finite set such that $A \subseteq A'$ and $b \da^{ns}_{A'} B$.
Then, by Lemma \ref{transitivity}, $b \da_A^{ns} B$.
\end{proof}
 
\begin{lemma}\label{eiaaretketjuja}
For all $a\in\M$, there is a number $n<\o$ such that
there are no models $\A_{0}\subseteq\A_{1}\subseteq ...\subseteq \A_{n}$
so that for all $i<n$, $a\nda^{ns}_{\A_{i}}\A_{i+1}$.
\end{lemma}

\begin{proof} 
Suppose models $\A_i$, $i \le n$, as in the statement of the lemma, exist.
Then, the same conditions hold also for some countable submodels, so we may assume each $\A_i$ is countable.
 We will show that player I can survive $n$ moves in $GI(a,\emptyset,\A_0 )$.
Then, the lemma will follow from AIV.

On the first move, player I chooses some finite $B_1 \subset \A_0$ so that $t(a/\A_0)$ does not split over $B_1$.
Then, there is some finite set $C_1 \subset \A_1$ so that $B_1 \subseteq C_1$ and $t(a/C_1)$ splits over $B_1$ and some $f_1 \in \textrm{Aut}(\M/B_1)$ such that $f(\A_1)=\A_0$.
Now player I plays $a_1=f(a)$ and $A_1'=f_1(C_1)$.
As $t(a/f_1(C_1))$ does not split over $B_1$ and $t(f_1(a)/f_1(C_1))$ splits over $B_1$, we have $t(a/f_1(C_1))\neq t(f_1(a)/f_1(C_1))$, and this is indeed a legitimate move.

On her move, player II chooses some finite $A_1 \subset \A_0$ such that $A_1' \subseteq A_1$.
On his second move, player I chooses some finite $B_2 \subset \A_0=f_1(\A_1)$ so that $A_1 \subset B_2$ and $t(a_1/B_2)$ does not split over $\A_0$.
Now there is some finite set $C_2 \subset f_1(\A_2)$ so that $t(a/C_2)$ splits over $B_2$ and some automorphism $f_2 \in \textrm{Aut}(\M/B_2)$ so that $f_2(f_1(\A_2))=\A_0$.
Player I plays $a_2=f_2(a_1)$ and $A_2'=f_2(C_2)$.
Continuing in this manner, he can survive $n$ many moves. 
\end{proof}

\begin{definition}\label{Urank}
For $a$ and a model $\A$, we define the \emph{$U$-rank} of $a$ over $\A$, denoted $U(a/\A )$, as follows:
\begin{itemize}
\item  $U(a/\A )\ge 0$ always;
\item $U(a/\A )\ge n+1$ if there is some model $\B$ so that $\A \subseteq \B$, 
$a\nda^{ns}_{\A}\B$ and $U(a/\B )\ge n$;
\item $U(a/\A )$ is the largest $n$ such that $U(a/\A )\ge n$.
\end{itemize}

For finite $A$ we write $U(a/A)$ for 
$\textrm{max}(\{U(a/\A) \, | \, \A \textrm{ is a model s.t. } A \subset \A \})$.
 \end{definition}

\begin{lemma}\label{Unonsplit}
Let  $\A\subseteq\B$ be models. Then
$a\da_{\A}^{ns}\B$ if and only if $U(a/\B )=U(a/\A )$.
\end{lemma}

\begin{proof}
From right to left the claim follows from the definition of $U$-rank.

For the other direction, suppose $a\da^{ns}_{\A}\B$.
It follows from the definition of $U$-rank that $U(a/\B) \le U(a/\A)$.
We will prove $U(a/\A) \le U(a/\B)$.

Let $n=U(a/\A )$, and choose models $\A_i'$, $i \le n$ so that $\A_0'=\A$ and for each $i<n$, $\A_i' \subseteq \A_{i+1}'$ and $a\nda^{ns}_{\A_{i}'}\A_{i+1}'$. 
Choose a model $\C$ so that $\A_n'  a \subseteq \C$.
By AVI, there is a model $\B'$ so that $t(\B'/\A)=t(\B/\A)$ and $\B' \da^{ns}_\A \C$.
Let $f \in \textrm{Aut}(\M/\A)$ be such that $f(\B')=\B$.
Denote $f(a)=b$ and $f(\A_i')=\A_i$ for $i \le n$.
Then, $\A_0=\A$, $t(b/\A)=t(a/\A)$ and
$b\nda^{ns}_{\A_{i}}\A_{i+1}$ for all $i<n$,
and $\B\da^{ns}_{\A} \A_n b$.

Let $\B_1$ be the unique $s$-prime model over $\B \cup \A_1$ (It exists by AIII since $\A \subseteq \B \cap \A_1$ and $\B \da_\A^{ns} \A_1$).
Suppose now that for $1 \le i <n$, $\B_{i-1} \da_{\A_{i-1}}^{ns} \A_i$, and  that we have defined $\B_i$ as the unique $s$-prime model over $\B_{i-1} \cup \A_i$ (taking $\B_0=\B$).
Then, we letÊ $\B_{i+1}$ be the unique $s$-prime model over $\B_i \cup \A_{i+1}$.
It exists, since from the "Furthermore" part in AIII it follows that $\B_i  \da_{\A_{i}}^{ns} \A_{i+1}.$
 
By Lemma \ref{vaplaajykskas}, $t(b/\B)=t(a/\B)$.
Thus, to show that $U(a/\B) \ge U(a/\A)$, it is enough that
$b\nda^{ns}_{\B_{i}}\B_{i+1}$ for all $i<n$.  
Suppose for the sake of contradiction that $b\da^{ns}_{\B_{i}}\B_{i+1}$ for some $i<n$.
Using induction and the "Furthermore" part in AIII, we get that $\B_i \da_{\A_i}^{ns} \A_n b$, and hence
by monotonicity and AV, $b \da_{\A_i}^{ns} \B_i$.
On the other hand, the counterassumption and monotonicity give $b\da^{ns}_{\B_{i}}\A_{i+1}$.
But from these two and Lemma \ref{transitivity}, it follows that
$b \da^{ns}_{\A_{i}}\A_{i+1}$,
a contradiction.
\end{proof}

\section{Indiscernible and Morley sequences}

\begin{definition}
We say that a sequence $(a_{i})_{i<\a}$ is \emph{indiscernible} over $A$
if every permutation of the sequence $\{ a_{i}\vert\ i<\a\}$
extends to an automorphism $f\in \textrm{Aut}(\M /A)$.

We say that a sequence $(a_{i})_{i<\a}$ is \emph{weakly indiscernible} over $A$
if every permutation of a finite subset of the sequence $\{ a_{i}\vert\ i<\a\}$
extends to an automorphism $f\in \textrm{Aut}(\M /A)$.

We say a sequence $(a_{i})_{i<\a}$ is \emph{strongly indiscernible} over $A$ if for all
cardinals $\kappa$, there are $a_{i}$, $\a\le i<\kappa$, such that
$(a_{i})_{i<\kappa}$ is indiscernible over $A$.

Let $\A$ be a model.
We say a sequence $(a_{i})_{i<\a}$ is \emph{Morley} over $\A$,  
if for all $i<\a$, $t(a_{i}/\A )=t(a_{0}/\A )$ and
$a_{i}\da_{\A}^{ns}\cup_{j<i}a_{j}$.
 \end{definition}

In the rest of this chapter, we will assume that all indiscernible sequences and Morley sequences that we consider are non-trivial,  i.e. they do not just repeat the same element.

\begin{lemma}\label{fodor}
Let $A$ be a finite set and $\kappa$ a cardinal such that $\kappa=\textrm{cf}(\kappa)>\omega$.
For every sequence $(a_i)_{i<\kappa}$, there is a model $\A \supset A$ and some $X \subset \kappa$ cofinal so that $(a_i)_{i \in X}$ is Morley over $\A$.
\end{lemma}

\begin{proof}
For $i<\kappa$, choose models $\A_i$ so that for each $i$, $A \subset \A_i$, $a_i \in \A_{i+1}$, $\A_j \subset \A_i$ for $j<i$, $\A_\g=\bigcup_{i<\g} \A_i$ for a limit $\g$, and $\vert \A_i \vert = \vert i \vert +\o$.
Then, for each limit $i$, there is some $\a_i<i$ so that $a_i \da_{\A_{\a_i}}^{ns} \A_i$ (By Lemma \ref{ansA}, there is some finite $A_i \subset \A_i$ so that $a_i \da_{A_{i}}^{ns} \A_i$; just choose $\a_i$ so that $A_i \subset \A_{\a_i}$).
By Fodor's Lemma, there is some $X' \subset \kappa$ cofinal and some $\a<\kappa$ so that $\a_i=\a$ for all $i \in X'$.
Choose $\A=\A_\a$.
By Lemma \ref{notypes}, there are at most $\vert \A\vert<\kappa$ many weak types over $\A$, and thus by the pigeonhole principle, there is some cofinal $X \subseteq X'$ so that $t(a_i/\A)=t(a_j/\A)$ for all $i,j \in X$.
\end{proof}

\begin{lemma}\label{morleyvapaa}
If $(a_{i})_{i<\a}$ is Morley over a countable model $\A$, then
for all $i<\a$, $a_{i}\da_{\A}^{ns}\cup\{ a_{j}\vert\ j<\a,\ j\ne i\}$.
\end{lemma}
 
\begin{proof}  
The claim holds if $a_i \da_\A^{ns} S$ for every finite $S \subset \cup\{ a_{j}\vert\ j<\a,\ j\ne i\}$.
Since we can always relabel the indices, 
it  thus suffices to show that for all $n<\o$, $a_i \da_\A^{ns} \{a_j \, | \, j \neq i, j \le n\}$.
We will prove that for any $n<\o$, if $n=I \cup J$, where $I \cap J= \emptyset$, then $\bigcup_{i \in I}  a_i \da^{ns}_\A \bigcup_{i \in J} a_i$, and the lemma will follow.
We do this by induction on $n$.
If $n=1$, the claim holds trivially, and if $n=2$, it follows directly from Remark \ref{symremark}.
Suppose now the claim holds for $n$, and consider the partition of $n+1$ into the sets $I$ and $J \cup \{n\}$.
Let $a_n'$ be such that $t(a_n'/\A)=t(a_n/\A)$ and
$$a_n' \da^{ns}_\A \A[a_i \, | \,i \in J]\cup \bigcup_{i<n}a_i.$$
Then, in particular, $a_n' \da_\A \bigcup_{i<n} a_i$, so 
 $t(a_n'/ \A \cup \bigcup_{i<n} a_i)=t(a_n/ \A \cup \bigcup_{i<n} a_i)$.
Now, 
$$a_n' \da_{\A[a_i \, | \,i \in J]} \bigcup_{i \in I} a_i,$$
and by Remark \ref{symremark} and monotonicity,
$$\bigcup_{i \in I} a_i \da_{ \A[a_i \, | \,i \in J]} a_n' \cup \bigcup_{i \in J} a_i.$$
By the inductive assumption, we have $\bigcup_{i \in I} a_i \da^{ns}_\A \bigcup_{i \in J} a_i,$ and thus, by Remark \ref{symremark} and Lemma \ref{domination},
$$\bigcup_{i \in I} a_i \da_\A^{ns} \A [a_i \, | \,i \in J].$$
Hence, by Lemma \ref{transitivity},
$$\bigcup_{i \in I} a_i \da_\A^{ns} a_n' \cup\bigcup_{i \in J} a_i,$$
and since $t(a_n'/ \A \cup \bigcup_{i<n} a_i)=t(a_n/ \A \cup \bigcup_{i<n} a_i)$, we have 
$$\bigcup_{i \in I} a_i \da_\A^{ns} \bigcup_{i \in J} a_i \cup \{a_n\},$$
as wanted.  
\end{proof}

\begin{lemma}
If $\A$ is a countable model, then Morley sequences over $\A$
are strongly indiscernible over $\A$.
\end{lemma} 
\begin{proof}
We show first that Morley sequences are weakly indiscernible.
If a sequence $(a_i)_{i<\a}$ is Morley over some model $\A$, then also every finite subsequence is Morley over $\A$.
Thus, as we may relabel any finite subsequence, it suffices 
to show that if a sequence $(a_i)_{i\le n}$, where $n \in \o$, is Morley over a model $\A$, then it is indiscernible over $\A$, i.e. that every permutation extends to an automorphism $f \in \textrm{Aut}(\M/\A)$.
We do this by induction on $n$.
The case $n=0$ is clear.
Suppose now $n=m+1$, where $m\ge0$.
We can obtain any permutation of the $a_i$, $i\le m+1$, by first permuting the $m$ first elements, then changing the place of the two last elements and permuting the $m$ first elements again.
Thus, it is enough to find some  $f \in \textrm{Aut}(\M / \A)$ so that $f(a_i)=a_i$ for $i<m$, $f(a_m)=a_{m+1}$ and $f(a_{m+1})=a_m$.

Since $t(a_m/\A)=t(a_{m+1}/\A)$, $a_m \da_\A^{ns} (a_i)_{i<m}$ and $a_{m+1} \da_\A^{ns} (a_i)_{i<m}$, we have by Lemma \ref{vaplaajykskas} that 
$t(a_m/ \A (a_i)_{i<m})=t(a_{m+1}/ \A (a_i)_{i<m})$ and thus there is some $g_1 \in \textrm{Aut}(\M/\A (a_i)_{i<m})$ such that $g_1(a_{m})=a_{m+1}$.
By Lemma \ref{morleyvapaa}, we have 
$$a_m \da_\A^{ns} (a_i)_{i < m} a_{m+1}.$$
Since $a_{m+1} \da_\A^{ns} (a_i)_{i\le m}$, we have 
$$g_1(a_{m+1}) \da_\A^{ns} (a_i)_{i< m} g_1(a_m),$$ 
so 
$$g_1(a_{m+1}) \da_\A^{ns} (a_i)_{i< m} a_{m+1}$$ since $g_1(a_m)=a_{m+1}$.
Thus, by Lemma \ref{vaplaajykskas}, there is some $g_2 \in \textrm{Aut}(\M/\A (a_i)_{i<m} a_{m+1})$ such that $g_2(g_1(a_{m+1}))=a_{m}$.
Then, $f= g_2 \circ g_1$ is the desired automorphism.
 
Next, we show that Morley sequences are indiscernible. 
Let $(a_i)_{i \in I}$ be a Morley sequence over $\A$, and let $\pi \in \textrm{Sym}(I)$ be a permutation.
We need to show that $\pi$ extends to some $F \in \textrm{Aut}(\M/\A)$.
This is done by constructing models $\A_i$ for $i<\kappa$ so that $\A_0=\A$, for each $i$, $\A_{i+1}$ is the unique $s$-prime model over $\A_i \cup \A[a_i]$, and unions are taken at limit steps.
For this we need to show that these $s$-prime models exist, i.e. that for each $i$, $\A_i \da^{ns}_\A \A[a_i]$.
 
By Lemmas \ref{symlemma} and \ref{domination}, it suffices to show that $a_i \da_{\A}^{ns} \A_i$.
For this, we will show that $a_{i_0}, a_{i_1}, \ldots, a_{i_n} \da_\A^{ns} \A_i$ for $i \le i_0< \ldots <i_n$ (the claim then clearly follows).
We prove this by induction on $i$.
The claim holds for $i=0$, since $\A_0=\A$.
Suppose now it holds for $j$. We show it holds for $j+1$.
For this, we will need two auxiliary claims.

\begin{claim}\label{domi}
The element $a_j$ dominates $\A[a_j]$ over $\A_j$.
\end{claim}

\begin{proof}
Let $c$ be such that $c \da_{\A_j}^{ns} a_j$.
By the inductive assumption, we have $a_j \da^{ns}_\A \A_j$, and thus, by symmetry and transitivity, $\A_j c \da^{ns}_\A a_j$.
By Lemma \ref{domination}, $\A_j c \da^{ns}_\A \A[a_j]$, and hence $c \da^{ns}_{\A_j} \A[a_j]$, as wanted.
\end{proof}

\begin{claim}\label{2}
Let $\B$ be a model such that $\A \subseteq \B$.
Suppose $a \da_\A^{ns} b$ and $ab \da_\A^{ns} \B$. 
Then, $a \da_\B^{ns} b$.
\end{claim}

\begin{proof}
Suppose not.
Choose some finite $A \subset \A$ such that $ab \da_A^{ns} \B$ and $a \da_A^{ns} \A b$.
Since $t(a/\B b)$ splits over $A$, there is some $c \in \B$ so that $a \not\da^{ns}_A bc$.
Choose $c' \in \A$ so that $t(c'/A)=t(c/A)$.
If we would have $t(c'/Aab) \neq t(c/Aab)$, then the pair $c,c'$ would witness that $ab \not\da_A \B$.
Hence, $t(c/Aab)=t(c'/Aab)$.
But now we have $a \not\da^{ns}_A bc'$, so $a \not\da^{ns}_A \A b$, a contradiction.
\end{proof}

Let  now $j<i_0< \ldots < i_n$.
By the inductive assumption, $a_j, a_{i_0}, \ldots, a_{i_n} \da^{ns}_\A \A_j$.
Thus, by Lemma \ref{morleyvapaa} and Claim \ref{2},
$$a_{i_0}, \ldots, a_{i_n} \da^{ns}_{\A_j} a_j,$$
and hence, by  Remark \ref{symremark} and Claim \ref{domi},
$$\A[a_j] \da^{ns}_{\A_j} \A_j a_{i_0}, \ldots, a_{i_n}.$$
By the inductive assumption we have  $ \A[a_j] \da_\A^{ns} \A_j$.
This, and Lemma \ref{transitivity} give
 $$\A[a_j] \da^{ns}_\A \A_j a_{i_0}, \ldots, a_{i_n},$$
and hence by the domination part in AIII,
$$\A_{j+1} \da^{ns}_{\A_j} \A_j a_{i_0}, \ldots, a_{i_n},$$
so
$$a_{i_0}, \ldots, a_{i_n} \da_{\A_j}^{ns} \A_{j+1}.$$ 
By applying the inductive assumption and transitivity, we get
$a_{i_0}, \ldots, a_{i_n} \da^{ns}_\A \A_{j+1},$
as wanted.

Let now $i$ be a limit ordinal.
Then, $\A_j \da^{ns}_{\A} a_i, a_{i_0}, \ldots, a_{i_n}$ for all successor ordinals $j<i$ and $i<i_0<\ldots<i_n$. 
Since $\A_i=\bigcup_{j<i} \A_i$, we have $\A_i \da^{ns}_{\A} a_i, a_{i_0}, \ldots, a_{i_n}$.

Thus, we have shown that the $s$-prime models required for the construction indeed exist.
Now, we construct models $\A_i^\pi$ so that $\A_0^\pi=\A$, for each $i$, $\A_{i+1}^\pi$ is the unique $s$-prime model over $\A_{i}^\pi \cup \A[a_{\pi(i)}]$, and at limit stages unions are taken.
We have already shown that any permutation of finitely many elements of the sequence $(a_i)_{i \in I}$ extends to an automorphism of $\M$ fixing $\A$.
Since being a Morley sequence is a local property (i.e. determined by finite subsequences of a sequence), also the sequence $(a_{\pi(i)})_{i \in I}$ is Morley.
Thus, the models $\A_i^\pi$ exist for each $i \in I$.

We claim that for each $i$, there is an isomorphism $F_i: \A_i \to \A_i^\pi$ fixing $\A$ pointwise.
Since $(a_{\pi(i)})_{i \in I}$ is a Morley sequence, we have $\A[a_{\pi(i)}] \da_\A \A_i^{\pi}$.

Clearly we may choose $F_0=id \raj \A$.
Suppose now the claim holds for $i$.
Now, $\A_i^\pi$ is isomorphic to $\A_i$ over $\A$, and by Lemma \ref{primary}, there is some mapping $f_i \in \textrm{Aut}(\M/\A)$ such that $f_i(\A[a_i])=\A[a_{\pi(i)}]$.
 
Now $\A[a_{\pi(i)}] \da^{ns}_\A f_i(\A_i)$, and similarly as in the proof of Lemma \ref{typegalois}, on sees that the map $F_i \cup f_i: \A_i \cup \A [a_i] \to \A_i^{\pi} \cup \A[\pi(i)]$ is weakly elementary.
Thus, it extends to an elementary map $F_{i+1}: \A_{i+1} \to \A_{i+1}^\pi$.
If $i$ is a limit, then we set $F_i=\bigcup_{j<i} F_j$.

Now $F=\bigcup_{i \in I} F_i$ is as wanted.

Clearly a Morley sequence can be extended to be arbitrarily long.
Thus, Morley sequences are strongly indiscernible.
\end{proof}

\section{Lascar types and the main independence notion}

In this section, we will present our main independence notion and prove that it has all the usual properties of non-forking. The notion will be based on independence in the sense of Lascar splitting.

\begin{definition}
Let $A$ be a finite set, and let $E$ be an equivalence relation on $M^{n}$, for some $n<\o$.
We say $E$ is \emph{$A$-invariant} if for all $f \in Aut(\M /A)$ and $a,b \in \M$, it holds that if $(a,b) \in  E$, then $(f(a),f(b)) \in E$.
We denote the set of all $A$-invariant equivalence 
relations that have only boundedly many equivalence classes by $E(A)$.
 
We say that $a$ and $b$ have the same \emph{Lascar type} over a set $B$, denoted  $Lt(a/B)=Lt(b/B)$, 
if for all finite
$A\subseteq B$ and all $E\in E(A)$, it holds that $(a,b)\in E$.
\end{definition}

\begin{lemma}\label{silt}
If $(a_{i})_{i<\o}$ is strongly indiscernible over $B$,
then $Lt(a_{i}/B)=Lt(a_{0}/B)$ for all $i<\o$
\end{lemma}

\begin{proof}
For each $\kappa$, there are $a_i$, $\o \le i < \kappa$, so that $(a_i)_{i<\kappa}$ is indiscernible over $B$.
If $E \in E(A)$ for some finite $A \subset B$, then $E$ has only boundedly many classes, and thus, for a large enough $\kappa$, there must be some indices $i<j<\kappa$ so that $(a_i, a_j) \in E$.
But this implies that $(a_i, a_j) \in E$ for all $i,j<\kappa$, and the lemma follows.
\end{proof}
 
\begin{lemma}\label{frown}
Let $\A$ be a model and let $t(a/\A)=t(b/\A)$.
Then, $Lt(a/\A)=Lt(b/\A)$.
\end{lemma}

\begin{proof}
Since the equality of Lascar types is determined locally (i.e. it depends on finite sets only), we may without loss assume that $\A$ is countable.

Since $t(a/\A)=t(b/\A)$, there is a sequence $(a_i)_{i<\o}$ such that $(a)\frown (a_{i})_{i<\o}$
and $(b)\frown (a_{i})_{i<\o}$ are Morley over $\A$. 
Because Morley sequences are strongly indiscernible, $Lt(a/\A)=Lt(b/\A)$ by Lemma \ref{silt}.
\end{proof}

In particular, by Lemma \ref{notypes}, for any finite set $A$, the number of Lascar types $Lt(a/A)$ is countable.
It follows that every equivalence relation $E \in E(A)$ has only countably many equivalence classes.

\begin{lemma}\label{lascsat}
Let $\A$ be a countable model, $A$ a finite set such that $A \subset \A$ and $b \in \M$.
Then, there is some $a \in \A$ such that $Lt(a/A)=Lt(b/A)$.
\end{lemma}

\begin{proof}
Since there are only countably many Lascar types over $A$, there is some countable model $\B$ containing
$A$ and realizing all Lascar types over $A$.
By AI, we can construct an automorphism $f \in \textrm{Aut}(\M/A)$ such that $f(\B)=\A$.
Let $b'=f^{-1}(b)$.
Then, there is some $a' \in \B$ such that $Lt(a'/A)=Lt(b'/A)$.
Let $a=f(a')$.
Then, $a \in \A$ and 
$Lt(a/A)=Lt(f(b')/A)=Lt(b/A).$
\end{proof}

\begin{lemma}\label{lasindisc}
Let $A$ be a finite set and let $a, b \in \M$.
Then, $Lt(a/A)=Lt(b/A)$ if and only if there are
$n<\o$ and
strongly indiscernible sequences $I_{i}$ over $A$, $i\le n$,
such that $a\in I_{0}$, $b\in I_{n}$ and for all
$i<n$, $I_{i}\cap I_{i+1}\ne\emptyset$.
\end{lemma}

\begin{proof}
The implication from right to left follows from Lemma \ref{silt} and the fact that all the strongly indiscernible sequences intersect each other.

For the other direction, we note that "there are
$n<\o$ and
strongly indiscernible sequences $I_{i}$ over $A$, $i\le n$,
such that $a\in I_{0}$, $b\in I_{n}$ and for all
$i<n$, $I_{i}\cap I_{i+1}\ne\emptyset$" is an $A$-invariant equivalence relation.
Since we assume that $Lt(a/A)=Lt(b/A)$, it is enough to prove that this equivalence relation has only boundedly many classes.  

Suppose, for the sake of contradiction, that it has unboundedly many classes.
Then, there is a sequence $(a_i)_{i<\o_1}$ where no two elements are in the same class. 
By Lemma \ref{fodor}, there is some $X \subseteq \o_1$, $\vert X \vert =\o_1$, and a model $\A \supset A$ such that $(a_i)_{i \in X}$ is a Morley sequence over $\A$ and thus strongly indiscernible over $A$.
But now by the definition of our equivalence relation, all the elements $a_i$, $i \in X$ are in the same equivalence class, a contradiction. 
 \end{proof}
 
 Now we are ready to introduce our main independence notion.

\begin{definition}\label{freedom}
Let $A \subset B$ be finite.
We say that $t(a/B)$ \emph{Lascar splits} over $A$,
if there are $b,c\in B$ such that
$Lt(b/A)=Lt(c/A)$ but $t(ab/A)\ne t(ac/A)$.

We say $a$ is free from $C$ over $B$, denoted $a\da_{B}C$, if there is some finite $A\subset B$ such that
for all $D\supseteq B \cup C$, there is some $b$ such that
$t(b/B \cup C)=t(a/B \cup C)$ and $t(b/D)$ does not Lascar split over $A$.
\end{definition}

\begin{remark}\label{freedomremark}
Note that it follows from the above definition that if $ab \da_A B$, then $a \da_A B$.

Also, the independence notion is monotone, i.e. if $A\subseteq B \subseteq C \subseteq D$ and $a \da_A D$, then $a \da_B C$.
Indeed, let $A_0 \subseteq A$ be a finite set witnessing that $a$ is free from $D$ over $A$.
We claim that it also witnesses that $a$ is free from $C$ over $B$.
Let $E$ be an arbitrary set such that $C \subseteq E$.
Since $D \subseteq D \cup E$, there is some 
$a'$ such that $t(a'/D)=t(a/D)$ and $t(a'/D \cup E)$ does not Lascar split over $A_0$.
In particular, $t(a'/C)=t(a/C)$ and $t(a'/E)$ does not Lascar split over $A_0$.
Thus, $a \da_B C$.
\end{remark}

\begin{lemma}\label{kolme}
Let $a \in \M$, let $\A$ be a model and let $B\supseteq\A$.  
The following are equivalent:
\begin{enumerate}[(i)]
\item $a\da_{\A}B$, 
\item $a\da^{ns}_{\A}B$ ,
\item $t(a/B)$ does not Lascar split over some finite $A\subseteq\A$.
\end{enumerate}
\end{lemma}

\begin{proof}
"$(i)\Rightarrow (iii)$" follows from Definition \ref{freedom} by choosing $D=B$.

For "$(ii) \Rightarrow (i)$", suppose $a\da^{ns}_{\A}B$.
Then, there is some finite $A \subset \A$ so that $t(a/B)$ does not split over $A$,
and in paricular $t(a/\A)$ does not split over $A$.
Let $D \supset B$ be arbitrary.
By Lemma \ref{vaplaajol2}, there is some $b$ such that $t(b/\A)=t(a/\A)$ and $t(b/D)$ does not split over $A$.
Since $a \da_\A^{ns} B$ and $b \da_\A^{ns} B$, we have by Lemma \ref{vaplaajykskas} that $t(b/B)=t(a/B)$.
Now, $t(b/D)$ does not Lascar split over $A$.
Indeed, if it would Lascar split, then we could find $c, d \in D$ so that $Lt(c/A)=Lt(d/A)$ 
but $t(bc/A) \ne t(bd/A)$.
By Lemma \ref{notypes}, the equivalence relation ``$t(x/A)=t(y/A)$" has only boundedly many classes, and thus Lascar types imply weak types, so $t(b/D)$ would split over $A$, a contradiction.

For $(iii)\Rightarrow (ii)$,
suppose that $t(a/B)$ does not Lascar split over $A$.
We may without loss assume that
$t(a/\A )$ does not split over $A$ (just enlarge $A$ if necessary).
We claim that $t(a/B)$ does not split over $A$.
If it does, then there are $b,c \in B$  witnessing the splitting. 
Let $\B \subseteq \A$ be a countable model containing $A$.
By Lemma \ref{lascsat}, we find $(b',c') \in \B$ so that $Lt(b',c'/A)=Lt(b,c/A)$.
Since $Lt(b/A)=Lt(b'/A)$ and $Lt(c/A)=Lt(c'/A)$, we must have $t(ab/A)=t(ab'/A)$ and $t(ac/A)=t(ac'/A)$ 
(otherwise $t(a/B)$ would Lascar split over $A$).
But since $t(ab/A) \neq t(ac/A)$, we have
$$t(ab'/A)=t(ab/A) \neq t(ac/A)=t(ac'/A),$$
which means that $t(a/\A)$ splits over $A$, a contradiction.
\end{proof}

\begin{remark}\label{kolmeremark}
Note that from the proof of  "(ii) $\Rightarrow$ (i)" for  Lemma \ref{kolme},
it follows that if $\A$ is a model such that $\A \subseteq B$ and $A \subset \A$ is a finite set so that
$a \da_A^{ns} B$, then $a \da_A B$.
In particular, for all models $\A$ and all $a \in \M$, there is some finite $A \subset \A$ such that $a \da_A \A$.
\end{remark}
 
\begin{lemma}\label{lemmanen}
Suppose $\A$ is a model, $A\subseteq\A$ is finite and
$U(a/A)=U(a/\A )$. Then $a\da_{A}\A$. 
\end{lemma}

\begin{proof}
We claim that it is enough to show that $t(a/\A)$ does not Lascar split over $A$.
Indeed, suppose so.
Choose a finite set $B$ such that $A \subseteq B \subset \A$ and $t(a/\A)$ does not split over $B$.
For an arbitrary $D \supseteq \A$,  
there is some $b$ so that $t(b/\A)=t(a/\A)$ and $t(b/D)$ does not split over $B$.
We will show that $t(b/D)$ does not Lascar split over $A$.
Suppose it does.
Then, we can find $c \in \A$ and $d \in D$ such that $Lt(c/A)=Lt(d/A)$ but $t(bc/A) \neq t(bd/A)$.
By Lemma \ref{lascsat}, there is some $d' \in \A$ such that $Lt(d'/B)=Lt(d/B)$.
Then, either $t(d'b/A) \neq t(cb/A)$ or $t(d'b/A) \neq t(db/A)$.
In the first case $t(b/\A)$ Lascar splits over $A$, and in the second case, $t(b/D)$ splits over $B$.
Both contradict our assumptions.

Suppose now, for the sake of contradiction, that 
$t(a/\A)$ does Lascar split over $A$. 
We enlarge the model $\A$ as follows.
First we go through all pairs $b,c \in \A$ so that $Lt(b/A)=Lt(c/A)$.
For each such pair, we find finitely many strongly indiscernible sequences over $A$ of length $\o_1$ as in Lemma \ref{lasindisc}.
We enlarge $\A$ to contain all these sequences.
After this, we repeat the process $\o$ many times.
Then, for every permutation of a sequence of length $\o_1$ that is strongly indiscernible over $A$ and contained in the model, we choose some automorphism $f \in \textrm{Aut}(\M/A)$ that extends the permutation.
We close the model under all the chosen automorphisms.
Next, we start looking again at pairs in the model that have same Lascar type over $A$ and adding $A$-indiscernible sequences of length $\o_1$ witnessing this.
After repeating the whole process sufficiently long, we have obtained a model $\A^*\supseteq \A$ such that for any $b,c \in \A^{*}$ with $Lt(b/A)=Lt(c/A)$, $\A^{*}$ contains $A$-indiscernible sequences witnessing this, and moreover
every permutation of a sequence of length $\o_1$ that is strongly indiscernible over $A$ and contained in $\A^{*}$ extends to an automorphism of $\A^{*}$.
 
Choose now an element $a^*$ so that $t(a^*/\A)=t(a/\A)$ and $a^* \da_\A^{ns} \A^*$.
Then, $U(a^*/\A^*)=U(a^*/A)$ by Lemma \ref{Unonsplit}.
Let $f \in \textrm{Aut}(\M/\A)$ be such that $f(a^*)=a$, and denote $\A'=f(\A^*)$.
Now, $U(a/\A')=U(a/A)$ and $t(a/\A')$ Lascar splits over $A$.

Let $b,c \in \A'$ witness the splitting.
Then, $Lt(b/A)=Lt(c/A)$ and inside $\A'$ there are for some $n<\o$, strongly indiscernible sequences $I_i$, $i \le n$, over $A$ of length $\o_1$ so that $b \in I_0$, $c \in I_n$ and $I_i \cap I_{i+1} \neq \emptyset$ for $i<n$.
Since $t(ab/A) \neq t(ac/A)$, in at least one of these sequences there must be two elements that have different weak types over $Aa$.
Since there are only countably many weak types over $Aa$, this implies that
there is inside $\A'$ a sequence $(a_{i})_{i<\o_{1}}$ strongly indiscernible over $A$
such that  
$t(aa_{0}/A)\ne t(aa_{1}/A)$ but $t(aa_{1}/A)=t(aa_{i}/A)$ for all $0<i<\o_{1}$.
Moreover, every permutation of $(a_{i})_{i<\o_{1}}$
extends to an automorphism $f\in Aut(\A' /A)$.
  
For each $i<\o_1$, let $f_i \in Aut(\M/A)$ be an automorphism permuting the sequence $(a_{i})_{i<\o_{1}}$ so that $f_i(a_0)=a_i$ and $f_i(\A')=\A'$.
Denote $b_i=f_i(a)$ for each $i<\o_1$.
Then, $U(b_i/\A')=U(b_i/A)$ and for all $j<i<\o_1$, $t(b_i/A)=t(b_j/A)$, but $t(b_i/\A') \neq t(b_j/\A')$ since  
$$t(b_i a_i/A)=t(f_i(a) f_i(a_0)/A)=t(a a_0 /A) \neq t(a f_j^{-1}(a_i) /A)=t(f_j(a) a_i/A)=t(b_j a_i /A).$$ 
  
Let $\B\subseteq\A$ be countable model such that $A\subseteq\B$.
Then for all $i<\o$, 
$$U(b_{i}/\A')=U(b_i/\B),$$
so $b_{i}\da^{ns}_{\B}\A'$ by Lemma \ref{Unonsplit}.
Thus, for all $i<j<\o_{1}$,
$t(b_{i}/\B )\ne t(b_{j}/\B )$, a contradiction by Lemma \ref{vaplaajykskas} since there are only countably many types over $\B$. 
\end{proof}

\begin{corollary}\label{typecoro}
For every $a \in \M$, every finite set $A$ and every $B \supseteq A$, there is some $b \in \M$ such that $t(a/A)=t(b/A)$ and $b \da B$.
\end{corollary}
  
\begin{proof}
Let $\A$ be a model such that $U(a/A)=U(a/\A).$
Let $\B$ be a model such that $\A \cup B \subseteq \B$, and let $b$ be such that $t(b/\A)=t(a/\A)$ and $b \da^{ns}_\A \B$.
Then, by Lemma \ref{Unonsplit}, 
$$U(b/\B)=U(b/\A)=U(a/\A)=U(a/A)=U(b/A).$$ 
By Lemma \ref{lemmanen}, $b \da_A \B$, and thus $b \da_A B$.
\end{proof}  
  
\begin{lemma}\label{lasctrans}
Suppose $A\subseteq\A\subseteq B$. Then $a\da_{A}B$ if and only if
$a\da_{A}\A$ and $a\da_{\A}B$.
\end{lemma}

\begin{proof}
"$\Rightarrow$": $a\da_{A}\A$ is clear and $a\da_{\A}B$ follows from Lemma \ref{kolme}.
  
"$\Leftarrow$": 
Since $a \da_A \A$, there is by definition some finite $A_0 \subseteq A$ and some $b$ such that
$t(b/\A )=t(a/\A )$ and $t(b/B)$ does not Lascar split over $A_0$.
By Lemma \ref{kolme}, $b\da^{ns}_{\A}B$ and $a \da^{ns}_\A B$.
Thus, by Lemma \ref{vaplaajykskas}, $t(b/B)=t(a/B)$.
Hence $a\da_{A}B$, as wanted.
\end{proof}

\begin{lemma}\label{lascsym}
Let $A$ be finite. Then, $a\da_{A}b$ if and only if $b\da_{A}a$. 
\end{lemma}

\begin{proof}
Suppose $a \da_A b$.
Let $\A_0$ be a model such that $A \subset \A_0$.
By Corollary \ref{typecoro}, there exists some $b'$ such that $t(b'/A)=t(b/A)$
and $b' \da_A \A_0$.
Let $f \in \textrm{Aut}(\M/A)$ be such that $f(b')=b$, and denote $\A=f(\A')$.
Then, $A \subset \A$ and $b \da_A \A$.
By Definition \ref{freedom},
there is some $a'$ such that
 $t(a'/Ab)=t(a/Ab)$
and $a'\da_{A}\A b$.
Then, $a' \da_\A b$, and by Lemma \ref{kolme} and Remark \ref{symremark}, $b \da_\A a'$.
By Lemma \ref{lasctrans}, $b \da_A a'$, and thus $bÊ\da_A a$.
\end{proof}
 
\begin{lemma}\label{lascvaplaaj}
For every $a$, every finite set $A$ and every  $B\supseteq A$, there is $b$
such that $Lt(b/A)=Lt(a/A)$ and $b\da_{A}B$.
\end{lemma}

\begin{proof}
Let $\A_0$ be a countable model such that $A \subset \A_0$.
By Corollary \ref{typecoro}, there is some element $a'$ so that $t(a'/A)=t(a/A)$ and $a' \da_A \A_0$.
Let $f \in \textrm{Aut}(\M/A)$ be such that $f(a')=a$.
Denote $\A=f(\A_0)$.
Now, $A \subset \A$ and $a \da_A \A$.
 
Choose now  $b$ so that $t(b/\A )=t(a/\A )$ and $b\da^{ns}_{\A}B$.  
Then, $b\da_{\A}B$.
By Lemma \ref{lasctrans}, $b\da_{A}B$.
Moreover, by Lemma \ref{frown}, $Lt(b/A)=Lt(a/A)$.
\end{proof}

\begin{lemma}\label{lascvaplaajykskas}
If $A \subseteq B$, $a\da_{A}B$, $b\da_{A}B$ and $Lt(a/A)=Lt(b/A)$, then
$Lt(a/B)=Lt(b/B)$.
\end{lemma}

\begin{proof}
Clearly it is enought to prove this under the assumption
that $A$ and $B$ are finite (if $Lt(a/B_0)=Lt(b/B_0)$ for every finite $B_0 \subset B$, then $Lt(a/B)=Lt(b/B)$).
Suppose the claim does not hold and choose countable models $\A_a$ and $\A_b$ so that $Aa \subset \A_a$, 
$Ab \subset \A_b$,  $B\da_{A}\A_{a}$ and $B\da_{A}\A_{b}$.  
By Lemma \ref{lascvaplaaj}, there is some $c$ such that $Lt(c/A)=Lt(a/A)$ and $c \da_A \A_a \cup \A_b \cup B$.
By monotonicity, we have $c \da_{\A_a} B$, and thus by Lemma \ref{lascsym}, $B \da_{\A_a} c$.
Hence, by Lemma \ref{lasctrans}, $B\da_{A}\A_{a}c$ and so
$ac\da_{A}B$.

By the counterassumption, we may without loss assume that $Lt(c/B) \neq Lt(a/B)$.
Choose a model $\B\supseteq B$  so that $ac\da_{A}\B$. 
By Lemma \ref{frown}, $t(a/\B )\ne t(c/\B )$. 
So there is some $b'\in\B$ that
withesses this, i.e. $t(ab'/A)\ne t(cb'/A)$. 
As $Lt(c/A)=Lt(a/A)$, this means $t(b'/Aac)$ Lascar splits over $A$, a contradiction since $b' \da_A ac$.
   
\end{proof}
 
\begin{lemma}\label{lasctrans2}
Suppose $A\subseteq B\subseteq C$, $a\da_{A}B$ and $a\da_{B}C$.
Then $a\da_{A}C$. 
\end{lemma}

\begin{proof}
Clearly it is enough to prove this for finite
$A$. 
Choose $b$ so that $Lt(b/A)=Lt(a/A)$ and
$b\da_{A}C$. 
Then, by monotonicity, $b \da_A B$, and thus by
Lemma \ref{lascvaplaajykskas}, 
$Lt(b/B)=Lt(a/B)$.
Again by monotonicity,
$b\da_{B}C$, and by Lemma \ref{lascvaplaajykskas}, 
$Lt(b/C)=Lt(a/C)$.
The claim follows.
\end{proof}

\begin{lemma}\label{local}
Suppose $A \subset B$ and $a \not\da_{A} B$.
Then there is some $b \in B$ such that $a \not\da_{A} b$.
\end{lemma}

\begin{proof}
Choose a finite $C \subseteq A$ such that $a \da_C A$ and an element $c$ such that
$Lt(c/C)=Lt(a/C)$ and $c \da_C A \cup B$ (they exist by Corollary \ref{coro} and  Lemma \ref{lascvaplaaj}).
Then, by Lemma \ref{lascvaplaajykskas}, $Lt(c/A)=Lt(a/A)$.
We have $a \not\da_C B$, and thus $t(c/B) \neq t(a/B)$.
Hence, there is some $b \in B$ so that $t(cb/C) \neq t(ab/C)$.
By monotonicity, we have $c \da_A B$, and in particular $c \da_A b$.
If $a \da_A b$, then $a \da_C b$ by Lemma \ref{lasctrans}.
Since  $t(a/Cb) \neq t(c/Cb)$, this contradicts Lemma \ref{lascvaplaajykskas}.
\end{proof}

\begin{definition}
Let $A$ be a finite set and let $f \in \textrm{Aut}(\M/A)$.
We say that $f$ is a \emph{strong automorphism} over $A$ if it preserves Lascar types over $A$, i.e. if 
for any $a$, $Lt(a/A)=Lt(f(a)/A)$.
We denote the set of strong automorphisms over $A$ by $\textrm{Saut}(\M/A)$.
\end{definition}

\begin{lemma}
Suppose $A$ is finite and $Lt(a/A)=Lt(b/A)$. Then there is
$f\in Saut(\M /A)$ such that $f(a)=b$.
\end{lemma}

\begin{proof}
Choose a countable model
$\A\supseteq A$ such that $ab\da_{A}\A$.
In particular, by Remark \ref{freedomremark}, 
$a\da_{A}\A$ and $b\da_{A}\A$.
By Lemma \ref{lascvaplaajykskas}, $Lt(a/\A)=Lt(b/\A)$. 
Thus, there is some $f\in \textrm{Aut}(\M /\A )$ such that $f(a)=b$.
By Lemma \ref{frown}, $f \in  \textrm{Saut}(\M /\A )$. 
\end{proof}

\begin{lemma}\label{nonsplitU}
Suppose $\A$ is a model, $A\subseteq\A$ is finite and $t(a/\A )$
does not split over $A$. Then $U(a/A)=U(a/\A )$.
\end{lemma}

\begin{proof}
 Suppose not. 
If we choose some countable model $\A'$ such that $A \subset \A' \subseteq \A$,
then $a \da_{\A'}^{ns} \A$, and thus, by Lemma \ref{Unonsplit}, $U(a/\A')=U(a/\A)$.
Hence, we may assume that $\A$ is countable.
 
Choose a countable model $\B$ such that $A \subset \B$ and $U(a/\B )=U(a/A)$. 
Now, there is some $f \in \textrm{Aut}(\M/A)$ so that $f(\B)=\A$.
Let $a'=f(a)$.
We have
$$U(a/\A) \neq U(a/A)=U(a/\B)=U(a'/\A),$$
and thus $t(a/\A) \neq t(a'/\A)$.
Hence there is some $c \in \A$ such that $t(ac/A) \neq t(a'c/A)$.
Let $b \in \B$ be such that $f(b)=c$ (and thus $t(b/A)=t(c/A)$).
Then, $t(a'c/A)=t(ab/A)$, so $t(ac/A) \neq t(ab/A)$.
Let $c' \in \A$ be such that $Lt(c'/A)=Lt(b/A)$, and thus $t(c'/A)=t(b/A)=t(c/A).$
Since $t(a/\A)$ does not split over $A$, we have $t(ac'/A) \neq t(ab/A).$
 
We note that since $a \da_{A}^{ns} \A$, we have by Remark \ref{kolmeremark}  $a \da_{A} \A$, and thus in particular $a \da_A c'$.
 Choose $g\in Saut(\M /A)$ so that
$g(b)=c'$. Let $a''=g(a)$. 
By Lemma \ref{lemmanen}, we have $a \da_A b$, and thus $a'' \da_A c'$.
But now $Lt(a''/A)=Lt(a/A)$,
$a\da_{A}c'$ and $a''\da_{A}c'$, yet 
$$t(ac'/A) \ne t(a'c'/A)=t(ab/A)=t(a''c'/A),$$
so in particular $t(a/Ac') \ne t(a''/Ac')$, 
a contradiction.
\end{proof}

\begin{corollary}\label{Ukoro}
Let $\A$ be a model.
Then, 
$$U(a/\A)=\textrm{min}(\{U(a/B) \, | \, B \subset \A \textrm{ finite } \}.$$
\end{corollary}

\begin{proof} 
By Definition \ref{Urank}, for each finite $B \subset \A$, it holds that $U(a/B) \ge U(a/\A)$. 
On the other hand, by Lemma \ref{ansA}, there is some finite $A \subset \A$ so that $t(a/\A)$ does not split over $A$.
By Lemma \ref{nonsplitU}, $U(a/\A)=U(a/A)$.
 \end{proof}

Corollary \ref{Ukoro} allows us to define $U(a/A)$ for arbitrary $A$ as follows.

\begin{definition}\label{U2}
Let $A$ be arbitrary.
We define  $U(a/A)$ to be the minimum of
$U(a/B)$, $B\subseteq A$ finite.
\end{definition}

\begin{lemma}\label{daiffU}
For all $A\subseteq B$ and $a$, $a\da_{A}B$ if and only if $U(a/A)=U(a/B)$.
\end{lemma}

\begin{proof} 
Suppose first $B$ is finite.

"$\Leftarrow$": Choose a model $\A\supseteq B$ such that
$U(a/\A )=U(a/B)$. 
Then $U(a/\A )=U(a/A )$ and thus by Lemma \ref{lemmanen},
$a\da_{A}\A$, and in particular $a\da_{A}B$.

 "$\Rightarrow$":
Choose a model $\A\supseteq A$ such that
$U(a/\A )=U(a/A)$ and a model  $\B\supseteq\A B$.
By Lemma \ref{vaplaajol}, there is some $a'$ 
such that $t(a'/\A )=t(a/\A )$ and $a'\da^{ns}_{\A}\B$. 
Then, by Lemma \ref{Unonsplit},
$$U(a'/A)=U(a'/\A )=U(a'/\B ),$$ 
so by Lemma \ref{lemmanen}, $a'\da_{A}\B$.
By Lemma \ref{frown}, $Lt(a'/A)=Lt(a/A)$, and thus by Lemma \ref{lascvaplaajykskas}, 
 $t(a'/B)=t(a/B)$. 
 Thus
 $$U(a/B)=U(a'/B)=U(a'/A)=U(a/A).$$
 
We now prove the general case.
Let $A, B$ be arbitrary such that $A \subseteq B$.

"$\Leftarrow$": 
Suppose $a \da_A B$.
There are some finite sets $A_0, A_0' \subseteq A$ such that $a \da_{A_0} B$ and some $U(a/A)=U(a/A_0'')$.
We may without loss assume that $A_0=A_0'$.
Indeed, this follows from monotonicity and the fact that if $A_0''$ is any set such that
$A_0' \subseteq A_0'' \subseteq A$, then $U(a/A_0'')=U(a/A)$.
Let $B_0 \subseteq B$ be a finite set such that $U(a/B_0)=U(a/B)$.
By similar argument as above, we may without loss suppose that $A_0 \subseteq B_0$.
Thus, since the result holds for finite sets, we have 
$$U(a/A)=U(a/A_0)=U(a/B_0)=U(a/B).$$

"$\Rightarrow$":
Suppose $U(a/A)=U(a/B)$, but $a \nda_A B$.
Let $A_0 \subseteq A$ and $B_0 \subseteq B$ be finite sets such that 
$$U(a/A_0)=U(a/A)=U(a/B)=U(a/B_0).$$
By monotonicity, we have $a \nda_{A_0} B$, and by Lemma \ref{local}, there is some $b \in B$ such that
$a \nda_{A_0} b$.
Then, also $a \nda_{A_0} B_0 b$.
But 
$$U(a/B_0b)=U(a/B)=U(a/A_0),$$
a contradiction.
\end{proof}

\begin{corollary}\label{coro}
For all $A$ and $a$ there is finite $B\subseteq A$ such that
$a\da_{B}A$. 
\end{corollary}

\begin{corollary}\label{lascvaplaaj2}
For all $a$ and all sets $A \subseteq B$, there is some $b$ such that $Lt(b/A)=Lt(a/A)$ and $b \da_A B$.
\end{corollary}

\begin{proof}
Let $A_0 \subseteq A$ be a finite set such that $U(a/A_0)=U(a/A)$.
Then, $a \da_{A_0} A$ by Lemma \ref{daiffU}.
By Lemma \ref{lascvaplaaj}, there is some $b$ such that $Lt(b/A_0)=Lt(a/A_0)$ and $b \da_{A_0} B$.
By Lemma \ref{lascvaplaajykskas}, $Lt(b/A)=Lt(a/A)$.
\end{proof}
  
\begin{lemma}\label{lascfullsym}
Let $A$ be arbitrary.
If $a\da_{A}b$, then $b\da_{A}a$. 
\end{lemma}

\begin{proof}
Suppose not.
Choose some finite $B\subseteq A$ such that
$a\da_{B}Ab$ and $b\da_{B}A$ (such a set can be found by Corollary \ref{coro}).  
Since $b \nda_A a$, we have $b \nda_B Aa$.
By Lemma \ref{local}, there is some finite set $C$ such that
$B \subseteq C \subseteq A$ and $b\not\da_{B}Ca$.  
By transitivity, $b\nda_{C}a$.
On the other hand, $a \da_B Ab$, and thus $a \da_B Cb$, so $a \da_C b$, which contradicts Lemma \ref{lascsym}.
\end{proof}

\begin{lemma}\label{Uaddit}
For any $a,b$ and $A$, it holds that
$$U(ab/A)=U(a/bA)+U(b/A).$$
\end{lemma}

\begin{proof}
We first note that it suffices to prove the lemma in case $A$ is finite.
Indeed, by definition $\ref{U2}$, we find finite $A_1, A_2, A_3 \subset A$ so that $U(ab/A)=U(ab/A_1)$,
$U(a/bA)=U(a/bA_2)$ and $U(b/A)=U(b/A_3)$.
Denote $A_0=A_1 \cup A_2 \cup A_3$.
Since the above ranks are minimal, we have $U(ab/A)=U(ab/A_0)$,
$U(a/bA)=U(a/bA_0)$ and $U(b/A)=U(b/A_0)$.
Thus it suffices to show that the lemma holds for $A_0$, a finite set.
 
Next, we show that for any $c$ and any finite set $B$, $U(c/B)$ is the maximal number $n$ such that there are
sets $B_i$, $i\le n$ so that $B_0=B$, and for all $i<n$, $B_i \subseteq B_{i+1}$ and $c \not\da_{B_i} B_{i+1}$.
By Lemma \ref{daiffU}, $U(c/B_i)>U(c/B_{i+1})$ for all $i<n$, and thus, $U(c/B) \ge n$.
On the other hand, by the definition of $U$-rank (Definition \ref{Urank}), there are models $\B_i$, $i \le m=U(c/B)$, so that 
$B \subset \B_0$, and for each $i<m$, $\B_i \subset \B_{i+1}$ and $c \not\da_{\B_i} \B_{i+1}$.
Write $B_0=B$.
By Lemma \ref{local}, for each $1 \le i<m$, we find some finite $B_i \subset \B_i$ so that $c \not\da_{B_{i-1}} B_i$.
Thus, $n \ge m=U(c/B)$.

To show $U(ab/A) \le U(a/bA)+U(b/A)$, we let $n=U(ab/A)$ and $A_i$, $i\le n$ be as above for $U(ab/A)$.
Then, for each $i <n$, we must have either $a \not\da_{bA_i} A_{i+1}$ or $b \not\da_{A_i} A_{i+1}$.
Indeed, if we would have both $a \da_{bA_i} A_{i+1}$ and $b \da_{A_i} A_{i+1}$, then by Lemma \ref{lascsym}, we would have $A_{i+1} \da_{A_i} b$ and $A_{i+1} \da_{bA_i} a$, and thus by applying first Lemma \ref{lasctrans2} and monotonicity, then Lemma \ref{lascsym} again, we would get $ab \da_{A_i} A_{i+1}$.
Thus, $U(a/bA)+U(b/A) \ge n$.
 
Let now $U(b/A)=m$ and let $A_i'$, $i \le m$ be the sets witnessing this (here $A_0'=A$).
Choose $a'$ so that $t(a'/Ab)=t(a/Ab)$ and $a' \da_{bA} A_m'$.
Using a suitable automorphism, we find $A_i$, $i \le m$, also witnessing $U(b/A)=m$ so that $a \da_{bA} A_m$.
Thus, by Lemma \ref{daiffU}, $U(a/bA_m)=U(a/bA)$.
Let $U(a/bA_m)=k$ and choose $B_i$, $i \le k$ witnessing this.
Now, $A=A_0, \ldots, A_{m-1}, B_0, \ldots, B_k$ witness that $U(ab/A) \ge m+k$ (note that we may without loss assume that $A_m=B_0$).
\end{proof}

\begin{lemma}\label{statilemma}
Suppose $\A$ is a model, $A$ a finite set such that $A \subset \A$, $B$ is such that $\A \subseteq B$,
$b \da_A^{ns} B$, $b' \da_A^{ns} B$ and $t(b/A)=t(b'/A)$.
Then, $t(b/B)=t(b'/B)$. 
\end{lemma}

\begin{proof}
By Lemma \ref{lascvaplaajykskas}, it suffices to show that $Lt(b/A)=Lt(b'/A)$. 
By Lemma \ref{frown}, this follows after we have shown that $t(b/\A)=t(b'/\A)$.
Suppose not. 
Then, there is some $a \in \A$ such that $t(ab/A) \neq t(ab'/A)$.
Let $f \in \textrm{Aut}(\M/A)$ be such that $f(b')=b$, and let $a'=f(a)$.
Since $b' \da_A^{ns} \A$, we have $b' \da_A \A$ by Remark \ref{kolmeremark}, and in particular, $b' \da_A a$.
Thus, $b \da_A a'$.
By Lemma \ref{lascsat} there is some $a'' \in \A$ such that $Lt(a''/A)=Lt(a'/A)$.
Since $b \da_A^{ns} \A$, we have $b \da_A a''$.
Then, $Lt(a''/Ab)=Lt(a'/Ab)$ by Lemmas \ref{lascfullsym} and \ref{lascvaplaajykskas}.
Now,
$$t(ab/A) \neq t(ab'/A)=t(a'b/A)=t(a''b/A),$$
a contradiction since $t(b/\A)$ does not split over $A$.
\end{proof}

\begin{definition}
We say that a set $A$ is \emph{bounded} if $\vert A \vert < \delta= \vert \M \vert$.
\end{definition}

\begin{definition}
We say $a$ is in the \emph{bounded closure} of $A$, denoted $a \in \textrm{bcl}(A)$, if $t(a/A)$ has only boundedly many realizations.
\end{definition}

\begin{lemma}\label{bclmalli}
Let $\A$ be a model.
Then, $\textrm{bcl}(\A)=\A$.
\end{lemma}

\begin{proof}
Clearly $\A \subseteq \textrm{bcl}(\A)$.
For the converse, suppose towards a contradiction that $a \in \textrm{bcl}(\A) \setminus \A$.
By Lemma \ref{ansA}, there is some finite $A \subset \A$ so that $a \da_A^{ns} \A$.  Choose now an element $a'$ such that $t(a'/\A)=t(a/\A)$ and $a' \da^{ns}_A \textrm{bcl}(\A)$.
Then, $a' \in \textrm{bcl}(\A)$.
By Axiom I, there is some $b \in \A$ such that $t(b/A)=t(a'/A)$ and thus $b \neq a'$.
In particular, $t(a'a'/A) \neq t(ba'/A)$.
Thus, $a'$ and $b$ witness that $t(a'/\textrm{bcl}(\A))$ splits over $A$, a contradiction.
\end{proof}

\begin{lemma}\label{bcllocal}
If $a \in \textrm{bcl}(A)$, then there is some finite $B \subseteq A$ so that $a \in \textrm{bcl}(B)$.
\end{lemma}

\begin{proof}
There is some finite $B \subseteq A$ such that $a \da_B A$.
We claim that $a \in \textrm{bcl}(B)$.
Suppose not.
Let $\A$ be a model such that $A \subseteq \A$.
Now there is some $a'$ so that $Lt(a'/A)=Lt(a/A)$ and $a' \da_B \A$.
By Lemma \ref{bclmalli}, $a' \in \textrm{bcl}(A) \subseteq \textrm{bcl}(\A)=\A$.
Since $a \notin \textrm{bcl}(B)$, the weak type $t(a/B)$ has unboundedly many realizations.
Hence, by Lemma \ref{fodor}, there is a Morley sequence $(a_i)_{i<\o}$ over some model $\B \supset B$ so that $a_0=a'$ (just use a suitable automorphism to obtain this).
By Axiom AI, there is an element $a'' \in \A$ so that $t(a''/a'B)=t(a_1/a'B)$,
and by Lemma \ref{silt} , $Lt(a_1/B)=Lt(a'/B)$.
Thus, there is an automorphism $f \in \textrm{Aut}(\M/B)$ such that $f(a'')=a_1$ and $f(a')=a'$.
Using Lemma \ref{lasindisc}, one sees that automorphisms preserve equality of Lascar types.
Hence, the fact that $Lt(a_1/B)=Lt(a'/B)$ implies 
$Lt(a''/B)=Lt(a'/B).$ 
But we have $a'=a_0 \neq a_1$, and thus also $a'' \neq a'$, so
$t(a'a'/B) \neq t(a'a''/B)$, which contradicts Lemma \ref{lascvaplaajykskas} since we assumed $a' \da_B \A$.
\end{proof} 

\begin{lemma}\label{xiv}
For every $A$, $\textrm{bcl}(\textrm{bcl}(A))=\textrm{bcl}(A)$.
\end{lemma}
 
\begin{proof}
By Lemma \ref{bcllocal}, we may assume that $A$ is finite.
Suppose now $a \in \textrm{bcl}(\textrm{bcl}(A)) \setminus \textrm{bcl}(A)$.
By Lemma \ref{bcllocal}, there is some $b \in \textrm{bcl}(A)$ so that $a \in \textrm{bcl}(Ab)$.
Let $\kappa$ be an uncountable cardinal such that $\kappa> \vert \textrm{bcl}(\textrm{bcl}(A)) \vert$.
Since $a \notin \textrm{bcl}(A)$, there are $a_i$, $i<\kappa$ so that $a_i \neq a_j$ when $i \neq j$ and $t(a_i/A)=t(a/A)$ for all $i<\kappa$.
For each $i$, there is some $b_i \in \textrm{bcl}(A)$ such that $t(b_ia_i/A)=t(ba/A)$.
By the pigeonhole principle, there is some $b' $ and some $X \subseteq \kappa$
so that $\vert X \vert = \kappa$ and $b_i=b'$ for $i \in X$.
Hence, for any $i \in X$, $t(a_i/Ab')$ has unboundedly many realizations, a contradiction since $a_i \in \textrm{bcl}(Ab')$.
\end{proof} 

\begin{lemma}\label{xi}
Let $A \subset B$.
If $a \in \textrm{bcl}(A)$, then $a \da_A B$.
\end{lemma}

\begin{proof}
By Lemma \ref{bcllocal}, we may assume that $A$ is finite.
Choose $a'$ so that $Lt(a'/A)=Lt(a/A)$ and $a' \da_A B$.
Then, $a' \in \textrm{bcl}(A)$.
Consider the equivalence relation $E$ defined so that $(x,y) \in E$ if either $x,y \notin \textrm{bcl}(A)$ or $x=y \in \textrm{bcl}(A)$.
This is an $A$-invariant equivalence relation.
Moreover, since $A$ is finite, we may choose a countable model $\A$ so that $A \subset \A$. 
By Lemma \ref{bclmalli}, $\textrm{bcl}(A) \subset \A$, so $E$ has boundedly many classes, and thus $(a,a') \in E$.
It follows that $a=a'$.
\end{proof}

\begin{lemma}\label{xii}
If $a \in \textrm{bcl}(B) \setminus \textrm{bcl}(A)$, then $a \not\da_A B$. 
\end{lemma} 

\begin{proof}
Suppose $a \da_A B$.
Choose a model $\A$ so that $B \subseteq \A$ and $a'$ so that $t(a'/B)=t(a/B)$ and $a' \da_A \A$.
By Lemma \ref{bclmalli}, $a' \in \A$.
Now we proceed as in the proof of Lemma \ref{bcllocal} to obtain a contradiction.
\end{proof}

Now we have shown that our main independence notion $\da$ has all the properties
of non-forking.

\begin{theorem}\label{main}
Suppose $A \subseteq B \subseteq C \subseteq D$.
Then, the following hold.
\begin{enumerate}[(i)]
\item For each $a$, there is some finite $A_0 \subseteq A$ such that $a \da_{A_0} A$.
\item If $a \not\da_{A} B$, then there is some $b \in B$ so that $a \not\da_{A} b$.
\item Suppose that $Lt(a/A)=Lt(b/A)$, $a \da_A B$ and $b \da_A B$.
Then, $Lt(a/B)=Lt(b/B)$.
\item For every $a$, there is some $b$ such that $Lt(b/A)=Lt(a/A)$ and $b \da_A B$.
\item If $a \da_A D$, then $a \da_B C$.
\item If $a \da_A B$ and $a \da_B C$, then $a \da_A C$.
\item If $a \da_A b$, then $b \da_A a$.
\item $a \da_A B$ if and only if $U(a/B)=U(a/A)$.
\item For all $a$, $U(a/\emptyset)<\o$.
\item For all $a,b$, $U(ab/A)=U(a/bA)+U(b/A)$.
\item If $a \in \textrm{bcl}(A)$, then $a \da_A B$.
\item If $a \in \textrm{bcl}(B) \setminus \textrm{bcl}(A)$, then $a \not\da_A B$. 
\item If $a \in \textrm{bcl}(A)$, then there is some finite $A_0 \subseteq A$ so that $a \in \textrm{bcl}(A_0)$.
\item $\textrm{bcl}(\textrm{bcl}(A))=\textrm{bcl}(A)$.
\item If $\A$ is a model, then $\textrm{bcl}(\A)=\A$.
\end{enumerate}
\end{theorem}

\begin{proof}
(i) This is Corollary \ref{coro}.

(ii) Lemma \ref{local}.

(iii) Lemma \ref{lascvaplaajykskas}.

(iv) Corollary \ref{lascvaplaaj2}.

(v) Remark \ref{freedomremark}. 

(vi) Lemma \ref{lasctrans2}

(vii) Lemma \ref{lascfullsym}.

(viii) Lemma \ref{daiffU}.

(ix) This follows from Definition \ref{Urank} and Lemmas \ref{Unonsplit} and \ref{eiaaretketjuja}.

(x) Lemma \ref{Uaddit}.
 
(xi) Lemma \ref{xi}.

(xii) Lemma \ref{xii}.

(xiii) Lemma \ref{bcllocal}.

(xiv) Lemma \ref{xiv}.

(xv) Lemma \ref{bclmalli}.

\end{proof}
  
\section{$\M^{eq}$ and canonical bases}

Let $\mathcal{E}$ be a countable collection of $\emptyset$-invariant
equivalence relations $E$ such that $E \subseteq \M^n \times \M^n$ for some $n$.
By this we mean that if $E \in \mathcal{E}$, then $E$ is an equivalence relation on some model in $\mathcal{K}$ (note that from this it follows that $E$ is an equivalence relation on every model in $\K$; indeed, it takes at most three tuples to prove that a relation is not an equivalence relation, and by axiom AI all models are $s$-saturated)
and there is some countable collection $G_E$ of Galois-types so that $(a,b) \in E$ if and only if $t^g(ab/\emptyset) \in G_E$.
We assume that the identity relation is in $\mathcal{E}$, $= \in \mathcal{E}$ (note that there are only countably many Galois types over $\emptyset$).
 For every
$\A\in\K$ we let $\A^{eq}$ be the set
$\{ a/E\vert\ a\in\A ,\ E\in\mathcal{E}\}$. 
We identify each element $a$ with $a/=$. 
For each $E \in \mathcal{E}$, we add to
our language a predicate $P_{E}$
with the interpretation $\{ a/E\vert\ a\in\A\}$
and a function $F_{E}: \A^n \to \A^{eq}$ (for a suitable $n$) such that $F_{E}(a)=a/E$.  
Then, we have all the structure of $\A$ on $P_{=}$.  
We let $\K^{eq}=\{\A^{eq}\vert\ \A\in\K\}$.
We write $\A^{eq}\preccurlyeq^{eq}\B^{eq}$ if
$\A^{eq}$ is a submodel of $\B^{eq}$ and $\A\preccurlyeq\B$.

We will show that $(\mathcal{K}^{eq},\preccurlyeq^{eq})$ is an AEC with AP, JEP and arbitrary
large models, that $LS(\mathcal{K}^{eq})=\o$ and that $\mathcal{K}^{eq}$
does not contain finite models.
Moreover, if $\mathcal{K}$ satisfies the axioms AI-AVI listed in the first section, then also $\mathcal{K}^{eq}$ satisfies them.

Notice first that for each model $\A$, the model $\A^{eq}$ is unique up to isomorphism over $\A$
and that every automorphism of $\A$ extends to an automorphism
of $\A^{eq}$. 
Thus it is easy to see that
$(\K^{eq},\preccurlyeq^{eq})$ is an AEC with AP, JEP and arbitrary
large models,  that $LS(\K^{eq})=\o$ and that $\K^{eq}$ does not contain finite models.
It is also easily seen that if the axioms AI, AIII, and AVI hold for $\mathcal{K}$, they hold also for $\mathcal{K}^{eq}$. 
 
We now show that also AII holds.

\begin{lemma}
Suppose the axioms AI-AVI hold for $\K$.
Then, axiom AII holds for $\K^{eq}$.
\end{lemma}

\begin{proof}
Let $\A^{eq} \in \K^{eq}$ be countable, and let $a$ be arbitrary.
We need to construct an $s$-primary model over $\A^{eq}a$.
Let $b_1, \ldots, b_n \in \M$ be such that $a=(F_{E_1}(b_1), \ldots, F_{E_m}(b_n))$ for some $E_1, \ldots, E_m \in \mathcal{E}$, and denote $b=(b_1, \ldots, b_n)$.
 We will first show that we may choose $b$ so that there is some finite $A \subset \A$ such that for all $b'$, $t(b'/Aa)=t(b/Aa)$ implies $t(b'/\A a)=t(b/\A a)$.
 
We note first that for this it suffices to find some $b=(b_1, \ldots, b_m)$ so that $a=(F_{E_1}(b_1), \ldots, F_{E_m}(b_n))$ and a finite set $A$ such that $t(b/\A)=t(b'/\A)$ whenever $t(b/A)=t(b'/A)$ and $(F_{E_1}(b_1'), \ldots, F_{E_m}(b_n'))=a$.
Indeed, suppose we have found such a tuple $b$ and such a set $A$.
Let $b'$ be such that $t(b'/Aa)=t(b/Aa)$.
Then, $a=(F_{E_1}(b_1'), \ldots, F_{E_m}(b_n'))$, and thus  $t(b/\A)=t(b'/\A)$.
We claim that moreover, $t(b' /\A a)=t(b / \A a)$.
If not, then there is some finite set $A' \subset \A$ such that $t(b'a/A') \neq t(ba/A')$.
Since $t(b'/A')=t(b/A')$, there is some $f \in \textrm{Aut}(\M/A')$ such that $f(b)=b'$.
But $f$ extends to an automorphism $f' \in \textrm{Aut}(\M^{eq}/A)$, and 
$$f'(a)=(F_{E_1}(f(b_1)), \ldots, F_{E_m}(f(b_m))=(F_{E_1}(b_1'), \ldots, F_{E_m}(b_n')=a,$$
where for each $i$, $f(b_i)$ denotes the relevant projection of the tuple $f(b)$.
Thus, $t(b'a/A')=t(ba/A')$, a contradiction.
 
To simplify notation, denote now by $F$ the function from $\M$ to $\M^{eq}$ that is given by $(F_{E_1}, \ldots, F_{E_m})$.
Let  $A_0 =\emptyset$ and let $b_0$ be such that $F(b_0)=a$.
If $b_0$ and $A_0$ are not as wanted, then there is some finite $A_1 \subset \A$ and some $b_1$ so that $F(b_1)=a$, $t(b_0/A_0)=t(b_1/A_0)$ and $t(b_1/A_1) \neq t(b_0/A_1)$.
Now we check if $A_1$ and $b_1$ are as wanted.
By AIV, we cannot continue this process infinitely, so at some step we have found $b=b_n$ and $A=A_n$ as wanted.  
  
Let now $\B=\A[b]=\A b \cup \bigcup_{i<\o} b_i \le \M$.
We claim that $\B^{eq}$ is $s$-primary over $\A^{eq}a$.
For this, we need to enumerate the elements of $\B^{eq}$ so that we may write $\B^{eq}=\A^{eq} a \cup \bigcup_{i<\o} c_i$.
Let $b=(b_0, \ldots, b_k)$, where each $b_i$ is a singleton.
For $0 \le i \le k$, we denote $c_i=b_i$.
By the above argument, the required isolation property is satisfied.
After this, we list the elements so that whenever $i<j$,  we have $b_i=c_{i'}$ and $b_j=c_{j'}$ for some $i'<j'$. 
Moreover, we take care that for each singleton $c \in \B^{eq}\setminus (\B \cup \A^{eq})$, the elements of some tuple $d \in \B$ such that $c=F_E(d)$ for some $E \in \mathcal{E}$ are listed before $c$ (i.e. if $d=(d_0, \ldots, d_k)$, then, $d_0=c_{i_0}, \ldots, d_k=c_{i_k}$ and $c=c_j$ for some $i_0< \ldots <i_k<j$).
Then, the required isolation properties are satisfied and we see that $\B^{eq}$ is indeed as wanted.
\end{proof}
  
\begin{lemma}
Suppose the axioms AI-AVI hold for $\K$.
Then, axiom AIV holds for $\K^{eq}$.
\end{lemma}
  
\begin{proof}
Let $a \in \M^{eq}$.
Again, there is some tuple $b \in \M$ so that $a=F(b)$ for some definable function $F$.
Then, there is a number $n<\o$ so that for any countable $\A \in \K$ and finite $A' \subset \A$, player II wins $GI(b, A', \A)$ in $n$ moves.
Let now $\A^{eq} \in \K^{eq}$ be countable and $A \subset \A^{eq}$ finite.
We claim that player  II will win $GI(a, A, \A^{eq})$ in $n$ moves.
Let $A' \subset \A$ be such that every element $x \in A$ can be written as $x=F(y)$ for some $y \in A'$ where $F$ is a definable function.
Now, player II wins $GI(b, A', \A)$ in $n$ moves.
If there are some tuples $a', a'' \in \M^{eq}$ and some finite sets $C \subset B \subset \A^{eq}$ such that $t(a'/C)=t(a''/C)$ but $t(a'/B) \neq t(a''/B)$, then there are tuples $b', b'' \in \M$ and  a definable function $F$ so that $a' =F(b')$, $a'' =F(b'')$, and some $B' \subset \A$ and a definable function $H$ so that $B \subseteq H(B')$ and $t(b'/B') \neq t(b''/B')$.
Thus, the claim follows.  
\end{proof}  

\begin{lemma}
Suppose the axioms AI-AVI hold for $\K$.
Then, axiom AV holds for $\K^{eq}$.
\end{lemma}

\begin{proof}
Suppose $\A^{eq}, \B^{eq} \in \K^{eq}$ are countable models, $\A^{eq} \subset \B^{eq}$, $a \in \M^{eq}$ and $\B^{eq} \da_{\A^{eq}}^{ns} a$.
We need to prove $a  \da_{\A^{eq}}^{ns} \B^{eq}$.
Suppose this does not hold.
Then, there is some $b \in \B$ such that
$$a  \nda_{\A^{eq}}^{ns} b.$$
Let $a' \in \M$ be such that $a =F(a')$ for some definable function $F$.
Choose some $b' \in \M$ such that $t(b'/\A^{eq})=t(b/\A^{eq})$ and $b' \da_{\A^{eq}}^{ns} aa'$
(note that we may apply Lemma \ref{vaplaajol} since we needed only axiom AI to prove it).
Since $\B^{eq} \da_{\A^{eq}}^{ns} a$, we have $b \da_{\A^{eq}}^{ns} a$, and thus 
$t(b'/\A^{eq}a)=t(b/\A^{eq}a)$ (note that also Lemma \ref{vaplaajykskas} requires only axiom AI).
Thus, to obtain a contradiction, it suffices to show that $a \da_{\A^{eq}}^{ns} b'$.
But $b' \da_{\A^{eq}}^{ns} aa'$ implies that $b' \da_\A^{ns} a'$ (in $\M$), and thus by Lemma \ref{symlemma}, 
$a' \da_\A^{ns} b'$.
It follows that $a' \da_{\A^{eq}}^{ns} b'$.

\end{proof}  

Let $\M' \in \K$ be a $\vert \M' \vert$ -model homogeneous and universal structure such that $\M \preccurlyeq \M'$ and $\vert \M' \vert > \vert \M \vert$.
 We call  $\M'$ the \emph{supermonster}. 
Then, every $f \in \textrm{Aut}(\M)$ extends to some $f' \in \textrm{Aut}(\M')$.
In the following, we will abuse notation and write just $f$ for both maps. 

By a \emph{global type} $p$, we mean a maximal collection $\{p_A \, | \, A \subset \M \textrm{ finite }\}$ such that $p_A$ is a Galois type over $A$, and whenever $A \subseteq B$ and $b \in \M$ realizes $p_B$, then $b$ realizes also $p_A$.
We denote the collection of global types by $S(\M)$.
Moreover, we require that global types are consistent, i.e. that for each $p \in S(\M)$, there is some $b \in \M'$ such that $b$ realizes $p_A$ for each finite set $A \subset \M$.
  
Let $f \in \textrm{Aut}(\M^{eq})$, $p \in S(\M)$.
We say that $f(p)=p$ if for all finite $A,B \subset \M$ such that $f(B)=A$ and all $b$ realizing $p_B$, it holds that $t(b/A)=p_A$.

\begin{definition}
Let $p \in S(\M)$.
We say that $\alpha \in \M^{eq}$ is a \emph{canonical base} for $p$ if it holds for every $f \in \textrm{Aut}(\M^{eq})$ that $f(p)=p$ if and only if $f(\alpha)=\alpha$.
\end{definition}

\begin{lemma}\label{cbexists}
Suppose $p \in S(\M)$. 
Then, there is some  $\alpha \in \M^{eq}$ so that $\alpha$ is a canonical base for $p$.
\end{lemma}
 
\begin{proof}
Let $p \in S(\M)$ be a global type that does not split over $a \in \M$.
Suppose $b \in \M'$ realizes $p$.
Consider an arbitrary $c \in \M$ and let $q=t(b,c/\emptyset)$.
Since $t(b/\M)$ does not split over $a$, there are types $q_i$, $i<\o$, over $\emptyset$, so that for all $d \in \M$ the following holds: $bd$ realizes $q$ if and only if $ad$ realizes $q_i$ for some $i<\o$.
Indeed, there are only countably many types over the empty set, and from the non-splitting it follows that if $t(d_1 a/ \emptyset)=t(d_2 a/\emptyset)$, then $t(bd_1a/\emptyset)=t(bd_2a/\emptyset)$.
Thus, we may choose the types  $q_i$ as wanted.

For $c \in \M$, denote $q_c=t(b,c/\emptyset)$, and let $q_i^c$, $i<\o$ be such that for all $d \in \M$, $bd$ realizes $q_c$ if and only if $ad$ realizes $q_i^c$ for some $i<\o$.
Define the equivalence relation $E$ as follows: $(a_0, a_1) \in E$ if  the following holds for all $c,d \in \M$:  $a_0 d$ realizes $\bigvee_{i<\o}q_i^c$ if and only if $a_1 d$ realizes $\bigvee_{i<\o}q_i^c$.

We claim that $a/E$ is a canonical base for $p$.
Suppose first that $f \in \textrm{Aut}(\M)$ is such that $f(p)=p$.
We will show that $(a,f(a)) \in E$, which implies that $f(a/E)=f(a)/E=a/E$.
 As $f$ is an automorphism, we have $t(a/\emptyset)=t(f(a)/\emptyset)$.
Let now $c \in \M$ be arbitrary. 
Since $f(p)=p$, we have $$q_c=t(bc/\emptyset)=t(bf(c)/\emptyset).$$ 
for every $c \in \M$
Thus, we may choose $q_i^c=q_i^{f(c)}$ for all $i<\o$, and moreover we have
\begin{eqnarray*}
f(a)f(d) \textrm{ realizes } \bigvee_{i<\o}q_i^c \iff ad \textrm{ realizes } \bigvee_{i<\o}q_i^c 
\iff af(d) \textrm{ realizes } \bigvee_{i<\o}q_i^c, 
\end{eqnarray*}
where the first equivalence follows from $f$ being an automorphism, and the second one from the fact that $t(bd/\emptyset)=t(bf(d)/\emptyset)$.
Since this holds for every $d$ and automorphisms are surjective, we may, for arbitrary $d' \in \M$, choose $d \in \M$ so that $d'=f(d)$ to obtain
\begin{eqnarray*}
f(a)d' \textrm{ realizes } \bigvee_{i<\o}q_i^{c} \iff  ad' \textrm{ realizes } \bigvee_{i<\o}q_i^{c}. 
\end{eqnarray*}

Suppose now $f(a/E)=a/E$.
Then, $(a,f(a)) \in E$.
We will show that $t(bc/\emptyset)=t(b f(c)/ \emptyset)$ for all $c \in \M$.
Then, clearly $f(p)=p$.
Let $c$ be arbitrary.
Denote $q^c=t(bc/\emptyset)$.
Then, $ac$, and thus also $f(a)f(c)$, realizes $\bigvee_{i<\o}q_i^c$.
But now, since $(a,f(a)) \in E$, we have that also $af(c)$ realizes $\bigvee_{i<\o}q_i^c$.
From this it follows that $bf(c)$ realizes $q^c$, i.e.
$$t(bc/\emptyset)=q_c=t(b f(c)/ \emptyset).$$

So, we have shown that $a/E$ is a canonical base for $p$, as wanted.
\end{proof} 

We note that from the proof of Lemma \ref{cbexists} it follows that there are only countably many equivalence relations needed to get the canonical bases of all global types.
Indeed, for each global type $p$ realized by an element $b \in \M'$, there is some tuple $a \in \M$ such that $t(b/\M)$ does not split over $a$.
Now the tuple $ba$ determines the equivalence relation $E$ from the proof of the lemma, and we claim that $E$ depends only on $t(ba/\emptyset)$.
Indeed, let $b' \in \M'$, $a' \in \M$ be such that $b' \da_{a'} \M$ and $t(b'a'/ \emptyset)=t(ba/\emptyset)$.
Then, there is some $f \in \textrm{Aut}(\M/\emptyset)$ such that  $f(a)=a'$, and some $F \in \textrm{Aut}(\M'/\emptyset)$ extending $f$ (and in particular, $F(m) \in \M$ for every $m \in \M$).
Let $b''=F(b)$.
Then, since $b \da_a^{ns} \M$, we have $b'' \da_{a'}^{ns} \M$.
On the other hand, we have 
$$t(b'a'/\emptyset)=t(ba/\emptyset)=t(b''a'/\emptyset),$$
so $t(b'/a')=t(b''/a')$ and by Lemma \ref{statilemma}, $t(b'/\M)=t(b/\M)$.
Thus, there is an automorphism $G$ of $\M'$ such that $G(ab)=a'b'$ and $G(\M)=\M$ (we first take $ab \mapsto a'b''$ and then fix $\M$ and take $b'' \mapsto b'$).
Then, the global type $p'=t(G(b)/\M)$ can be determined from $f(a)=a'$ in the same way as $p$ was determined from $a$, and the definition of $E$ stays the same.  
 
\begin{definition}
Let $a \in \M$ and let $A \subset \M$.
By Theorem \ref{main}, (iv), there is some $b \in \M'$ such that $Lt(a/A)=Lt(b/A)$ and $b \da_A \M$.
Let $p=t(b/\M)$.
By a \emph{canonical base} for $a$ over $A$, we mean a canonical base of $p$.
We write $\alpha=Cb(a/A)$ to denote that $\alpha$ is a canonical base of $a$ over $A$.
\end{definition}

Next, we prove some important properties of canonical bases.

\begin{lemma}\label{cbbcl}
Let $a \in \M$ and let $A \subset \M$ be a finite set.
Then, $Cb(a/A) \in \textrm{bcl}(A)$.
\end{lemma}

\begin{proof}
Let $b \in \M'$ be such that $Lt(a/A)=Lt(b/A)$ and $b \da_A \M$,
and let $p=t(b/\M)$.
Then, there is some finite $A_0 \subseteq A$ such that $b \da_{A_0} \M$.
We may without loss assume $A=A_0$. 
Denote $\alpha=Cb(a/A)$, and suppose $\alpha \notin \textrm{bcl}(A)$.
Then, $t(\alpha/A)$ has unboundedly many realizations.
By the proof of Lemma \ref{cbexists}, each one of them defines a global type, and by the definition of a canonical base, the global types defined by these unboundedly many elements are pairwise distinct. 
Let $f \in \textrm{Aut}(\M/A)$ and let $\alpha'=f(\alpha)$. 
Then $f$ extends to an automorphism $g$ of $\M'$, and we have $g(b)=b'$ for some $b' \in \M' \setminus \M$.
Since $g(\M)=\M$, we have $b' \da_A \M$.
Let $\A \subset \M$ be a countable model such that $A \subset \A$.
Then, by (v) in Theorem \ref{main}, we have $b \da_\A \M$ and $b' \da_\A \M$.
Since $t(b/\M) \neq t(b'/\M)$, by Lemmas \ref{kolme} and \ref{vaplaajykskas}, we must have 
$t(b/\A)Ê\neq t(b'/\A)$.
This means that we have uncountably many distinct types over the countable model $\A$, a contradiction against
Lemma \ref{notypes}.
\end{proof}

\begin{remark}
Let $a \in \M$, and let $A$ and $B$ be sets such that $A \subsetneq B$, and let $\alpha \in \M^{eq}$
If $a \da_A B$, then $\alpha=Cb(a/A)$ if and only if $\alpha=Cb(a/B)$.
\end{remark}
 
\begin{lemma}
Let  $a \in \M$ and let $\alpha=Cb(a/A)$.
Then, $a \da_\alpha A$.
\end{lemma}

\begin{proof}
Let $b \in \M'$ be such that $Lt(a/A)=Lt(b/B)$ and $b \da_A \M$, and let $p=t(b/\M)$.
Then, $\alpha$ is a canonical base of $p$.

We note first that $b \da^{ns}_\alpha \M$.
Indeed, if there were some $b,c \in \M$ that would witness the splitting, then there would be some automorphism $f$ fixing $\alpha$ such that $f(b)=c$.
Let $b'=f(b)$.
Now, we have
$$t(bd/\alpha) \neq t(bc/\alpha)=t(b'd/\alpha),$$
so $t(b/d\alpha) \neq t(b'/d\alpha)$, which is a contradiction since $f$ fixes the type $p$ (since it fixes $\alpha$).

In particular, this implies $b \da_\alpha A$.
Since $\alpha \in \textrm{bcl}(A)$, we have $b \da_A \alpha$ and $a \da_A \alpha$, so $t(a/A\alpha)=t(b/A \alpha)$, and thus $a \da_\alpha A$.

\end{proof}

\begin{definition}
Let $p=t(d/B)$ for some $d \in \M$ and $B \subset \M$.
Suppose $B \subset C$.
We say that $p'=t(d'/C)$ is a \emph{free extension} of $p$ into $C$ if $t(d/B)=t(d'/B)$ and $d' \da_B C$.
We call a type \emph{stationary} if it has a unique free extension to any set.
If $p$ is a stationary type, we will denote the free extension of $p$ into $C$ by $p \vert_C$. 
\end{definition}

\begin{lemma}
Let $\alpha=Cb(a/A)$.
Then, $t(a/\alpha)$ is stationary.
\end{lemma}

\begin{proof}
Let $b \in \M'$ be such that $Lt(a/A)=Lt(b/A)$ and $b \da_A \M$, and let $p=t(b/\M)$.
Then, $\alpha$ is a canonical base of $p$.
As in the proof of the previous Lemma, $b \da^{ns}_\alpha \M$, so $b \da^{ns}_\alpha \M^{eq}$ and thus $b \da_\alpha \M^{eq}$.
Also, we have (as seen in the proof of the previous lemma) $t(b/\alpha)=t(a/\alpha)$.
Suppose now $\alpha \in B \subset \M^{eq}$ and there is some $c \in \M$ such that $t(c/\alpha)=t(a/\alpha)=t(b/\alpha)$, $c \da_\alpha B$ but $t(c/B) \neq t(b/B)$.
Let $b' \in \M'$ be such that $Lt(b'/\alpha)=Lt(c/\alpha)$ and $b' \da_\alpha \M$.
Then, $b' \da^{ns} \M^{eq}$, so
by Lemma \ref{statilemma}, $t(b'/\M)=t(b/\M)$, which is a contradiction, since 
$$t(b/B)\neq t(c/B)=t(b'/B).$$
\end{proof}
 
\section{The axioms in quasiminimal classes}

We have seen that the model class $\K$ of Example \ref{equivalence} satisfies the axioms AI-AVI.
In this section, we will show that something more general is true: If
 $\K$ is a quasiminimal class in the sense of \cite{monet} (which is the same as in \cite{kir}, but without finite-dimensional structures), then $\K$ satisfies the axioms, given that $\K$ only contains infinite-dimensional models.
We first present some definitions.

\begin{definition}
Let $X$ be a set and let $\textrm{cl}: \mathcal{P}(X) \to \mathcal{P}(X)$ be an operator on the power set of $X$.
We say $(X, \textrm{cl})$ is a \emph{pregeometry} if the following hold:
\begin{enumerate}[(i)]
\item If $A \subseteq X$, then $A \subseteq \textrm{cl}(A)$ and $\textrm{cl}(\textrm{cl}(A))=\textrm{cl}(A)$;
\item If $A\subseteq B \subseteq X$, then $\textrm{cl}(A) \subseteq \textrm{cl}(B)$;
\item (exchange) If $A \subseteq X$, $a,b \in X$, and $a \in \textrm{cl}(A \cup \{b\})$, then $a \in \textrm{cl}(A)$ or $b \in \textrm{cl}(A \cup \{a\})$;
\item (local character) If $A \subset X$ and $a \in \textrm{cl}(A)$, then there is a finite $A_0 \subseteq A$ such that $a \in \textrm{cl}(A_0)$.
\end{enumerate}

We call $\textrm{cl}$ the \emph{closure operator} on $X$ and we say a set $A \subseteq X$ is \emph{closed} if $\textrm{cl}(A)=A$.
\end{definition}
 
\begin{definition}
If $(X, \textrm{cl})$ is a pregeometry, we say that $A$ is \emph{independent} if $a \notin \textrm{cl}(A\setminus \{a\})$ for any $a \in A$.
We say $B$ is a \emph{basis} for $Y$ if $B \subseteq Y$ is independent and $Y \subseteq \textrm{cl}(Y)$.
\end{definition} 
 
If $(X,cl)$ is a pregeometry, $Y \subseteq X$, $B_1,B_2 \subseteq Y$ and both $B_1$ and $B_2$ are bases for $Y$, then $\vert B_1 \vert= \vert B_2 \vert$ (see e.g. \cite{MarIntro}).

\begin{definition}
Let $B$ be a basis for $Y$.
We say that $\vert B \vert$ is the \emph{dimension} of $Y$ with respect to cl, and write $\textrm{dim}_{\textrm{cl}}(Y)=\vert B \vert$. 
\end{definition}

\begin{definition}
We say that $a$ and $b$ have the same \emph{quantifier-free type} over the set $A$ if they satisfy exactly the same quantifier-free first-order formulae with parameters from $A$.
This is denoted $\textrm{tp}(a/A)=\textrm{tp}(b/A)$.
We write just $\textrm{tp}(a)$ for $\textrm{tp}(a/\emptyset)$.
For two sets $A$ and $B$ of the same cardinality, we say $\textrm{tp}(A)=\textrm{tp}(B)$ if the elements satisfy exactly the same quantifier-free first order formulae given a suitable enumeration.
\end{definition}
 
In \cite{monet}, a quasiminimal pregeometry structure and a quasiminimal class are defined as follows.  

\begin{definition}\label{quasiminclass}
Let $M$ be an $L$-structure for a countable language $L$, equipped with a pregeometry $cl$ (or $\textrm{cl}_M$ if it is necessary to specify $M$).
We say that $M$ is a \emph{quasiminimal pregeometry structure} if the following hold:
\begin{enumerate}
\item (QM1) The pregeometry  is determined by the language.
That is, if $a$ and $a'$ are singletons and
$\textrm{tp}(a, b)=\textrm{tp}(a', b')$, then $a \in \textrm{cl}(b)$ if and only if $a' \in \textrm{cl}(b')$.
\item (QM2) $M$ is infinite-dimensional with respect to cl.
\item (QM3) (Countable closure property) If $A \subseteq M$ is finite, then $\textrm{cl}(A)$ is countable.
\item (QM4) (Uniqueness of the generic type) Suppose that $H,H' \subseteq M$ are countable closed subsets, enumerated so that $\textrm{tp}(H)=\textrm{tp}(H')$.
If $a \in M \setminus H$ and $a' \in M \setminus H'$ are singletons, then $\textrm{tp}(H,a)=\textrm{tp}(H',a')$ (with respect to the same enumerations for $H$ and $H'$).
\item (QM5) ($\aleph_0$-homogeneity over closed sets and the empty set)
Let $H,H' \subseteq M$ be countable closed subsets or empty, enumerated so that $\textrm{tp}(H)=\textrm{tp}(H')$,
and let $b, b'$ be finite tuples from $M$ such that $\textrm{tp}(H, b)=\textrm{tp}(H',b')$, and let $a$ be a singleton such that $a \in \textrm{cl}(H, b)$.
Then there is some singleton $a' \in M$ such that $\textrm{tp}(H, b,a)=\textrm{tp}(H',b',a').$
\end{enumerate}

We say $M$ is a \emph{weakly quasiminimal pregeometry structure} if it satisfies all the above axioms except possibly QM2.
\end{definition}

It is easy to see that the class $\K$ from Example \ref{equivalence} satisfies the axioms.

\begin{definition}
Suppose $M_1$ and $M_2$ are weakly quasiminimal pregeometry  $L$-structures.
Let $\theta$ be an isomorphism from $M_1$ to  some substructure of $M_2$.
We say that  $\theta$ is a \emph{closed embedding} if $\theta(M_1)$ is closed in $M_2$ with respect to $\textrm{cl}_{M_2}$, and $\textrm{cl}_{M_1}$ is the restriction of $\textrm{cl}_{M_2}$ to $M_1$.
\end{definition}

Given a quasiminimal pregeometry structure $M$, let $\K^-(M)$ be the smallest class of $L$-structures which contains $M$ and all its closed substructures and is closed under isomorphisms, and let $\K(M)$ be the smallest class containing $\K^-(M)$ which is also closed under taking unions of chains of closed embeddings.

From now on, we suppose that $\K=\K(\M)$ for some  quasiminimal pregeometry structure $\M$, and that we have discarded all the finite-dimensional structures from $\K$.
For $\A, \B \in \K$, we define $\A \preccurlyeq \B$ if $\A$ is a closed submodel of $\B$. 
It is well known that $(\K, \preccurlyeq)$ is an AEC with $LS(\K)=\o$.
We may without loss assume that $\M$ is a monster model for $\K$.
In \cite{monet}, it is shown that $\K$ is totally categorical and has arbitrarily large models (Theorem 2.2).
We will show that $\K$ has AP and JEP and satisfies the axioms AI-AVI.
  
\begin{lemma}
$\K$ has AP and  JEP. 
\end{lemma}

\begin{proof}
Let $\A, \B \in \K$, let $\A' \preccurlyeq \A$, and let $f:\A' \to \B$ be an elementary embedding.
We may without loss assume that $\textrm{dim}(\B) \ge \textrm{dim}(\A)$.
 Let $B$ be a basis for $\A'=\textrm{dom}(f)$.
Then, $B'=f(B)$ is a basis for $\textrm{ran}(f)$.
We may extend $B$ to $C$, a basis for $\A$,  and $B'$ to $C'$, 
a basis for $\B$.
Let $\psi: C \setminus B \to C' \setminus B'$ be an injection.  
By Theorem 3.3. in \cite{kir}, $f \cup \psi$ extends to an embedding of $\A$ into $\B$.
Thus, $\K$ has AP.
  
For JEP, let $\A, \B \in \K$, and let $B$ be a basis for $\A$.
Again, we may without loss assume $\textrm{dim}(\B)\ge \textrm{dim}(\A)$.
Let $B'$ be a basis for $\B$.
By Theorem 3.3. in \cite{kir}, $\A$ embeds into $\B$.
\end{proof}

We now note that we may reformulate the conditions QM4 and QM5 so that the concept of Galois type is used instead of the concept of quantifier-free type. This will be useful in the arguments we later present.
Indeed, for QM4, let $H, H' \subset \M$ be countable and closed, let $t^g(H)=t^g(H')$, and let $a, a'$ be singletons such that $a \notin \textrm{cl}(H)$ and $a' \notin \textrm{cl}(H')$.
As $H$ and $H'$ are closed, they are models.
Since $H$ and $H'$ are countable, there is some isomorphism $f: H \to H'$.
Using QM4, we may extend $f$ to a map $g_0: H a \to H' a'$ that preserves quantifier-free formulae.
Let $\A=\textrm{cl}(Ha)$ and $\B=\textrm{cl}(H'a')$.
We will extend $g_0$ to an isomorphism $g: \A \to \B$.
Indeed, if $b \in \A=\textrm{cl}(Ha)$, then by QM5 and QM1, there is some $b' \in \B= \textrm{cl}(H'a')$ such that 
$\textrm{tp}(H,a,b)=\textrm{tp}(H,a',b')$, so $f_0$ extends to a map $f_1: H,a,b \to H',a',b'$ preserving quantifier-free formulae.
Since both $\A$ and $\B$ are countable, we can do a back-and-forth construction to obtain an isomorphism $g:  \A \to  \B$.
Then, $g(H,a)=(H',a')$ and 
 $g$ extends to an automorphism of $\M$, so 
 $t^g(H,a)=t^g(H',a')$, as wanted.
 
For QM5,  suppose $H, H' \subset \M$ are either countable and closed or empty, let $t^g(H)=t^g(H')$, and let $b,b' \in \M$ be such that $t^g(H,b)=t^g(H',b')$ and let $a \in \textrm{cl}(H,b)$.
Again, there is a map $f$ such that $f(H)=H'$, $f(b)=b'$ and $f$ preserves quantifier-free formulae.
As in the case of QM4, we may extend $f$ to an isomorphism $g: \textrm{cl}(Hb) \to \textrm{cl}(H'b')$.
If $a \in \textrm{cl}(Hb)$, then $t^g(H, b, a)=t^g(H',b',g(a))$.
 
\begin{lemma}
$\K$ satisfies the axioms AI-AVI.
\end{lemma}

\begin{proof}
For AI, suppose $a \in \M$, $\A \in \K$ and $A$ is a finite set such that $A \subset \A$.
Since $\A$ is closed and infinite-dimensional, we have $\textrm{cl}(A) \subsetneq \A$.
If $a \notin \A$, then $a \notin \textrm{cl}(A)$.
Let $b \in \A \setminus \textrm{cl}(A)$.
By QM4, $t^g(\textrm{cl}(A),b)=t^g(\textrm{cl}(A),a)$, and hence there is some $f \in \textrm{Aut}(\M/A)$ such that $f(a)=b$. 

AII follows directly from Proposition 5.2 in \cite{monet}.

For AIII, suppose $\A \da_\B^{ns} \C$, where $\B=\A \cap \C$.
Then, $\D=\textrm{cl}(\A \cup \C)$ is the desired $s$-prime model over $\A \cup \C$.
Indeed, let $\mathcal{E}$ be a model and let $g: \A \cap \C \to \mathcal{E}$ be weakly elementary.
To extend $g$ to an elementary map $f: \D \to \mathcal{E}$, we do the same construction as in
the proof of Theorem 3.3 in \cite{kir}, with $\A$ in place of $G$ and taking care that for every finite $X \subseteq B$, where $B$ is a pregeometry basis for $\D$ over $\A$, we have $f_X \raj \A \cup \C=g$.

It is also easy to see that there are no other $s$-prime models over $\A \cup \C$.
Indeed, any model containing $\A \cup \C$ must contain its closure, $\D$.
Suppose $\D \subsetneq \mathcal{E}$ for some model $\mathcal{E}$.
The identity map $\textrm{id}: \A \cup \C \to \D$ is clearly weakly elementary, but it does not extend to an elementary map from $\mathcal{E} \to \D$.
Thus, $\mathcal{E}$ is not $s$-prime over $\A \cup \C$.
 
Suppose now $\C'$ is a model such that $\C \subseteq \C'$ and $\A \da_\B^{ns} \C'$.
We need to prove that $\D \da_\C^{ns} \C'$.
Let $d \in \D$ be arbitrary, and let $a \in \A$, $c \in \C$ be such that $d \in \textrm{cl}(ac)$.
Let $B \subset \B$ be a finite set such that $a \da_B^{ns} \C'$.
Then, we have also $a \da_{Bc}^{ns} \C'$.
We claim that $d \da_{Bc}^{ns} \C'$ and thus $d \da_\C^{ns} \C'$.
By the proof of Proposition 5.2 in \cite{monet}, $t(d/Bac)$ determines $t(d/\C'a)$.
Suppose $d \nda_{Bc}^{ns} \C'$.
Then, there are some $e,f \in \C'$ such that $t(e/Bc)=t(f/Bc)$ but $t(de/Bc) \neq t(df/Bc)$.
We note that we have $t(e/Bac)=t(f/Bac)$, since otherwise $e$ and $f$ would witness the splitting of $t(a/\C')$ over $Bc$.
Let $\sigma$ be an automorphism fixing $Bac$ such that $\sigma(f)=e$.
Then, $t(\sigma(d)/Bce) \neq t(d/Bce)$ although $t(\sigma(d)/Bac)=t(d/Bac)$, a contradiction.
 
 For AIV, on move $n$ in $GI(a, A, \A)$,
 let player II choose the set $A_{n+1}$ so that $t(a_{n+1}/\A)$ does not split over $A_{n+1}$ (such a set exists by Proposition 4.2 in \cite{monet}).
Then, player $I$ plays some tuple $a_{n+2}$ and some set $A_{n+2}'=A_{n+1}b$ so that $t(a_{n+2}/A_{n+1})=t(a_{n+1}/A_{n+1})$ but  $t(a_{n+2}/A_{n+2}')\neq t(a_{n+1}/A_{n+2}')$.
We claim that $\textrm{dim}_{cl}(a_{n+2}/\A)<\textrm{dim}_{cl}(a_{n+1}/\A)$, which means that player I can only move $\textrm{dim}_{cl}(a/\A)$ many times.

Let $m=\textrm{dim}_{cl}(a_{n+1}/\A)$.
By Lemma 4.3 in \cite{monet}, $\textrm{dim}_{cl}(a_{n+1}/A_{n+1})=m$, and thus  $\textrm{dim}_{cl}(a_{n+2}/A_{n+1})=m$, so $\textrm{dim}_{cl}(a_{n+2}/\A) \le m$.
Suppose $\textrm{dim}_{cl}(a_{n+2}/\A) = m$.
Then, in particular, $\textrm{dim}_{cl}(a_{n+2}/A_{n+1}b)=\textrm{dim}_{cl}(a_{n+1}/A_{n+1}b)=m$.
For $i=1,2$, write $a_{n+i}=a_{n+i}' a_{n+i}''$, where $a_{n+i}'$ is an $m$-tuple free over $A_{n+1}b$ and $a_{n+i}'' \in \textrm{cl}(A_{n+1} b a_{n+i}')$.
By QM4, there is some automorphism $\sigma$ fixing $A_{n+1} b$ pointwise so that $\sigma(a_{n+2}')=a_{n+1}'$.
Since $t(a_{n+2}/A_{n+1})=t(a_{n+1}/A_{n+1})$, we have
$$t(\sigma(a_{n+2}'')/A_{n+1} a_{n+1}')=t(a_{n+1}''/A_{n+1} a_{n+1}').$$
On the other hand, $t(a_{n+2}/A_{n+2}')\neq t(a_{n+1}/A_{n+2}')$, so
$$t(a_{n+1}''/A a_{n+1}' b)\neq t(\sigma(a_{n+2}'')/A a_{n+1}' b),$$
which contradicts the fact that $t(a_{n+1}''/A_{n+1} a_{n+1}')$ determines $t(a_{n+1}''/\A a_{n+1}')$ by the proof of Proposition 5.2 in \cite{monet}.
 
Before proving AV and AVI, we first show that if $\A \subseteq \B$ for some model $\B$, then $a \da^{ns}_\A \B$ if and only if $\textrm{dim}_{cl}(a/\B)=\textrm{dim}_{cl}(a/\A)$.

Suppose first $a \da^{ns}_\A \B$.
Let $A \subset \A$ be such that $t(a/\B)$ does not split over $A$.
By Lemma 4.3 in \cite{monet}, $\textrm{dim}(a/A)=\textrm{dim}(a/\B)$, so in particular $\textrm{dim}(a/\A)=\textrm{dim}(a/\B)$.

Suppose now $\textrm{dim}_{cl}(a/\A)=\textrm{dim}_{cl}(a/\B)$.
Write $a=a'a''$, where $a'$ is a tuple independent over $\A$ and $a'' \in \textrm{cl}(\A a').$
By the proof of Proposition 5.2. in \cite{monet} there is a finite set $A \subset \A$ such that $t(a''/Aa')$ determines $t(a''/\A a')$.  
We now show that $t(a''/Aa')$ also determines $t(a''/\B a')$.
Suppose not.
Then, there is some $c''$ and some $b \in \B$ so that $t(a''/Aa')=t(c''/Aa')$ but $t(a''/Aa'b) \neq t(c''/Aa'b)$.
Let $b' \in \A$ be such that $t(b'/A)=t(b/A)$, and let $\sigma$ be an automorphism fixing $A$ such that $\sigma(b)=b'$.
Since $a'$ is independent over $\A$, QM4 implies that $t(\sigma(a')/Ab')=t(a'/Ab')$, and thus we may without loss assume that $\sigma(a')=a'$.
Then,
$t(\sigma(a'')/Aa')=t(\sigma(c'')/Aa')=t(a''/Aa')$,
but $t(\sigma(a'')/Aa'b') \neq t(\sigma(c'')/Aa'b')$, a contradiction.
 
Next, we show that $a' \da^{ns}_A \B$.
Suppose not.
Then, there are $b, c \in \B$ such that $t(b/A)=t(c/A)$ but $t(a'b/A) \neq t(a'c/A)$.
Let $\sigma$ be an automorphism fixing $A$ and mapping $b$ to $c$.
Then, $t(\sigma(a')/Ac) \neq t(a'/Ac)$.
Since we assumed $\textrm{dim}_{cl}(a/\A)=\textrm{dim}_{cl}(a/\B)$, the tuple $a'$ is independent over $\B$ (in the pregeometry sense), so in particular it is independent over $Abc$.
Thus, $\sigma(a')$ is independent over $A \sigma(b)=Ac$.
But then, by QM4, $t(\sigma(a')/Ac) = t(a'/Ac)$, a contradiction.
 
Assume now $a \nda^{ns}_\A \B$, and let $b,c \in \B$ witness the splitting over $A$.
Let $\sigma$ be an automorphism fixing $A$ so that $\sigma(c)=b$.
Since $a' \da_A^{ns} \B$, we may assume $\sigma(a')=a'$.
Thus, $t(a''/Aa')=t(\sigma(a'')/Aa')$ but $t(a''/Aa'b) \neq t(\sigma(a'')/Aa'b)$, which contradicts the fact that $t(a''/\B a')$ is determined by $t(a''/Aa')$.
Thus, we have seen that if $\A \subseteq \B$, then $a \da^{ns}_\A \B$ if and only if $\textrm{dim}_{cl}(a/\B)=\textrm{dim}_{cl}(a/\A)$.

We now prove AV.
 Suppose $\B \da_\A^{ns} a$
but $a \nda_\A^{ns} \B$.
Then, $\textrm{dim}(a/\B)<\textrm{dim}(a/\A)$, so there is some $b \in \B$ such that $\textrm{dim}(a/\A b)< \textrm{dim} (a /\A)$.
But now, applying exchange,
we see that $\textrm{dim}_{cl}(b/\A a) < \textrm{dim}_{cl}(b/\A)$. 
Choose a model $\A'$ such that $\A a \subseteq \A'$ and $b \da_\A^{ns} \A'$.
Now, $\textrm{dim}_{cl}(b/\A)=\textrm{dim}_{cl}(b/\A')$, a contradiction.

For AVI,  
let $B_\A$ be a pregeometry basis for $\A$, $B_\B$ a basis for $\B$ over $\A$ and $B_\D$ a basis for $\D$ over $\A$.
Let $B_\B'$ be a set of the same cardinality as $B_\B$, independent over $B_\A \cup B_\D$.
Then, by Theorem 3.3 in \cite{kir}, there is a model $\C$ with basis $B_\A \cup B_\B'$ and $t(\C/\A)=t(\B/\A)$.
We are left to show that $\C \da_\A^{ns} \D$.
Let $c \in \C$ be arbitrary. 
Since the pregeometry is determined by the language, it follows from the theory of pregeometries that  $\textrm{dim}_{cl}(c/\A)=\textrm{dim}_{cl}(c/\D)$, and thus  $c \da_\A^{ns} \D$.
\end{proof}

\begin{remark}
Note that from QM1, QM3 and QM4 it follows that in a quasiminimal pregeometry structure, $\textrm{cl}(A)=\textrm{bcl}(A)$ for any set $A$.
Thus, from now on we will use bcl for cl.
\end{remark}

We next remark that in a quasiminimal pregeometry structure, $U(a/A)=\textrm{dim}_{\textrm{bcl}}(a/A)$.

\begin{lemma}
$U(a/A)=\textrm{dim}_{\textrm{bcl}}(a/A)$.
\end{lemma}
 
\begin{proof}
We prove the lemma by induction on $\textrm{dim}_{\textrm{bcl}}(a/A)$.
Suppose first $\textrm{dim}_{\textrm{bcl}}(a/A)=0$, i.e. $a \in \textrm{bcl}(A)$.
Then, by Theorem \ref{main}, (xi), $a \da_A \A$ for every model $\A \supset A$.
By (viii) of the same theorem, $U(a/A)=U(a/\A)$ for every such model.
Assume towards a contradiction that $U(a/A)=n>0$.
Then, by the definition of $U$-rank, there are some models $\B, \C$ such that $A \subseteq \B \subseteq \C$,
$U(a/A)=U(a/\B)$, $a \nda_\B^{ns} \C$ and $U(a/\C) \ge n-1$.
But $U(a/\C)=n$, which implies $U(a/\B)=n+1$, a contradiction.

Suppose next that $\textrm{dim}_\textrm{bcl}(a/A)=1$.
We will show that $U(a/A)=1$.
We may assume $a=(a_1, \ldots, a_m)$, where $a_1 \notin \textrm{bcl}(A)$ and $a_2, \ldots, a_m \in \textrm{bcl}(Aa_1)$.
By Theorem \ref{main} (x),
$$U(a/A)=U(a_1/A)+U(a_2, \ldots, a_m/Aa_1)=U(a_1/A),$$
so it suffices to show that $U(a_1/A)=1$.
Let $\A$ be a model such that $U(a_1/\A)=U(a_1/A)$, and $\B$ a model containing $\A a$.
Then, $U(a/\B)=0$ and $a \nda_\A^{ns} \B$, so $U(a_1/A) \ge 1$.
Suppose now $U(a_1/A)=n>1$.
Then, there are some models  $\B, \C$ such that $A \subseteq \B \subseteq \C$,
$U(a/A)=U(a/\B)$, $a \nda_\B^{ns} \C$ and $U(a/\C) \ge n-1>0$.
By Lemmas  \ref{vaplaajol2} and \ref{Unonsplit}, there is some $a'$ such that $t(a'/\B)=t(a/\B)$ and $a' \da^{ns}_\B \C$.
Since $a, a' \notin \C$, we have $t(a/\C)=t(a'/\C)$ by QM4, a contradiction.
 
Now suppose $n\ge 2$ and for all $m\le n$ it holds that if $a'$ and $A'$ are such that $\textrm{dim}_{\textrm{bcl}}(a'/A')=m$, then $U(a'/A')=m$.
Assume $\textrm{dim}_{bcl}(a/A)=n+1$.
We may without loss write $a=(a_1, \ldots, a_{n+1}, \ldots, a_r)$, where $\textrm{dim}_{\textrm{bcl}}(a_1, \ldots, a_{n+1})=n+1$ and $a_{n+2}, \ldots, a_r \in \textrm{bcl}(A, a_1, \ldots, a_{n+1})$.
By the inductive assumption, $U(a_1, \ldots, a_n/A)=n$ and $U(a_{n+1}, \ldots, a_r/A, a_1, \ldots, a_n)=1$, so by Theorem \ref{main} (x), we have $U(a/A)=n+1$, as wanted.
\end{proof}

\chapter{The Group Configuration}
  
In this chapter, we adapt E. Hrushovski's group configuration for the setting of quasiminimal classes.
We assume that $\K=\K(\M)$ for some quasiminimal pregeometry structure, as in the last section of Chapter 2.
We may without loss of generality assume that $\M$ is a monster model for the class $\K$.
The group configuration was originally presented for stable first-order theories in E. Hrushovski's Ph. D. Thesis.  
There, he proved that if a certain kind of configuration of tuples can be found in a model, then there is a group interpretable there.
The proof can be found in e.g. \cite{pi}.  
In the original context, properties of algebraic closures and the forking calculus for stable theories were used.
Doing a similar construction in our setting, we will make use of the independence calculus developed in Chapter 2.
We will prove that if a certain kind of configuration of tuples can be found in $\mathbb{M}$, then there is
a group interpretable in $(\mathbb{M}^{eq})^{eq}$ such that its generic elements have $U$-rank $1$.

We will be working in $\mathbb{M}^{eq}$ and occasionally in $(\mathbb{M}^{eq})^{eq}$.
To avoid confusion, we will write $\textrm{bcl}^{eq}(A)$ for the bounded closure of $A$ in $\mathbb{M}^{eq}$.
In this case, $A$ might contain some element $a \in \mathbb{M}^{eq} \setminus \mathbb{M}$.
  
We will say a set $A$ is independent over $B$ if $a \downarrow_B  (A \setminus \{a\})$ for each $a \in A$. 
 
\begin{definition}
We say $x$ and $y$ are \emph{interbounded} over a set $A$ if $x \in \textrm{bcl}(Ay)$ and $y \in \textrm{bcl}(Ax)$.
\end{definition}

We are now ready to present the configuration that will yield a group.
 
 \begin{pdfpic}
\psset{unit=9mm,dotscale=1}
\begin{pspicture}
\pspicture*(-2,-2)(10,10)
\psline(2,0)(2,6)
\psline(2,6)(8,6)
\psline(2,3)(8,6)
\psline(2,0)(5,6)
\psdots(2,0)(2,6)(8,6)(2,3)(5,6) 
\uput{3mm}[280](1.5,0.5){$c$}
\uput{3mm}[280](1.5,3.6){$b$}
\uput{3mm}[280](1.5,6.8){$a$}
\uput{3mm}[280](4.2,4.2){$z$}
\uput{3mm}[280](8.3,6.8){$y$}
\uput{3mm}[280](5,6.8){$x$}
\end{pspicture}
\end{pdfpic}

\begin{definition}\label{conf}
By a \emph{strict bounded partial quadrangle} over  a finite set $A$ we mean a $6$-tuple of elements $(a,b,c,x,y,z)$ in $\M^{eq}$, each of $U$-rank $1$ over $A$, such that 
\begin{enumerate}[(i)]
\item any triple of non-collinear points is independent over $A$ (see the picture), i.e. has $U$-rank $3$ over $A$;
\item every line has $U$-rank $2$ over $A$ (see the picture).
\end{enumerate}
\end{definition}

\begin{remark}
If each of $a$,$b$,$c$,$x$,$y$,$z$ is replaced by an element interbounded with it over $A$, then the new $6$-tuple $(a',b',c',x',y',z')$ is also a strict bounded partial quadrangle over $A$.
We see this as follows:
First of all, $a' \notin \textrm{bcl}(A)$, since otherwise we would have $a \in \textrm{bcl}(A)$ which would mean $U(a/A)=0$.
Thus, $U(a'/A) \ge 1$.
Suppose $U(a'/A)>1$.
Then, 
$$U(aa'/A)=U(a/Aa')+U(a'/A)=0+U(a'/A)>1.$$
On the other hand,
$$U(aa'/A)=U(a'/Aa)+U(a/A)=0+1=1,$$
a contradiction.
Thus, $U(a'/A)=1$.
Similarly one shows that
$$U(b'/A)=U(c'/A)=U(x'/A)=U(y'/A)=U(z'/A)=1.$$

For (i), we show that $U(a',b',x'/A)=3$.
The rest of the non-collinear triples are treated similarly.
We note first that $b' \notin \textrm{bcl}(Aa')$, since then we would have
$$b \in \textrm{bcl}(b'A) \subseteq \textrm{bcl}(Aa') \subseteq \textrm{bcl}(Aa).$$
Thus, $U(b'/Aa') \ge 1$. 
Since $U(b'/A)=1$, we have $U(b'/Aa')=1$.
Similarly, one shows that $U(x'/Aa'b')=1$, and it follows that $U(a',b',x'/A)=3$.
 
For (ii), it suffices to show that $\textrm{bcl}(A,a,b)=\textrm{bcl}(A,a',b')$, $\textrm{bcl}(A,a,c)=\textrm{bcl}(A,a',c')$ and $\textrm{bcl}(A,b,c)=\textrm{bcl}(A,b',c')$.
Since $a',b' \in \textrm{bcl}(A,a,b)$, we have $\textrm{bcl}(A,a',b') \subseteq \textrm{bcl}(A,a,b)$.
On the other hand, $a, b \in \textrm{bcl}(A,a',b')$, so $\textrm{bcl}(A,a,b) \subseteq \textrm{bcl}(A,a',b')$.
The other equalities are similar.

We say that this new partial quadrangle is \emph{boundedly equivalent} to the first one.
\end{remark}

\begin{remark}\label{cbinterbd}
If we have a strict bounded partial quadrangle, as in Definition \ref{conf}, then
$a$ is interbounded with $\textrm{Cb}(xy/Aa)$, $b$ is interbounded with $\textrm{Cb}(yz/Ab)$, and $c$ is interbounded with $\textrm{Cb}(zx/Ac)$.
We show that $a$ is interbounded with $\textrm{Cb}(xy/Aa)$.
The other statements are similar.
Denote $\alpha=\textrm{Cb}(xy/Aa)$.
Clearly $\alpha \in \textrm{bcl}(Aa)$.  

For the other part, we first note that $\alpha \notin \textrm{bcl}(A)$.
Indeed, we have 
$xy \da_\alpha Aa$, and thus
$$U(xy/\alpha)=U(xy/ Aa \alpha)=U(xy/Aa)=1.$$
Hence $\alpha \in \textrm{bcl}(A)$ would imply $U(xy/A)=1$, a contradiction.
Thus, $\alpha \in \textrm{bcl}(Aa) \setminus \textrm{bcl}(A)$, and hence $\alpha \nda_A Aa$.
This implies $a \nda_A \alpha$, so $U(a/A \alpha)<U(a/A)=1$, and thus $a \in \textrm{bcl}(A \alpha)$.
 
\end{remark}

\begin{definition}
We say that a tuple $a$ is \emph{Galois definable} from a set $A$, if it holds for every $f \in \textrm{Aut}(\mathbb{M}/A)$
that $f(a)=a$.
We write $a \in \textrm{dcl}(A)$, and say that $a$ is in the \emph{definable closure} of $A$.

We say that $a$ and $b$ are \emph{interdefinable} if $a \in \textrm{dcl}(b)$ and $b \in \textrm{dcl}(a)$.
We say that they are interdefinable over $A$ if  $a \in \textrm{dcl}(Ab)$ and $b \in \textrm{dcl}(Aa)$.

\end{definition}

\begin{definition}
We say that a set $B$ is \emph{Galois definable} over a set $A$, if $f(B)=B$ for all $f \in \textrm{Aut}(\M/A)$.
\end{definition}

\begin{definition}
We say that a group $G$ is \emph{Galois definable} over $A$ if $G$ and the group operation on $G$ are both Galois definable over $A$ as sets.
\end{definition}
 
 \begin{definition}\label{genericeka}
 Let $B \subset \M$.
 We say an element $b \in B$ is \emph{generic} over some set $A$ if $U(b/A)$ is maximal (among the elements of $B$). 
 The set $A$ is not mentioned if it is clear from the context.
 For instance, if $B$ is assumed to be Galois definable over some set $D$, then we usually assume $A=D$. 
  
Let $p=t(a/A)$ for some $a \in \M$ and $A \subset A'$.
We say $b \in \M$ is a \emph{generic} realization of $p$ (over $A'$) if $U(b/A')$ is maximal among the realizations of $p$.
 
 \end{definition}
 
We are now ready to state the main theorem of this chapter.
We will prove it as a series of lemmas.

\begin{theorem}\label{groupmain}
Suppose $A$ is a finite set, $(a,b,c,x,y,z)$ is a strict bounded partial quadrangle over $A$ and $t(a,b,c,x,y,z/A)$ is stationary.
Then, there is a group $G$ in $(\M^{eq})^{eq}$, Galois definable over some finite set $A' \subset \M$.
Moreover, a generic element of $G$ has $U$-rank $1$.
\end{theorem}
   
\begin{proof} 
We note first that if we replace the closure operator cl with the closure operator $\textrm{cl}_A$ defined by $\textrm{cl}_A(B)=\textrm{cl}(A\cup B)$, we get from $\M$ a new quasiminimal class that is closed under isomorphisms and consists of models containing the set $A$.
We may think of this new class as obtained by adding the elements of $A$ as parameters to our language.
Then, $A \subseteq \textrm{cl}(\emptyset)$.
Thus, to simplify notation, we assume from now on that $A=\emptyset$.
When using the independence calculus developed in Chapter 2, we will write $A \da B$ for $A \da_\emptyset B$.

We begin our proof by replacing 
the tuple $(a,b,c,x,y,z)$ with one boundedly equivalent with it so that $z$ and $y$ become interdefinable over $b$.
For each $n$ we first define an equivalence relation $E^n$ on $\M$ so that $xE^ny$ if and only if $\textrm{bcl}(x)=\textrm{bcl}(y)$.
Similarly, define an equivalence relation $E^{*}$ on $\M^{eq}$ so that $xE^{*}y$ if and only if $\textrm{bcl}^{eq}(x)=\textrm{bcl}^{eq}(y)$.

\begin{lemma}\label{luokat}
For each $u \in \M^n$, the element $u/E^n$ is interdefinable with $(u/E^n)/E^*$ in $(\mathbb{M}^{eq})^{eq}.$
\end{lemma}

\begin{proof}
Clearly $(u/E^n)/E^* \in \textrm{dcl}(u/E^n)$.
Suppose now $a \in \M^n$ is such that $(a/E^n)/E^*=(u/E^n)/E^*$.
We note that for each $x \in \M^n$, $(\textrm{bcl}^{eq}(x/E)) \cap \M=\textrm{bcl}(x)$ and $((\textrm{bcl}^{eq})^{eq}((x/E)/E^*)) \cap \M^{eq}=\textrm{bcl}^{eq}(x/E)$.
Thus, $\textrm({bcl}^{eq})^{eq}((x/E)/E^*) \cap \M=\textrm{bcl}(x)$.
Hence, $\textrm{bcl}(a)=\textrm{bcl}(u)$, so $a/E^n=u/E^n$.
We have thus seen that $(u/E^n)/E^{*}$ determines $u/E^n$, so $u/E^n \in \textrm{dcl}((u/E^n)/E^{*})$. 
\end{proof}

We also note that if $U(u)=1$, then $U(u/E^n)=U((u/E^n)/E^*)=1$.
Indeed,  $u/E^n$ is interbounded with $u$ and thus has $U$-rank $1$. 
As $((u/E^n)/E^*)$ is interbounded with $u/E^n$, it also has $U$-rank $1$.

Replace now $x$ with $x/E^n$, $y$ with $y/E^n$ and $z$ with $z/E^n$.
The new elements are interbounded with the old ones, so we still have a strict bounded partial quadrangle over $A$. 
From now on, denote this new $6$-tuple by $(a,b,c,x,y,z)$.
  
Let $a' \in \M$ be such that $Lt(a'/b,z,y)=Lt(a/b,z,y)$
and $a' \da  abcxyz$.  
Then, there are tuples $c',x'$ such that $\textrm{Lt}(a',c',x'/ b,z,y)=Lt(a,c,x/b,z,y)$ and in particular, $t(a',b,c',x',y,z/\emptyset)=t(a,b,c,x,y,z/\emptyset)$.
Thus, $(a',b,c',x',y,z)$ is a strict bounded partial quadrangle over $\emptyset$.
Similarly, we find an element $c'' \in \M$ such that $c'' \da abcxyza'c'x'$  and  elements $a'', x''$ so that $(a'',b,c'',x'',y,z)$ is a strict bounded partial quadrangle over $\emptyset$.
The below picture may help the reader.
\begin{figure}[htbp] \label{fig:multiscale}
\begin{center}
   \includegraphics[scale=0.9]{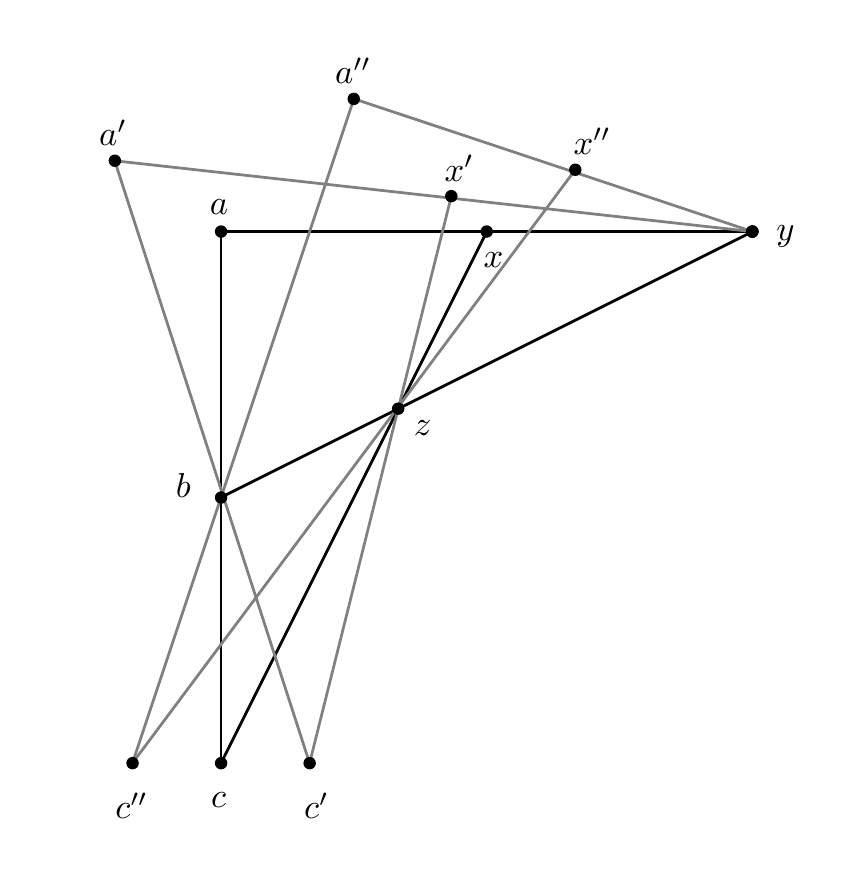}
\end{center}
\end{figure}

We will add the elements $a'$ and $c''$ as parameters in our language, but this will affect the closure operator and the independence notion.
In our arguments, we will be doing calculations both in the set-up we have before adding these parameters and the one after adding them.
We will use the notation cl and $\da$ for the setup before adding the parameters, and $\textrm{cl}^*$ and $\da^*$ for the setup after adding the parameters, i.e. for any sets $B,C,D$, $\textrm{cl}^*(B)=\textrm{cl}(B,a',c'')$ and 
$B \da^{*}_C D$ if and only if $B \da_{Ca' c''} D$.
Similary, we write $u \in \textrm{dcl}^*(B)$ if and only if $u \in \textrm{dcl}(Ba'c'')$ and use the notation $\textrm{Cb}^*(u/B)$ for $\textrm{Cb}(u/Ba'c'')$.
 
\begin{lemma}\label{interdefinable}
The tuples $yx'$ and $zx''$ are interdefinable over $a''bc'$ in $\mathbb{M}^{eq}$ after adding the parameters $a'$ and $c''$ to the language.
\end{lemma}
 
\begin{proof} 
We first prove an auxiliary claim.
 
\begin{claim}\label{tyyppi}
If $t(z'/byc'x')=t(z/byc'x')$, then $\textrm{bcl}^{eq}(z)=\textrm{bcl}^{eq}(z')$ in $\M^{eq}$.
\end{claim}

\begin{proof}
We will show that both $z$ and $z'$ are interbounded (with respect to bcl) with $\textrm{Cb}(b,y/c',x')$ and thus interbounded with each other.
Denote $\alpha=\textrm{Cb}(b,y/c',x')$.
The set $\{b,c',z\}$ is independent. 
In particular, $b \downarrow_z c'$.
But $y \in \textrm{bcl}(b,z)$ and $x' \in \textrm{bcl}(c',z)$.
Thus, $by \downarrow_z c' x'$, so $\alpha \in \textrm{bcl}(z)$.
We also have $by \da_\alpha c'x'$, so
$$U(by/\alpha)=U(by/\alpha c'x')=U(by/c'x')=1,$$
where the second equality follows from the fact that $\alpha \in \textrm{bcl}(z)$ and $z \in \textrm{bcl}(c',x')$.
Now $\alpha \notin \textrm{bcl}(\emptyset)$, since then we would have $U(by/\emptyset)=1$, contradicting our assumptions.
Thus, $\alpha \nda z$, and hence $z \nda \alpha$, so $z \in \textrm{bcl}(\alpha)$.
Hence we have seen that  $z$ is interbounded with $\alpha=\textrm{Cb}(b,y/c',x')$.
Since $t(z'/byc'x')=t(z/byc'x')$, the same holds for $z'$.
Thus, $z$ and $z'$ are interbounded.
\end{proof}

By Claim \ref{tyyppi}, 
$u=z/E^*$ if and only if there is some $w$ such that $t(w/byc'x')=t(z/byc'x')$ and $w/E^*=u$.
From this, it follows that $z/E^* \in \textrm{dcl}(byc'x')$. 
Thus, by Lemma \ref{luokat}, $z \in \textrm{dcl}(byc'x') \subseteq \textrm{dcl}(a''bc'yx')$.  

For $zx'' \in \textrm{dcl}^{*}(a''bc'yx')$, it suffices to show that $x'' \in \textrm{dcl}^{*}(a''bc'yx'z)$.
If $t(x^*/a''yzc'')=t(x''/a''yzc'')$, then $\textrm{bcl}(x^*)=\textrm{bcl}(x'')$ (this is proved like Claim \ref{tyyppi}), and thus $x''/E^* \in \textrm{dcl}(a''c''yz) \subseteq \textrm{dcl}^{*}(a''bc'x'yz)$ (note that $\textrm{dcl}^*$ is defined with $c''$ as a parameter).
By Lemma \ref{luokat}, $x'' \in \textrm{dcl}^*(a''bc'x'yz)$.

Similarly, one proves that $yx' \in \textrm{dcl}^*(a''bc'zx'')$.
\end{proof}

Let $q_1=t(yx'/a'c'')$, $q_2=t(zx''/a'c'')$.
We will consider $Cb(yx',zx''/a''bc')$ as a function from $q_1$ to $q_2$.
To see precisely how this is done, we need to introduce some concepts.
 
Suppose $p$ and $q$ are stationary types over some set $B$.
By a \emph{germ of an invertible definable function} from $p$ to $q$, we mean a Lascar type $r(u,v)$ over some finite set $C$, such that 
\begin{itemize}
\item $Lt(u,v)=r$ implies $t(u/B)=p$ and $t(v/B)=q$;
\item Suppose $Lt(u,v)=r$ and $D \subset \M$ is such that $C \subseteq D$. 
If $(u',v')$ realizes $r \vert_D$, then $u' \da_B D$, $v' \da_B D$, $v' \in \textrm{dcl}(u',D)$ and $u' \in \textrm{dcl}(v',D)$.
\end{itemize}

We will denote germs of functions by the Greek letters $\sigma, \tau$, etc.
We note that the germs can be represented by elements in $\M^{eq}$. Just represent the germ determined by some Lascar type $r$, as above, by some canonical base of $r$.
If $\sigma$ is this germ and $u$ realizes $p \vert_\sigma$, then $\sigma(u)$ is the unique element $v$ such that $(u,v)$ realizes $r \vert_\sigma$.
Note that if $a$ realizes $p\vert_B$ and $\sigma \in B$, then $\sigma(a)$ realizes $q \vert_\sigma$, and as $\sigma(a) \da_\sigma B$, the element $\sigma(a)$ realizes $q \vert_B$.

We note that the germs can be composed.
Suppose $q'$ is another stationary type over $B$, $\sigma$ is a germ from $p$ to $q$ and $\tau$ is a germ from $q$ to $q'$.
Then, by $\tau.\sigma$ we denote a germ from $p$ to $q'$ determined as follows.
Let $u$ realize $p \vert_{\sigma, \tau}$.
Then, we may think of $\tau.\sigma$ as some canonical base of $Lt((u, \tau(\sigma(u)))/\sigma, \tau)$.
We note that $t(u, \tau(\sigma(a)))/\sigma, \tau)$ is stationary since $t(u/B \sigma \tau)$ is stationary as a free extension of a stationary type and since $\tau(\sigma(u))$ is definable from $u$, $\sigma$ and $\tau$.
Thus, $\tau.\sigma \in \textrm{dcl}(\sigma, \tau)$ (see the proof of Lemma \ref{cbexists}), and the notation is meaningful.
 
We wish to apply the above methods to  the types $q_1$ and $q_2$, and thus we will do a small trick to make them stationary.
To simplify notation, denote for a while $d=(a,b,c,x,y,z,c',x',a'', x'')$.
Choose now a tuple $d' \in \M^{eq}$ such that $Lt(d'/a'c'')=Lt(d/a'c'')$ and $d' \da_{a'c''} d$.
Now, there is some $d'' \in \M$ such that $d'=F(d'')$ for some definable function $F$ and $d''  \da_{a'c''} d$.
We claim that for any subsequence $e \subseteq d$, the type $t(e/a'c''d'')$ is stationary.
Indeed, there is some subsequence $e' \subset d''$ such that $Lt(F(e')/a'c'')=Lt(e/a'c'')$ for some definable function $F$.
Thus, $t(e/a'c''e')$ (and hence $t(e/a'c''d'')$) determines $Lt(e/a'c'')$.
Let $a'c''d'' \subseteq B$ and  $f_1, f_2$ are such that $t(f_1/a'c''d'')=t(f_2/a'c''d'')=t(e/a'c''d'')$ and $f_i \da_{a'c''d''} B$ for $i=1,2$. 
Then, $Lt(f_1/a'c'')=Lt(f_2/a'c'')$.
Since $f_i \da_{a'c''} d''$ for $i=1,2$, we have by transitivity $f_i \da_{a'c''} B$, and thus $Lt(f_1/B)=Lt(f_2/B)$.
So the type is indeed stationary.

Now, we add the tuple $d''$ as parameters to our language.
Since it is independent over $a'c''$ from everything that we will need in the independence calculations that will follow, the calculations won't depend on whether we have added $d''$ or not.
Thus, we may from now without loss assume $d''=\emptyset$ to simplify notation.

Now, we may assume $q_1$ and $q_2$ are stationary. 
We will consider $Cb(yx',zx''/a''bc')$ as a germ of an invertible definable function from $q_1$ to $q_2$, and show 
that we may without loss suppose that $b=Cb(yx',zx''/a''bc')$.
Then, we will prove that for independent  $b_1, b_2$ realizing $\textrm{tp}(b/a'c'')$, $b_1^{-1}.b_2$ is a germ of an invertible definable function from $q_1$ to $q_1$.  

Note first that as $a'' \in \textrm{bcl}(bc'') \subseteq \textrm{bcl}^{*}(b)$ and $c' \in \textrm{bcl}(a'b) \subseteq \textrm{bcl}^{*}(b)$,
 we have $\textrm{Cb}(yx',zx''/a''bc')=\textrm{Cb}^*(yx',zx''/b)$.
Thus, from Lemma \ref{interdefinable}, it follows that the tuples $yx'$ and $zx''$ are interdefinable over $\textrm{Cb}^*(yx',zx''/b)$ after adding the parameters.
We will then view $\textrm{Cb}^*(yx',zx''/b)$ as a germ of a function taking $yx' \mapsto zx''$.

We claim that after adding the parameters, $b$ is interbounded with $\textrm{Cb}^*(yx',zx''/b)$.
Denote now $\alpha=\textrm{Cb}^*(yx',zx''/b)$.
Clearly, $\alpha \in \textrm{bcl}^*(b)$.
We have $yx' zx'' \da_\alpha^* b$ and thus $b \da_\alpha^* yx'zx''$.
Since $b \in \textrm{bcl}^*(y,z)$, we have $b \in \textrm{bcl}^*(\alpha)$ by Theorem \ref{main}, (xii).
Thus, we may without loss assume that $b=\textrm{Cb}^*(yx',zx''/b)$.
 
 Let $r=t(b/a',c'')$.
If $b_1, b_2$ realize $r$, then by $b_1^{-1}.b_2$ we mean the germ of the invertible definable function from $q_1$ to $q_1$ obtained by first applying $b_2$, then $b_1^{-1}$.
In other words, let $y_1 x_1'$ realize $q_1 |_{b_1b_2}$, and let $z_1x_1''=b_2.(y_1 x_1')$.
So $z_1x_1''$ realizes $q_2|_{ b_1b_2}$.
Let $y_2 x_2'=b_1^{-1}.(z_1x_1'')$ (i.e. $z_1x_1''=b_1.(y_2 x_2')$).
We may code the germ $b_1^{-1}.b_2$ by some canonical base of $t(y_1 x_1',y_2 x_2'/b_1,b_2,a',c'')$, i.e.
we will have $b_1^{-1}.b_2=\textrm{Cb}^*(y_1 x_1',y_2 x_2'/b_1,b_2)$.
At this point, we fix the type of this canonical base.
As noted before, we have
$b_1^{-1}.b_2 \in \textrm{dcl}^*(b_1,b_2)$.

\begin{lemma}\label{vapaus}
Let $b_1$, $b_2$ realize $r$ ($=\textrm{tp}(b/a'c'')$), and let $b_1 \da^* b_2$. 
Then, $b_1^{-1}.b_2 \da^* b_i$ for $i=1,2$.
In particular, $U(b_1^{-1}.b_2/a'c'')=1.$
\end{lemma}

\begin{proof}
Without loss of generality, $b_2=b$ and $b_1 \da^* a,b,c,x,y,z,c',x',a'',x''$.
Indeed, if the lemma holds for these tuples and $b_1', b_2'$ are arbitrary realizations of $r$ such that $b_1' \da^* b_2'$, then there is some automorphism $f \in \textrm{Aut}(\M/a'c'')$ such that $f(b_1)=b_1'$.
Then, $b_1' \da^* f(b_2)$. 
The type $r$ is stationary due to the trick we have done above, so this implies
$t(f(b_2)/a' c'' b_1')=t(b_2'/a'c'' b_1')$.
Hence there is an automorphism $g\in \textrm{Aut}(\M/a'c'')$ such that $g(b_1', f(b_2))=(b_1',b_2')$.
Then, $g \circ f$ is an automorphism taking $(b_1, b_2)$ to $(b_1',b_2')$ so the claim holds also for $b_1',b_2'$. 
 
We have $a' \da bzx$, and thus $b \da_{zx} a'$.
Since $b \da zx$, we get $b \da a'zx$.
On the other hand, $c'' \da a'bzx$, and thus (since $b \da a'zx$) $b \da a'c''zx$.
This implies $b \da^* zx$.
Since, $x'' \in \textrm{bcl}^*(z)$ and $c \in \textrm{bcl}^*(zx)$, we have $b \downarrow^* cxzx''$.
 Hence, $t(b/a'c''cxzx'')=t(b_1/a'c''cxzx'')$ (remember that $r$ is stationary), and there are elements $a_1, y_1, c_1',x_1',a_1''$ so that
$$t(a_1,b_1,c,x,y_1,z,c_1',x_1',a_1'',x''/a'c'')=t(a,b,c,x,y,z,c',x',a'',x''/a'c'').$$
To visualize this, think of the picture just before Lemma \ref{interdefinable}.
In the picture, keep the lines $(c,x,z)$ and $(c'',z,x'')$ fixed pointwise and move $b$ to $b_1$ by an automorphism fixing $a'c''$.
As a result, we get another similar picture drawn on top of the first one, with new elements $a_1, y_1, c_1'$ and $a_1''$ in the same configuration with respect to the fixed points as $a,y,c$ and $a''$ in the original picture.
 
\begin{claim}\label{i}
$a a_1 b b_1 \downarrow^* yx'.$ 
\end{claim}
\begin{proof}
By similar arguments as before, one sees that $y \da_{a'c''abc} b_1$ and $y \da a'c''abc$, so $y \da_{a' c''} abcb_1$
by transitivity.
As $a_1 \in \textrm{bcl}^*(b_1, c)$, we have  (by symmertry) $aa_1b b_1c \downarrow^* y$ and thus $a a_1 b b_1 \downarrow^* y.$
As $x' \in \textrm{bcl}^*(y)$, we have $a a_1 b b_1 \downarrow^* yx'.$
\end{proof}

\begin{claim}\label{ii}
$y_1 x_1' \in \textrm{bcl}^*(a,a_1,y)$
\end{claim}
\begin{proof}
$x \in \textrm{bcl}^*(a,y)$, $y_1 \in \textrm{bcl}^*(a_1,x)$ and $x_1' \in \textrm{bcl}^*(y_1)$.
\end{proof}

\begin{claim}\label{iii}
$y_1x_1'=(b_1^{-1}.b)(yx')$. 
\end{claim}
\begin{proof}
By Claim \ref{i}, $yx' \downarrow^* b b_1$, so it realizes $q_1 \vert_{bb_1}$.
On the other hand, $t(b_1y_1x_1'/a'c'')=t(byx'/a'c'')$ so $y_1 x_1' \da^* b_1$.
By similar arguments that were used to show that we may assume $b=Cb^*(yx',zx''/b)$, we also
see that we may assume $b_1=Cb^*(y_1x_1',zx''/b_1)$.
Thus, $b: yx' \mapsto zx''$ and $b_1:y_1x_1' \mapsto zx''$.
\end{proof}

\begin{claim}\label{iv}
$a a_1 \downarrow b.$
\end{claim}
\begin{proof}
$abc \downarrow^* b_1$, and thus $ab \downarrow_c^* b_1$.
As $a_1 \in \textrm{bcl}^*(b_1,c)$, we have $ab \downarrow_c^* a_1$.
By similar type of calculations that we have done before, we see that $a \da a' c'' c$.
Since $t(a_1/a'c''c)=t(a/a'c''c)$, we have that $a_1 \da a' c'' c$.
Together with $ab \downarrow_c^* a_1$, this implies $ab \da^* a_1$.
Using the independence calculus, one can verify that $b \da a'c''$ and $a \da^* b$, and thus 
$$U(aa_1b/a'c'')=U(b/a'c'')+U(a/ba'c'')+U(a_1/aba'c'')=1+1+1=3,$$
so $aa_1 \da^* b$, as wanted.
 \end{proof}

\begin{claim}\label{v}
$aa_1 \downarrow b_1.$
\end{claim}
\begin{proof}
Like Claim \ref{iv}.
\end{proof}

Denote $\sigma=b_1^{-1}.b$.
Now by Claim \ref{i}, $yx' \downarrow_{aa_1}^* a a_1 b b_1$.
Thus, by Claim \ref{ii}, $yx'y_1x_1' \downarrow _{aa_1}^* a a_1 b b_1$.
On the other hand, by Claim \ref{i}, $yx' \downarrow_{bb_1}^* a a_1 b b_1$.
By Claim \ref{iii}, $y_1x_1' \in \textrm{bcl}^*(yx',b,b_1)$, so $yx'y_1x_1' \downarrow_{bb_1}^* a a_1bb_1$.
Since $\sigma=\textrm{Cb}^*(yx',y_1x_1'/b,b_1)$, we also have 
$$\sigma=\textrm{Cb}^*(yx',y_1x_1'/a,a_1,b,b_1).$$
So, $\sigma \in \textrm{bcl}^*(a,a_1)$ since $yx'y_1x_1' \downarrow _{aa_1}^* a a_1 b b_1$.
By Claims \ref{iv} and \ref{v}, $\sigma \downarrow^* b$ and $\sigma \downarrow^* b_1$.
\end{proof}
 
Denote now $\sigma=b_1^{-1}.b_2$ (from Lemma \ref{vapaus}) and let $s=t(\sigma/a'c'')$ (note that $t(\sigma^{-1}/a'c'')=s$ also).

\begin{lemma}\label{ykkonen}
Let $\sigma_1, \sigma_2$ be realizations of $s$ such that $\sigma_1 \da^* \sigma_2$.
Then, $\sigma_1.\sigma_2$ realizes $s|_{\sigma_i}$ for $i=1,2$.
\end{lemma}

\begin{proof}
Choose some $\beta$ so that $\beta$ realizes $r|_{\sigma_1 \sigma_2}$ (remember that $r=\textrm{tp}(b/a'c'')$).
As $\sigma_1$ realizes $s$, there are some $\beta_1, \beta_2$ realizing $r$ such that $\sigma_1=\beta_1^{-1}.\beta_2$. 
By Lemma \ref{vapaus}, $\sigma_1 \downarrow^* \beta_i$ for $i=1,2$.
Thus, there is an automorphism fixing $\sigma_1$ and mapping $\beta_2 \mapsto \beta$.
Hence, there is some $\tau$ realizing $r$ so that $\sigma_1=\tau^{-1}.\beta$.
Similarly, we find $\tau'$ realizing $r$ so that $\sigma_2=\beta^{-1}.\tau'$.
Since $\beta \downarrow^* \sigma_1, \sigma_2$ and $\sigma_1 \downarrow^* \sigma_2$, we have $\sigma_1 \downarrow^*_\beta \sigma_2$.
Thus, $\tau \downarrow^*_\beta \tau'$. (as $\tau \in \textrm{bcl}^*(\sigma_1, \beta)$, $\tau' \in \textrm{bcl}^*(\sigma_2,\beta))$.
On the other hand, $\tau \da^{*} \beta$,  and thus $\tau \downarrow^* \tau'$.
Hence, $\sigma_1.\sigma_2=\tau^{-1}.\tau'$ which realizes $s$.

We still need to prove $\tau^{-1}.\tau'  \downarrow \sigma_i$ for $i=1,2$.
By Lemma \ref{vapaus}, $\tau \downarrow^* \sigma_1$.
Since $\sigma_1 \in \textrm{bcl}^*(\tau, \beta)$, we have $\tau \downarrow^* \beta$.
Thus, one sees easily the set $\{\tau, \tau',\beta\}$ is independent over $a'c''$.
By Lemma \ref{vapaus}, $\tau^{-1}.\tau' \downarrow^* \tau$, so $\tau^{-1}.\tau' \downarrow^* \tau \beta$.
Indeed, if we would have $\tau^{-1}.\tau' \in \textrm{bcl}^*(\tau, \beta)$, then it would hold that $\tau' \in \textrm{bcl}^*(\tau^{-1}.\tau', \tau) \subseteq \textrm{bcl}^*(\tau, \beta)$, which contradicts the independence of the set $\{\tau, \tau',\beta\}$.
As $\sigma_1 \in \textrm{bcl}^*(\tau, \beta)$, it follows that $\tau^{-1}.\tau' \downarrow^* \sigma_1$.
Similarly, $\tau^{-1}.\tau' \downarrow^* \sigma_2$.
\end{proof}

Let $G$ be the group of germs of functions from $q_1$ to $q_1$ generated by $\{\sigma \, | \, \sigma \textrm{ realizes } s\}$ (note that this set is closed under inverses and thus indeed a group). 

\begin{lemma}
For any $\tau \in G$, there are $\sigma_1, \sigma_2$ realizing $s$ such that $\tau=\sigma_1.\sigma_2$. 
\end{lemma}

\begin{proof}
It is enough to show that if $\tau_i$ realize $s$ for $i=1,2,3$ then there are $\sigma_1, \sigma_2$ realizing $s$ so that $\tau_1.\tau_2.\tau_3=\sigma_1.\sigma_2$.
Let $\sigma$ realize $s|_{\tau_1 \tau_2 \tau_3}$.
By Lemma \ref{ykkonen}, $\sigma^{-1}.\tau_2$ realizes $s |_{\tau_2}$.
Now, $\sigma \tau_2 \downarrow_{\tau_2}^* \tau_1 \tau_2 \tau_3$, and thus $\sigma^{-1}.\tau_2 \downarrow_{\tau_2}^* \tau_1 \tau_2 \tau_3$.
As $\sigma^{-1}.\tau_2 \downarrow^* \tau_2$, we get $\sigma^{-1}.\tau_2 \downarrow^* \tau_1 \tau_2 \tau_3$.
Thus, by Lemma \ref{ykkonen}, $(\sigma^{-1}.\tau_2).\tau_3$ realizes $s$.
Again by Lemma \ref{ykkonen}, $\tau_1.\sigma$ realizes $s$.
By choosing $\sigma_1=\tau_1. \sigma$ and $\sigma_2=\sigma^{-1}.\tau_2.\tau_3$, we get $\sigma_1.\sigma_2=\tau_1.\tau_2.\tau_3$.
\end{proof}

Consider the set 
$$G'=\{(\sigma_1.\sigma_2) \, | \, \textrm{$\sigma_1, \sigma_2$ are realizations of $s$}\}.$$
It is clearly Galois definable over $a'c''$.
Let $E$ be the equivalence relation such that for $\gamma_1, \gamma_2 \in G'$, $(\gamma_1, \gamma_2) \in E$ if and only if $\gamma_1(u)=\gamma_2(u)$ for all $u$ realizing $q_1 \vert_{\gamma_1 \gamma_2}$.
Then, $G=G'/E$, and $G$ is  Galois definable over $a'c''$.

It remains to prove that  for a generic $\sigma_1. \sigma$ it holds that $U(\sigma_1.\sigma_2/a'c'')=1$.
We note first that for $\sigma=b_1^{-1}.b_2$, we have $U(\sigma/a'c'')=1$.
Indeed, since $\sigma \da^* b_1$, we have
$$U(\sigma/a'c'')=U(\sigma/a'c''b_1) \le U(b_2/a'c''b_1)=1,$$
where the inequality follows from the fact that $\sigma \in \textrm{dcl}^*(b_1, b_2)$, and the last equality from the fact that $b_1 \da^* b_2$.
On the other hand, we cannot have $\sigma \in \textrm{bcl}(a'c'')$, since $b_2 \in \textrm{dcl}(b_1, \sigma)$.
Thus, $U(\sigma/a'c'')=1$.
If $\sigma_1 \da^* \sigma_2$, then by Lemma \ref{ykkonen}, $\sigma_1.\sigma_2$ realizes $s$, and thus $U(\sigma_1.\sigma_2/a'c'')=1$.
If $\sigma_1 \nda^* \sigma_2$, then $U(\sigma_1.\sigma_2/a'c'') \le 1$.
This proves the theorem.
\end{proof}

From now on, we will use the term \emph{group configuration} for a configuration as in Definition \ref{conf}.
We will next give an example of a situation where the configuration arises - that of a non-trivial locally modular pregeometry.
For this, we need some definitions.

\begin{definition}
Let $(S, \textrm{cl})$ be a pregeometry.

We say it is modular, if it holds for all closed sets  $A, B \subseteq S$ that
$$\textrm{dim}(A \cup B)=\textrm{dim}(A)+\textrm{dim}(B)-\textrm{dim}(\textrm{cl}(A \cap B)).$$

If  there exists some tuple $a \in S$ such that the pregeometry $(S, \textrm{cl}_{a})$ is modular, where the operator $\textrm{cl}_a$ is defined so that $\textrm{cl}_a(A)=\textrm{cl}(Aa)$ for any $A \subseteq S$, then we say the pregeometry $(S, \textrm{cl})$ is \emph{locally modular}.

If a pregeometry is not locally modular, we say it is \emph{non locally modular}.
\end{definition}

It is easy to see that if $V$ is a vector space and span is the linear span operator, then $(V, \textrm{span})$ is modular.
On the other hand, affine geometry is not modular, but once you add the point of origin, it becomes a vector space.
Thus, it is locally modular.
An algebraically closed field together with the algebraic closure operator provides an example of a non locally modular pregeometry.
For more details on these, see e.g. \cite{MarIntro}.
 
\begin{definition}
Let $(S, \textrm{cl})$ be a pregeometry.
We say it is \emph{trivial} if it holds for every $A \subseteq S$ that 
$$\textrm{cl}(A)=\bigcup_{a \in A} \textrm{cl}(a).$$
\end{definition}
 
\begin{lemma}\label{locmodgroup}
Suppose $(\M, \textrm{bcl})$ is a non-trivial locally modular pregeometry.
Then, there exists a group configuration in $\M$.
\end{lemma}

\begin{proof}
We may without loss assume that we have added the necessary parameters in our language so that $(\M, \textrm{bcl})$ is modular.
Let $a_1, \ldots, a_{n}$ be such that $\textrm{dim}(a_1, \ldots, a_n)=n-1$ and every $n-1$ -element subset also has dimension $n-1$.
Such elements exist by non-triviality.
Suppose moreover that $n$ is the least number so that such elements can be found.
The modularity of the pregeometry is preserved in further localizations, so we now localize at $(a_3, \ldots, a_n)$.
Again, we may simplify notation by assuming that these elements are parameters in our language.
 
Write now $a=a_1$ and $b=a_2$.
Then, $a$ and $b$ are independent of each other and $\textrm{bcl}(a) \cup \textrm{bcl}(b) \subsetneq \textrm{bcl}(a,b)$.
Next, pick some element $x$ such that it is independent from $\{a,b\}$ and $t(x/a)=t(b/a)$.

Choose now some $c \in \textrm{bcl}(a,b) \setminus (\textrm{bcl}(a) \cup \textrm{bcl}(b))$.
We note that then, $c \notin \textrm{bcl}(a,x) \cup \textrm{bcl}(b,x)$.
Indeed, if we had $c \in \textrm{bcl}(a,x)$, we would have $x \in \textrm{bcl}(a,c) \subseteq \textrm{bcl}(a,b)$, a contradiction. 
Similarly, one sees that $c \notin \textrm{bcl}(b,x)$.

Next, we choose some 
$y \in \textrm{bcl}(a,x) \setminus (\textrm{bcl}(a) \cup \textrm{bcl}(x)).$
Then, it will hold that
$$y \notin \textrm{bcl}(a,b) \cup \textrm{bcl}(b,c) \cup \textrm{bcl}(a,c) \cup \textrm{bcl}(x,b) \cup \textrm{bcl}(x,c).$$
This is again easily seen by using the exchange property of the pregeometry.

Since the pregeometry is modular, we have
$$\textrm{dim}(\textrm{bcl}(x,c) \cup \textrm{bcl}(b,y))=\textrm{dim}(\textrm{bcl}(x,c))+\textrm{dim}(\textrm{bcl}(b,y))-\textrm{dim}(\textrm{bcl}(x,c) \cap \textrm{bcl}(b,y)).$$
But
$$\textrm{dim}(\textrm{bcl}(x,c) \cup \textrm{bcl}(b,y)) \le \textrm{dim}(x,c,b,y)=\textrm{dim}(x,c,b)=3,$$
so $\textrm{dim}(\textrm{bcl}(x,c) \cap \textrm{bcl}(b,y)) \ge 1$.
Let $z \in \textrm{bcl}(x,c) \cap \textrm{bcl}(b,y)) \setminus \textrm{bcl}(\emptyset)$.
Using exchange, one shows first that
$$z \notin \textrm{bcl}(x) \cup \textrm{bcl}(c) \cup \textrm{bcl}(b) \cup \textrm{bcl}(y).$$
Then, again using exchange, one shows that
$$z \notin \textrm{bcl}(a,x) \cup \textrm{bcl}(a,y) \cup \textrm{bcl}(x,y) \cup \textrm{bcl}(a,b) \cup \textrm{bcl}(a,c) \cup \textrm{bcl}(b,c) \cup \textrm{bcl}(b,x) \cup \textrm{bcl}(c,y).$$

Now, $(a,b,c,x,y,z)$ is a group configuration.
\end{proof}
 
\chapter{Groups in Zariski-like structures}

In this chapter, we suppose that $\M$ is a monster model for a quasiminimal class as introduced in Chapter 2.
As an attempt to generalize Zariski geometries to this context, we will present axioms for a Zariski-like structure.
These axioms capture some of the properties of the irreducible closed sets in Zariski geometries that are needed for finding a group in that context.
We then apply the group configuration theorem from Chapter 3 to show that if $\M$ satisfies these axioms and the pregeometry obtained from the bounded closure operator is non-trivial, then a $1$-dimensional group can be found in $(\M^{eq})^{eq}$.
The argument is a modification of the one presented for Zariski geometries in \cite{HrZi}.

To simplify notation, we often write $a \da b$ for $a \da_\emptyset b$ and $U(a)$ for $U(a/\emptyset)$.
In the following definition, when speaking about indiscernible sequences, we don't assume that they are non-trivial.

\begin{definition}\label{zariskilike}
We say that an infinite-dimensional quasi-minimal pregeometry structure (in the sense of \cite{monet} and \cite{kir}) $\M$ is \emph{Zariski-like} if for each $n$, there is a countable collection of subsets of $\M^n$, called the \emph{irreducible sets} satisfying the following nine axioms:  

(ZL1)
The irreducible sets are Galois definable, i.e. if $C \subset \M^n$ is irreducible, then $f(C)=C$ for every $f \in \textrm{Aut}(\M/ \emptyset)$.

(ZL2) 
For each $n$ and each $a \in \M^n$, there is some irreducible $C \subset \M^n$ such that $a$ is a generic point of $C$ (over $\emptyset$).  
 
(ZL3) 
The generic elements (i.e. elements of maximal $U$-rank over $\emptyset$) of an irreducible set have the same Galois type.

(ZL4) 
If $C_1, C_2$ are irreducible, $a \in C_1$ generic and $a \in C_2$, then $C_1 \subseteq C_2.$

(ZL5) 
If $C_1, C_2$ are irreducible, $(a,b) \in C_1$ is generic, $a$ is a generic element of $C_2$ and $(a',b') \in C_1$, then $a' \in C_2$. 

(ZL6)
If $C \subset \M^n$ is irreducible and $f$ is a coordinate permutation on $\M^n$, then $f(C)$ is irreducible.

\vspace{4mm}

Before we can continue listing the axioms, we have to adapt the definition of specialization from the Zariski geometry context (Definition \ref{SpecHrZi}) to our setting.
 
\begin{definition}
If $A \subset \M$, we say that a function $f: A \to \M$ is a \emph{specialization} if for any $a_1, \ldots, a_n \in A$ and for any irreducible set $C \subseteq \M^n$, it holds that if $(a_1, \ldots, a_n) \in C$, then $(f(a_1), \ldots, f(a_n)) \in C$.
 If $A=(a_i : i \in I)$, $B=(b_i:i\in I)$ and the indexing is clear from the context, we write $A \to B$ if the map $a_i \mapsto b_i$, $i \in I$, is a specialization.

If $a$ and $b$ are finite tuples and $a \to b$, we denote $\textrm{rk}(a \to b)=U(a/\emptyset)-U(b/\emptyset)$.
\end{definition}
  
We also present the definitions of strongly regular and strongly good specializations as generalizations of the regular and good specializations of Definition \ref{goodHrZi}.

\begin{definition}
We define a \emph{strongly regular} specialization as follows:
\begin{itemize}
\item Isomorphisms are strongly regular;
\item If $a \to a'$ is a specialization and $a \in \M$ is generic over $\emptyset$, then $a \to a'$ is strongly regular;
\item $aa' \to bb'$ is strongly regular if $a \downarrow_\emptyset a'$ and the specializations $a \to b$ and $a' \to b'$ are strongly regular.
\end{itemize}
\end{definition}

\begin{definition}\label{good}
We define \emph{strongly good} specializations recursively as follows.
Strongly regular specializations are strongly good.
Let $a=(a_1, a_2, a_3)$, $a'=(a_1', a_2', a_3')$, and $a \to a'$.
Suppose:
\begin{enumerate}[(i)]
\item $(a_1, a_2) \to (a_1', a_2')$ is strongly good.
\item $a_1 \to a_1'$ is an isomorphism.
\item $a_3 \in \textrm{bcl}(a_1)$.
\end{enumerate}
Then, $a \to a'$ is strongly good.
\end{definition}

(ZL7)
Let $a \to a'$ be a strongly good specialization such that $U(a)-U(a') \le 1$.
Then any specializations $ab \to a'b'$, $ac \to a'c'$ can be amalgamated: there exists $b^{*}$, independent from $c$ over $a$, such that $\textrm{t}^g(b^{*}/a)=\textrm{t}^g(b/a)$, and $ab^{*}c \to a'b'c'$.
 
(ZL8)
Let $(a_i:i\in I)$ be independent and strongly indiscernible over $b$.
Suppose $(a_i':i\in I)$ is strongly indiscernible over $b'$, and $a_ib \to a_i'b'$ for each $i \in I$.
Further suppose $b \to b'$ is a strongly good specialization and $U(b)-U(b') \le 1$.
Then, $(ba_i:i \in I) \to (b'a_i':i\in I)$.

\vspace{4mm}

To be able to state the last axiom, we need to recall the concept of an unbounded set.

\begin{definition}
Denote by $\mathcal{P}_{<\omega}(I)$ the set of finite subsets of $I$.

We say that $S \subseteq \mathcal{P}_{<\omega}(I)$ is \emph{unbounded} if for every $A \in \mathcal{P}_{<\omega}(I)$, there is some $B \in S$ such that $A \subseteq B$.
\end{definition}
 
(ZL9) 
Let $\kappa$ be a (possibly finite) cardinal.
Let  $a_i, b_i \in \M$ with $i < \kappa$, such that $a_0 \neq a_1$ and $b_0=b_1$.
Suppose $(a_i)_{i<\kappa} \to (b_i)_{i <\kappa}$ is a specialization.
Assume there is some unbounded  $S \subset \mathcal{P}_{<\omega}(\kappa)$ satisfying the following conditions:
\begin{enumerate}[(i)]
\item  $0,1 \in X$ for all $X \in S$;
\item For all $X,Y \in S$ such that $X \subseteq Y$,  and for all sequences $(c_i)_{i \in Y}$ from $\M$, the following holds: 
If $c_0=c_1$, $ (a_i)_{i \in Y} \to (c_i)_{i \in Y} \to (b_i)_{i \in Y},$
and $\textrm{rk}((a_i)_{i \in Y} \to (c_i)_{i \in Y}) \le 1$, then $\textrm{rk}((a_i)_{i \in X} \to (c_i)_{i \in X}) \le 1$.
\end{enumerate}

Then, there are $(c_i)_{i<\kappa}$ such that   
$$(a_i)_{i \in \kappa} \to (c_i)_{i \in \kappa} \to (b_i)_{i \in \kappa},$$
$c_0=c_1$ and $\textrm{rk}((a_i)_{i \in X} \to (c_i)_{i \in X}) \le 1$ for all $X \in S$.
\end{definition}

The axioms (ZL1)-(ZL6) state some general properties of the irreducible closed sets in the Zariski geometry context.
Axioms (ZL7) and (ZL8) restate Lemmas \ref{amalgHrZi} and \ref{yhdistHrZi} in our context.
Axiom (ZL9) captures Lemma \ref{dimthmspecHrZi} and the traces of Compactness needed in finding the group configuration.
 
\begin{remark}\label{ZL9implies}
We note that (ZL9) implies Lemma \ref{dimthmspecHrZi}, i.e.
the usual dimension theorem of Zariski geometry.
Indeed, suppose $\kappa=n$, a finite cardinal and $a=(a_0, \ldots, a_{n-1})$ and $b=(b_0, \ldots, b_{n-1})$ are such that
$a \to b$, $a_0 \neq a_1$ and $b_0=b_1$.
Let $S=\{n\}$.
Then, conditions (i) and (ii) in (ZL9) hold, so we find an $n$-tuple $c$ such that
$a \to c \to b$, $U(a)-U(c) \le 1$ and $c_0=c_1$.
\end{remark}

In the following, we note that Zariski-like structures are indeed generalizations of Zariski geometries.
 
\begin{example}\label{ZariskiZariskiLike} 
Let $D$ be a Zariski geometry.
Since $D$ is strongly minimal, it is also quasiminimal.
Consider the collection of closed sets in the language.
Then, the irreducible (in the topological sense) ones among them satisfy the axioms (ZL1)-(ZL9).
Indeed, the axioms (ZL1)-(ZL6) are clearly satisfied.
It is well known that on a strongly minimal structure, $U$-ranks and Morley ranks coincide.
On a Zariski geometry, first-order types imply Galois types.
Moreover, every strongly regular specialization is regular, and every strongly good specialization is good.
Hence, (ZL7) is Lemma \ref{amalgHrZi} and (ZL8) is Lemma \ref{yhdistHrZi}.
(ZL9) holds by Compactness.
 \end{example}

\begin{example}
Consider the model class from Example \ref{equivalence}.
For each $n$, define the irreducible sets of $\M^n$ to be those definable with finite conjunctions of formulae of the form $x_i=x_j$ or $E(x_i, x_j) \wedge \neg x_i =x_j$.
In addition, we require that   
if $E(x_i, x_j) \wedge \neg x_i =x_j$ belongs to the conjunction, then also $E(x_j, x_i) \wedge \neg x_i =x_j$ belongs there, that if 
both $E(x_i, x_j)  \wedge \neg x_i =x_j$ and $E(x_j, x_k) \wedge \neg x_j =x_k$ belong to the conjunction, then either $x_i=x_k$ or $E(x_i,x_k) \wedge \neg x_i=x_k$ belongs there,
and that if both $E(x_i, x_j) \wedge \neg x_i=x_j$ and $x_i=x_k$ belong to the conjunction, then also $E(x_k, x_j) \wedge \neg x_k=x_j$ belongs there.

Now, it is quite easy to verify that the class satisfies the axioms (ZL1)-(ZL9).
\end{example} 
 
\section{Families of plane curves}
 
We will show that in a quasiminimal structure with a non-trivial pregeometry, satisfying the axioms (ZL1)-(ZL9), we can find the group configuration from Chapter 3.
When doing this for non locally modular structures, families of plane curves will play an important role.

\begin{definition}\label{goodfor}
Let $C \subset \M^{n+m}$ be an irreducible set.
We say an element $a \in \M^n$ is \emph{good} for $C$ if there is some $b \in \M^m$ so that $(a,b)$ is a generic element of $C$.
\end{definition}

\begin{definition}
Let $\M$ be a Zariski-like structure, let $E \subseteq \M^n$ be irreducible, and let $C \subseteq \M^2 \times E$ be an irreducible set.
For each $e \in E$, denote $C(e)=\{(x,y) \in \M^2 \, | \, (x,y,e) \in C\}$.
Suppose now $e \in E$ is a generic point.
If $e$ is good for $C$ and the generic point of $C(e)$ has $U$-rank $1$ over $e$, then we say that $C(e)$ is a \emph{plane curve}.
We say $C$ is a \emph{family of plane curves} parametrized by $E$.  

We say that $\alpha$ is the \emph{canonical parameter} of the plane curve $C(e)$ if $\alpha=\textrm{Cb}(x,y/e)$ for a generic element $(x,y) \in C(e)$.
We define the $\emph{rank}$ of the family to be the $U$-rank of $\textrm{Cb}(x,y/e)$ over $\emptyset$, where $e \in E$ is generic, and $(x,y)$ is a generic point of $C(e)$.
\end{definition}

\begin{definition}
We say a family of plane curves $C \subset \M^2 \times E$ is \emph{relevant} if for a generic  $e \in E$ and a generic point $(x,y) \in C(e)$ it holds that $x,y \notin \textrm{bcl}(e)$.
\end{definition}

When proving that a one-dimensional group can be found from a Zariski-like structure, the non locally modular case will be the difficult one.
In this case, finding the group configuration will lean heavily on the fact that not being locally modular implies the existence of a relevant family of plane curves of rank at least $2$.
This fact is stated in the following lemma that is essentially the same as (ii) $\Rightarrow$ (iii) of Lemma 3.4. in \cite{Hy}.

\begin{lemma}\label{relevantfamily}
Suppose $\M$ is a Zariski-like structure, and every relevant family of plane curves on $\M$ has rank $1$.
Then, $\M$ is locally modular.
\end{lemma}

\begin{proof}
We prove first the following claim and then show that local modularity follows.

\begin{claim}\label{auxaux}
If $B$ is a finite subset of $\M$, $e$ is a tuple of elements of $\M$, $b,c \in \M \setminus \textrm{bcl}(\emptyset)$,   $U((b,c)/e)=2$  and $U((b,c)/Be)=1$, then there is some $d \in (\textrm{bcl}(bce) \cap \textrm{bcl}(Be)) \setminus \textrm{bcl}(e)$.
\end{claim}

\begin{proof}
We suppose first that $e$ is a singleton.
If we have $b \in \textrm{bcl}(Be)$ or $c \in \textrm{bcl}(Be)$, then we may choose $d=b$ or $d=c$, respectively.
Assume now $b, c \notin \textrm{bcl}(Be)$.
Let $f=Cb((b,c)/Be)$.
Then, $f \in \textrm{bcl}(Be)$ and $U(bc/f)=1$.
We claim that  $U(f/\emptyset)=1$.
Suppose for the sake of contradiction that $U(f/\emptyset)>1$.
Let $C$ be the locus of $(b,c,Be)$ and let $E$ be the locus of $Be$.
Then, applying (ZL5), we see that $C$ is a relevant family of plane curves parametrized by $E$.
This family has rank greater than $1$, which contradicts our assumptions.
 
Since $U(b,c/\emptyset)=2$, we have $bc \nda f$, and thus, by symmetry, $f \nda bc$.
Hence, $f \in \textrm{bcl}(bc)$.
Since we have $bc \da e$, it follows that $f \da e$.
We note that we also have $f \da b$.
Indeed, otherwise we would have $f \in \textrm{bcl}(b)$, and thus $c \in \textrm{bcl}(f,b) \subset \textrm{bcl}(b)$, contradicting the fact that $U(b,c/\emptyset)=2$.
Hence, by the uniqueness of the generic type, we have $t(e/f)=t(b/f)$, and thus there is some $d \in \M$ such that
$t(ed/f)=t(bc/f)$.
Since $c \in \textrm{bcl}(bf)$, we have $d \in \textrm{bcl}(ef) \subseteq \textrm{bcl}(Be)$.
Since $f \in \textrm{bcl}(bc)$, we also have $d \in \textrm{bcl}(bce)$.
Moreover,  $t(ed/f)=t(bc/f)$ implies $t(ed/\emptyset)=t(bc/\emptyset)$, so $d \notin \textrm{bcl}(e)$.
So, the claim holds in the case that $e$ is a singleton.

Suppose now $e=(e_0, \ldots, e_n)$ is a finite tuple from $\M$, $U((b,c)/e)=2$ and $U((b,c)/Be)=1$.
 We now have $bc \nda_{e_0} Be$, and thus, by what we have proved above, there is some
$d \in (\textrm{bcl}(bce_0) \cap \textrm{bcl}(Be)) \setminus \textrm{bcl}(e_0)$.
Since $U((b,c)/e)=2$, we have $bc \da_{e_{0}} e$.
Since $d \in \textrm{bcl}(bce_0)$, this implies $d \da_{e_{0}} e$, so $d \notin \textrm{bcl}(e)$.
\end{proof}

We now claim that $\M$ becomes modular after a localization. 
For any finite tuple $e$ of elements of $\M$, denote $\textrm{bcl}_e(X)=\textrm{bcl}(Xe)$ and
$\textrm{dim}_e(X)=\textrm{dim}_{\textrm{bcl}_e}(X)$.
We need to prove that there is some $e$ such that for any finite sets $A$ and $B$, it holds that 
$$\textrm{dim}_e(A \cup B)=\textrm{dim}_e(A)+\textrm{dim}_e B- \textrm{dim}_e(\textrm{bcl}_e(A) \cap \textrm{bcl}_e(B)).$$

Our auxiliary claim expresses that if  we have $\textrm{dim}_e(b,c)=2$, $\textrm{dim}_e(b/B)=\textrm{dim}_e(c/B)=1$ and $\textrm{dim}_e(b,c/B)=1$, then there is some $d \in \textrm{bcl}_e(b,c) \cap \textrm{bcl}_e(B)$ such that $\textrm{dim}_e(d)=1$. 
In other words, $d$ proves that modularity holds in case that $A$ is a two-element set, $A=\{b,c\}$.
Indeed, in this case, we have $\textrm{dim}_e(\textrm{bcl}_e(b,c) \cap \textrm{bcl}_e(B)) \le 1$ since $\textrm{dim}_e(b,c)=2$ and $\textrm{dim}_e(b,c/B)=1$, so at most one of these elements is in $\textrm{bcl}_e(B)$. 
The existence of $d$ proves that this dimension is actually $1$, as needed for the local modularity.

In Claim \ref{auxaux}, the tuple $e$ was arbitrarily chosen, so we know  that the local modularity condition holds whenever $A$ is a two-element set, no matter where we localize.
Suppose now the pregeometry is not locally modular.
Let now $A=\{a_1, \ldots, a_n\}$, independent over $e$, and suppose $e, B$ are such that $A \nda_e B$ forms a counterexample which shows that the pregeometry does not become modular when localizing at $e$.
Moreover, assume that $n$ is the least possible number for which we can find sets $A$, $B$ and $e$ forming a counterexample (then, of course, $n>2$).
We note that if we localize again at some finite tuple, then our assumption concerning one-dimensionality of all plane curves still holds.
Thus, we may assume that $\textrm{bcl}(Ae) \cap \textrm{bcl}(Be)=\textrm{bcl}(e)$ (if this does not hold, then enlarge $e$ so that it holds).
 
It now follows that $a_1 \ldots a_{n-1} \nda_{ea_n} B$.
Indeed, if we had $a_1 \ldots a_{n-1} \da_{ea_n} B$, then we would have $a_n \nda_e B$ (otherwise, we would get $A \da_e B$ by transitivity), so 
$$a_n \in \textrm{bcl}(Ae) \cap \textrm{bcl}(Be) \setminus \textrm{bcl}(e).$$
But this contradicts the assumption that $A$ and $B$ form a counterexample for the modularity when localizing at $e$
So, we have  $a_1 \ldots a_{n-1} \nda_{ea_n} B$.

Localize at $a_n$.
But then, since the local modularity condition holds for sets of size $n-1$, 
there is some $k \in (\textrm{bcl}(a_1, \ldots, a_n e) \cap \textrm{bcl}(B,a_n,e)) \setminus \textrm{bcl}(a_n,e)$.
Since $k \in \textrm{bcl}(B,a_n,e) \setminus \textrm{bcl}(a_n,e)$, we get $ka_n \nda_e B$ , which implies that there is some
 $$d \in (\textrm{bcl}(ka_ne) \cap \textrm{bcl}(Be))\setminus \textrm{bcl}(e)\subseteq (\textrm{bcl}(Ae) \cap (Be))\setminus \textrm{bcl}(e)= \emptyset,$$ 
a contradiction.
\end{proof}
 
\section{Groups from indiscernible arrays} 
 
In the non locally modular case, we are going to use a relevant family of plane curves of rank at least $2$ to build the group configuration from Chapter 3.
In our setting, it will be useful to reformulate this configuration in terms of indiscernible arrays.

 \begin{definition}
We say that $f=(f_{ij} : i \in I, j \in J)$, where $I$ and $J$ are ordered sets, is an \emph{indiscernible array} over $A$ if whenever $i_1, \ldots, i_n \in I$, $j_1, \ldots, j_m \in J$, $i_1 < \ldots < i_n$, $j_1 <\ldots < j_m$, then $t((f_{i_{\nu} j_{\mu}}: 1 \le \nu \le n, 1 \le \mu, \le m)/A)$ depends only on the numbers $n$ and $m$.

If at least the $U$-rank of the above sequence depends only on $m,n$, and $U((f_{i_{\nu} j_{\mu}}: 1 \le \nu \le n, 1 \le \mu, \le m)/A)= \alpha(m,n)$, where $\alpha$ is some polynomial of $m$ and $n$, we say that $f$ is \emph{rank-indiscernible}  over $A$, of type $\alpha$, and write $U(f;n,m/A)=\alpha(n,m)$.
\end{definition}  
 
If $(c_{ij}: i \in I, j \in J)$ is an array and $I' \subseteq I$, $J' \subseteq J$ , we write $c_{I'J'}$ for $(c_{ij}: i \in I', j \in J'$).
If $|I'|=m$ and $|J'|=n$, we call $c_{I'J'}$ an $m \times n$ -\emph{rectangle} from $c_{ij}$.
 
\begin{lemma}
Let $f=(f_{ij} : i, j \in \kappa)$ be an indiscernible array over $A$, and let $\kappa \ge \o_1$.
Then, for all $m,n$, all the $m \times n$ rectangles of $f$ have the same Lascar type over $A$. 
\end{lemma}

\begin{proof}
Suppose not.
Let $m, n$ be such that all the $m \times n$ -rectangles don't have the same Lascar type over $A$.
Let $(B_k)_{k< \kappa}$ be a sequence of disjoint $m \times n$ -rectangles such that if $f_{ij} \in B_{k_1}$ and $f_{i'j'} \in B_{k_2}$, where $k_1 <k_2$, then $i<i'$ and $j<j'$.
There is some $I \subset \kappa$, $\vert I \vert =\kappa$ such that $(B_k)_{k \in I}$ is Morley over some model $\A \supset A$.
Relabel the indices so that $I=\kappa$.
By the counterassumption, there is some $m \times n$ rectangle $B$ such that $t(B/A) \neq t(B_0/A)$.
Let $0<\lambda < \kappa$ be such that whenever $f_{ij} \in B$ and $f_{i'j'} \in B_\lambda$, then $i<i'$ and $j<j'$.
Now, $B_0 B$ and $B_0 B_\lambda$ are both $2m \times 2n$ -rectangles, so $t(B_0 B/A)=t(B_0 B_\lambda /A)$.
This is a contradiction, since $Lt(B_0/A) \neq Lt(B/A)$, $Lt(B_0/A)=Lt(B_\lambda/A)$ and automorphisms preserve the equality of Lascar types.
\end{proof} 
 
The following lemma will yield the connection between the indiscernible arrays and the group configuration from Chapter 3.

\begin{lemma}\label{groupexists}
 Let $(f_{ij} : i,j< \omega_1)$  be an indiscernible array of elements of $\M$, of type  $\alpha(m,n)=m+n-1$ over some finite parameter set $B$.
Then there exists a Galois definable $1$-dimensional group in $(\M^{eq})^{eq}$. 
 \end{lemma}

\begin{proof}
We will show that there is in $\M$ a group configuration as in Definition \ref{conf}, and thus a Galois definable $1$-dimensional group by Theorem \ref{groupmain}. 
Let $\A$ be a countable model such that $B \subset \A$ and $f \da_B \A$
(note that we can find such a model by constructing a sequence $(a_i)_{i<\o}$ independent from $f$ over $B$, and then taking $\A=\textrm{bcl}(B, (a_i)_{i<\o})$).
We write $\textrm{bcl}_\A(X)$ for $\textrm{bcl}(\A \cup X)$.
To simplify notation, we assume that $B=\emptyset$.

We prove first an auxiliary claim. 
 
\begin{claim}\label{claimura}
Suppose $U(c/ d_1 d_2 \A)=U(c/d_1 \A)=U(c/d_2 \A)$.
Then there exists $e \in \textrm{bcl}_\A(d_1) \cap \textrm{bcl}_\A(d_2)$ such that $U(c/e \A)=U(c/d_1 d_2 \A)$.
\end{claim}

\begin{proof}
 
Let $E= \textrm{bcl}_\A(d_1) \cap \textrm{bcl}_\A(d_2)$. 
Let $c=(c_1, \ldots, c_m)$ and suppose $c_1, \ldots, c_k$ are independent over $\A$ from $d_1 d_2$ while $c \in \textrm{bcl}_\A(c_1, \ldots, c_k, d_1, d_2)$.
Then
\begin{displaymath}
c \in \textrm{bcl}_\A(c_1, \ldots, c_k, d_i, E)
\end{displaymath}
for $i=1,2$.
We will show that $c \in \textrm{bcl}_\A(c_1, \ldots, c_k, E)$.
Let $e \in E$ be a finite tuple such that $E=\textrm{bcl}_\A(e)$.
Suppose $U(d_2/\A d_1)=r$ and $U(d_2/\A \cup E)=r+l$.
We may assume without loss of generality that $d_2=e \cup \{d_{2,1}, \ldots, d_{2,r}, d_{2,r+1}, \ldots, d_{2,r+l}\}$ where $d_{2,1}, \ldots, d_{2,r}$ are independent over $\A d_1$ and $d_{2,r+1}, \ldots, d_{2,r+l} \in \textrm{bcl}_\A(d_1, d_{2,1}, \ldots, d_{2,r})$.
Now  
\begin{displaymath}
c \in \textrm{bcl}(c_1, \ldots, c_k, d_{2,1}, \ldots, d_{2,r}, \ldots, d_{2,r+l}, e,a)
\end{displaymath}
for some $a \in \A$ such that $d_{2,r+1}, \ldots, d_{2,r+l} \in \textrm{bcl}(d_1, d_{2,1}, \ldots, d_{2,r},a)$.
We will show that we can move the parameters $d_{2,1}, \ldots, d_{2,r}, \ldots, d_{2,r+l}$ one by one to $E$ using automorphisms.
We do this first for $d_{2,1}$.

We note first that $c \da_{\A d_1, d_{2,2}, \ldots, d_{2,r}} d_{2,1}$.
Indeed,
\begin{displaymath}
U (c/ \A d_1) \ge U(c/ \A, d_1, d_{2,2}, \ldots, d_{2,r}) \ge U(c/ \A, d_1, d_{2,1}, d_{2,2}, \ldots, d_{2,r})=U(c/\A, d_1, d_2)=U(c/\A d_1),
\end{displaymath}
so $U(c/ \A, d_1, d_{2,2}, \ldots, d_{2,r})=U(c/ \A, d_1, d_{2,1}, d_{2,2}, \ldots, d_{2,r})$.

We have $d_{2,1} \notin \textrm{bcl}_\A(d_1, d_{2,2}, \ldots, d_{2,r})$, and thus, $d_{2,1} \da_\A d_1 d_{2,2},\ldots d_{2,r}$.
Hence, by transitivity,
$$d_{2,1} \da_\A d_1 d_{2,2} \ldots, d_{2,r} c.$$
Then, there is some finite set $A \subset \A$ such that $a \in A$, $d_{2,1} \da_A d_1 d_{2,2} \ldots, d_{2,r} c$ and $d_1 d_{2,2} \ldots, d_{2,r} c \da_A \A$.
Let $d_{2,1}' \in \A$ be such that $Lt(d_{2,1}'/A)=Lt(d_{2,1}/A)$.
Now, there is some $f \in \textrm{Aut}(\M/A d_1 d_{2,2} \ldots, d_{2,r} c)$ such that $f(d_{2,1})=d_{2,1}'$.
For $1 \le i \le l$, denote $d_{2,r+i}'=f(d_{2,r+i})$.
Then, we have
\begin{displaymath}
c \in \textrm{bcl} (c_1, \ldots, c_k, d_{2,2}, \ldots, d_{2,r}, d_{2,r+1}' \ldots, d_{2,r+l}', e, d_{2,1}', a).
\end{displaymath}

We now repeat the above argument with $d_{2,2}$ in place of $d_{2,1}$.
When choosing a finite set $A' \subset \A$ such that $d_{2,2} \da_{A'} d_1 d_{2,3} \ldots, d_{2,r} c$ and $d_1 d_{2,3} \ldots, d_{2,r} c \da_{A'} \A$, we take care that $a,d_{2,1}' \in A$.
After doing the argument $r$ times, we have obtained elements $d_{2,r+1}^* \ldots, d_{2,r+l}'^* \in \textrm{bcl}_\A(d_1)$ such that 
\begin{displaymath}
c \in \textrm{bcl}_\A(c_1, \ldots, c_k, d_{2,r+1}^* \ldots, d_{2,r+l}^*, e).
\end{displaymath} 
If $d_{2,r+1}^* \ldots, d_{2,r+l}^* \in E$, we are done.

If not, there are some numbers $0<n \le m \le l$ such that (after renaming the elements in $\{d_{2, r+1}^{*}, \ldots, d_{2, r+l}^{*} \} \setminus (E \cup \A)$)  we have $U(d_{2,1}^{*}, \ldots, d_{2,n}^{*}, d_{2,n+1}^{*}, \ldots, d_{2,m}^{*}/E \cup \A)=m$, $U(d_{2,1}^{*}, \ldots, d_{2,n}^{*}/\A d_2)=n$, and $d_{2,n+1}^{*}, \ldots, d_{2,m}^{*} \in \textrm{bcl}_\A(d_2, d_{2,1}^{*}, \ldots, d_{2,n}^{*})$.
As $d_{2,1}^{*}, \ldots, d_{2,n}^{*} \in \textrm{bcl}_\A(d_1)$, we have $n \le U(d_1/E \cup \A)$.
Thus
\begin{displaymath}
U(c/\A d_2) \ge U(c/\A,d_2, d_{2,1}^{*}, \ldots, d_{2,n}^{*}) \ge U(c/ \A, d_1,d_2)=U(c/ \A d_2),
\end{displaymath}
so 
$$c \da_{ \A d_2} d_{2,1}^{*} \ldots d_{2,n}^{*},$$
and thus e.g. $d_{2,1}^{*} \da_{ \A, d_2, d_{2,2}^{*}, \ldots, d_{2,n}^{*}} c.$
Hence we may move $d_{2,1}^{*}, \ldots, d_{2,n}^{*}$ to $E$ with the same process as before with $d_2$ in place of $d_1$.
We keep repeating the process, and as at every step we move one element to $E$, we will eventually have moved them all, so we get $c \in \textrm{bcl}_\A (c_1, \ldots, c_k, E)$ as wanted. 

\end{proof}

From now on, we will simplify the notation by assuming that the elements of $\A$ are symbols in our language.

Let $a=f_{1,2}$, $c=f_{2, 2}$, $y=f_{1,3}$, $z=f_{2,3}$.
We will find elements $x$ and $b$ so that $\{a,b,c,x,y,z\}$ will form a group configuration.

Let $d=(f_{3, 2}, f_{3,3})$.
One can compute using the type $\alpha$ of the array that 
\begin{displaymath}
U(d/a y)=U(d/c z)=U(d/a c z y)=1.
\end{displaymath}
Thus, by Claim \ref{claimura}, there exists $x \in \textrm{bcl}(a y) \cap \textrm{bcl}(c z)$ such that $U(d/x)=1$.
We prove that $U(x)=1$.
We have $U(x) \geq 1$ since $U(d)=2$.
Now
\begin{displaymath}
3-U(x)=U(a y c z/x) \leq U(a y/x)+ U(c z/ x)=U(a y)+U(c z) -2U(x)=4- 2 U(x),
\end{displaymath}
where we use the type of the array and the fact that $x \in \textrm{bcl}(a y) \cap \textrm{bcl}(c z)$.
Thus, $U(x) \leq 1$.

Let $a'=f_{1, 1}$, $c'=f_{2,1}$.
By the type of the array, 
\begin{displaymath}
U(yz/a c)= U(yz /a' c')=U(yz /a c a' c') =1.
\end{displaymath}
By Claim \ref{claimura}, there exists $b \in \textrm{bcl}(a c) \cap \textrm{bcl}(a' c')$ such that $U(yz/b)=1$.
We prove that $U(b)=1$.
By the type of the array, $U(yz)=2$, and thus we must have $U(b) \ge 1$.
On the other hand, we have
\begin{displaymath}
3-U(b)=U(a ca'c'/b) \leq U(ac)+U(a'c')-2U(b)=4-2U(b),
\end{displaymath}
so $U(b) \leq 1$.

It is clear from the type of the array that 
\begin{displaymath}
U(z)=U(y)=U(c)=U(a)=1,
\end{displaymath}
and
\begin{displaymath}
U(z, y)=U(a, c) =U(a, y) =U(c, z)=2.
\end{displaymath}
Also, 
\begin{displaymath}
U(a, b, c)=U(a, y, x)=U(z, y, b)=U(z, c, x)=2,
\end{displaymath}
and
\begin{displaymath}
U(z,x,y,a,b,c)=U(z,y, a,c)=3
\end{displaymath}
by the type of the array and the choice of $x$ and $y$.
Thus, we are left to prove that the rest of the pairs have $U$-rank $2$ and that the rest of the triples have $U$-rank $3$.

We prove first that $U(a, c, y)=U(a, c, z)=3$ (and it of course follows that $U(y,c)=U(z, a)=2$).
Suppose that $y \in \textrm{bcl}(a, c)$.
Consider the concatenated sequence $(f_{i,2}f_{i,3})_{i < \omega_1}$.
Now, there is some stationary set $S \subseteq \omega_1$ and some model $\B$ such that the sequence $(f_{i,2}f_{i,3})_{i \in S}$ is Morley over $\B$.
Let $j, k \in S$ be such that $j<k$.
Since the sequence $(f_{i,2}f_{i,3})_{i < \omega_1}$ is order indiscernible, there is some automorphism $g$ of $\M$ such that $g(f_{1,2} f_{1,3})=f_{j,2} f_{j,3}$ and $g(f_{2,2} f_{2,3})=f_{k,2} f_{k,3}$.
Since $(f_{i,2}f_{i,3})_{i \in S}$ is Morley over $\B$, there is an automorphism $\pi \in \textrm{Aut}(\M/\B)$ such that
$\pi(f_{j,2} f_{j,3})=f_{k,2} f_{k,3}$ and $\pi(f_{k,2} f_{k,3})=f_{j,2} f_{j,3}$.
The map $g^{-1} \circ \pi \circ g$ is an automorphism taking  $f_{1,2} f_{1,3} \mapsto f_{2,2} f_{2,3}$,  and $f_{2,2} f_{2,3} \mapsto f_{1,2} f_{1,3}$.
Hence
\begin{displaymath}
t(f_{1,3}f_{1,2} f_{2,2}/\emptyset)=t(f_{2,3}f_{2,2} f_{1,2}/\emptyset).
\end{displaymath}
So, $z \in \textrm{bcl}(a, c)$ and $U( a, c,y,z)=2$ which is a contradiction (by the type of the array it should be $3$).
One proves similarly that $z \notin \textrm{bcl}(a, c)$.

Now we prove $U(c, y, z)=3$.
Suppose that $y \in \textrm{bcl}(z, c)$.
Considering the concatenated sequence $(f_{1,j}f_{2,j})_{j < \omega_1}$, we notice that there is an automorphism mapping $f_{1,2} f_{2,2} \mapsto f_{1,3} f_{2,3}$ and $f_{1,3} f_{2,3} \mapsto f_{1,2} f_{2,2}$, and thus
\begin{displaymath}
t(f_{1,2}f_{2,2}f_{2,3}/\emptyset)=t(f_{1,3}f_{2,3} f_{2,2}/\emptyset).
\end{displaymath}
Hence $a \in \textrm{bcl}(c, z)$ and we get a contradiction.
Similarly, $U(a, y, z)=3$.

Now we prove $U(x, z)=2$.
Suppose not.
Then, $x \in \textrm{bcl}(z).$
We chose $x$ so that $U(d/x)=1$.
As $U(d)=2$, we have $U(d/z)=1$ and thus $U(d, z)=2$.
So, $z \in \textrm{bcl}(d)= \textrm{bcl}(f_{3,2}, f_{3,3})$.
By the indiscernibility of the array, $y \in \textrm{bcl}(c, z)$, and we already proved this is not the case.
Similarly, $U(x, y)=2$.

Next we prove that $U(x,y,z)=3$.
If not, then $z \in \textrm{bcl}(x, y) \subseteq \textrm{bcl}(y, a)$, and we already proved this is not the case.
By similar arguments, $U(z, x, a)=U(x, y, c)=U(x,a,c)=3$, and it follows that $U(x,a)=U(x,c)=2$.

Now we prove that
\begin{displaymath}
U(a, b)=U(b,c)=U(z,b)=U(y,b)=2.
\end{displaymath}
If $U(a, b) \neq 2$, we would have $a \in \textrm{bcl}(b) \subseteq \textrm{bcl}(z,y)$ (note that $U(z,y,b)=2$), contradicting the fact that $U(a, y,z)=3$.
The fact that the $U$-rank of the three other pairs is also $2$ is proved similarly.

We are left to prove that the rest of the triples have $U$-rank $3$.
Suppose $U(a, b,z)=2$.
Then $z \in \textrm{bcl}(a,b) \subseteq \textrm{bcl}(a,c)$, and we have already proved this is not the case.
One proves similarly that $U(z, b, c)=U(y, a, b)=U(y, b,c)=3$.

Suppose $U(b,z,x)=2$.
Then, $b \in \textrm{bcl}(x,z) \subseteq \textrm{bcl}(c,z)$ which is again a contradiction.
Thus, $U(b,x,z)=3$ and it follows that $U(x,b)=2$.
One proves similarly that $U(x, a, b)=U(x, y, b)=U(x, b,c)=3$.

Thus, we have obtained a group configuration over $\A$.
There is some finite $A \subset \A$ so that the configuration is over $A$.
Hence, we may apply Theorem \ref{groupmain} to see that the group exists.
\end{proof}

\section{Finding the group}

In this section we will prove that the group exists.
Before doing it, we still need to present the following technical lemma.

\begin{lemma}\label{tekn1}
Let $(A_{ij}: 1 \le i \le M, 1 \le j \le N)$, $M, N \ge 2$, be a subarray of an indiscernible array of size $\o_1 \times \o_1$ over some finite tuple $b$.  
Assume $U(A;m,n/b)=m+n$ for any $m \le M$, $n \le N$, and that $\textrm{dim}_{\textrm{bcl}_b}(\textrm{dcl}(A_{12} A_{22}b) \cap \textrm{dcl}(A_{11} A_{21}b))=2$.
 
Let $b(A_{ij}) \to b(a_{ij})$ be a rank-$1$ specialization.
Suppose $Lt(a_{ij}/ b)$ is constant with $i$, $j$,   $U(a_{ij}/b)=1$ for each pair $i$, $j$, and $U(a;2,1/b)=2$.
Also assume $bA_{ij} A_{i'j} \to ba_{ij} a_{i'j}$ is strongly good for any $i$, $i'$, $j$.
Then $a$ is a rank-indiscernible array of type $m+n-1$ over $b$.
\end{lemma}

\begin{proof}
To simplify notation, we assume $b=\emptyset$.
All the arguments are similar in the more general case.

We prove the lemma as a series of auxiliary claims.

\begin{claim}\label{apu}  
Let $c_{ij}$ ($i,j=1,2$) be a $2 \times 2$ -rectangle from $a$. Assume
\begin{displaymath}
(A_{ij} : i,j=1,2) \to (c_{ij} : i,j=1,2)
\end{displaymath}
is a rank-$1$ specialization.  
Then, $U(c_{21} c_{22}/c_{11} c_{12}) \le 1$.
\end{claim}

\begin{proof}
Suppose for the sake of contradiction that 
\begin{equation}\label{vol}
U(c_{21} c_{22}/c_{11} c_{12})=2.
\end{equation}
By our assumptions on $a$ and $A$, we have $U(c_{11} c_{21})=2$ and $U(c_{12} c_{22})=2$, and $U((A_{ij} : i,j=1,2))=4$.
Thus, as we have a rank-$1$ specialization, $U((c_{ij} : i,j=1,2))=3$, so $U(c_{12} c_{22}/c_{11} c_{21})=1$.

By our assumptions, $c_{11} \da c_{21}$ and $Lt(c_{11})=Lt(c_{21})$.
We claim that $Lt(c_{11} c_{21})=Lt(c_{21} c_{11})$
Indeed, there is some strong automorphism $f_1$ such that $f_1(c_{11})=c_{21}$.
Let $c_{21}'=f_1(c_{21})$.
We have $c_{21}' \da c_{21}$, and thus $Lt(c_{21}'/c_{21})=Lt(c_{11}/c_{21})$, so there is some strong automorphism $f_2$ such that $f_2(c_{21}'c_{21})=c_{11} c_{21}$.
Then, $f=f_2 \circ f_1$ is a strong automorphism, and  $f(c_{11} c_{21})=c_{21} c_{11}$.
    
Denote $c_{22}'=f(c_{22})$ and $c_{12}'=f(c_{12})$.
Let $c_{10}$ and $c_{20}$  be such that $Lt(c_{10}c_{20}/c_{11}c_{21})=Lt(c_{22}'c_{12}'/c_{11}c_{21})$ and 
$$c_{10} c_{20} \da_{c_{11} c_{21}} c_{12} c_{22}.$$ 
Then, there is some strong automorphism $g$ such that $g(c_{22}',c_{12}',c_{11},c_{21})=(c_{10},c_{20},c_{11},c_{21})$.
Now $g \circ f (c_{22},c_{12},c_{21},c_{11})=(c_{10},c_{20},c_{11},c_{21})$, and thus
$$Lt(c_{10}c_{20}c_{11}c_{21}/\emptyset)=Lt(c_{22}c_{12}c_{21}c_{11}/\emptyset).$$
 By ZL1, we have 
\begin{displaymath}
c_{21} c_{11} c_{22} c_{12} \to c_{11} c_{21} c_{10} c_{20}.
\end{displaymath}
Since $A$ is an indiscernible array, we have
$t(A_{11} A_{21} A_{12} A_{22})=t(A_{21} A_{11} A_{22} A_{12})$ (see the proof of Lemma \ref{groupexists}), so
\begin{displaymath}
A_{11} A_{21} A_{12} A_{22} \to A_{21} A_{11} A_{22} A_{12}.
\end{displaymath}
By our assumptions, we also have
\begin{displaymath}
A_{21} A_{11} A_{22} A_{12} \to c_{21} c_{11} c_{22} c_{12}.
\end{displaymath}
By composing these three specializations, we get
\begin{equation}\label{spec1}
A_{11} A_{21} A_{12} A_{22} \to c_{11} c_{21} c_{10} c_{20}.
\end{equation}
Of course, we also have
\begin{equation}\label{spec2}
A_{11} A_{21} A_{12} A_{22} \to c_{11} c_{21} c_{12} c_{22}.
\end{equation}
By our assumptions, $A_{11} A_{21} \to  c_{11} c_{21}$ is a strongly good specialization, and since $U(c_{11} c_{21})=2$ and $U(A; 2,1)=3$, it has rank $1$.
Applying ZL7 to the specializations \ref{spec1} and \ref{spec2}, we find $A_{10}$ and $A_{20}$ with
\begin{displaymath}
t(A_{10} A_{20}/A_{11} A_{21})=t(A_{12} A_{22}/A_{11} A_{21})
\end{displaymath}
and
\begin{displaymath}
A_{10} A_{20} \downarrow_{A_{11} A_{21}} A_{12} A_{22}
\end{displaymath}
such that $(A_{ij} : i=1,2, j=0,1,2) \to (c_{ij} : i=1,2, j=0,1,2)$, and in particular
\begin{equation}\label{spec3}
(A_{ij} : i=1,2, j=0,2) \to (c_{ij} : i=1,2, j=0,2).
\end{equation}

Next, we prove that the specialization \ref{spec3} is actually an isomorphism.
This will lead to a contradiction, since $U(A_{12} A_{22})=3$ but $U(c_{12} c_{22}) \le 2$.

We prove first that $U(A_{ij} : i=1,2, j=0,2)=4$.
Denote 
\begin{displaymath}
X= \textrm{dcl}(A_{12} A_{22}) \cap \textrm{dcl}(A_{11} A_{21}).
\end{displaymath}
By our assumptions, $dim_{\textrm{bcl}}(X)=2$, and thus $U(A_{12} A_{22} /X)=1$.
As $t(A_{10} A_{20}/A_{11} A_{21})=t(A_{12} A_{22}/A_{11} A_{21})$ and $X \subseteq \textrm{dcl}(A_{11} A_{21})$, we have also $U(A_{10} A_{20}/X)=1$.
Moreover, $U(A_{10} A_{20}/A_{11} A_{21})=U(A_{12} A_{22} /A_{11} A_{21})=1$, so
$$A_{10} A_{20} \downarrow_{X} A_{11} A_{21}.$$
On the other hand,
$$A_{10} A_{20} \da_{A_{11} A_{21}} A_{12} A_{22},$$ 
so by transitivity,
\begin{displaymath}
A_{10} A_{20} \downarrow_{X} A_{11} A_{21} A_{22} A_{12},
\end{displaymath}
and therefore $U(A_{10} A_{20} /A_{22} A_{12})=U(A_{10} A_{20} /X)=1$, so $U(A_{ij} : i=1,2, j=0,2)=4$.
 
To get the contradiction, we have to prove that also $U(c_{ij} : i=1,2, j=0,2)=4$.
Now, we have chosen $c_{10}$ and $c_{20}$ so that $c_{10} c_{20} \downarrow_{c_{11} c_{21}} c_{12} c_{22}$.
By the counterassumption (\ref{vol}), we have $U(c_{21}/c_{11} c_{12} c_{22})=1=U(c_{21})$, so
\begin{displaymath}
c_{12} c_{22} \downarrow_{c_{11}} c_{21},
\end{displaymath}
and thus,
\begin{equation}\label{fork2}
c_{10} c_{20} \downarrow_{c_{11}} c_{12} c_{22}.
\end{equation}
From (\ref{vol}), it follows that $U(c_{21}/c_{12})=1$.
Thus, $c_{21} \da c_{12}$, so by symmetry $c_{12} \da c_{21}$ and $U(c_{12}/c_{21})=1$.
As $Lt(c_{11} c_{21} c_{10} c_{20})=Lt(c_{21} c_{11} c_{22} c_{12})$, we have $U(c_{20}/c_{11})=U(c_{12}/c_{21})=1$. 
Thus $c_{20} \downarrow c_{11}$, and from this together with (\ref{fork2}) it follows that $U(c_{20}/c_{12} c_{22})=1$.
As $U(c_{12} c_{22})=2$, we are left to prove that $U(c_{10}/c_{20} c_{12} c_{22})=1$, i.e. $c_{10} \downarrow c_{12} c_{22} c_{20}$.
From (\ref{fork2}) it follows that
\begin{equation}\label{fork3}
c_{10} \downarrow_{c_{11} c_{20}} c_{12} c_{22}.
\end{equation}
Now 
\begin{equation}\label{ranks}
U(c_{10}/c_{20} c_{11} c_{21})=U(c_{22}/c_{12} c_{21} c_{11})=1,
\end{equation}
where the second equality follows from (\ref{vol}), so 
\begin{equation}\label{fork4}
c_{10} \downarrow_{c_{11}} c_{20}.
\end{equation}
From (\ref{fork3}) and (\ref{fork4}) it follows by transitivity that
\begin{displaymath}
c_{10} \downarrow_{c_{11}} c_{12} c_{22} c_{20}.
\end{displaymath}
But now by (\ref{ranks}), $c_{11} \downarrow c_{10}$, and hence $c_{10} \downarrow c_{12} c_{22} c_{20}$ as wanted. 
Thus, $U(c_{ij} : i=1,2, j=0,2)=4$, so we get the contradiction and the claim is proved.
\end{proof}

\begin{claim}\label{iclaim}
For any set $*$ of $i$-indices and $j \ge 2$, $U(a_{*,j}/a_{*,<j}) \le 1$.
\end{claim}
 
\begin{proof}
Suppose not.
Then, there are some indices $i \neq i'$ in the index set $*$ such that $U(a_{ij} a_{i' j}/a_{*, <j})=2$.
In particular, $U(a_{ij} a_{i' j}/a_{i1} a_{i'1})=2$.
Since $U(a; 2, 1)=2$, we have $U(a_{ij} a_{i'j} a_{i1} a_{i' 1})=4=U(A;2,2)$.
Thus, $A_{ij} A_{ij'} A_{i1} A_{i' 1} \to a_{ij} a_{ij'} a_{i1} a_{i' 1}$ is an isomorphism.
But this contradicts the fact that $U(A;2,1)=3.$
\end{proof}

\begin{claim}\label{jclaim}
For any set $*$ of $j$-indices and $i \ge 2$, $U(a_{i,*}/a_{<i,*}) \le 1$.
\end{claim}

\begin{proof}
Suppose otherwise.
Then there exist $j < j'$ such that  $U(a_{ij} a_{ij'}/ a_{1j} a_{1j'})=2$.
Write $c_{11}=a_{1j}$, $c_{12}=a_{1j'}$, $c_{21}=a_{ij}$ and $c_{22}=a_{ij'}$.
As $U(c_{11})=1$,  we have $U(c_{21} c_{22}c_{11} c_{12}) \ge 3$.
Thus, $(A_{ij} : i,j=1,2) \to (c_{ij}: i,j=1,2)$ is either an isomorphism or a rank-$1$ specialization.
In the first case we get a contradiction because $U(c_{11} c_{21})=2$ but $U(A_{11} A_{21})=3$.
In the second case, Claim \ref{apu} gives us $U(c_{21} c_{22}/ c_{11} c_{12}) \le 1$ which is also a contradiction.
\end{proof}
 
Using the Claims \ref{iclaim} and \ref{jclaim} one proves by induction on $m$ and $n$ that any $m \times n$ -rectangle from $a$ has $U$-rank at most $m+n-1$.  
Suppose, for the sake of contradiction, that the inequality is strict for some $m \times n$ -rectangle.
Then, by Claim \ref{iclaim} the inequality remains strict for the $m \times (n+1)$ rectangle obtained by adjoining an   $m \times 1$ -array. Similarly, by Claim \ref{jclaim} it remains strict for the $(m+1) \times (n+1)$ -rectangle obtained by further adjoining a $1 \times (n+1)$ array.
Continuing this way one finds that the inequality is strict for $m=M$, $n=N$.
But we assumed $U(A;M,N)=M+N$ and that the specialization $A \to a$ has rank $1$, so this is a contradiction.
\end{proof}

Now we are ready to present our main theorem.

\begin{theorem}\label{grouptheorem}
Let $\M$ be a Zariski-like structure with a non-trivial pregeometry.
Then, there exists a Galois definable one-dimensional group in $(\M^{eq})^{eq}$.
\end{theorem}

\begin{proof}  
If $\M$ is locally modular, then the theorem follows from Lemma \ref{locmodgroup}.

So suppose $\M$ is non locally modular.
By Lemma \ref{relevantfamily}, there exists a relevant family of plane curves that has rank $r \ge 2$.
Let $\alpha$ be the canonical parameter for one of the curves in this family, and suppose $U(\alpha)=r$.
Let $(x,y)$ be a generic point on this curve, i.e. $\alpha=\textrm{Cb}(x,y/\alpha)$. 
Since the family is relevant, we have $x \da \alpha$, $y \da \alpha$.
We also have $x \da y$, because otherwise $1=U(xy/\emptyset)=U(xy/\alpha)$, so $xy \da_\emptyset \alpha$, which would imply $\alpha \in \textrm{bcl}(\emptyset)$ since $\alpha$ is a canonical parameter.
 
Let $c_1, \ldots, c_r, d_1, \ldots, d_r$  be such that $t(c_i,d_i/\alpha)=t(x,y/\alpha)$ for $1 \le i \le r$,
and the sequence $x, c_1, \ldots, c_r$ is independent over $\alpha$.
We claim that $U(c_1,d_1, \ldots, c_r,d_r)=2r$.
For this, we first show that $\alpha \in \textrm{bcl}(c_1,d_1, \ldots, c_r,d_r)$.
 We have $c_1 \in \textrm{bcl}(d_1, \alpha)$, so $c_1 \nda_\alpha d_1$, and thus
$$U(c_1/d_1)=U(c_1)=U(c_1/\alpha)>U(c_1/d_1\alpha),$$
so $c_1 \nda_{d_1} \alpha$.
Hence, $U(\alpha/d_1)>U(\alpha/c_1d_1)$.
Thus, $U(\alpha/ c_1 d_1) \le r-1$.
 
Suppose now $0<U(\alpha/c_1 d_1, \ldots, c_k d_k) \le r-k$ for some $k<r$.
We claim that $U(\alpha/c_1 d_1, \ldots, c_k d_k, c_{k+1} d_{k+1}) \le r-k-1$.
Suppose towards a contradiction that
$$U(\alpha/c_1 d_1, \ldots, c_k d_k)=U(\alpha/c_1 d_1, \ldots, c_k d_k, c_{k+1} d_{k+1}).$$
Then, 
\begin{eqnarray}\label{rairai}
c_{k+1} d_{k+1} \da_{c_1 d_1, \ldots, c_k d_k}  \alpha.
\end{eqnarray}
We have
\begin{eqnarray}\label{rairai2}
\alpha=\textrm{Cb}(c_{k+1} d_{k+1}/\alpha)=\textrm{Cb}(c_{k+1} d_{k+1}/\alpha, c_1 d_1, \ldots, c_k d_k),
\end{eqnarray}
where the second equality follows from the fact that $c_{k+1} d_{k+1}Ê\da_\alpha c_1 d_1, \ldots, c_k d_k.$
From (\ref{rairai}) and and (\ref{rairai2}), it follows that $\alpha \in \textrm{bcl}( c_1 d_1, \ldots, c_k d_k),$ a contradiction.

Thus, $\alpha \in \textrm{bcl}(c_1, d_1, \ldots, c_r,d_r)$.
On the other hand,
$$U(\alpha)+U(c_1 d_1/\alpha)+ \ldots + U(c_r d_r/\alpha)=r+r=2r,$$
so $U(c_1, d_1, \ldots, c_r,d_r)=2r$. 

Next, we show that for $1 \le k \le r$, $U(\alpha/c_1, \ldots, c_r, d_1, \ldots, d_k)=r-k$.
Indeed, 
\begin{eqnarray*}
2r&=&U(\alpha, c_1, \ldots, c_r, d_1, \ldots, d_k)=U(c_1, \ldots, c_r, d_1, \ldots, d_k)+U(\alpha/c_1, \ldots, c_r, d_1, \ldots, d_k) \\ &=&(r+k)+U(\alpha/c_1, \ldots, c_r, d_1, \ldots, d_k). 
\end{eqnarray*} 

Let now $C$ be the locus of $(x,y,c_1, \ldots, c_r, d_1, \ldots, d_r)$ and $E$ the locus of $(c_1, \ldots, c_r, d_1, \ldots, d_r)$.
Then, $C$ is a family of plane curves parametrized by $E$, and 
 $C(c_1, \ldots, c_r, d_1, \ldots, d_r)$ is a curve in this family.
Denote $d=(c_1, \ldots, c_r, d_1, \ldots, d_{r-2})$ and $e_0=(d_{r-1}, d_r)$.
Since $xy \da_\alpha c_i d_i$ for each $i$, we have $\alpha=\textrm{Cb}(x,y/c_1, \ldots, c_n, d_1, \ldots, d_n)$. 
It is interbounded with $e_0$ over $d$, and $U(e_0/d)=2$.

Let $e \in E(d)$ be a generic element.
We now write  $C(e; a,b)$ for "$(a,b)$ is a generic point of $C(ed)$". 
We write $C^2(e; ab, a'b')$ if the following hold:
\begin{enumerate}
\item $C(e; a, b)$ and $C(e; a', b')$; 
\item $ab \downarrow_{de} a' b'$; 
\item $Lt(ab/de)=Lt(a'b'/de).$  
\end{enumerate}

\begin{claim}\label{aclspec}
\begin{enumerate}[(i)]
\item If $a \neq a'$ and $b \neq b'$, then $C^2(e; ab, a'b')$ implies that $U(aba'b'/d)=4$.
\item If $a \neq a'$ and $b \neq b'$, and $C^2(e; ab, a'b')$, then  $deaba' b' \to deabab$. 
\end{enumerate}
\end{claim}

\begin{proof}
Since $e$ is interbounded with $\textrm{Cb}(a,b/de)$, we see using similar arguments as above, that (i) holds.
 
 For (ii), we will apply (ZL8).
Let $\A$ be a model such that $de \in \A$ and $aba'b' \da_{de} \A$.
Then, $ab \da_{dea'b'} \A$, and since $ab \da_{de} a'b'$, we get by transitivity that $ab \da_{de} a'b' \A$, which implies $ab \da_\A a'b'$.
On the other hand, we have $Lt(ab/de)=Lt(a'b'/de)$, $ab \da_{de} \A$ and $a'b' \da_{de} \A$, so $Lt(ab/\A)=Lt(a'b'/\A)$.
Thus, we may extend $(ab, a'b')$ to a Morley sequence over $\A$.
It follows that $a'b'$ and $ab$ are strongly indiscernible over $de$.
Of course also $ab$ and $ab$ are strongly indiscernible over $de$ (just repeat $ab$ arbitrarily many times to extend the sequence).
Clearly $dea'b' \to deab$ and $deab \to deab$, $\textrm{rk}(de \to de)=0 \le 1$, and $de \to de$ is strongly good.   
Hence, we may apply (ZL8) to get $deaba' b' \to deabab$.
\end{proof}
 
Pick some generic point $e \in E(d)$, and independent generics $a_0, b_0, a, b \in \M$ such that $C^{2}(e; ab, a_0 b_0)$. 
Let $\kappa$ be some cardinal large enough and let $a_i, b_j$, $i,j <\kappa$ be a sequence of  generic elements of $\M$ independent over $d$ such that $Lt(a_i b_j /d a_0 b_0)=Lt(ab/ d a_0 b_0)$ for all $i,j$.  
For each pair $i,j$, let $f_{ij}$ be an automorphism fixing $a_0, b_0, d$ such that $f_{ij}(a,b)=(a_i,b_j)$.
Denote $e_{ij}=f_{ij}(e)$.
Then, $C^2(e_{ij}; a_i b_j, a_0 b_0)$ holds  for each pair $i,j$.
Let $A_{ij}=(a_i, b_j, e_{ij})$, $A=(A_{ij})_{i,j \ge 1}$.
We will next show that if we choose $\kappa$ to be large enough, then we can find an indiscernible array of size $\o_1 \times \o_1$ such that each one of its finite subarrays is isomorphic to some finite subarray of $A$.
 
Let  $\lambda< \kappa$ be a cardinal large enough (but not too large) for the argument that follows. 
For each $i< \kappa$, denote $A_{i, <\lambda}=(A_{ij} | j< \lambda)$.
Using Erd\"os-Rado and an Ehrenfeucht-Mostowski construction, one finds a  sequence $(A'_{i, < \lambda})_{i<\o_1}$ such that every finite permutation of the sequence preserving the order of the indices $i$ extends to 
some $f \in \textrm{Aut}(\M/d a_0 b_0)$.
Moreover, an isomorphic copy of every finite subsequence can be found in the original sequence $(A_{i, < \lambda})_{i<\kappa}$.
This construction is due to Shelah, and the details can be found in e.g. \cite{CatTran}, Proposition 2.13.
There it is done for a sequence of finite tuples (whereas we have a sequence of sequences of length $\lambda$), but the proof is similar in our case.
 
We may now without loss  assume that $(A_{i,< \lambda}')_{i<\o_1}$ are the $\o_1$ first elements in the sequence 
$(A_{i,< \lambda})_{i<\kappa}$.
Since we have chosen $\lambda$ to be large enough, we may now apply the same argument to $(A_{<\o_1, j}')_{j<\lambda}$ to obtain an array $(A_{<\o_1, j}'')_{j<\o_1}$. 
This is an array of size $\o_1 \times \o_1$, indiscernible over $d a_0 b_0$, and we may assume it is a subarray of the original array $A$.
From now on, we will use $A$ to denote this indiscernible array of size $\o_1 \times \o_1$.
 
We write $x \to *y$ for $(x, d, a_0, b_0) \to (y, d, a_0, b_0)$.

\begin{claim}\label{ekac}
Let $A'_{ij}=A_{i1}$ for $j \ge 1$.
Then, $A \to *A'$.
\end{claim} 
\begin{proof}
For each $i<\o_1$, consider the the sequence $(A_{ij})_{j<\o_1}$.
Now, there is some cofinal set $X_i \subset \o_1$ such that $(A_{ij})_{j \in X_i}$ is Morley, and thus strongly indiscernible, over $d a_0 b_0$.
For each $j$, we have $A_{ij} d a_0 b_0 \to A_{i1} d a_0 b_0$.
Moreover, $\textrm{rk}(d a_0 b_0 \to d a_0 b_0)=0 \le 1$ and $d a_0 b_0 \to d a_0 b_0$ is a strongly good specialization.
Also, $(A'_{ij})_{j \in X_i}$ is strongly indiscernible (since it just repeats the same entry).
Thus, by (ZL8), there is, for each $i$, a specialization 
$(A_{ij})_{j \in X_i} \to * (A'_{ij})_{j \in X_i}.$
If we enumerate the set $X_i$ again, using the order type of $\o_1$, then we get (we still use the notation with the index set $X_i$ to denote that the sequence so indexed is the Morley one)
$$(A_{ij})_{j \in \o_1}\to *(A_{ij})_{j \in X_i} \to *(A'_{ij})_{j \in X_i}\to *(A_{ij}')_{j \in \o_1},$$
so in particular $(A_{ij})_{j \in \o_1} \to *(A'_{ij})_{j \in \o_1}$.

To prove that $A \to *A'$, it suffices to show $(A_{ij})_{i<\o_1, j \in J} \to *(A_{ij}')_{i<\o_1, j \in J}$ for all finite $J \subset \o_1$.
So, let $J \subset \o_1$ be finite.
Since $(A_{ij})_{j \in \o_1} \to * (A'_{ij})_{j \in \o_1}$ holds for every $i$, we have  $(A_{ij})_{ j \in J} \to *(A_{ij}')_{ j \in J}$ for every $i \in \o_1$.
Thus, applying (ZL8) similarly as we did above, we obtain $(A_{ij})_{i<\o_1, j \in J} \to *(A_{ij}')_{i<\o_1, j \in J}$, as wanted.
It then follows that $A \to *A'$.  
\end{proof} 

\begin{claim}\label{tokac}
Let $A''_{ij}=(a_0, b_0, e_{i1})$.
Then, $A' \to *A''$.
\end{claim}
\begin{proof} 
As $A'_{ij}$ and $A''_{ij}$ do not depend on $j$ and as specializations respect repeated entries, it suffices to show that $(A'_{i1}:i) \to *(A''_{i1}:i)$.
By Claim \ref{aclspec} (ii), $d a_0 b_0 b_1 \to d a_0 b_0 b_0$.
By Claim \ref{aclspec} (i), we have $U(a_0 b_0 a_1 b_1/d)=4$, so $(d, a_0, b_0, b_1)$ is a generic point of $\M^{2r+1}$, and this is a strongly good specialization.
It is also clearly of rank $1$.
By Claim \ref{aclspec} (ii), $(a_i, b_1, e_{i1}, d, a_0, b_0) \to (a_0, b_0, e_{i1}, d, a_0, b_0)$ for every given $i$.
Thus, we may apply (ZL8) similarly as in the proof of the previous claim.
\end{proof}

\begin{claim}\label{aputulos}
If $(i,j) \neq (i',j')$, then $e_{ij} \downarrow_{d a_0 b_0} e_{i'j'}$.
\end{claim}
\begin{proof}
Suppose not.
By the same arguments that we used to prove Claim \ref{aclspec} (i), $U(e_{ij}/d a_0 b_0)=1$, and $e_{ij} \in \textrm{bcl}(d a_0 b_0 a_i b_j)$.
From the first of these statements it follows that $e_{ij} \in \textrm{bcl}(d a_0 b_0 e_{i'j'})$, since $U(e_{ij}/da_0b_0e_{i'j'})<U(e_{ij}/da_0b_0)$ by the counterassumption.
From the second statement it follows that
$a_i b_j$ dominates $e_{ij}$ over $d a_0 b_0$.
Similarly, $a_{i'} b_{j'}$ dominates $e_{i'j'}$ over $d a_0 b_0$.
Suppose first $i \neq i'$ and $j \neq j'$.
Then, $a_i b_j \downarrow_{d a_0 b_0} a_{i'} b_{j'}$ (the sequence was chosen to consist of elements independent over $d$),  and by domination $e_{ij} \downarrow_{d a_0 b_0} e_{i'j'}$, a contradiction.
 
Suppose now $i=i'$ and $j \neq j'$ (the other case is symmetric).
Similarly as before, we get that $b_j$ dominates $e_{ij}$ over $d a_0 b_0 a_i$ and $b_{j'}$ dominates $e_{ij'}$ over $d a_0 b_0 a_i$.
As $b_j \downarrow_{d a_0 b_0 a_i} b_{j'}$, we get that $e_{ij} \downarrow_{d a_0 b_0 a_i} e_{ij'}$.
Thus, to get a contradiction it suffices to show that $a_i \downarrow_{a_0 b_0} e_{ij}$, since $e_{ij} \downarrow_{d a_0 b_0} e_{i'j'}$ then follows by transitivity.  
Suppose not.
As $U(e_{ij}/d a_0 b_0)=1$, we must now have $e_{ij} \in \textrm{bcl}(d a_0 b_0 a_i)$.
But then we have $b_j \in \textrm{bcl}(d a_i e_{ij}) \subseteq \textrm{bcl}(d a_0 b_0 a_i)$ which is a contradiction since the sequence $a_0, b_0, a_i, b_j$ was chosen to be independent over $d$.
\end{proof}
 
By claims (\ref{ekac}) and (\ref{tokac}), $A \to *A''$.
We will apply (ZL9) to this specialization and eventually obtain an infinite rank-indiscernible array $A^*$ such that
$A \to A^* \to A''$.
The array $A^*$ will be of type $m+n-1$ over the parameters $d a_0 b_0$, as desired.

Let now $A_0$ be a finite subarray of $A$ containing the entry $A_{11}$, and let $A''_0$ be the corresponding finite subarray of $A''$.
Then, there is a specialization $A_0 \to * A_0''$.
After suitably rearranging the indices, we may assume that the tuple on the left begins with ``$a_0 a_1 \ldots$", whereas the tuple on the right begins with ``$a_0 a_0 \ldots$".
By Remark \ref{ZL9implies}, the dimension theorem holds, and thus there is a finite array $A_0^*$
such that $d a_0 b_0 A_0 \to  d' a_0' b_0'  A^*_0 \to d a_0 b_0 A''_0$ for some $d', a_0', b_0'$, ${A^*_0}_{11}=a_0'b_1^*e_{11}^*$ for some $b_1^*$, $e_{11}^*$, and $U(A_0)-U(A_0^*) \le 1$.
In particular, we have $d a_0 b_0 \to  d' a_0' b_0' \to d a_0 b_0$.
By (ZL3) this implies that $t^g (d a_0 b_0/\emptyset)=t^g (d' a_0' b_0'/\emptyset)$.
Thus, we may assume that $d' a_0' b_0'=d a_0 b_0$ (if it is not, then just apply to the array $A^*_0$ an automorphism taking $d' a_0' b_0' \mapsto d a_0 b_0$).
In particular, we may assume $A_0 \to * A^*_0 \to *A''_0$ and ${A^*_0}_{11}=a_0b_1^*e_{11}^*$ for some $b_1^*$ and $e_{11}^*$.
We next show that  the assumptions of Lemma \ref{tekn1} hold for $A_0$ and $A_0^*$ over the parameters $d a_0 b_0$.
  
By the calculations made for Claim \ref{aclspec} (i), $e_{ij} \in \textrm{bcl}(d a_0 b_0 a_i b_j)$.
Thus, as the elements $a_i, b_j$ were chosen to be independent over $d$ for $i,j \ge 0$, we have $U(A; 1, 1/d a_0 b_0)=2$, and it is easy to show by induction that $U(A; m, n/ d a_0  b_0)=m+n$. 
Write $C=A_{11} A_{21}$ and $C'=A_{12} A_{22}$.
Now $U(C/d a_0 b_0)=U(C'/d a_0 b_0)=3$, and $U(C \cap C'/d a_0 b_0)=2$.
Thus, $2 \le \textrm{dim}_{\textrm{bcl}}(\textrm{dcl}(Cda_0b_0) \cap \textrm{dcl}(C'da_0b_0)/d a_0 b_0) \le 3$.
Denote $X=\textrm{dcl}(Cda_0b_0) \cap \textrm{dcl}(C'da_0b_0)$, and suppose $\textrm{dim}_{\textrm{bcl}}(X/d a_0 b_0)=3$.
Since $X \subseteq \textrm{bcl}(Cda_0b_0)$,  we must have $\textrm{bcl}(X)=\textrm{bcl}(Cda_0b_0)$.
But this is impossible since $b_1 \in \textrm{bcl}(Cda_0b_0) \setminus \textrm{bcl}(X)$.
Thus, $\textrm{dim}_{\textrm{bcl}}(X/d a_0 b_0)=2$.
 
Consider now $A^*_0$.
We will show that the assuptions posed for the array $a$ in the statement of Lemma \ref{tekn1} hold for $A_0^*$ over the parameters $d a_0 b_0$, and it will then follow that $A_0^*$ is of type $m+n-1$ over the parameters.
We prove first that for any indices $i,j$, ${A^{*}_0}_{ij}=(a_0 b_0 e_{ij}^{*})$ for some $e_{ij}^{*}$.

Denote  ${A^{*}_0}_{ij}=(a_i^* b_j^* e_{ij}^{*})$.
Since $A_0 \to *A^{*}_0$, we have $a_1 e_{1j} d \to a^{*}_1 e^{*}_{1j} d$ for each $j$, and thus $3=U(a_1 e_{1j}/d) \ge U(a^{*}_1 e^{*}_{1j}/d)$ for each $j$.  
On the other hand, we have $a_1^{*} b_j^{*} e^{*}_{1j}d a_0 b_0 \to a_0 b_0 e_{11} da_0 b_0$, and thus  $U(a^{*}_1 e_{1j}^{*}/d) \ge U(a_0 e_{11}/d)=3$, so $U(a^{*}_1 e_{1j}^{*}/d)=3$.
Similarly one shows that $U(a_1^{*} b_j^{*} e_{1j}^{*}/d)=3$, so $b_j^{*} \in \textrm{bcl}(da_1^{*} e_{1j}^{*})$.
Hence, 
$$U(a_1^{*}b_j^{*} e_{1j}^{*} a_0 b_0/d)=U(a_1^{*} e_{1j}^{*} a_0 b_0/d)=U(a_0 e_{1j}^{*} a_0 b_0/d)=U(e_{1j}^{*} a_0 b_0/d),$$
where the second equality follows from the fact that $a_1^*=a_0$ (this holds by the choice of our enumeration for the specialization). 
As we have  $e_{1j} da_0 b_0 \to e_{1j}^{*} da_0 b_0 \to e_{11}d a_0 b_0$, we get $U(e_{1j}^{*} a_0 b_0/d)=3$.
Thus, $a_1^{*} b_j^{*} e_{1j}^{*}d a_0 b_0 \to a_0 b_0 e_{11} da_0 b_0$ is an isomorphism for each $j$.
Hence, for each $j$, $b_j^{*}=b_0$.

In particular, $b_1^{*}=b_0$.
By applying similar arguments as above to the specialization $a_i^*b_1^*e_{i1}^*da_0 b_0 \to a_0 b_0 e_{i1} da_0 b_0$, we get that $U(a_i^*b_1^*e_{i1}^*a_0 b_0/d)=U(e_{i1}^*a_0b_0/d)=3$ for each $i$, so the specialization is an isomorphism for each $i$.
Thus, $a_i^*=a_0$ for each $i$.
 
We next show that $Lt({A_0^{*}}_{ij}/d a_0 b_0)$ does not depend on $i,j$ and has $U$-rank $1$.
The specialization $A^{*} \to *A''$ also gives 
\begin{equation}\label{eijei1}
e_{ij}^{*} e_{i'j'}^{*} \to *e_{i1} e_{i'1}.
\end{equation}
Suppose  $i \neq i'$. 
We have $U(e_{ij}^{*} e_{i'j'}^{*}/d a_0 b_0) \ge U(e_{i1} e_{i'1}/d a_0 b_0)$.
By Claim \ref{aputulos}, $U(e_{i1} e_{i'1}/d a_0 b_0)=2$. 
As $e_{ij} e_{i'j'} \to *e_{ij}^{*} e_{i'j'}^{*}$, we also have $U(e_{ij}^{*} e_{i'j'}^{*}/d a_0 b_0) \le 2$.
Thus, equality holds, and the specialization (\ref{eijei1}) is an isomorphism, so $t(e_{ij}^{*} e_{i'j'}^{*}/d a_0 b_0)=t(e_{i1}e_{i'1}/d a_0 b_0)$.  

We note that $Lt(e_{i1}/d a_0 b_0)=Lt(e_{i'1}/d a_0 b_0)$.
Indeed, there is some cofinal subset $X \subset \o_1$ such that $(A_{k1})_{k \in X}$ is Morley over $da_0b_0$.
We may without loss assume $i <i'$.
Let $k, k' \in X$ be such that $k<k'$.
As $(A_{k1})_{k \in X}$ is Morley over $da_0b_0$, we have $Lt(e_{k1}/da_0 b_0)=Lt(e_{k'1}/da_0 b_0)$.
By the indiscernibility of the array $A$ over $da_0b_0$, there is some automorphism fixing $da_0b_0$ and taking $(e_{i1}, e_{i'1})$ to $(e_{k1}, e_{k'1})$.
Thus, $Lt(e_{i1}/a_0 b_0)=Lt(e_{i'1}/a_0 b_0)$.
 
Hence $Lt(e_{ij}^{*}/da_0 b_0)=Lt(e_{i'j'}^{*}/da_0 b_0)$, and this of course remains true if $i=i'$.  
It follows that $Lt({A_0^{*}}_{ij}/da_0 b_0)$ does not depend on $i,j$. 
It has $U$-rank $1$ since $t(e_{ij}^{*}/da_0 b_0)=t(e_{i1}/da_0 b_0)$. 
   
From the above calculations we see that $U(A^{*}_0; 2, 1/d a_0 b_0)=2$.
Thus, the fact that $A^{*}_0$ is of type $m+n-1$ follows from Lemma \ref{tekn1} as soon as we verify that the specialization ${A_0}_{ij} {A_0}_{i'j} \to *{A_0^{*}}_{ij} {A_0^{*}}_{i'j}$ is strongly good for any $i, i', j$.
In other words, we must show that 
\begin{equation}\label{pitkaspec}
(d, a_0, b_0, a_{i}, a_{i'}, b_j, e_{ij}, e_{i'j}) \to (d, a_0, b_0, a_0, a_0, b_0, e_{ij}^{*}, e_{i'j}^{*})
\end{equation}
is strongly good.
Now $U(a_0 b_0 b_j e_{ij} e_{i'j}/d) \ge 5$ since $a_i \in \textrm{bcl}(b_j e_{ij}d)$ and $a_{i'} \in \textrm{bcl} (b_j e_{i'j}d)$, and the elements $a_0,b_0,a_i, a_{i'}, b_j$ form an independent sequence over $d$ and are each independent from $d$.
By Claim \ref{aputulos}, $U(a_0 b_0 e_{ij} e_{i'j}/d)=4$.
Thus, $b_j$ is independent from $(d, a_0, b_0, e_{ij}, e_{i'j})$ and
\begin{displaymath}
(a_0, b_0, b_j, e_{ij}, e_{i'j}) \to (a_0, b_0,  b_0, e_{ij}^{*}, e_{i'j}^{*})
\end{displaymath}
is strongly regular  because $(a_0, b_0, e_{ij}, e_{i'j}) \to (a_0, b_0, e_{ij}^{*}, e_{i'j}^{*})$ is an isomorphism and $b_j \to b_0$ is a strongly regular specialization ($b_j$ is a generic element of $\M$).
Also $(d,b_j, e_{ij}) \to (d,b_0, e_{ij}^{*})$ is an isomorphism, and $a_i \in \textrm{bcl}(d,b_j, e_{ij})$.
Thus,
\begin{displaymath}
(d,a_0, b_0, a_i, b_j, e_{ij}, e_{i'j}) \to (d,a_0, b_0, a_0, b_0, e_{ij}^{*}, e_{i'j}^{*})
\end{displaymath}
is strongly good.
Similarly, $(d, b_j, e_{i'j}) \to (d,b_0, e_{i'j}^{*})$ is an isomorphism, and $a_{i'} \in \textrm{bcl}(d,b_j, e_{i'j})$, so the specialization \ref{pitkaspec} is strongly good by the recursive definition.
 
Hence, by Lemma \ref{tekn1}, $A_0^*$ is of type $m+n-1$ over $da_0b_0$.

Next, we apply (ZL9) to the specialization $A \to *A''$ to eventually obtain an infinite indiscernible array of type $m+n-1$ over $da_0b_0$.
Enumerate the elements on the left side of the specialization so that $a_0$ is the element enumerated by $0$ and $a_1$ the element enumerated by $1$, and use a corresponding enumeration on the right side (there, both the element enumerated by $0$ and the element enumerated by $1$ will be $a_0$).
Let  $S$ be a collection of index sets corresponding to all $m \times n$ subarrays of $A$ containing the entry $A_{11}$ for all natural numbers $m,n$.
Moreover, we add $0$ to every $X \in S$.
The set $S$ is unbounded and directed, and by what we just proved,  every $X \in S$ corresponds to an array $A^*_X$ of type $m+n-1$ over $da_0b_0$ (we get the correspondence by removing the element indexed by $0$ from each $X$).
Thus, the conditions of (ZL9) hold for the set $S$, and hence we 
obtain an infinite array $A^*$ where each $m \times n$ -subarray containing the entry $A_{11}^*$ has $U$-rank $m+n-1$ over $da_0 b_0$.

We claim that $A^*$ is actually of type $m+n-1$ over $d a_0 b_0$.
To prove this, let $A_0^*$ be an arbitrary $m_0 \times n_0$ subarray of $A^*$.
Then, there is some $(m_0+1) \times (n_0+1)$ subarray $A_1^{*}$ of $A^*$ such that $A_1^*$ contains the entry $A_{11}^*$ and $A_0^*$ is a subarray of $A_1^*$.
We have already shown that $A_1^*$ is of type $m+n-1$ over $d a_0 b_0$.
Hence, $U(A_0^*/d a_0 b_0)=m_0+n_0-1$, as wanted.

If we have chosen the cardinals $\kappa$ and $\lambda$ large enough when starting to construct the array $A$, we may assume that $A$ and thus $A^*$ is big enough that we may apply the Shelah trick again.
Thus, we may without loss suppose that $A^*$ is indiscernible.
By Lemma \ref{groupexists}, there is a $1$-dimensional Galois-definable group in $(\M^{eq})^{eq}$.

\end{proof}

\chapter{An example: covers of the multiplicative group of an algebraically closed field}
 
In this chapter, we show that curves on a cover of the multiplicative group of an algebraically closed field satisfy the axioms
for a Zariski-like structure.
In sections 1-4 we give some results about the cover structures equipped with a topology obtained by taking positive quantifier-free definable sets as basic closed sets.
In section 5, we develop dimension theory for these sets and discuss the connection with the dimension on the Zariski topology of the field sort. 
Most of the results in these sections have been presented previously in \cite{Lucy}.

In section 6, we discuss bounded closures and show that the dimension given by the pregeometry obtained from the bounded closure coincides with the dimension obtained from the closed sets. 
In section 7, we apply the results to show that a cover of the multiplicative group of an algebraically closed field is Zariski-like.
Most of the arguments are similar as those in \cite{HrZi}.
   
\begin{definition}
Let $V$ be a vector space over $\mathbb{Q}$ and let $F$ be an algebraically closed field of characteristic $0$.
A \emph{cover} of the multiplicative group of $F$ is a $2$-sorted structure $(V, F^{*})$ represented by an exact sequence 
\begin{displaymath}
0 \to K \to V \to F^{*} \to 1,
\end{displaymath}
where the map $V \to F^{*}$ is given by exp, a surjective group homomorphism from $(V,+)$ onto $(F^{*},\cdot)$ with kernel $K$.
\end{definition}

We will consider a cover as a structure $V$ in the language $\mathcal{L}=\{0, +, f_q, R_{+}, R_0\}_{q \in \mathbb{Q}}$, where $V$ consists of the elements in the vector space, $0$ is a constant symbol denoting the zero element of the vector space $V$, $+$ is a binary function symbol denoting addition on $V$, and for each $q \in \mathbb{Q}$, $f_q$ is a unary function symbol denoting scalar multiplication by the number $q$.
The symbol $R_+$ is a ternary relation symbol interpreted so that $R_+(v_1, v_2, v_3)$ if and only if $\textrm{exp}(v_1)+\textrm{exp}(v_2)=\textrm{exp}(v_3)$, and $R_0$ is a binary relation symbol interpreted so that $R_0(v_1, v_2)$ if and only if $\textrm{exp}(v_1)+\textrm{exp}(v_2)=0$.
Note that field multiplication is definable using vector space addition.
 
However, for the sake of readability, we will be using the concepts of a vector space $V$ (called the \emph{cover}) and a field $F$ together with the usual algebraic notation when expressing statements about the structure.
If $v=(v_1, \ldots, v_n) \in V^n$, we write $\textrm{exp}(v)$ for $(\textrm{exp}(v_1), \ldots, \textrm{exp}(v_n)) \in F^n$.

The first-order theory of this type of cover structures is complete, submodel complete, superstable and admits elimination of quantifiers (\cite{ZilUusi}, \cite{Zilber}).
Moreover, with an additional axiom (in $L_{\omega_1 \omega}$) stating $K \cong \mathbb{Z}$, the class is categorical in uncountable cardinalities. 
This was originally proved in \cite{Zilber} but an error was later found in the proof and corrected in \cite{monet}.
Throughout this presentation, we will make the assumption $K \cong \mathbb{Z}$.

\section{Varieties and tori}
 
We will eventually define a topology on the cover and show that the irreducible sets of that topology satisfy our axioms.
To be able to do this, we first need to look at some properties of varieties.
When using the word variety, we always mean a Zariski closed subset of $F^n$ for some $n$, defined as the zero locus of some set of polynomials (as in Chapter 1).
That is, we only consider affine varieties, and we don't require them to necessarily be irreducible.

To be able to understand the behaviour of first-order types on the cover structure, we need to understand some properties of roots.
 
\begin{definition}
Let $W$ be an irreducible variety.
For any natural number $n$, we say that an irreducible variety $X$ is an $n$:th root of $W$ if $X^n=W$. 
 
Suppose now $W$ is an arbitrary variety with a decomposition $W=W_1 \cup \ldots \cup W_r$ into irreducible components.
Then, we define the $n$:th roots of $W$ to be all the unions of the form $\bigcup_{i=1}^r {{W_i}^{\frac{1}{n}}}_{(j)}$, where each ${{W_i}^{\frac{1}{n}}}_{(j)}$ is an $n$:th root of $W_i$.
\end{definition}

We note that every variety has only finitely many $n$:th roots.
Also, if $W$ is any variety and $X$ is a $n$:th root of $W$, then $X^n=W$. 

\begin{remark}
For any variety $W$ and any nonzero natural number $m$, we have 
$$\textrm{log }W^m= \bigcup_i m \textrm{log}(W^m)^{\frac{1}{m}}_{(i)}$$
and
$$\frac{1}{m} \textrm{log}(W)=\bigcup_{i} \textrm{log }W^{\frac{1}{m}}_{(i)}.$$

To see that the first equation holds, suppose first $u \in m \textrm{log}(W^m)^{\frac{1}{m}}_{(i)}$ for some $i$.
Then, $u=mv$ for some $v \in \textrm{log}(W^m)^{\frac{1}{m}}_{(i)}$ and $\textrm{exp}(u)=\textrm{exp}(mv) \in W^m$.
On the other hand, suppose $u \in \textrm{log }W^m$.
Then, $\textrm{exp}(u)=x$ for some $x \in W^m$, and thus $\textrm{exp}(\frac{u}{m})=x^{\frac{1}{m}}_{(i)} \in (W^m)^{\frac{1}{m}}_{(j)}$ for some $m$:th roots $x^{\frac{1}{m}}_{(i)}$ and $(W^m)^{\frac{1}{m}}_{(j)}$ of $x$ and $W^m$, respectively.  
Thus, $u \in m \textrm{log}((W^m)^{\frac{1}{m}}_{(j)})$.

For the second equation, suppose first $u=\frac{v}{m}$ for some $v \in \textrm{log}(W)$.
Then, $\textrm{exp}(u)$ is an $m$:th root of $x$ for some $x \in W$.
Thus, $\textrm{exp}(u) \in W^{\frac{1}{m}}_{(i)}$ for some $i$.
On the other hand, suppose $u \in \textrm{log }W^{\frac{1}{m}}_{(i)}$ for some $i$.
Then, $\textrm{exp}(mu) \in W$,  so $u \in \frac{1}{m} \textrm{log }W$.
\end{remark}

Let $W\subset F^n$ be an irreducible variety with $m$:th roots $W^{\frac{1}{m}}_{(i)}$. 
We say that an element $x \in F^n$ is an $m$:th \emph{root of unity} if each of its coordinates is an $m$:th root of unity in $F$.
We note that multiplication by $m$:th roots of unity permutes the $m$:th roots of $W$.
Suppose $x \in W^{\frac{1}{m}}_{(i)}$ for some $i$.
Then, $x^m \in W$.
If $\zeta \in F^n$ is an $m$:th root of unity, then $(\zeta x)^m=x^m \in W$.
Hence, $\zeta x \in W^{\frac{1}{m}}_{(j)}$ for some $j$.
Now we must have $\zeta W^{\frac{1}{m}}_{(i)}= W^{\frac{1}{m}}_{(j)}$ because $W^{\frac{1}{m}}_{(i)}$ is irreducible (if different elements were mapped into different roots in multiplication by $\zeta$, then we could write $W^{\frac{1}{m}}_{(i)}=\bigcup_j \{x \, | \, \zeta x \in W^{\frac{1}{m}}_{(j)}\})$.
Note also that the image consists of the whole of $W^{\frac{1}{m}}_{(j)}$ as multiplication by $\zeta$ is an injection and the different $m$:th roots have the same dimension.
  
Suppose now $W^{\frac{1}{m}}_{(i)}$ and $W^{\frac{1}{m}}_{(j)}$ are two different $m$:th roots of $W$.
Let $x \in W^{\frac{1}{m}}_{(i)} \setminus W^{\frac{1}{m}}_{(j)}$.
Then, there is some $y \in W^{\frac{1}{m}}_{(j)}$ such that $x^m=y^m$.
Hence, $y=\zeta x$ for some root of unity $\zeta$.
Now, $\zeta W^{\frac{1}{m}}_{(i)}=W^{\frac{1}{m}}_{(j)}$.
 
We can now determine quantifier-free types on the cover. 
One easily sees that the following lemma holds. 
 
\begin{lemma}\label{tyypit}[\cite{ZilUusi}]
Let $(V,F)$ be a cover and $A \subset V$.
Let $v \in V^n$ with linearly independent coordinates.
Then, the quantifier free type of $v$ over $A$ is determined by the formulae
\begin{eqnarray*}
\textrm{exp}\left(\frac{v}{l}\right) &\in& W^{\frac{1}{l}} \quad l \in \mathbb{N}, \\
\textrm{exp}(v) &\notin& Y \quad Y \subset W, \textrm{dim}(Y)<\textrm{dim}(W), \\
mv &\neq& 0 \quad m \in \mathbb{Z}^n, m \neq 0,
\end{eqnarray*}
where $W$ is the locus of $\textrm{exp}(v)$ over $\mathbb{Q}(\textrm{exp}(A))$ (the smallest field containing $\textrm{exp}(A))$ and each $W^{\frac{1}{l}}$ is an $l$:th root of $W$.
\end{lemma}

\subsection{Linear sets and tori}

\begin{definition}
Let a subset $L$ of the cover sort be called \emph{linear} if it can be defined by $\mathbb{Q}$-linear equations only.
\end{definition}

\begin{remark}\label{linearintersect}
We note that for any linear set $L \subset V^n$ and any $k=(k_1, \ldots, k_n) \in K^n$, it holds that $L \cap (L+k) \neq \emptyset$ if and only if $L+k=L$.
Indeed, if $L+k \neq L$, then one of the equations defining $L$ is of the form
 $$q_1 v_1 +\ldots +q_n v_n+b=0, \quad q_i \in \mathbb{Q}, \quad b \in V,$$
 where for some $1 \le i \le n$, $q_i \neq 0$ and $k_i \neq 0$.
Now every element of $L+k$ satisfies
$$q_1 v_1 +\ldots +q_n v_n+b=(q_1 k_1, \ldots, q_n k_n) \neq 0.$$
No element can satisfy both equations.
\end{remark}

In our analysis of the definable sets on the cover it will be very useful that every linear set will correspond to a torus on the field sort and that we can thus use linear sets to analyze tori.
Many of the definable sets on the cover are obtained as inverse images of varieties under the map exp.
Since the field element $0$ does not have an inverse image, we can do the same analysis by thinking of our varieties as the Zariski closed subsets of the Zariski open set $({F^*})^n$.
Thus, every time we will be considering a variety, we will mean the Zariski closure of some such set.
For instance, we would not consider the variety $W$ given by the polynomial $xy-x=0$, since $W \cap (F^*)^2$ will 
already be given by the Zariski closure $\overline{W \cap (F^*)^2}$ which is 
given by the polynomial $y-1=0$. 
 
We give $({F^*})^n$ group structure by taking coordinate-wise multiplication as the group operation.
Then, we will think of any coset of a subgroup as a torus.
This gives the following definition.

\begin{definition}
 Call a set $T \subseteq (F^*)^m$ a \emph{torus} if it can be defined using equations of the form
\begin{displaymath}
\prod_i {x_i}^{z_i}=c \qquad z_i \in \mathbb{Z} \quad c \in F^{*} 
\end{displaymath}

If $T \subseteq (F^*)^m$ is a torus such that $T \neq (F^*)^m$, we say that $T$ is a \emph{proper torus}.
\end{definition}

We will sometimes view a torus as a variety.
Then, we mean the Zariski closure of a set that is defined as above.
The ideal corresponding to this kind of a variety is generated by polynomials of the form
$$\prod_{i \in I} {x_i}^{n_i}-c \prod_{i \in J} {x_i}^{n_i},$$
where $c \in F^{*}$, $n_i \in \mathbb{N} \setminus \{0\}$ for each $i$, and $I, J$ are finite index sets such that $I \cap J=0$.
 
We will say that a torus $T$ is \emph{irreducible}, if the Zariski closure of $T$ is irreducible as a variety (in the usual sense).

A basic property of irreducible tori is that they can be transformed in a canonical form by a birational coordinate change on $(F^{*})^n$.
Namely, if $T \subseteq {F^*}^n$ is an irreducible torus, then there is a natural number $k$ such that $T$ can be expressed as $x_i=c_i$ for $1\le i \le k$, where $c_i \in F^{*}$ for each $i$.
To show this, we first prove the following auxiliary lemma.

\begin{lemma}\label{matriisi}
Let $z_1, \ldots, z_n$ be integers and suppose $1$ is the greatest integer that divides each one of the $z_i$.
Then, there exists an $n \times n$ integer matrix $A$ such that the first row of $A$ is $(z_1, \ldots, z_n)$ and $\textrm{det}(A)= \pm 1$. 
\end{lemma}

\begin{proof}
We use induction on $n$ to show that the matrix $A$ exists.
If $n=1$, this is clear as we must have $z_1= \pm 1$.
Suppose now $n=2.$  
Since $\textrm{gcd}(z_1, z_2)=1$, there are $d_1, d_2 \in \mathbb{Z} \setminus \{0\}$ such that $d_1 z_1+d_2 z_2=1$, and we may choose 
\begin{displaymath} 
A=\left( \begin{array}{cc}
 z_1 & z_2 \\
-d_2 & d_1 
\end{array} \right).
\end{displaymath}

Suppose now the claim holds for $n$, and consider $n+1$.
Let $m=\textrm{gcd}(z_{n}, z_{n+1})$.
Write $z_{n}'=\frac{z_{n}}{m}$ and $z_{n+1}'=\frac{z_{n+1}}{m}$. 
Let $d_n, d_{n+1} \in \mathbb{Z} \setminus \{0\}$ be such that $d_n z_{n}'+d_{n+1} z_{n+1}'=1$,
and let
\begin{displaymath}
M=\left( \begin{array}{ccccc}
 1 & 0 && \ldots & 0\\
0 &1 && \ldots & 0 \\
\vdots && \ddots && \vdots \\
0&0&\ldots &z_{n}' &z_{n+1}'\\
0&0&\ldots &-d_{n+1}& d_n
\end{array} \right).
\end{displaymath}
Then, $\textrm{det}(M)=1$. 
 Now $1$ is the greatest integer dividing each one of $z_{1}, \ldots, z_{n-1}, m$, and thus by the inductive assumption, there is an integer matrix $(A')_{ij}$ with $\textrm{det}(A')= 1$ and first row $(z_{1}, \ldots, z_{n-1}, m)$. 
Let $(B)_{ij}$ be the $(n+1) \times (n+1)$ -matrix such that $B_{ij}=A'_{ij}$ for $1 \le i,j \le n$ and $B_{ij}=\delta_{ij}$ otherwise.
Then, 
$\textrm{det}(B)=\textrm{det}(A')=1$, and the matrix $A=BM$ is as wanted. 
\end{proof}

\begin{lemma}\label{torusbirat}
Let $T \subset( {F^{*}})^m$ be an irreducible torus given in the coordinates $x_1, \ldots, x_m$.
Then, there is a birational coordinate change given by 
\begin{eqnarray}\label{chaange}
y_i=\prod_{j=1}^m x_j^{z_{ij}}, \quad x_i=\prod_{j=1}^m y_j^{z_{ij}'}, \quad 1 \le i \le m, \quad z_{ij}, z_{ij}' \in \mathbb{Z},
\end{eqnarray}
such that in the new coordinates, $T$ is of the form 
$$y_i=c_i, \quad d_i \in \mathbb{Z}, \quad c_i \in F^*, \quad 1 \le i \le k,$$
where $1 \le k \le m$, and $c_i \in F^*$ for each $i$.
 \end{lemma}
 
\begin{proof} 
Suppose the torus $T$ is given by the equations 
\begin{eqnarray}\label{tallaa}
\prod_{i=1}^m x_i^{n_{ji}}=c_j, 
\end{eqnarray}
where $1\le j \le k$ for some $k$, $n_{ij} \in \mathbb{Z}$ and $c_j \in F^*$.
 
To prove the lemma, we will view the multiplicative group $(F^{*})^m$ as a $\mathbb{Z}$-module where $\mathbb{Z}$ acts by exponentiation. 
Then we  look for invertible endomorphims that would give a suitable coordinate change. 
 
We start looking at the first one of the equations (\ref{tallaa}). 
Since $T$ is irreducible, $1$ is the 
greatest integer that divides each one of the $n_{1i}$ for $1 \le i \le m$.
By Lemma \ref{matriisi}, there exists an integer matrix $A$ such that the first row of $A$ is $(n_{11}, \ldots, n_{1m})$ and $\textrm{det}(A)= \pm 1$.
Then, the coordinate change given by $A$ is of the form (\ref{chaange}) and transforms our equation into $y_1=c_1$.
Using Cramer's rule, we see that  $A^{-1}$ is also an integer matrix, 
and thus the reverse coordinate change is also given in the form (\ref{chaange}).
 
Since the coordinate change we have done is given by equations of the form (\ref{chaange}), 
all the equations in (\ref{tallaa}) are still in the torus form after the transformation.
Consider the second equation.
After substituting $y_1=c_1$, it will be given by 
\begin{eqnarray*} 
y_2^{z_2} \cdots y_m^{z_m}=c_2'
\end{eqnarray*}
for some $z_2, \ldots, z_m \in \mathbb{Z}$, $c_2' \in F^*$. 
Let $d$ be the greatest integer dividing each one of the numbers $z_2, \ldots, z_m$.
Then, an integer matrix with determinant $ \pm 1$ and first row $(\frac{z_2}{d}, \ldots, \frac{z_m}{d})$ transforms the equation into
$$y_2^{d}=c_2',$$
which gives us 
$$y_2=\zeta^{i} a, \quad i=0, \ldots, d-1,$$
where $a$ is a number such that $a^{d}=c_2'$ and $\zeta$ is a primitive $d$:th root of unity.

We substitute these values to the third equation to get at most $d$ distinct equations.
Then, we deal with each one of them as we did with the second equation above.
Proceeding this way and going through all the equations, we will get $T$ in the new coordinates as a union of smaller tori, each given by equations of the form
$$y_i=c_i, \quad 1 \le i \le k$$
for $c_i \in F^*$ and some $k \le n$.
Since our coordinate change and its inverse are both given by rational functions, it is a homeomorphism in the Zariski topology, and thus maps irreducible sets to irreducible sets.
Since we assumed $T$ to be irreducible, only one of the components listed is nonempty. 
This proves the lemma.
\end{proof}
  
\begin{remark}
We note that since the coordinate change in Lemma \ref{torusbirat} and its inverse are both given by equations of the form (\ref{chaange}), it maps a variety $W$ to a torus if and only if $W$ is a torus.
\end{remark}

Now it is easy to prove the following properties of tori.

\begin{lemma}\label{torusom}
The following hold:
\begin{enumerate}[(a)]
\item If $T_1$, $T_2$ are tori, then $T_1 \cap T_2$ is a torus.
\item If $T$ is a torus, then every irreducible component of $T$ is a torus.
\item If $T$ is a torus, then $T$ has distinct $m$:th roots for any $m$.
Moreover, any $m$:th root of $T$ is a torus.
\item If $T$ is an irreducible torus, then $T^m$ is a torus for every natural number $m$.
\end{enumerate}
\end{lemma}

\begin{proof}
(a) is clear from the definition.

(b) 
Consider an equation of the form
$$x_1^{z_1} \cdots x_n^{z_n}=c, \quad c \in F^*, z_i \in \mathbb{Z} \textrm{ for } 1 \le i \le n.$$
If the greatest integer dividing each one of the numbers $z_1, \ldots, z_n$ is $1$, then we may birationally
transform the equation into $y_1-c=0$ as in the proof of Lemma \ref{torusbirat}.
On the other hand, if it is some $d>1$, then we get
$$y_1^{d}-c=\prod_{i=0}^{d-1} (y_1-\zeta^i a)=0,$$
where $\zeta$ is a primitive $d$:th root of unity and $a$ is a number such that $a^d=c$.
From this, we see that the corresponding polynomial is irreducible if and only if the greatest number dividing each one of the numbers $z_1, \ldots, z_n$ is $1$.
So every polynomial in the torus form divides into irreducible factors that are also in the torus form.
This proves (b).
 
(c) It is enough to show this for irreducible $T$.
Let $m$ be a non-zero natural number. 
By Lemma \ref{torusbirat}, the variety $T^{\frac{1}{m}}$ (union of all roots) is defined by equations of the form
$x_i^m-c_i=0,$ where $c_i \in F^*$, and is thus clearly reducible.
Also, $T^{\frac{1}{m}}$ is clearly a torus, so the $m$:th roots are tori by (b).

(d)  $T$ is (without loss of generality) given by equations $x_i=c_i$ for $1\le i \le k$ where $k \le n$.
Now, $T^m$ is given by the equations $x_i=c_i^m$ for $1\le i \le k$, and is clearly a torus.
\end{proof} 

\begin{remark}
We note that if $L \subset V^n$ is linear, then $\textrm{exp}(L)\subset F^n$ is a torus.
Also, using Lemma \ref{torusbirat}, it is easy to see that 
any irreducible torus $T\subset F^n$ can be written as $T=\textrm{exp}(L)$ for some linear set $L\subset V^n$
(note that the matrix giving the coordinate change on $F^*$ can also be applied on $V$).
\end{remark}

Now we can state the following theorem that 
follows from Theorem 2.3. in \cite{bays}.
 \begin{theorem}\label{zilbertheorem}
Let $(V,F)$ be a cover with ($K=\mathbb{Z}$).
Let $(G, \textrm{exp}(G))$ be a countable submodel such that $G=\textrm{log}(\textrm{exp}(G))$, and let $h \in V^m$.
 Let $W$ be the locus of $\textrm{exp}(h)$ over $\textrm{exp}(G)$.
Suppose $W$ is not contained in a torus. 
Then the subtype of $h$ over $G$ consisting of formulae $\textrm{exp}(\frac{h}{l}) \in W^\frac{1}{l}$ is finitely determined.
 \end{theorem}

\section{PQF-topology}

By positive, quantifier free formulae, we mean first-order formulae that don't contain any negation symbols or any quantifiers.
In our context this means that we obtain all the sets definable by positive, quantifier-free formulae by first taking all the sets defined by equations of the form
$\sum_{i}q_i v_i=a$, where $q_i \in \mathbb{Q}$ and $a \in V,$ or of the form $\textrm{exp}(\frac{v}{l}) \in W$, where
$W$ is a variety and $l \in \mathbb{N}$, and then closing this collection with respect to finite unions and finite intersections.

\begin{definition}
Define a topology on our structure by taking the sets definable by positive quantifier-free first-order formulae as the basic closed sets.
Call this the \emph{PQF}-topology.
\end{definition}
 
We define the notion of an irreducible set in the usual way.

\begin{definition}
We say a nonempty closed set is \emph{irreducible} if it cannot be written as the union of two proper closed subsets.
\end{definition}

\begin{remark}
We note that the PQF-topology is not Noetherian.
Indeed, let $C_0=\{ u \in V \, | \, \textrm{exp}(u)=1\}$.
For $i=1, 2, \ldots$, denote 
$$C_i=\left\{ u \in V \, | \, \textrm{exp}\left(\frac{u}{2^i} \right)=1 \right\}.$$
Then, $C_0 \supsetneq C_1 \supsetneq \ldots$ is an infinite descending chain of PQF-closed sets.
\end{remark}

The proof of the following lemma is similar to the one presented in \cite{HrZi} for Zariski geometries.

\begin{lemma}\label{PQFcartesianirred}
Let $C_1, C_2$ be irreducible closed sets in the PQF-topology.
Then, $C_1 \times C_2$ is irreducible.
\end{lemma}

\begin{proof}
Suppose $C_1 \times C_2=F_1 \cup F_2$ where $F_1, F_2$ are closed.
For $i=1,2$, let 
\begin{displaymath}
F_i^{*}=\{a \in C_1 \, | \, (a,x) \in F_i \textrm{ for every } x \in C_2\}.
\end{displaymath}
We note that $F_1^*$ is PQF-closed.
Indeed, for each $b \in C_2$, the set $D_b=\{a \in C_1 \, | \, (a,b) \in F_1\}$ is PQF-closed, and thus
$F_1^{*}= \bigcap_{b \in C_2} D_b$ is closed.
Similarly, $F_2^{*}$ is closed.
Let $a \in C_1$.
 For $i=1,2$, the set $F_i(a)=\{x \in C_2 \, | \, (a,x) \in F_i\}$ is PQF-closed.
Now $C_2=F_1(a) \cup F_2(a)$, and thus, as $C_2$ is irreducible, we have either $C_2=F_1(a)$ or $C_2=F_2(a)$.
Hence, for each $a \in C_1$, there is an $i \in \{1,2\}$ such that $(a,x) \in F_i$ for every $x \in C_2$.
So, $C_1=F_1^{*} \cup F_2^{*}$, and as $C_1$ is irreducible,  $C_1=F_i^{*}$ for some $i$.
Thus, $C_1 \times C_2=F_i$.
\end{proof}

\begin{definition}
Let $W$ be a variety.
If $W$ has distinct $n$:th roots for some natural number $n$, we say that $W$ \emph{branches}.
We say that $W$ \emph{stops branching at the finite level} if there is a natural number $l$ such that the $l$:th roots $W^{\frac{1}{l}}$ no longer branch.
We say $W$ \emph{branches infinitely} if it does not stop branching at the finite level.
 \end{definition}
 
Suppose $W$ is a variety, $v \in \textrm{log }W$.
 For any $l$, denote by $W^{\frac{1}{l}}_{(v)}$ the $l$:th root of $W$ such that $\textrm{exp}(\frac{v}{l}) \in W^{\frac{1}{l}}_{(v)}$.
If $W$ is a variety not contained in any torus, then, Theorem \ref{zilbertheorem} implies that there is some number $m$ such that for any $m'>m$, the $m'$:th root $W^{\frac{1}{m'}}_{(v)}$ is determined by the $m$:th root $W^{\frac{1}{m}}_{(v)}$.

\begin{lemma}[\cite{Lucy}]\label{PQFbasicform}
Any set definable by a positive quantifier free formula is a finite union of sets of the form
\begin{displaymath}
m \cdot (L \cap \textrm{log } W)
\end{displaymath}
for some linear set $L$, a variety $W$ and some $m \in \mathbb{N}$.
\end{lemma}

\begin{proof} 
On the cover sort, a basic PQF-closed subset of $V^m$ is defined by some positive boolean combination of equations
\begin{displaymath}
\sum_i q_i v_i=a, \qquad \textrm{exp}(\frac{v}{l}) \in W,
\end{displaymath}
where $a \in V$, $W$ is a variety and $l \in \mathbb{N}$.

To see this, suppose that $v=(v_1, \ldots, v_m) \in V^m$ satisfies  
\begin{displaymath}
(\textrm{exp}(q_{11} v_1+ \ldots + q_{1m}v_m +a_1), \ldots, \textrm{exp}(q_{n1} v_1+\ldots +q_{nm}v_m+a_n)) \in W_0 
\end{displaymath}
for $q_{ij} \in \mathbb{Q}$ and some variety $W_0 \subset F^n$.
Now, $q_{ij}=\frac{k_{ij}}{l_{ij}}$ for some $k_{ij} \in \mathbb{Z}$, $l_{ij} \in \mathbb{N}$, and we may write for each $i$
\begin{displaymath}
\textrm{exp}(q_{i1} v_1+ \ldots + q_{im}v_m +a_i)=\textrm{exp}\left(\frac{v_1}{l_{i1}}\right)^{k_{i1}} \cdots \textrm{exp}\left(\frac{v_1}{l_{im}}\right)^{k_{im}} \textrm{exp}(a_i).
\end{displaymath}
By suitably changing the $k_{ij}$ and $l_{ij}$ (by expanding the fractions) we may assume that $l_{ij}=l$ for each $i,j$.
When we substitute these values in the equations of the variety $W_0$ and clear the denominators (note that some of the $k_{ij}$ might be negative), we get equations for a new variety $W$ such that $(\textrm{exp}(\frac{v_1}{l}), \ldots, \textrm{exp}(\frac{v_m}{l}))= \textrm{exp}(\frac{v}{l}) \in W$.

 For any $l$,
 \begin{displaymath}
 \textrm{exp}(v) \in W \Leftrightarrow \textrm{exp}(\frac{v}{l}) \in \bigcup_i W^{\frac{1}{l}}_{(i)},
 \end{displaymath}
 where the union is taken over all possible choices of the $l$:th root $W^{\frac{1}{l}}$.
 Hence,
 \begin{displaymath}
 \left(\textrm{exp}\left(\frac{v}{l_1}\right) \in W_1\right) \wedge \left(\textrm{exp}\left(\frac{v}{l_2}\right) \in W_2\right) \Leftrightarrow v \in l_1 \cdot l_2 \textrm{log} W,
 \end{displaymath}
 where $W=\bigcup (W_1^{\frac{1}{l_2}} \cap W_2^{\frac{1}{l_1}})$ and the union is again over all possible roots.
 Since we also have $v \in L \Leftrightarrow \frac{v}{l} \in \frac{1}{l} L$, this proves the lemma (note that if we have two linear sets $L_1$ and $L_2$, then $L_1 \cap L_2$ is linear). 
 \end{proof}

\begin{corollary}\label{variety}[\cite{Lucy}]
Let $C$ be a set on the cover sort, definable by a positive, quantifier-free first-order formula (i.e. a basic closed set in the PQF-topology).
Then $\textrm{exp}(C)$ is a Zariski closed set on the field sort.
\end{corollary}
\begin{proof}
By Lemma \ref{PQFbasicform}, it suffices to consider sets of the form $C=m(L \cap \textrm{log }W)$, where $L$ is a linear set and $W$ is a variety.
Let $T=\textrm{exp}(L)$.
Then, $\textrm{exp}(C)=(T \cap W)^m$.
But $(T \cap W)^m$ is the image of the Zariski closed set $T \cap W$ under the finite map $x \mapsto x^m$.
Hence, it is Zariski closed by the Corollary of Lemma I.5.2 in \cite{shafa}.
\end{proof}

\section{Irreducible Sets}

In this section we present some basic properties of the sets irreducible on the cover.
We will show that all irreducible sets are actually definable by positive quantifier-free formulae.
First, we give a canonical way to write any irreducible  variety $W$ as $W=T \cap W'$, where $T$ is a torus and $W'$ is not contained in any proper torus.
  
\begin{lemma}\label{leikkaus}[\cite{Lucy}]
Any irreducible variety $W \subset F^n$ can be written as $W' \cap T$ where $T$ is the minimal torus containing $W$ (note that this could be $F^n$) and $W'$ is a  variety not contained in any proper torus. 
\end{lemma}

\begin{proof}  
By Lemma \ref{torusbirat}, we may assume $T$ is (without loss of generality) given by equations $x_i=c_i$ for $1\le i \le k$ where $k \le n$.
Let $a$ be a generic point of $W$ and let $I \subseteq F[x_{k+1}, \ldots, x_{k_n}]$ be the ideal consisting of all polynomials $f$ such that $f(a)=0$.
Let $J=\langle I \rangle \subseteq F[x_1, \ldots, x_n]$, the ideal generated by $I$ in $F[x_1, \ldots, x_n]$.
Let $W'$ be the variety associated to $J$.
Since the ideal $J$ does not contain any of the polynomials $x_i-c_i$ for $1 \le i \le k$ and since $T$ is the minimal torus containing $W$, the variety $W'$ is not contained in any torus.
\end{proof}

\begin{remark}
The variety $W'$ given in the proof of Lemma \ref{leikkaus} is irreducible.
Indeed, since $W$ is irreducible, also the variety $V(I)$ given by the ideal $I$ is irreducible.
Now, $W'=F^k \times V(I)$ which is irreducible as a Cartesian product of two irreducible varieties.
\end{remark}
 
We still need two lemmas before being able to show that irreducible sets are definable by positive, quantifier-free formulae.

\begin{lemma}\label{singlemth}
 Let $T$ be an irreducible torus and let $L\subset V^n$ be a linear set such that $\textrm{exp}(L)=T$.
Then, for each $m$, $\textrm{exp}(\frac{L}{m})=T^{\frac{1}{m}}_{(i)}$ for a single $m$:th root of $T$.
\end{lemma}

\begin{proof}
As $T$ is a torus, we may, by Lemma \ref{torusbirat}, assume it is given by a finite set of equations of the form
\begin{eqnarray*}  
x_i-c_i=0,
\end{eqnarray*}
where $c_i \in F^*$.
Consider the $m$:th roots of $T$ for some arbitrary $m$.
The torus equations give us equations of the form
\begin{eqnarray*} 
x_i-\zeta_{ij}=0,
\end{eqnarray*}
where each $\zeta_{ij}$ is an $m$:th root of $c_i$ ($j=1, \ldots, m$).
Then, each $m$:th root of $T$ satisfies exactly one of these equations for each $i$.

The set $\frac{L}{m}$ is given by a set of linear equations in the variables $u_1, \ldots, u_n$.
If it were to contain some elements $\frac{a}{m}$ and $\frac{b}{m}$ that would map into distinct roots, then we would have for some $i$ that 
$\textrm{exp}\left(\frac{a_i}{m}\right) =\zeta_{ij}$ 
but 
$\textrm{exp}\left(\frac{b_i}{m}\right)=\zeta_{ij'}$
where  $j \neq j'$.
This is impossible, as the linear equations defining $\frac{L}{m}$ cannot imply both 
$u_i=d_j$ 
and 
$u_i=d_{j'}$, 
where $d_j \neq d_{j'}$ (if we choose $d_j, d_{j'}$ so that $\textrm{exp}(d_j)=\zeta_{ij}$ and $\textrm{exp}(d_{ij'})=\zeta_{j'}$, then clearly $d_j \neq d_{j'}$). 
\end{proof}

\begin{lemma}\label{linearirre}
Let $C \subset V^n$ be irreducible, and let $L \subset V^n$ be linear.
Suppose $C \subset \bigcup_{k \in K^n} L+k$.
Then, $C \subset L+k$ for a single $k$.
 \end{lemma}

\begin{proof}
Suppose there are  $a,b \in C$ such that $a \in L+h_1$, $b \in L+h_2$ for some $h_1, h_2 \in K^n$ such that $L+h_1 \neq L+h_2$.
For simplicity of notation, we denote $L+h_1$ by $L'$ and set $k=h_2-h_1$ which allows us to denote $L+h_2$ by $L'+k$.
Denote $T=\textrm{exp}(L')$.
Now $T$ is an irreducible torus and thus, by lemma \ref{torusbirat}, we may assume it is given by equations of the form
\begin{eqnarray*} 
x_i-c_i=0, \quad 1 \le i \le m
\end{eqnarray*}
for some $m \le n$.

Write $k=(k_1 z, \ldots, k_n z)$, where $z$ is the generator of the kernel and $k_i \in \mathbb{Z}$ for each $i$.
We may without loss suppose that $k_1 \neq 0$.
Let $M>1$ be a natural number such that $\textrm{gcd}(k_1, M)=1$.
By Lemma \ref{singlemth}, $\frac{L'}{M}$ maps to a single $M$:th root of $T$, and so does
$\frac{L'+k}{M}$.
Use coordinates $u_1, \ldots, u_n$ for the cover sort.
Then, there is an element $\zeta \in V$ such that $\textrm{exp}(\zeta)=c_1$,  
every point of $\frac{L'}{M}$ satisfies the equation 
\begin{eqnarray}\label{jakamatta}
u_1-\frac{\zeta}{M}=0
\end{eqnarray}
 and every point of $\frac{L'+k}{M}$ satisfies the equation 
\begin{eqnarray}\label{jaettu}
 u_1-\frac{\zeta+k_1z}{M}=0.
\end{eqnarray}
By the equation (\ref{jakamatta}), every point of $\textrm{exp}(\frac{L'}{M})$ satisfies the equation
$$x_1=\textrm{exp}\left(\frac{\zeta}{M}\right),$$
and by the equation (\ref{jaettu}), every point of  $\textrm{exp}(\frac{L'+k}{M})$ satisfies the equation
$$x_1=\textrm{exp}\left(\frac{\zeta+k_1z}{M}\right).$$
Since $k_1$ is not divisible by $M$, we have $\frac{k_1z}{M} \notin K$ and thus
$\textrm{exp}(\frac{\zeta}{M}) \neq \textrm{exp}(\frac{\zeta+k_1z}{M}).$
Hence, $\frac{L'}{M}$ and $\frac{L'+k}{M}$ map to distinct $M$:th roots of $T$.
  
Now, we may write
\begin{displaymath}
C \subset \bigcup_i \left\{x \, | \, \textrm{exp}\left(\frac{x}{M}\right) \in T^{\frac{1}{M}}_i \right\},
\end{displaymath}
where $T^{\frac{1}{M}}_i$ are the distinct $M$:th roots of $T$.
Since $a$ and $b$ are in different members of the union, $C$ is not contained in any single one of them.
This contradicts the irreducibility of $C$.
\end{proof}

The following lemma gives a canonical form for the irreducible sets.
It also implies that in particular, they are definable.
 
\begin{lemma}\label{irredform}[\cite{Lucy}]
An irreducible PQF-closed subset of $V^n$ has the form 
\begin{displaymath}
L \cap m \cdot \textrm{log } W,
\end{displaymath}
for a linear $L$, a variety $W$ which does not branch and $m \in \mathbb{N}$. 
\end{lemma}

\begin{proof}
A general PQF-closed set $C$ is an intersection of basic PQF-closed sets.
If $C$ is to be irreducible, we may, by Lemma \ref{PQFbasicform}, assume that each of these basic PQF-closed sets is of the form $L \cap m \textrm{log }W$ for some linear set $L$ and some variety $W$.
Thus, we may write
\begin{displaymath}
C=\bigcap_{i<\kappa} (L_i \cap m_i \textrm{log } W_i)
\end{displaymath}
for some cardinal $\kappa$.
By Noetherianity of the linear topology on $V^n$, the linear part stabilizes, so writing $L=\bigcap_{i < \kappa} L_i$, we get
\begin{eqnarray}\label{cee}
C= L \cap \bigcap_{i <\kappa} m_i \textrm{log } W_i.
\end{eqnarray}
By Lemma \ref{leikkaus}, each $W_i$ can be written as $W_i=W_i'\cap T_i$ where $T_i$ is a torus and $W_i'$ is not contained in any torus.
As $\textrm{log}(W_1 \cap W_2)=\textrm{log } W_1 \cap \textrm{log } W_2$ for any varieties $W_1, W_2$, we may assume that each $W_i$ is itself either a torus or contained in no torus.

If $W_i$ is a torus for some $i$, then $\textrm{log }W_i=\bigcup_{k \in K^n} (L'+k)$, where $K$ is the kernel of the map exp and $L'$ is a linear set such that $\textrm{exp}(L')=W_i$. 
Now $m_i \textrm{log }W_i=\bigcup_{k \in K^n} (m_iL'+m_ik)$, and by Lemma \ref{linearirre}, $C \subset m_iL'+m_ik$ for a single $k$.
The set $m_iL'+m_ik$ is linear, so we may without loss assume it is contained in the intersection $L=\bigcap_{i \in I} L_i$.
Thus, we only need to consider the $W_i$ that are not contained in any torus.
 
If $W_i$ is contained in no torus, we may assume that $W_i$ is irreducible and that $W_i$ does not branch at all.
Indeed, if $W_i$ would be reducible, we could write
\begin{displaymath}
\textrm{log } W_i= \bigcup_j \{x \, | \, \textrm{exp}(x)  \in W_i^j\},
\end{displaymath}
where the $W_i^j$ are the irreducible components of $W_i$, and if it would branch for some $m$, we could write
\begin{displaymath}
\textrm{log } W_i= \bigcup_j \{x \, | \,\textrm{exp}( \frac{x}{m}) \in {W_i^{\frac{1}{m}}}_{(j)}\},
\end{displaymath}
where the ${W_i^{\frac{1}{m}}}_{(j)}$ are the $m$:th roots of $W$.
For any varieties $W_1, W_2$ which don't branch at all, we have that both $\textrm{exp}(\frac{v}{m_1}) \in W_1$ and $\textrm{exp}(\frac{v}{m_2}) \in W_2$ if and only if
\begin{displaymath}
\textrm{exp}\left(\frac{v}{m_1 m_2}\right) \in W_1^{\frac{1}{m_2}} \cap W_2^{\frac{1}{m_1}}
\end{displaymath}
for the unique roots $W_1^{\frac{1}{m_2}}$, $W_2^{\frac{1}{m_1}}$ of $W_1$ and $W_2$, respectively.

Consider now the representation (\ref{cee}).
If we have $W_1^{\frac{1}{m_2}}=W_1^{\frac{1}{m_2}} \cap W_2^{\frac{1}{m_1}}$, then 
$$m_1 \textrm{log }W_1 \cap m_2 \textrm{log }W_2=m_1 \textrm{log }W_1,$$
and we may drop $m_2 \textrm{log }W_2$ from the representation (\ref{cee}).
If not, then there is some irreducible component $W_{1,2}$ of $W_1^{\frac{1}{m_2}} \cap W_2^{\frac{1}{m_1}}$ such that $\textrm{exp}(\frac{C}{m_1 m_2}) \subseteq W_{1,2}$.
Write $W_{1,2}=W_{1,2}' \cap T_{1,2}$ where $T_{1,2}$ is a torus and $W_{1,2}'$ is not contained in any torus.
Then, there is a linear set $L_{1,2}$ such that $C \subset L_{1,2}$ and $\textrm{exp}(L_{1,2})=T_{1,2}$.
Now, we may replace $L$ by $L \cap L_{1,2}$ in the representation (\ref{cee}) and only consider $W_{1,2}'$ from now on.
As $W_{1,2}'$ is not contained in any torus, there is a number $n$ such that the $n$:th roots of $W_{1,2}'$ no longer branch.
Moreover, for one of these $n$:th roots, say $W_{1,2}^{*}$, we have that $\textrm{exp}(\frac{C}{n m_1 m_2}) \subseteq W_{1,2}^{*}$.
Thus, we may write
\begin{eqnarray*}
C=L \cap L_{1,2} \cap n m_1 m_2 \textrm{log } W_{1,2}^{*} \cap \bigcap_{3 \le i <\kappa} m_i \textrm{log }W_i.
\end{eqnarray*}
We note that since $W_1$ is irreducible and does not branch, also $W_1^{\frac{1}{m_2}}$ is irreducible.
Thus, since definable finite-to-one maps  preserve Morley ranks,
$$\textrm{MR}(W_1^{\frac{1}{m_2}} \cap W_2^{\frac{1}{m_1}})<\textrm{MR}(W_1^{\frac{1}{m_2}})=\textrm{MR}(W_1),$$
and in particular, $$\textrm{MR}(W_{1,2}^{*})<\textrm{MR}(W_1).$$ 
Now, we repeat the above process for $W_{1,2}^{*}$ and $W_3$ to obtain $W_{1,2,3}^{*}$ (in case we ended up discarding  $m_2 \textrm{log }W_2$ from the representation (\ref{cee}), we do this for $W_1$ and $W_3$). 

We claim that we can go on this way for at most finitely many steps (meaning that we may discard all but finitely many of the $m_i \textrm{log }W_i$).
Suppose not. 
Then, we may define $W_{1, \ldots, n}^*$ for arbitrary large $n$ (after relabeling the indices to account for sets that were discarded from the representation (\ref{cee})),
and we always have $\textrm{MR}(W_{1, \ldots, n}^*)>\textrm{MR}(W_{1, \ldots, n+1}^*)$.
Moreover, for each $W_{1, \ldots, n}^*$, there is some number $M_n$ so that
$$W_1 \supseteq (W_{1,2}^*)^{M_2} \supseteq \ldots \supseteq (W_{1, \ldots, n}^*)^{M_n} \supseteq (W_{1, \ldots, n+1}^*)^{M_{n+1}} \supseteq \ldots$$
Since definable finite-to-one maps preserve Morley ranks, we see that the rank drops at each inclusion above, and thus all the inclusions must be proper.
This contradicts the Noetherianity of the Zariski topology on the field sort $F$.

So the process eventually terminates, and we get for $C$ a representation that is as wanted.
\end{proof}

\begin{corollary}
If $C$ is a closed irreducible subset of the cover sort, then $\textrm{exp}(C)$ is an irreducible variety.
\end{corollary}
\begin{proof}
By Corollary \ref{variety}, $\textrm{exp}(C)$ is a variety.
If it were reducible, then also $C$ would be reducible.
\end{proof}

\section{Irreducible Components}

Since the PQF-topology is not Noetherian, we cannot speak about irreducible components in the classical sense.
However, we give a more general definition of irreducible components that makes sense in the context of PQF-closed sets.
We then prove some basic properties of the irreducible components of sets of the form $\textrm{log}(W)$, where $W$ is a variety.
The results in this section were presented in \cite{Lucy}.

\begin{definition}
If $C$ is a PQF-closed set, we say that the \emph{irreducible components} of $C$ are the maximal irreducible subsets of $C$.
\end{definition}

We will eventually prove that for any irreducible $PQF$-closed $C$ and any variety $W$, it holds that $C$ is an irreducible component of $\textrm{log}(W)$ if and only if $\textrm{exp}(C)$ is an irreducible component of $W$.
We first show that for an irreducible variety $W$, the irreducible components of $\textrm{log}(W)$ are of a certain form (Theorem \ref{komponentit}).
When we study a variety $W$, we will from now on write it as $W=T \cap W'$ where $T$ is the minimal torus containing $W$ and $W'$ is a variety not contained in a torus.
Moreover, we will always assume the variety $W'$ is obtained as in the proof of Lemma \ref{leikkaus}.

We will now prove a couple of auxiliary results before proving Theorem \ref{komponentit}, which will give us more information of the irreducible components of $\textrm{log } W$ for any variety $W$.

\begin{lemma}\label{singlek}
Let $C, D \subset V^n$ be irreducible sets.
Suppose $C \subset \bigcup_{k \in K^n} D+k$.
Then, $C \subset D+k$ for a single $k$.
 \end{lemma}

\begin{proof}
By Lemma \ref{irredform}, $D=L\cap m \textrm{log }W$ for some linear set and some variety $W$ that does not branch.
If  $D=L$ for some linear set $L$, then this is Lemma \ref{linearirre}.
 
Suppose now $D=m \textrm{log } W$ for some variety $W$.
Let $k \in K^n$.
We claim that
$$\textrm{log }W + \frac{k}{m}=\textrm{log}(\zeta W),$$
where $\zeta=\textrm{exp}(\frac{k}{m})$ is an $m$:th root of unity.
Indeed, $\textrm{exp}(\textrm{log}(W)+\frac{k}{m})=\zeta W$, so $\textrm{log }W + \frac{k}{m} \subseteq \textrm{log}(\zeta W)$.
On the other hand, suppose $u \in \textrm {log}(\zeta W)$.
Then, $\textrm{exp}(u)=\zeta x$ for some $x \in W$.
Let $u' \in \textrm{log}(W)$ be such that $\textrm{exp}(u')=x$.
Then, there is some $k' \in K^n$ such that
$$u=u'+\frac{k}{m}+k'=(u'+k')+\frac{k}{m} \in \textrm{log}(W)+\frac{k}{m},$$
as wanted.
Thus, $\textrm{log }W + \frac{k}{m}=\textrm{log}(\zeta W)$, and
$$m \textrm{log }W + k=m \textrm{log}(\zeta W).$$
Since there are only finitely many distinct $m$:th roots of unity $\zeta$, the union $\bigcup_{k \in K^n} D+k$ has only 
finitely many distinct members.
Since $C$ is irreducible, it must be contained in one of them.

Suppose next $D= L \cap m \textrm{log } W$ for a linear set $L$ and a variety $W$.
We may without loss assume the union $\bigcup_{k \in K^n} D+k$ does not contain any members $D+k$ such that $C \cap (D+k)=\emptyset$.
By the two above results, there are some $k_1, k_2 \in K^n$ such that
$C \subseteq L+k_1$ and $C \subseteq  m \textrm{log } W +k_2$.
If $L+k_1 \neq L+k_2$, then $(L+k_1) \cap (L+k_2)=\emptyset$ and hence $C \cap (D+k_2)=\emptyset$, a contradiction.
Thus, 
$$C \subseteq (L+k_1) \cap (m \textrm{log } W +k_2)= (L+k_2) \cap (m \textrm{log } W +k_2)=D+k_2,$$
which proves the lemma.
 \end{proof}
 
\begin{lemma}\label{rootsintersect}
Let $W=T\cap W' \subset F^n$ be an irreducible variety, and let $m$ be a natural number.
Then, the $m$:th roots of $W$ are exactly the varieties $T^{\frac{1}{m}}_{(i)} \cap W'^{\frac{1}{m}}_{(j)}$, where
 $T^{\frac{1}{m}}_{(i)}$ goes through the $m$:th roots of $T$ and $W'^{\frac{1}{m}}_{(j)}$ goes through the distinct $m$:th roots of $W'$.
\end{lemma}

\begin{proof}
Let $X$ be a $m$:th root of $W$.
Write $X=T' \cap Y$, where $T'$ is the minimal torus containing $X$ and $Y$ is a variety not contained in any torus, obtained as in the proof of Lemma \ref{leikkaus}.

Since we have $X^m \subset T$, there is some $m$:th root $T^{\frac{1}{m}}_{(i)}$ of $T$ such that $X \cap T^{\frac{1}{m}}_{(i)} \neq \emptyset$.
Then,
$$(T^{\frac{1}{m}}_{(i)} \cap X)^m =T \cap W=T \cap(T \cap W')=T \cap W'=W.$$
Hence, 
$$\textrm{MR}(T^{\frac{1}{m}}_{(i)} \cap X)=\textrm{MR}(W)=\textrm{MR}(X),$$  
so $X \cap T^{\frac{1}{m}}_{(i)}=X$ by irreducibility of $X$.
By minimality of $T'$, we have $T' \subseteq T^{\frac{1}{m}}_{(i)}$.
On the other hand,
$$T'^m \cap Y^m=(T' \cap Y)^m=X^m=W,$$
so $W \subseteq T'^m$ and hence $T \subseteq T'^m$.
Thus, we must have $T' = T^{\frac{1}{m}}_{(i)}$.

Next, we prove that $Y$ is a $m$:th root of $W'$.
Since $Y$ is irreducible, it suffices to show that $Y^m=W'$. 
By Lemma \ref{torusbirat}, we may assume $T$ is given by the equations
$$x_i=c_i, \quad 1\le i \le k$$
for some $k \le n$, where $c_i \in F^{*}$ for each $i$.
Then, as in the proof of Lemma \ref{leikkaus}, $W'$ is given by the ideal $J_W'=\langle I_{W'} \rangle$, where $I_{W'}$ consists of all the polynomials $f \in  F[x_{k+1}, \ldots, x_n]$ such that $f(a)=0$ for a generic point $a$ of $W$.
Let $J_Y$ be the ideal corresponding to $Y$, obtained similarly.

Suppose $b=(b_1, \ldots, b_n) \in T' \cap Y$.
Then, $b^m \in T \cap W'$, and in particular $b^m \in W'$.
Thus, the tuple $(b_{k+1}^m, \ldots, b_n^m)$ is a zero of every polynomial in $I_{W'}$. 
Hence, for every $f \in I_{W'}$, the ideal $I_{Y}$ contains the polynomial $f^*=f(x_{k+1}^m, \ldots, x_n^m)$,
so $Y^m \subseteq W'$.  

On the other hand, suppose $Y^m \subsetneq W'$, and denote by $I(Y^m)$ the ideal corresponding to $Y^m$.
Then, there is some polynomial $g \in I(Y^m) \setminus J_{W'}$.
We may assume $g$ is a generator of $I(Y^m)$ and thus $g \in F[x_{k+1}, \ldots, x_n]$ (note that since the projection of $Y$ on the first $k$ coordinates is $F^k$, also the projection of $Y^m$ is $F^k$).
For each $b \in Y$, $g(b^m)=0$, and thus $g(x_{k+1}^m, \ldots, x_n^m) \in I(Y) \subseteq I(X)$.
Since $W=X^m$, we have $g(c)=0$ for a generic point $c \in W$, and hence by the construction of $W'$, 
 $g \in I_{W'} \subseteq J_{W'}$, a contradiction.
 
Let now $T^{\frac{1}{m}}_{(i)}$ be an $m$:th root of $T$ and $W'^{\frac{1}{m}}_{(j)}$ an $m$:th root of $W'$.
We will show that $T^{\frac{1}{m}}_{(i)} \cap W'^{\frac{1}{m}}_{(j)}$ is an $m$:th root of $W$.
As
$(T^{\frac{1}{m}}_{(i)} \cap W'^{\frac{1}{m}}_{(j)})^m=T \cap W'=W,$
there is some $m$:th root $X$ of $W$ such that $X \cap (T^{\frac{1}{m}}_{(i)} \cap W'^{\frac{1}{m}}_{(j)}) \neq \emptyset.$
Since
$$(X \cap T^{\frac{1}{m}}_{(i)} \cap W'^{\frac{1}{m}}_{(j)})^m=W \cap T \cap W'=W,$$
we have 
$$\textrm{MR}(X \cap T^{\frac{1}{m}}_{(i)} \cap W'^{\frac{1}{m}}_{(j)})=\textrm{MR}(W)=\textrm{MR}(X).$$
As $X$ is irreducible, we must have $X \subseteq T^{\frac{1}{m}}_{(i)} \cap W'^{\frac{1}{m}}_{(j)}.$
We have already proved that $X=T^{\frac{1}{m}}_{(i')} \cap W'^{\frac{1}{m}}_{(j')}$ for some $m$:th roots $T^{\frac{1}{m}}_{(i')}$ and $W'^{\frac{1}{m}}_{(j')}$ of $T$ and $W'$, respectively.
If $T^{\frac{1}{m}}_{(i')} \neq T^{\frac{1}{m}}_{(i)}$, then $T^{\frac{1}{m}}_{(i')} \cap T^{\frac{1}{m}}_{(i)}= \emptyset$.
Thus, we must have $T^{\frac{1}{m}}_{(i')} =T^{\frac{1}{m}}_{(i)}$.

Suppose now $X \subsetneq T^{\frac{1}{m}}_{(i)} \cap W'^{\frac{1}{m}}_{(j)}.$ 
There is some $m$:th root of unity $\zeta=(\zeta_1, \ldots, \zeta_n)$ such that $W'^{\frac{1}{m}}_{(j)}=\zeta W'^{\frac{1}{m}}_{(j')}$.
By the choice of $W'$, we know that the first $k$ coordinates of both $W'^{\frac{1}{m}}_{(j)}$ and $W'^{\frac{1}{m}}_{(j')}$ get all possible values in $F$.
Thus, we may assume $\zeta_1=\ldots=\zeta_k=1$.
On the other hand, the $n-k$ last coordinates of $T^{\frac{1}{m}}_{(i)}$ get all possible values in $F^*$, and hence
$\zeta T^{\frac{1}{m}}_{(i)}=T^{\frac{1}{m}}_{(i)}$.
Thus, by multiplying the equation
$$T^{\frac{1}{m}}_{(i)} \cap W'^{\frac{1}{m}}_{(j')} \subsetneq  T^{\frac{1}{m}}_{(i)} \cap \zeta W'^{\frac{1}{m}}_{(j')}$$
successively by powers of $\zeta$, we get
\begin{eqnarray*}
T^{\frac{1}{m}}_{(i)} \cap W'^{\frac{1}{m}}_{(j')} \subsetneq  T^{\frac{1}{m}}_{(i)} \cap \zeta W'^{\frac{1}{m}}_{(j')}
\subsetneq T^{\frac{1}{m}}_{(i)} \cap \zeta^2 W'^{\frac{1}{m}}_{(j')} \subsetneq
\ldots \subsetneq T^{\frac{1}{m}}_{(i)} \cap \zeta^m W'^{\frac{1}{m}}_{(j')}= T^{\frac{1}{m}}_{(i)} \cap W'^{\frac{1}{m}}_{(j')},
\end{eqnarray*}
a contradiction.
\end{proof}

\begin{theorem}\label{komponentit}
Let $W$ be an irreducible variety in the field sort, and let $W=W' \cap T$ where $T$ is the minimal torus containing $W$, and $W'$ is contained in no torus.
Let $m$ be the level at which $W'$ stops branching and let $W'^{\frac{1}{m}}_{(i)}$ be the $m$:th roots of $W'$.
Let $L$ be linear such that $\textrm{exp}(L)=T$.
Then, the irreducible components of $\textrm{log } W \subset V^n$ are
\begin{displaymath}
(L+k) \cap m \cdot (\textrm{log } W'^{\frac{1}{m}}_{(i)}) \qquad k \in K^n.
\end{displaymath}
\end{theorem}

\begin{proof}
To show that a set $X$ is irreducible, it suffices that if $X \subset C_1 \cup C_2$ for some closed $C_1, C_2$, then $X \subseteq C_1$ or $X \subseteq C_2$. 
We will show that if this holds whenever $C_1, C_2$ are basic closed sets, then $X$ is irreducible.
Indeed, suppose this holds, and let  $C_1, C_2$ be arbitrary closed sets.
Let $C_1=\bigcap_i A_i$ and $C_2= \bigcap_j B_j$ where each $A_i$ and $B_j$ is a basic closed set.
Then, $C_1 \cup C_2= \bigcap_{i,j} (A_i \cup B_j)$.
Thus, if $X \subseteq C_1 \cup C_2$, we must have $X \subseteq A_i \cup B_j$ for each pair $(i,j)$.
Suppose now $X \not\subseteq C_1$.
Then, there is some index $i_0$ such that $X \not\subseteq A_{i_0}$.
But we have $X \subseteq A_{i_0} \cup B_j$ for each $j$, and thus, as the claim holds for basic closed sets, we must have $X \subseteq B_j$ for each $j$, and hence $X \subseteq C_2$.

We divide the analysis into three cases: \\
(a) $W=T$, a torus.\\
(b) $W$ is not contained in any torus\\
(c) $W=T \cap W'$, where $T$ is the minimal torus containing $W$ and $W'$ is not contained in any torus.
 
(a) 

Suppose $W=T$, a torus.
We must show that $L$ is irreducible and that it is maximal irreducible in $\textrm{log } T$.
By Lemma \ref{PQFbasicform}, basic closed sets are finite unions of sets of the form $L'\cap m \textrm{log } W$ where $L'$ is a linear set and $W$ is a variety. 
So suppose  $L \subset \bigcup_{i=1}^r L_i \cap m_i \textrm{log } W_i$ for some linear sets $L_i$ and varieties $W_i$.
We may assume that $L \cap (L_i \cap m_i \textrm{log } W_i) \neq \emptyset$ for each $i$.
 
Denote $m= \prod_{i=1}^r m_i$.
Now, for each $i$,
\begin{displaymath}
m_i \textrm{log }W_i=m \textrm{log }W_i^{\frac{1}{\prod_{j \neq i} m_j}},
\end{displaymath}
where on the right hand side we have the union of all possible $(\prod_{j \neq i} m_j)$:th roots.
Denote now $W_i^{\frac{1}{\prod_{j \neq i} m_j}}=W_i'$ (note that again we have the union of all roots).
Now we have
\begin{displaymath}
L \subset \bigcup_{i=1}^r L_i \cap m \textrm{log }W_i',
\end{displaymath} 
and so
\begin{displaymath}
\frac{L}{m} \subset \bigcup_{i=1}^r \frac{L_i}{m} \cap  \textrm{log }W_i',
\end{displaymath}
which gives us
\begin{displaymath}
\textrm{exp}\left(\frac{L}{m}\right) \subset \bigcup_{i=1}^r \textrm{exp}\left(\frac{L_i}{m}\right) \cap W_i'.
\end{displaymath}
By Lemma \ref{singlemth}, $\textrm{exp}(\frac{L}{m})=T^{\frac{1}{m}}_{(j)}$ for some single $m$:th root $T^{\frac{1}{m}}_{(j)}$ of $T$.
Denote from now on this $m$:th root simply by $T^{\frac{1}{m}}$.
For each $L_i$, there is a torus $T_i$ such that $\textrm{exp}(L_i)=T_i$.
Again, $\textrm{exp}(\frac{L_i}{m})$ is a single $m$:th root of this torus, and we denote it simply by $T_i^{\frac{1}{m}}$. 

Since $T$ is irreducible, $T^{\frac{1}{m}}$ is also irreducible, and thus we must have
\begin{displaymath}
\textrm{exp}\left(\frac{L}{m}\right)=T^{\frac{1}{m}} \subset T_i^{\frac{1}{m}} \cap W_i'
\end{displaymath}
for some single $i$.
This means that 
\begin{displaymath}
\frac{L}{m} \subset \textrm{log } T^{\frac{1}{m}} \subset \left(\bigcup_{k \in K^n} \frac{L_i}{m}+k\right) \cap \textrm{log }W_i',
 \end{displaymath}
so $\frac{L}{m} \subset \bigcup_{k \in K^n}  \frac{L_i}{m}+k$ and $\frac{L}{m} \subset  \textrm{log }W_i'$.
 
Now, as $\frac{L}{m} \subset \bigcup_{k \in K^n}  \frac{L_i}{m}+k$ and $\frac{L}{m}$ is irreducible in the linear topology of the vector space, we must have $\frac{L}{m} \subset \frac{L_i}{m}+k$ for a single $k$ (note that the vector space is a compact structure).  
But now $\frac{L_i}{m}+k=\frac{L_i}{m}$.
Indeed, if $\frac{L_i}{m}+k \neq \frac{L_i}{m}$, then by Remark \ref{linearintersect}, $(\frac{L_i}{m}+k) \cap \frac{L_i}{m}=\emptyset$. 
But then, as $\frac{L}{m} \subset \frac{L_i}{m}+k$, we must have $\frac{L}{m} \cap \frac{L_i}{m}=\emptyset$, which is a contradiction since we assumed that $L \cap (L_i \cap m_i \textrm{log } W_i) \neq \emptyset$. 
Thus, we must have $\frac{L_i}{m}+k=\frac{L_i}{m}$ and so $L \subset L_i$.
Moreover, we have 
\begin{displaymath}
L \subseteq m \textrm{log } W_i'=m_i \textrm{log } W_i,
\end{displaymath}
and hence $L \subseteq L_i \cap m_i \textrm{log } W_i$, as wanted.
 
As for the maximality of $L$,  suppose 
\begin{displaymath}
L \subsetneq I \subset \textrm{log }T= \bigcup_{k \in K^n} L+k,
\end{displaymath}
where $I$ is a closed set.
Now there is some $k$ such that $L+k \neq L$ and $I \cap (L+k) \neq \emptyset$.
By Lemma \ref{singlek}, $I$ is reducible.
 
(b)

Suppose now $W$ is not contained in any torus. 
We need to prove that $m \textrm{log}(W^{\frac{1}{m}})$ is a maximal irreducible subset of $\textrm{log } W$, where $W^{\frac{1}{m}}$ is a choice of the $m$:th root of $W$ for $m$ the level where $W$ stops branching.
Since $m \cdot \textrm{log}(W^{\frac{1}{m}})$ is irreducible if and only if $\textrm{log }W^{\frac{1}{m}}$ is irreducible, we only need to show that the latter is irreducible.
From now on, denote $X=W^{\frac{1}{m}}$. 
We note that $X$ is not contained in any torus as it does not branch.

Suppose now $\textrm{log } X \subseteq \bigcup_{i=1}^r L_i \cap m_i \textrm{log }W_i$, for linear $L_i$ and varieties $W_i$.
Again, we find a number $m'$ and varieties $W'_i$ such that $m' \textrm{log } W_i'=m_i \textrm{log } W_i$ for each $i$.
Then, we get
\begin{displaymath}
X^{\frac{1}{m'}} \subseteq \bigcup_{i=1}^r T_i^{\frac{1}{m'}} \cap W_i',
\end{displaymath}
where $T_i^{\frac{1}{m'}}=\textrm{exp}(\frac{L_i}{m'})$ (note that by Lemma \ref{singlemth} this is just one single $m$:th root).
As $X^{\frac{1}{m'}}$ is irreducible, we must have $X^{\frac{1}{m'}} \subset T_i^{\frac{1}{m'}} \cap W_i'$ for a single $i$.
But now, as $X$ does not branch, $X^{\frac{1}{m'}}$ cannot branch either so it cannot be contained in any (proper) torus.
Hence, we must have $T_i^{\frac{1}{m'}}=({F^{*}})^n$ which means $L_i=V^n$.
Thus, we have $X^{\frac{1}{m'}} \subset W_i'$, so 
\begin{displaymath}
\frac{1}{m'} \textrm{log } X=\textrm{log }X^{\frac{1}{m'}} \subseteq \textrm{log }W_i',
\end{displaymath}
where the equality holds because $X$ does not branch.
This gives us
\begin{displaymath}
\textrm{log } X \subseteq m' \textrm{log } W_i' = m_i \textrm{log } W_i=V^n \cap m_i \textrm{log } W_i= L_i \cap m_i \textrm{log } W_i.
\end{displaymath}

For maximality, suppose
\begin{displaymath}
W^{\frac{1}{m}}_{(i)} \subsetneq I \subset \textrm{log }W = \bigcup_i m \textrm{log } W_{(i)}^{\frac{1}{m}}. 
\end{displaymath}
Now  $I \cap (W^{\frac{1}{m}}_{(j)} \setminus W^{\frac{1}{m}}_{(i)}) \neq \emptyset$ for some $j \neq i$.
Thus, 
\begin{displaymath}
I=\bigcup_i \left\{x \in I \, | \, \textrm{exp}\left(\frac{x}{m}\right) \in W^{\frac{1}{m}}_{(i)}\right\},
\end{displaymath}
where at least two of the sets are distinct, so $I$ is not irreducible.

(c)

Suppose now $W=T \cap W'$ where $T$ is the minimal torus containing $W$ and $W'$ is a variety not contained in any torus, obtained as in the proof of Lemma \ref{leikkaus}. 
Let $L$ be such that $\textrm{exp}(L)=T$, and let $W'^{\frac{1}{m}}$ be a choice of the $m$:th root of $W'$ for $m$ the level where $W'$ stops branching.
We will prove that $L \cap m \textrm{log} W'^{\frac{1}{m}}$ is a maximal irreducible subset of $\textrm{log }W$.

We note first that $L \cap m \textrm{log} W'^{\frac{1}{m}}$ is irreducible if and only if $\frac{L}{m} \cap \textrm{log} W'^{\frac{1}{m}}$ is.
We denote $X=W'^{\frac{1}{m}}$ and prove that $\frac{L}{m} \cap \textrm{log X}$ is irreducible. (Note that $X$ does not branch.)

Suppose $\frac{L}{m} \cap \textrm{log} X \subset \bigcup_{i=1}^r L_i \cap m_i \textrm{log }W_i$ for linear $L_i$ and varieties $W_i$.
As before, we find a number $m'$ and varieties $W_i'$ so that $m_i \textrm{log }W_i=m' \textrm{log }W_i'$ for each $i$.
We get that 
\begin{eqnarray}\label{esityss}
\frac{L}{mm'} \cap \frac{1}{m'} \textrm{log} X \subset \bigcup_{i=1}^r \frac{L_i}{m'} \cap \textrm{log }W_i',
\end{eqnarray}
and thus
\begin{displaymath}
T^{\frac{1}{mm'}} \cap X^{\frac{1}{m'}} \subset \bigcup_{i=1}^r T_i^{\frac{1}{m'}} \cap W_i',
\end{displaymath}
where $T^{\frac{1}{mm'}}=\textrm{exp}(\frac{L}{mm'})$ is a single $mm'$:th root of $T$ and for each 
 $i$, $T_i^{\frac{1}{m'}}=\textrm{exp}(\frac{L_i}{m'})$ is a single $m'$:th root of $T_i$ (by Lemma \ref{singlemth}).
As $X$ does not branch, also $X^{\frac{1}{m'}}$ is a single $m'$:th root of $X$, and thus a single $mm'$:th root of $W'$.
By Lemma \ref{rootsintersect},  $T^{\frac{1}{mm'}} \cap X^{\frac{1}{m'}}$ is a single $mm'$:th root of $T \cap W'$,
and thus irreducible.
Hence,
\begin{eqnarray}\label{singlei}
T^{\frac{1}{mm'}} \cap X^{\frac{1}{m'}} \subseteq T_i^{\frac{1}{m'}} \cap W_i'
\end{eqnarray}
for a single $i$.

As $W'$ is not contained in any torus, $X^{\frac{1}{m'}}$ is not contained in any torus either (otherwise it, and thus $W'$, would branch infinitely).
Moreover, $T^{\frac{1}{mm'}}$ is a torus.
We claim that it is the minimal torus containing $T^{\frac{1}{mm'}} \cap X^{\frac{1}{m'}}$.
Suppose $T'$ is some torus such that $T^{\frac{1}{mm'}}  \cap X^{\frac{1}{m'}} \subseteq T'$.
Then, $T \cap W' \subseteq T'^{mm'}$, and thus $T \subseteq T'^{mm'}$ as $T$ is the minimal torus containing $W=T \cap W'$.
Hence, $T^{\frac{1}{mm'}}  \subseteq (T'^{mm'})^{\frac{1}{mm'}}_{(i)}$ for some $mm'$:th root $(T'^{mm'})^{\frac{1}{mm'}}_{(i)}$ of $T'^{mm'}$.
The torus $T'$ is one of the $mm'$:th roots of $T'^{mm'}$, and as $T^{\frac{1}{mm'}}\cap X^{\frac{1}{m'}} \subseteq T'$, we have $T' \cap T^{\frac{1}{mm'}} \neq \emptyset$.
Since the distinct $mm'$:th roots of the torus $T'$ do not intersect, we must have $T^{\frac{1}{mm'}} \subseteq T'$.
Thus, $T^{\frac{1}{mm'}}$ is the minimal torus containing $T^{\frac{1}{mm'}} \cap X^{\frac{1}{m'}}$.

Since $X$ does not branch, we have $\textrm{log } X^{\frac{1}{m'}}=\frac{1}{m'} \textrm{log } X$, and thus by taking
logarithms we get from (\ref{singlei})
\begin{displaymath} 
\frac{L}{mm'} \cap \frac{1}{m'}\textrm{log}{X} \subseteq \left( \bigcup_{k \in K^n}\frac{L_i}{m'}+k\right) \cap \textrm{log}W_i'.
\end{displaymath}
Since $T^{\frac{1}{mm'}}$ is the minimal torus containing $T^{\frac{1}{mm'}} \cap X^{\frac{1}{m'}}$, (\ref{singlei}) also gives $T^{\frac{1}{mm'}} \subseteq T_i^{\frac{1}{m'}}$ and hence
\begin{displaymath}
\frac{L}{mm'} \subseteq   \bigcup_{k \in K^n} \frac{L_i}{m'}+k.
\end{displaymath}
Since the vector space is a compact structure and linear sets are irreducible there, we must have
$\frac{L}{mm'} \subset \frac{L_i}{m'}+k$ for a single $k$.
We may assume that the union in the representation (\ref{esityss}) does not contain any redundant members and thus $\frac{L}{mm'} \cap \frac{L_i}{m'} \neq \emptyset$.
Then, $\frac{L_i}{m'}+k=\frac{L_i}{m'}$ and thus $\frac{L}{mm'} \subset \frac{L_i}{m'}$.
It follows that
\begin{displaymath}
\frac{L}{mm'} \cap \frac{1}{m'}\textrm{log}{X} \subset \frac{L_i}{m'} \cap \textrm{log}W_i',
\end{displaymath}
and hence 
\begin{displaymath}
\frac{L}{m} \cap \textrm{log}{X} \subset L_i \cap m_i \textrm{log}W_i,
\end{displaymath}
as wanted.

For maximality, let $L \cap m\textrm{log} W'^{\frac{1}{m}}_{(i)} \subsetneq I \subseteq \textrm{log } W$, where $I$ is a closed set.
Then, there either is some $k \in K^n$ such that $L+k \neq L$ and $(L+k) \cap I \neq \emptyset$, or $I \cap (m\textrm{log }W'^{\frac{1}{m}}_{(j)} \setminus m \textrm{log }W'^{\frac{1}{m}}_{(i)}) \neq \emptyset$ for some $j \neq i$. 
In both cases, $I$ is reducible.
\end{proof}

\begin{corollary}\label{olennainen}
If $W\subseteq F^n$ is a variety and $C\subseteq V^n$ is an irreducible subset of the cover sort, then $C$ is an irreducible component of $\textrm{log }W$ if and only if $\textrm{exp}(C)$ is an irreducible component of $W$.
\end{corollary}

\begin{proof}
Suppose first $W$ is irreducible.
If $\{C_i \, | \, i \in \omega\}$ are the irreducible components of $\textrm{log W}$, then by Theorem \ref{komponentit}, $\textrm{exp}(C_i)=W$ for each $C_i$.

Assume now $\textrm{exp}(C)=W$.
Then, $\textrm{log }W=\bigcup_{k \in K^n} C+k$.
Let $X$ be an irreducible component of $\textrm{log}(W)$.
By Lemma \ref{singlek}, $X \subseteq C+k$ for some $k \in K^n$.
Since $X$ is a maximal irreducible subset of $\textrm{log}(W)$, we must have $X=C+k$.
If there were some irreducible set $Y$ such that $C \subsetneq Y \subseteq \textrm{log}(W)$, then $X \subsetneq Y+k \subseteq \textrm{log}(W)$, a contradiction.
Thus, $C$ is an irreducible component of $\textrm{log}(W)$.

Suppose now
$$W=W_1 \cup \ldots \cup W_r,$$
written as the union of its irreducible components.
Then, 
$$\textrm{log }W=\textrm{log }W_1 \cup \ldots \cup \textrm{log }W_r.$$
We will show that $C$ is an irreducible component of $\textrm{log }W$ if and only if it is an irreducible component of $\textrm{log } W_i$ for some $i$.
Indeed, if $C$ is an irreducible component of $\textrm{log }W_i$ for some $i$ and not an irreducible component of $\textrm{log }W$, then there is some irreducible $C' \subset \textrm{log }W$ so that $C \subsetneq C'$ and $C' \not\subseteq \textrm{log }W_i$.
But now we may write $C'=\bigcup_{i=1}^r (C' \cap \textrm{log }W_i)$, where at least two members of the union are distinct. This contradicts the irreducibility of $C'$.

On the other hand, any irreducible component of $\textrm{log } W$ must be contained in $\textrm{log }W_i$ for some $i$, and is hence an irreducible component of $\textrm{log }W_i$.
Thus, $C$ is an irreducible component of $\textrm{log }W$ if and only if $C$ is an irreducible component of $\textrm{log } W_i$ for some $i$.
Since we assumed the claim holds for irreducible varieties, this is equivalent to the statement that $\textrm{exp}(C)=W_i$ for some $i$.
\end{proof}

\begin{remark}\label{expcomps}
Note that from Corollary \ref{olennainen} it follows that any irreducible PQF-closed set $C$ on the cover sort is an irreducible component of $\textrm{log}(\textrm{exp}(C))$.  
\end{remark}

\section{Dimension Analysis}

In this section, we define dimensions for the PQF-closed sets much in the same way that they are defined in the Noetherian case.
We then prove that the dimension of a PQF-closed set is equal to the dimension of its image under the map exp (where the latter is defined as in Chapter 1), and that the axiom (Z3) of Zariski geometries holds for PQF-closed sets on the cover.

\begin{lemma}\label{infchains}
There are no infinite descending chains of irreducible closed subsets of the cover sort $V$.
\end{lemma}

\begin{proof}
Suppose $C_1 \supsetneq C_2 \supsetneq \ldots$ is an infinite descending chain of irreducible closed sets on $V$.
Denote $W_i=\textrm{exp}(C_i)$.
Then, $W_i$ is an irreducible variety, and we get that
\begin{displaymath}
W_1 \supseteq W_2 \supseteq \ldots
\end{displaymath}
Since there are no infinite descending chains of closed sets on the field sort $F$, there is a number $m$ such that
$$W_m=W_{m+1}=W_{m+2}= \ldots.$$
Thus,  $\textrm{exp}(C_m)=\textrm{exp}(C_{m+1})$.
Both $C_m$ and $C_{m+1}$ are irreducible, and thus, by Remark \ref{expcomps}, they both are irreducible components of $\textrm{log}(W_m)$.
But this is impossible as $C_{m+1} \subsetneq C_m$ and irreducible components were defined to be maximal irreducible subsets.
\end{proof}
 
\begin{definition}
If $C$ is an irreducible, closed and nonempty set on the cover sort, we define the \emph{dimension} of $C$ inductively as follows: 
\begin{itemize}
\item $\textrm{dim}(C) \ge 0$,
\item $\textrm{dim}(C)=\textrm{sup }\{\textrm{dim}(F)+1 \, | \, F \subsetneq C, \textrm{$F$ closed, irreducible and nonempty } \}$.
\end{itemize}
\end{definition}

We define the concepts of locus and rank the same way as it was done for Zariski geometries in Chapter 1.

\begin{definition}
Let $a \in V^n$, $A \subset V$.
By the \emph{locus} of $a$ over $A$, we mean the smallest PQF-closed subset definable over $A$ containing $a$.
When not specified, we assume the set $A$ to be empty. 

We define $\textrm{rk}(a/A)=\textrm{dim}(C)$, where $C$ is the locus of $a$ over $A$.
We write $\textrm{rk}(a)$ for $\textrm{rk}(a/\emptyset)$.
\end{definition}

We now prove that the dimension of an irreducible PQF-closed set equals the dimension of its image under the map exp (where the latter is calculated as in Chapter 1).

\begin{lemma}\label{dimensio}
Let $C \subset V^n$ be an irreducible PQF-closed subset of the cover sort.
Then,
\begin{displaymath}
\textrm{dim}(C)=\textrm{dim}(\textrm{exp}(C)).
\end{displaymath} 
\end{lemma}

\begin{proof} 
Let $\textrm{dim}(C)=n$, and let $C_0 \subsetneq C_1 \subsetneq \ldots \subsetneq C_n=C$ be a maximal chain of irreducible closed sets in $C$.
Then, 
\begin{displaymath}
\textrm{exp}(C_1) \subsetneq \ldots \subsetneq \textrm{exp}(C_n)=\textrm{exp}(C),
\end{displaymath}
where the fact that the inclusions are proper is shown similarly as in the proof of Lemma \ref{infchains}.
Thus, $\textrm{dim}(C) \le \textrm{dim}(\textrm{exp}(C))$.

We prove the other inequality by induction on $\textrm{dim}(\textrm{exp}(C))$.
Suppose first $\textrm{dim}(\textrm{exp}(C))=0$.
Then, $\textrm{exp}(C)=\{x\}$ for some $x \in ({F^*})^n$.
Let $v \in V^n$ be such that $\textrm{exp}(v)=x$.
By Corollary \ref{olennainen},
$C$ is an irreducible component of $\textrm{log}(\{x\})$.
Hence, $C=\{v+k\}$ for some $k \in K$, so $\textrm{dim}(C)=0$.

Suppose now $\textrm{dim}(\textrm{exp}(C)) \le \textrm{dim}(C)$ whenever 
$\textrm{dim}(\textrm{exp}(C))=n$.
Assume $\textrm{dim}(\textrm{exp}(C))=n+1$.
Denote $\textrm{exp}(C)=W$.
There is some irreducible $W' \subsetneq W$ such that $\textrm{dim}(W')=n$.
Since $W' \subsetneq W$, there is some irreducible $Y \subsetneq C$ such that $\textrm{exp}(Y)=W'$.
By the inductive hypothesis, 
$$\textrm{dim}(Y) \ge \textrm{dim}(W')=n.$$
As $Y \subsetneq C$, we have $\textrm{dim}(C) \ge n+1$.
 \end{proof}

In a Noetherian topology, closed sets have finitely many irreducible components.
We now prove that in our context the PQF-closed sets have countably many irreducible components.

\begin{lemma}
Any basic PQF-closed set has countably many irreducible components.
\end{lemma}

\begin{proof}
Let $D$ be basic a PQF-closed set.
We will do a construction that yields all the irreducible components of $D$.
Let $W_0=\textrm{exp}(D)$.
For each irreducible component $W_0'$ of $W_0$, let $C_{0i}$, $i<\omega$, be the irreducible components of $\textrm{log }W_0'$ (by Theorem \ref{komponentit}, there are countably many of these).
For each $i$, consider $C_{0i} \cap D$.
If $C_{0i} \cap D=C_{0i}$, then $C_{0i}$ is an irreducible component of $D$.
Indeed, $C_{0i}$ is irreducible, and if we have $C_{0i} \subseteq C_{0i}'$ for some irreducible set $C_{0i}' \subset D$, then we must have $C_{0i}' \subseteq \textrm{log }W_0'$ (otherwise we could write $C_{0i}'=\bigcup_{j=1}^n \{x \, | \, \textrm{exp}(x) \in W_j'\}$, where $W_1', \ldots, W_n'$ are the irreducible components of $W_0$).
Since $C_{0i}$ is a maximal irreducible set in $\textrm{log }W_0'$, we have $C_{0i}=C_{0i}'$.
 
If $C_{0i} \cap D \neq C_{0i}$, then $W_1=\textrm{exp}(C_{0i} \cap D) \subsetneq W_0'$, and hence $\textrm{dim}(W_1) < \textrm{dim}(W')$.
We now repeat the process with $C_{0i} \cap D$ in place of $D$, looking at the irreducible components $W_1'$ of $W_1$.
At each iteration step, the dimension of the image under the exponential map drops, so eventually the process must terminate.
Thus, at some point we will have found an irreducible variety $W_n'$ such that $D \cap C_{ni}=C_{ni}$ for some 
irreducible component $C_{ni}$ of $\textrm{log }W_n'$.
We claim that $C_{ni}$ is an irreducible component of $D$.

Suppose there is some irreducible set $Y$ such that $C_{ni} \subseteq Y \subseteq D$.
Since $Y$ is irreducible, we have $Y \subset (D \cap C_{0i})$ for one of the irreducible components $C_{0i}$ of $\textrm{log }W_0'$.
Since $C_{ni} \subseteq Y$, the component $C_{0i}$ must be the one containing $C_{ni}$.
By the construction, we have either $Y=C_{0i}$, in which case $n=0$ and $D=C_{0i}$,  
or $Y \subseteq C_{0i} \cap D \subsetneq C_{0i}$. 
In the latter case, $Y \subseteq \textrm{log }W_1'$, where $W_1'$ is the irreducible component of $W_1$ containing $\textrm{exp}(C_{ni})$, and thus $Y \subseteq C_{1i}$ for some irreducible component $C_{1i}$ of $\textrm{log }W_1'$.
Now we repeat the process.
Going down like this, we eventually get $Y=C_{ni}$. 
 
At each step, we need to consider the finitely many irreducible components of some variety $X$ and the countably many irreducible components of $\textrm{log }X$.
Therefore, in the process, we have obtained countably many irreducible components of $D$.
We will show that these are all the irreducible components.
For this, it suffices that an arbitrary irreducible set $C \subseteq D$ is included in one of the components obtained in the construction.
We have $C \subseteq \textrm{log } W_r'$ for some $r$.
Assume moreover that $r$ is the largest possible such number.
By Theorem \ref{komponentit}, $\textrm{log } W_r'$ can be written as the union of its irreducible components, which are of the form $(L+k) \cap m \textrm{log}(W^{\frac{1}{m}}_{(i)})$, where $L$ is a linear set, $k \in K$, and $W$ is a variety not contained in any torus.
By Lemma \ref{singlek}, $C \subseteq L+k$ for a single $k$.
Since $W$ has only finitely many $m$:th roots, we also have $C \subset m \textrm{log}(W^{\frac{1}{m}}_{(i)})$ for a single $i$.
Thus, $C$ is contained in a single irreducible component $C_{rj}$ of $\textrm{log }W_r'$.
Since $r$ was assumed to be maximal, $C_{rj} \cap D=C_{rj}$, and thus $C_{rj}$ is an irreducible component of $D$.  
\end{proof}

\begin{remark}\label{ffff}
Let $D$ be a PQF-closed set, and let $C$ be an irreducible component of $D$.
Then, $\textrm{exp}(C) \subseteq \textrm{exp}(D)$, so 
$$\textrm{dim}(C)=\textrm{dim}(\textrm{exp}(C)) \le \textrm{dim}(\textrm{exp}(D)).$$ 
\end{remark}

\begin{definition}\label{dimgendef}
For an arbitrary PQF-closed set $C$, we define $\textrm{dim}(C)$ to be the maximum dimension of the irreducible components of $C$ (by Remark \ref{ffff}, this is a finite number). 

For any PQF-closed set $C$, we say that an element $a \in C$ is \emph{generic} if $\textrm{rk}(a/ \textrm{log }F_0)=\textrm{dim}(C)$, where $F_0$ is the smallest algebraically closed subfield of $F$ such that $C$ is definable over $\textrm{log}(F_0)$.
\end{definition}
 
It is now easy to see that for an arbitrary PQF-closed set $D$, it holds that $\textrm{dim}(D)=\textrm{dim}(\textrm{exp}(D))$.
 
In the context of Chapter 3, we gave a different definition for generic elements (in the context of a quasiminimal class).
We will later show that we can assume that our structure $V$ is a monster model for a quasiminimal class and that the above definition actually coincides with the one given in Chapter 3. 

Now we can prove that the Dimension Theorem (i.e. axiom (Z3) of Zariski geometries) holds on the cover sort.
 
\begin{theorem}[Dimension Theorem, \cite{Lucy}]\label{dimthm}
Let $C_1, C_2 \subset V^n$ be closed and irreducible.
Let $X$ be a non-empty irreducible component of $C_1 \cap C_2$.
Then, $\textrm{dim}(X) \ge \textrm{dim}(C_1)+\textrm{dim}(C_2)-n$.
\end{theorem}

\begin{proof}
Let $C_1, C_2 \subset V^n$ be closed and irreducible.
Let $X$ be a non-empty irreducible component of $C_1 \cap C_2$.
The theorem holds on the field sort, and thus, by Lemma \ref{dimensio}, it suffices to show that $\textrm{exp}(X)$ is an irreducible component of $\textrm{exp}(C_1) \cap \textrm{exp}(C_2)$.
For $i=1,2$, denote
$C_i=L_i \cap m_i \textrm{log } W_i$, where $L_i$ is linear, $W_i$ does not branch.
By Corollary \ref{olennainen} and Theorem \ref{komponentit}, we may assume that $T_i=\textrm{exp}(L_i)$ is the minimal torus containing $\textrm{exp}(C_i)=T_i \cap W_i^{m_i}$ and that $W_i^{m_i}$ is chosen as in the proof of Lemma \ref{leikkaus}.
Now,
$$\textrm{exp} (C_1) \cap \textrm{exp}(C_2)=T_1 \cap T_2 \cap W_1^{m_1} \cap W_2^{m_2}.$$
 
Let $W_1^{m_1} \cap W_2^{m_2}=T_3 \cap W$ where $T_3$ is a torus and $W$ is a variety not contained in any torus, chosen as in the proof of Lemma \ref{leikkaus}.
Set $T=T_1 \cap T_2 \cap T_3$ so that
 $$T_1 \cap T_2 \cap W_1^{m_1} \cap W_2^{m_2}=T \cap W.$$
 Let $X_1, \ldots, X_r$ be the irreducible components of $T\cap W$.
 For each $i$, we may write $X_i= T \cap T_i' \cap Y_i$, where $T_i'$ is the minimal torus containing $X_i$ and $Y_i$ is not contained in any torus, again chosen as in the proof of Lemma \ref{leikkaus}.
 Let $L_3$ be a linear set so that $\textrm{exp}(L_3)=T_3$ and $L_1 \cap L_2 \cap  L_3 \neq \emptyset$.  
 By Theorem \ref{komponentit}, the irreducible components of $\textrm{log}(T\cap W)$ are the sets
 $$((L \cap L_i')+k) \cap m_i' \textrm{log } {Y_i^{\frac{1}{m_i'}}}_{(j)},$$
 where $L=L_1 \cap L_2 \cap L_3$, and $L_i'$ is a linear set  such that $\textrm{exp}(L_i')=T_i'$ and $L_i' \cap L \neq \emptyset$, and $m_i'$ is such that the $m_i'$:th roots of $Y_i$ no longer branch.
 
Let $W_1^{\frac{1}{m_2}}$ and $W_2^{\frac{1}{m_1}}$ be the unique $m_2$:th and $m_1$:th roots of $W_1$ and $W_2$, respectively.
Then,
\begin{eqnarray*}
C_1 \cap C_2 &=& L_1 \cap L_2 \cap m_1 m_2 \textrm{log}(W_1^{\frac{1}{m_2}} \cap W_2^{\frac{1}{m_1}}),
\end{eqnarray*}
where
$$W_1^{\frac{1}{m_2}} \cap W_2^{\frac{1}{m_1}}=(W_1^{m_1} \cap W_2^{m_2})^{\frac{1}{m_1 m_2}}=(T_3 \cap W)^{\frac{1}{m_1 m_2}}$$
for a suitable choice of the $m_1 m_2$:th root.
Thus, for suitable choices of the $m_1 m_2$:th roots,
\begin{eqnarray*}
C_1 \cap C_2 &=& L_1 \cap L_2 \cap m_1 m_2 \textrm{log}((T_3 \cap W)^{\frac{1}{m_1 m_2}}) \\
&=& L_1 \cap L_2 \cap m_1 m_2 \textrm{log}((T_1 \cap T_2 \cap T_3 \cap W)^{\frac{1}{m_1 m_2}}) \\
&=&L_1 \cap L_2 \cap m_1 m_2 \textrm{log}(T^{\frac{1}{m_1 m_2} } \cap \bigcup_{i=1}^r (T_i'^{\frac{1}{m_1 m_2}}\cap Y_i^{\frac{1}{m_1 m_2}})).
\end{eqnarray*}
Hence, by Theorem \ref{komponentit}, irreducible components of $C_1 \cap C_2$ are of the form
$$((L_1 \cap L_2 \cap L_3 \cap L_i') +k) \cap m_i' \textrm{log } {Y_i ^{\frac{1}{m_i'}}}_{(j)},$$
and thus each one of them is an irreducible component of $\textrm{log}(T\cap W)$.
(Note that for each $i$, we have $T_i' \cap Y_i \subseteq T \cap T_i' \cap Y_i$.)
By Corollary  \ref{olennainen}, $\textrm{exp}(X)$ is an irreducible component of $T\cap W=\textrm{exp}(C_1 \cap C_2)$, as wanted.
\end{proof}

\section{Bounded closures}
 
\begin{definition}
Define a closure operator, cl, on $\mathcal{P}(V)$ so that for any $A \subset V$, 
$$\textrm{cl}(A)=\textrm{log}(\textrm{acl}(\textrm{exp}(A))).$$
\end{definition} 

It is easy to see that $\mathcal{P}(V)$ forms a pregeometry with respected to cl.
In this section we will we see that after adding countably many symbols to our language, $(V, \textrm{cl})$ will be a quasiminimal pregeometry structure in the sense of Definition \ref{quasiminclass} and  \cite{monet}. 
Then, we will show that the closure operator defined will coincide with the bounded closure,
i.e. that for any $A \subset V$,
$$\textrm{cl}(A)=\textrm{bcl}(A).$$
  
\begin{remark}\label{pregstruct}
Now it is easy to see that $(V, \textrm{cl})$ forms a pregeometry determined by the language (i.e. if $a \in V$, $\overline{b} \in V^n$ for some $n$ and $(a, \overline{b})$ and $(a', \overline{b'})$ have the same quantifier-free type, then $a \in \textrm{cl}(\overline{b})$ if and only if $a' \in \textrm{cl}(\overline{b'})$) and that $V$ is infinite-dimensional with respect to cl. 
Also, if $A \subset V$ is finite, then $\textrm{cl}(A)$ is countable.

Moreover, $V$ has a unique generic (quantifier-free) type with respect to this pregeometry.
Suppose that $H, H' \subset V$ are countable subsets closed with respect to the pregeometry (we write ``cl-closed" for this, to avoid confusion with ``PQF-closed"), enumerated so that $\textrm{tp}_{q.f}(H)=\textrm{tp}_{q.f} (H')$, where $\textrm{tp}_{q.f.}$ denotes the quantifier free type.  
Suppose $a \in V \setminus H$ and $a' \in V \setminus H'$.
As $H$ is cl-closed, $a$ cannot satisfy any linear dependencies over $H$ and $\textrm{exp}(a)$ must be transcendental over $\textrm{exp}(H)$.
Also, $\textrm{exp}(\frac{a}{n})$ must be transcendental over $\textrm{exp}(H)$ for every $n$, as otherwise we would have $\textrm{exp}(\frac{a}{n}) \in \textrm{exp}(H)$, and thus $\textrm{exp}(a)=(\textrm{exp}(\frac{a}{n}))^n \in \textrm{exp}(H)$.
The same holds with respect to $a'$ and $H'$.
Thus, $\textrm{tp}_{q.f.}(H,a)=\textrm{tp}_{q.f.}(H',a')$.

Also, $V$ is $\aleph_0$-homogeneous over closed sets. 
Let $H, H' \subset V$ be countable cl-closed sets enumerated so that $\textrm{tp}_{q.f}(H)=\textrm{tp}_{q.f.}(H')$, and let $\bar{b}$, $\bar{b'}$ be finite tuples from $V$ such that $\textrm{tp}_{q.f.}(H, \bar{b})=\textrm{tp}_{q.f.}(H', \bar{b'})$, and let $a \in \textrm{cl}(H, \bar{b})$.
Then, in particular, $\textrm{exp}(H)$ and $\textrm{exp}(H')$ are algebraically closed fields.
Suppose first $a$ satisfies some $\mathbb{Q}$-linear equation over $H\bar{b}$.
Then, $a$ is the unique solution of that equation and we may find a unique $a' \in V$ so that $a'$ satisfies the same equation over $H' \bar{b'}$.
Then, $\textrm{tp}_{q.f.}(H, \bar{b}, a)=\textrm{tp}_{q.f.}(H', \bar{b'}, a')$.

If $a$ does not satisfy any $\mathbb{Q}$-linear equation over $H\bar{b}$, then $W$, the locus of $\textrm{exp}(\bar{b},a)$ over $\textrm{exp}(H)$ is not contained in any torus definable over $\textrm{exp}(H)$.
Since $\textrm{exp}(H)$ is an algebraically closed field, it follows from Model Completeness that $W$ is not contained in any torus.
Thus, by Theorem \ref{zilbertheorem}, there is some natural number $m$ such that the $m$:th roots of $W$ (over $\textrm{exp}(H)$) no longer branch.
Let $W^{\frac{1}{m}}_{(i)}$ be such that $\textrm{exp}(\frac{1}{m}(\bar{b},a)) \in W^{\frac{1}{m}}_{(i)}$.
Then, there is a variety $W'$, a $m$:th root $W'^{\frac{1}{m}}_{(i)}$ defined over $\textrm{exp}(H')$ using the same equations that define $W$ and $W^{\frac{1}{m}}_{(i)}$ over $\textrm{exp}(H)$ such that $W'^{\frac{1}{m}}_{(i)}$ doesn't branch, and there is some $a' \in V$ such that $\textrm{exp}(\bar{b'}, a')$ is a generic point of $W'$ and $\textrm{exp}(\frac{1}{m}(\bar{b'},a')) \in W'^{\frac{1}{m}}_{(i)}$.
Now,  $\textrm{tp}_{q.f.}(H, \bar{b}, a)=\textrm{tp}_{q.f.}(H', \bar{b'}, a')$.

If we add constants for the elements of $\textrm{log}(\overline{\mathbb{Q}})$ to our language (here $\overline{\mathbb{Q}}$ stands for the field of algebraic numbers), we can use similar arguments to show that $V$ is $\aleph_0$-homogeneous over the empty set (using the $\emptyset$-definable algebraically closed field $\overline{\mathbb{Q}}$ in place of $\textrm{exp}(H)$). 
Thus, after adding countably many symbols to our language, $V$ will be a quasiminimal pregeometry structure as defined in \cite{monet}.
From now on we will assume we have added these symbols.
\end{remark}
 
As in Chapter 2, section 5, we can construct an AEC $\K(V)$ from some model $V$ of the theory of the covers.
Moreover, we may view $V$ as the monster model for the class.
Then, by the results of Chapter 2, we have an independence calculus for $\mathcal{K}(V)$.
From now on, when we write $V$, we always assume it to be the monster model for an AEC constructed this way.
  
\begin{lemma}\label{bcl} 
For any $A \subset V$, 
$$\textrm{bcl}(A)=\textrm{cl}(A).$$
\end{lemma}

\begin{proof}
Clearly $\textrm{cl}(A)=\textrm{log}(\textrm{acl}(\textrm{exp}(A))) \subseteq \textrm{bcl}(A)$.
 
On the other hand, suppose $v \notin \textrm{cl}(A).$
As in Remark \ref{pregstruct}, we see that $v$ has the same quantifier-free type with any other element not in $\textrm{cl}(A)$.
From quantifier elimination and \cite{monet} it follows that if $v' \notin \textrm{cl}(A)$, then $t^g(v/A)=t^g(v'/A)$.
Since there are uncountably many such $v'$, we have $v \notin \textrm{bcl}(A)$.
\end{proof}

From now on we will always write ``bcl" for ``cl" and ``bcl-closed" for ``cl-closed".
   
\begin{remark}\label{ranksama}
The dimension obtained from the pregeometry agrees with the one defined topologically with respect to the PQF-closed sets, i.e. for any $a \in V^n$, $A \subset V$, it holds that $$\textrm{rk}(a/A)=\textrm{dim}_{\textrm{bcl}}(a/A).$$ 
To see this, let $X \subset V^n$ be the smallest $A$-definable PQF-closed set containing $a$.
Then, $\textrm{exp}(X)$ is definable by some first-order formula $\phi(x_1, \ldots, x_n)$ in the field language over $\textrm{exp}(A)$.
By Corollary \ref{variety},  the set $\textrm{exp}(X)$ is Zariski closed, and thus there is some set of polynomials $S=\{p_1, \ldots, p_m\}$ such that $\textrm{exp}(X)$ is the zero locus of $S$.

Let $a_1, \ldots, a_r \in F$ be the coefficients of the polynomials $p_1, \ldots, p_m$ that are in $F \setminus \textrm{acl}(\textrm{exp}(A))$.
Replacing these by the variables $y_1, \ldots, y_r$, we may write a formula in the field  language with parameters from $\textrm{exp}(A)$ expressing
$$\exists y_1 \cdots \exists y_r (\phi(x_1, \ldots, x_n) \leftrightarrow \bigwedge_{i=1}^m p_i(y_1, \ldots, y_r, x_1, \ldots, x_n)=0).$$
This formula holds in $F$.
By Model Completeness in algebraically closed fields, it holds already in 
$\textrm{acl}(\textrm{exp}(A))$.  
Thus there are some polynomials $p_1', \ldots, p_m'$ over $\textrm{acl}(\textrm{exp}(A))$ such that $\textrm{exp}(X)$ is the zero locus of the set consisting of these polynomials.

Hence, we have seen that $\textrm{exp}(X)$ is definable as a variety over $\textrm{acl}(\textrm{exp}(A))$, the smallest algebraically closed field containing $A$.
This implies that $\textrm{exp}(X)$ is the locus (in the field sense) of $\textrm{exp}(a)$ over $\textrm{acl}(\textrm{exp}(A))$, so
$$\textrm{dim}(X)=\textrm{dim}(\textrm{exp}(X))=\textrm{dim}_{\textrm{acl}}(\textrm{exp}(a)/\textrm{acl}(\textrm{exp}(A)))=\textrm{dim}_{bcl}(a/A),$$
where the second equality follows form the fact that in algebraically closed fields, the dimension in Zariski topology, Morley rank, and the dimension with respect to the pregeometry defined by $acl$ are all equal (see e.g. \cite{MarIntro}, \cite{HrZi}).
The last equality follows from the definition of the closure operator cl and from  Lemma \ref{bcl}.

Note that since $U(a/A)=\textrm{dim}_{bcl}(a/A)$, this also gives us
$$\textrm{rk}(a/A)=U(a/A).$$

It also follows that the definition of generic elements of a PQF-closed set given in \ref{dimgendef} coincides with the one given in Definition \ref{genericeka}
\end{remark}
  
\section{Axioms for irreducible sets in the general framework}

In this section we show that the cover satisfies the axioms for Zariski-like structures if we take the irreducible (in the topological sense) PQF-closed sets that are PQF-closed already over $\emptyset$ to be the irreducible sets in the definition of a Zariski-like structure (Definition \ref{zariskilike}). 
Note that we have added the elements of $\textrm{log}(\overline{\mathbb{Q}})$ as constants in our language, so in practice this means that we are allowed to use parameters from $\textrm{bcl}(\emptyset)$ when defining these sets with positive, quantifier free formulae.
In the rest of the section, when we say ``irreducible", we mean irreducible in the sense of the PQF-topology (and not necessary PQF-closed over $\emptyset$).
If we mean irreducible and PQF-closed over $\emptyset$ (i.e. ``irreducible" in the sense of Definition \ref{zariskilike}), we always specify it.
Since all irreducible sets PQF-closed over $\emptyset$ are irreducible (in the topological sense), all the results proved in the more general framework of irreducible sets hold for them.
 
We recall the axioms:
\begin{enumerate}[(1)]
\item The irreducible sets PQF-closed over $\emptyset$ are Galois definable.
\item For each $n$ and each $v \in V^n$, there is some irreducible $C \subset \M^n$, PQF-closed over $\emptyset$, such that $v$ is generic in $C$, i.e. of all elements in $C$, $v$ has the maximal $U$-rank over $\emptyset$.  
\item The generic elements of an irreducible set that is PQF-closed over $\emptyset$ have the same Galois type.
\item If $C, D$ are irreducible and PQF-closed over $\emptyset$, $a \in C$ generic and $a \in D$, then $C \subseteq D.$
\item If $C_1, C_2$ are irreducible and PQF-closed over $\emptyset$, $(a,b) \in C_1$ is generic, $a$ is a generic element of $C_2$ and $(a',b') \in C_1$, then $a' \in C_2$.
\item If $C \subset D^n$ is irreducible and PQF-closed over $\emptyset$, and $f$ is a coordinate permutation on $V^n$, then $f(C)$ is irreducible and PQF-closed over $\emptyset$.
\item Let $a \to a'$ be a strongly good specialization and let $U(a)-U(a') \le 1$.
Then any specializations $ab \to a'b'$, $ac \to a'c'$ can be amalgamated: there exists $b^{*}$, independent from $c$ over $a$ such that $\textrm{t}^g(b^{*}/a)=\textrm{t}^g(b/a)$, and $ab^{*}c \to a'b'c'$.
\item Let $(a_i:i\in I)$ be independent and indiscernible over $b$.
Suppose $(a_i':i\in I)$ is indiscernible over $b'$, and $a_ib \to a_i'b'$ for each $i \in I$.
Further suppose $(b \to b')$ is a strongly good specialization and $U(b)-U(b') \le 1$.
Then, $(ba_i:i \in I) \to (b'a_i':i\in I)$.
\item Let $\kappa$ be a (possibly finite) cardinal and let  $a_i, b_i \in V$ with $i < \kappa$, such that $a_0 \neq a_1$ and $b_0=b_1$. 
Suppose $(a_i)_{i<\kappa} \to (b_i)_{i <\kappa}$ is a specialization.
Assume there is some unbounded and directed $S \subset \mathcal{P}_{<\omega}(\kappa)$ satisfying the following conditions:
\begin{enumerate}[(i)]
\item  $0,1 \in X$ for all $X \in S$;
\item For all $X,Y \in S$ such that $X \subseteq Y$,  and for all sequences $(c_i)_{i \in Y}$ from $V$, the following holds: 
If $c_0=c_1$, $ (a_i)_{i \in Y} \to (c_i)_{i \in Y} \to (b_i)_{i \in Y},$
and $\textrm{rk}((a_i)_{i \in Y} \to (c_i)_{i \in Y}) \le 1$, then $\textrm{rk}((a_i)_{i \in X} \to (c_i)_{i \in X}) \le 1$.
\end{enumerate}
 
Then, there are $(c_i)_{i<\kappa}$ such that   
$$(a_i)_{i \in \kappa} \to (c_i)_{i \in \kappa} \to (b_i)_{i \in \kappa},$$
$c_0=c_1$ and $\textrm{rk}((a_i)_{i \in X} \to (c_i)_{i \in X}) \le 1$ for all $X \in S$.
 \end{enumerate}
 
\vspace{3mm}

\begin{remark}\label{irredykskas}
Let $W$ be a Zariski closed set on the field sort $F$.
Note that if $W$ is irreducible on some algebraically closed subfield $F' \subset F$ (e.g. over $\overline{\mathbb{Q}}$), then it is irreducible also on $F$.
This follows from Model Completeness for algebraically closed fields or from the fact that if an ideal $\langle f_1, \ldots, f_r \rangle \subset F'[x_1, \ldots, x_n]$ is prime, then also the ideal $\langle f_1, \ldots, f_r \rangle \subset F[x_1, \ldots, x_n]$ is prime.
\end{remark}

\begin{remark}
Let $C$ be a PQF-closed subset of the cover sort, definable (by a positive quantifier-free formula) over $A \subset V$.
By Corollary \ref{variety}, $\textrm{exp}(C)$ is a variety.
By Corollary \ref{variety} and Model Completeness, $\textrm{exp}(C)$ is definable by polynomials with coefficients in the smallest algebraically closed field containing $A$ (see the arguments in Remark \ref{ranksama}).
Thus, in particular, if $C$ is definable over $\emptyset$, then $\textrm{exp}(C)$ is definable by polynomials with coefficients in $\overline{\mathbb{Q}}$.
\end{remark}

\begin{remark}
Note that from Remark \ref{ranksama}, it follows that any element is generic on its locus.
\end{remark}
 
\begin{lemma}\label{generic}
Let $v \in V^n$.
Then there is some irreducible set $C$, PQF-closed over $\emptyset$, such that $C$ is the locus of $v$.
\end{lemma}

\begin{proof}
Let $W$ be the smallest variety definable over $\overline{\mathbb{Q}}$ containing $\textrm{exp}(v)$.
Then, $W$ is irreducible on $\overline{\mathbb{Q}}^n$ and thus on $F^n$ by Remark \ref{irredykskas}.
There is a torus $T$ and a variety $W'$ not contained in any torus, both definable (as varieties) over $\overline{\mathbb{Q}}$ so that $W=T \cap W'$ (this is proved similarly as Lemma \ref{leikkaus}).
There is a linear $L \subset V$ definable over $\textrm{log}(\overline{\mathbb{Q}})$ (and thus over $\emptyset$)
such that $T=\textrm{exp}(L)$.
Also, the $m$:th roots of $W'$ are definable over $\overline{\mathbb{Q}}$ for all $m$.
Now the irreducible components of $\textrm{log }W$ are the sets
$$(L+k) \cap m \textrm{log }W'^{\frac{1}{m}}_{(i)}$$
for $m$ such that the $m$:th roots of $W'$ no longer branch (this can be proved completely similarly as Lemma \ref{komponentit}). 
They are clearly PQF-closed over $\emptyset$.

The element $v$ must lie on one of these components.
Denote this component by $C$.
Now $C$ is the locus of $v$.
Indeed, $v$ cannot be contained in any irreducible set $C' \subsetneq C$, PQF-closed over $\emptyset$, as then we would have $\textrm{exp}(v) \in \textrm{exp}(C') \subsetneq \textrm{exp}(C)=W$, where the equality holds by Corollary \ref{olennainen}.
This is a contradiction, as $W$ was taken to be the locus of $\textrm{exp}(v)$.
\end{proof}
 
Also, (4) is satisfied:

Suppose $C$ is irreducible and PQF closed over $\emptyset$, and $a \in C$ is generic, i.e. $U(a/\emptyset)$ is maximal.
Then, by Remark \ref{ranksama}, $U(a/\emptyset)$ is the dimension of the locus of $a$ over $\emptyset$, so $C$ must be the locus of $a$.
Thus, if $a \in D$, we must have $C \subseteq D$.

(5) is satisfied:

As $C_2$ is PQF-closed, also the set $\{(x,y) \, | \, x \in C_2\}$ is PQF-closed.
Since the intersection of two PQF-closed sets is PQF-closed, the set $D=\{(x,y) \, | \, (x,y) \in C_1 \textrm{ and } x \in C_2\}$ is PQF-closed.
Since $C_1, C_2$ are PQF-closed over $\emptyset$, also $D$ is PQF-closed over $\emptyset$.
As $(a,b)$ is a generic element of $C_1$ and $(a,b) \in D$, we have that $C_1 \subseteq D$.
Thus, $a' \in C_2$.

Also (1) and (6) clearly hold.
Now, we prove (3)

\begin{lemma}\label{samegaloistype}
The generic elements of an irreducible set PQF-closed over $\emptyset$ have the same Galois type over $\emptyset$.
\end{lemma}

\begin{proof}
Let $C \subset V^n$ be irreducible and PQF-closed over $\emptyset$, $a=(a_1, \ldots, a_n), b=(b_1, \ldots, b_n) \in C$ generic.
We may without loss assume that there is some $k \le n$ such that $a_1, \ldots, a_k$ are linearly independent and $a_{k+1}, \ldots, a_n \in \textrm{span}(a_1, \ldots, a_k)$.
Let $L$ be the linear set given by the equations expressing this.
Then, we have $a \in L$, and thus, as $a$ is generic on $C$, also $C \subseteq L$.
Hence, $b \in L$, so $b_{k+1}, \ldots, b_n \in \textrm{span}(b_1, \ldots, b_k)$.
In particular, $a_{k+1}, \ldots, a_n \in \textrm{dcl}(a_1, \ldots, a_k)$ and $b_{k+1}, \ldots, b_n \in \textrm{dcl}(b_1, \ldots, b_k)$.
This implies that 
$$U(a_1, \ldots, a_k)=U(a)=U(b)=U(b_1, \ldots, b_k).$$

Denote $a'=(a_1, \ldots, a_k)$ and $b'=(b_1, \ldots, b_k)$.
Let $D$ be the locus of $a'$.
Since $a$ is generic on $C$, we have $b' \in D$, and thus also $b'$ is generic on $D$ as it has the same $U$-rank as $a'$.
Since there are no linear dependencies between the elements $a_1, \ldots, a_k$, we have $D=m \textrm{log } W$ for some variety $W$ that does not branch.
Then, $W^m$ is the locus of $\textrm{exp}(a')$.
Let $l \in \mathbb{N}$ be arbitrary.
We have $\textrm{exp}(\frac{a'}{m}) \in W.$
Since $W$ does not branch, it has a unique $l$:th root, $W^{\frac{1}{l}}$.
Now, $\frac{a'}{ml} \in \frac{1}{l} \textrm{log } W$, and thus $\textrm{exp}(\frac{a'}{ml}) \in W^{\frac{1}{l}}$.
Then,
$$\textrm{exp}\left(\frac{a'}{l} \right) = \left( \textrm{exp} \left( \frac{a'}{lm}\right)\right)^m \in (W^{\frac{1}{l}})^m,$$
where $(W^{\frac{1}{l}})^m$ is a single $l$:th root of $W^m$.
Thus, for each $l$, we are able to determine in which $l$:th root of $W^m$ the element $\textrm{exp}(\frac{a'}{l})$ lies. 
Since $W^m$ is the locus of $\textrm{exp}(a')$, we have $\textrm{exp}(a') \notin Y$ for all $Y \subset W$ such that $\textrm{dim}(Y)<\textrm{dim}(W)$, and since $a_1, \ldots, a_k$ are linearly independent, $za' \neq 0$ for every $z \in \mathbb{Z}^n \setminus \{0\}$.
Hence, by Lemma \ref{tyypit}, the set $D$ determines the quantifier-free type of $a'$ over $\emptyset$.

The element $b'$ is also generic on $D$, so the same argument applies to it (note that $b_1, \ldots, b_k$ must be linearly independent; if they weren't, there would be some set $L'$ determined by the equations giving the dependencies, and by the genericity of $b'$ we would then have $a' \in L'$).
By \cite{monet} and \cite{kir}, quantifier free types determine Galois types. 
Hence $t^g(a'/\emptyset)=t^g(b'/\emptyset)$.
 
There is an automorphism taking $a'=(a_1, \ldots, a_k)$ to $b'=(b_1, \ldots, b_k)$.
Since the elements $b_{k+1}, \ldots, b_n$ are determined from the elements $b_1, \ldots, b_k$ by exactly the same linear equations that determine the elements $a_{k+1}, \ldots, a_n$ from $a_1, \ldots, a_k$, it must map $(a_{k+1}, \ldots, a_n)$ to $(b_{k+1}, \ldots, b_n)$, and hence $t^g(a/\emptyset)=t^g(b/\emptyset)$. 
 \end{proof}

Next, we prove (7).
This can be done essentially in the same way as for Zariski geometries (\cite{HrZi}). 

We note first that  if $u \to v$ is a strongly regular specialization on the cover sort, then $\textrm{exp}(u) \to \textrm{exp}(v)$ is a strongly regular specialization on the field sort (we consider the field as a Zariski-like structure where we take the irreducible, $\emptyset$-closed sets to be the irreducible sets of Definition \ref{zariskilike}).
It follows that $\textrm{exp}(u) \to \textrm{exp}(v)$ is a regular specialization on the field sort (as defined in \cite{HrZi}), in particular that $\textrm{exp}(v)$ is regular on the locus of $\textrm{exp}(u)$.
  
In Chapter 3, Definition \ref{goodfor}, we defined what it means for an element to be good for an irreducible set.
Although in the context of that chapter, we meant irreducible sets in the sense of the definition of Zariski-like, the concept can be defined similarly for all sets that are irreducible in the sense of the PQF-topology.

\begin{definition}
Let $C \subset V^{n+m}$ be an irreducible set.
We say an element $a \in V^n$ is \emph{good} for $C$ if there is some $b \in V^m$ so that $(a,b)$ is a generic element of $C$.
\end{definition}

We now prove a couple of lemmas needed for (7).

\begin{lemma}\label{poikkeuspisteet}
Suppose $C \subset V^{n+m}$ is irreducible, $a \in V^n$ is good for $C$, $a \to a'$ and $U(a)-U(a') \le 1$.
Then, $\textrm{dim}(C(a')) \le \textrm{dim}(C(a))$.
\end{lemma}

\begin{proof}
We present the argument in the case that $C$ is PQF-closed over $\emptyset$.
If it is not, then there are some parameters needed in defining $C$, and we have to take them into account in the calculations that will follow.
However, the calculations will remain similar to the ones we present here.  

Suppose $r_1=\textrm{dim}(C(a'))> \textrm{dim} (C(a))=r_2$.
It follows from the assumptions that $\textrm{exp}(a) \to \textrm{exp}(a')$ and $\textrm{MR}(\textrm{exp}(a))-\textrm{MR}(\textrm{exp}(a')) \le 1$.
By \cite{HrZi}, Lemma 4.12, the result holds for Zariski geometries, in particular algebraically closed fields.  
Thus, applying the result to $\textrm{exp}(C)$, $\textrm{exp}(a)$ and $\textrm{exp}(a')$, we get
\begin{eqnarray}\label{aboveinequality}
\textrm{dim}(\textrm{exp}(C)(\textrm{exp}(a'))) \le \textrm{dim}(\textrm{exp}(C)(\textrm{exp}(a))).
\end{eqnarray}

Let $b' \in C(a')$ be such that $U(b'/a')=r_1$.
As $(a',b') \in C$, we have $(\textrm{exp}(a'), \textrm{exp}(b')) \in \textrm{exp}(C)$, so $\textrm{exp}(b') \in \textrm{exp}(C)(\textrm{exp}(a'))$.
By Remark \ref{ranksama} and Lemma \ref{dimensio}, $\textrm{dim}(\textrm{exp}(b')/\textrm{exp}(a'))=r_1$, and thus $\textrm{dim}(\textrm{exp}(C)(\textrm{exp}(a'))) \ge r_1$.
By the inequality (\ref{aboveinequality}), there is some element $x \in \textrm{exp}(C)(\textrm{exp}(a))$ such that $\textrm{MR}(x/\textrm{exp}(a)) \ge r_1$.
Hence there is some $k \in K^n$ and some $b \in V^m$ such that $(a+k, b) \in C$ and $x=\textrm{exp}(b)$.
Since $U(a+k)=U(a)$, we have
$$U(b/a+k) \le \textrm{dim}(C)-U(a)=r_2,$$
a contradiction since we also have 
$$U(b/a+k)=\textrm{MR}(\textrm{exp}(b)/\textrm{exp}(a))\ge r_1>r_2.$$
\end{proof}
 
\begin{lemma}\label{henkitore}
Suppose $C$ is irreducible and $a$ is good for $C$.
Suppose $(C_i)_{i<\omega}$ is a collection of irreducible sets such that it is permuted by all automorphisms.
Let $b$ be an element such that $ab$ is good for each $C_i$ and that for any $c \in C(a)$ generic over $b$
it holds that $c \in C_i(ab)$ for some $i$.
Assume $(a, b) \to (a', b')$, $a \to a'$ is a strongly regular specialization, $U(a)-U(a') \le 1$, and $c'$ is such that $(a',c') \in C$.

Then, there is some $i <\omega$ so that $(a',b',c') \in C_i$.
\end{lemma}

\begin{proof} 
Let $D$ be the locus of $(a, b)$, and let $m=\textrm{lg}(c)$.
Denote 
$$C^{*}=\{(x,y,z) \in D \times V^m \, | \, (x,z) \in C\}.$$
Let $E$ be the locus of $a$.%

We first show that the irreducible components of $C^*$ not contained in
$$X=\{(x,y,z) \in D \times V^m \, | \, U(a)-U(x)>1 \textrm{ or $\textrm{exp}(x)$ is not regular on $\textrm{exp}(E)$}\}$$
all have same dimension $\textrm{dim}(D)+U(c/a)$.
(Note that the set $X$ is not definable.)
Suppose $(d,e,f)$ is a generic element in some such irreducible component.
Then, we have $U(d) \ge U(a)-1$, and thus by Lemma \ref{poikkeuspisteet}, $U(f/d) \le U(c/a)$.
 Hence, 
 $$U(d,e,f) \le \textrm{dim}(D)+U(c/a).$$
 
On the other hand, let $Y$ be an irreducible component of $C^*$ not contained in $X$.
Denote $\Delta_E=\{(x,y) \in E \times E \, | \, x=y \}$.
We note first that in the PQF-topology, $C^*$ is isomorphic to $(D \times C) \cap (\Delta_E \times V^{n+m})$, where $n=\textrm{lg}(b)$.
Both $D \times C$ and $\Delta_E \times V^{n+m}$ are Cartesian products of irreducible sets and thus irreducible.
As in the proof of Theorem \ref{dimthm}, one sees that if $Y'$ is an irreducible component of $(D \times C) \cap (\Delta_E \times V^{n+m})$, then $\textrm{exp}(Y')$ is an irreducible component of $\textrm{exp}((D \times C) \cap (\Delta_E \times V^{n+m}))$.
Hence also $\textrm{exp}(Y)$ is an irreducible component of $\textrm{exp}(C^*)$.

On the field sort $F$, we may obtain a copy of $\textrm{exp}(C^*)$ by intersecting $\textrm{exp}(D) \times \textrm{exp}(C)$ with a suitable $\textrm{exp}(E)$-diagonal.
Indeed, if we have $(x, y, c_1, c_2) \in \textrm{exp}(D) \times \textrm{exp}(C)$, where $(x, y) \in \textrm{exp}(D)$ and $(c_1, c_2) \in \textrm{exp}(C)$, then the diagonal $\Delta$ expressing ``$x=c_1$" is as wanted.  
As $Y$ is not contained in $X$, the isomorphic copy of $\textrm{exp}(Y)$ obtained in the above procedure contains some point $(c_1, y, c_1, c_2)$ where $c_1$ is a regular point of $\textrm{exp}(E)$.
Thus, using Lemma 5.4 in \cite{HrZi} to calculate dimensions in $W=\textrm{exp}(E) \times F^n \times \textrm{exp}(E) \times F^m$ and keeping in mind that for any PQF-closed set $C$, $\textrm{dim}(C)=\textrm{dim}(\textrm{exp}(C))$ (Lemma \ref{dimensio}), we obtain
\begin{eqnarray*}
\textrm{dim}(Y) &\ge& \textrm{dim}(D)+\textrm{dim}(C)+\textrm{dim}(\Delta)-\textrm{dim}(W) \\
&=& \textrm{dim}(D)+\textrm{dim}(C)-\textrm{dim}(E)\\
&=& \textrm{dim}(D)+U(c/a),
\end{eqnarray*}
since $\Delta$ is of codimension $\textrm{dim}(E)$ in $W$ and $\textrm{dim}(C)-\textrm{dim}(E)=U(c/a)$.
 
Let $C'$ be the irreducible component of $C^*$ containing $(a', b',c')$.
Then, $C'$ is not included in $X$.  
Indeed, $U(a') \ge U(a)-1$ and the fact that $a \to a'$ is strongly regular implies that $\textrm{exp}(a) \to \textrm{exp}(a')$ is regular, i.e. $\textrm{exp}(a')$ is regular on the locus of $\textrm{exp}(a)$ which is
$\textrm{exp}(E)$. 
Let $(d,e,f)$ be a generic point of $C'$.
Then, $(d,f) \in C$ and $U(d) \ge U(a)-1$. 
By Lemma \ref{poikkeuspisteet}, $U(f/d) \le U(c/a)$, and hence $U(f/d,e) \le U(c/a)$.
As $U(d,e,f)=\textrm{dim}(D)+U(c/a)$, this implies that $(d,e)$ is a generic point of $D$.
It also follows that $U(c/a)=U(f/d,e)=U(f/d)$, so in particular $f \da_d e$.
There is some automorphism taking $(d,e) \mapsto (a,b)$.
This automorphism then takes $f$ to some element $c \in C(a)$ such that $c \da_a b$.
By our assumptions, $(a,b,c) \in C_i$ for some $i<\omega$.
Since automorphisms permute the collection of the sets $C_i$, we have that $(d,e,f) \in C_j$ for some $j< \omega$.
Hence, $C' \subset C_j$, so in particular $(a',b',c') \in C_j$.
\end{proof}

\begin{lemma}\label{prop4.7}
Let $C$, $D$ be irreducible, $a \in D$ generic, and suppose $a$ is good for $C$.
 Let $r=\textrm{dim}(C(a))\ge 0$.
Let $c$ and $c'$ be such that $(a,c) \to (a',c')$ and $a \to a'$ is strongly regular and $U(a)-U(a') \le 1$.
Let $b' \in C(a')$.

Then, there exists $b \in C(a)$ such that $U(b/ac)=r$ and $(a,b,c) \to (a',b',c')$.
 \end{lemma}

\begin{proof}
Let $E$ be the locus of $(a,c)$, $C^{*}=\{(x,z,y) \, | \, (x,z) \in E, (x,y) \in C\}$.
Let $C_i$, $i<\omega$ be the irreducible components of $C^*$ satisfying $\textrm{dim}(C_i(a,c))=r$.
We claim that
$$C(a) \subseteq \bigcup_{i<\omega} C_i(a,c).$$

We note first that every irreducible component of $C(a)$ is of dimension $r$. 
Let $m=\textrm{lg}(b')$.
Then, $C \subset D \times V^m$ and $C(a)=C \cap (\{a\} \times V^m)$
As $a$ is generic in $D$, $\textrm{exp}(a)$ is regular on $\textrm{exp}(D)$.
Hence, applying the dimension theorem on $\textrm{exp}(D \times V^m)$
similarly as in the proof of Lemma \ref{henkitore},
we get that every nonempty irreducible component of $C(a)$ has dimension at least
$\textrm{dim}(C)+m-\textrm{dim}(D)-m=r.$

So in every such component, there is an element $x$ such that $U(x/ca)=r$.
Then, $(a,c,x) \in C'$ for some irreducible component $C'$ of $C^*$, and thus $x \in \bigcup_{i<\omega} C_i(a,c)$.
Hence, the irreducible component of $C(a)$ containing $x$ is included in $\bigcup_{i<\omega} C_i(a,c)$.
We conclude that $C(a) \subseteq \bigcup_{i<\omega} C_i(a,c)$.
 
Now, by Lemma \ref{henkitore}, $b' \in C_i(a',c')$ for some $i$.
Let $b$ be a generic point of $C_i(a,c)$.
Since $C_i$ is irreducible and $(a,b,c)$ is a generic point of $C_i$, we have $(a,b,c) \to (a',b',c')$.
\end{proof} 

Now we are ready to prove Axiom (7).
 
\begin{theorem}\label{amalg} 
Let $a \to a'$ be a strongly good specialization and let $U(a)-U(a') \le 1$.
Then any specializations $ab \to a'b'$, $ac \to a'c'$ can be amalgamated: there exists $b^{*}$, independent from $c$ over $a$ such that $\textrm{t}^g(b^{*}/a)=\textrm{t}^g(b/a)$, and $ab^{*}c \to a'b'c'$.
\end{theorem}
 
\begin{proof}
We prove the Lemma by induction, using Definition \ref{good}.
If $a \to a'$ is strongly regular, this follows from Lemma \ref{prop4.7}:
Let $D_1$ be the locus of $a$, $m=\textrm{lg}(b)$ and $C$ the locus of $(a,b)$ in $D_1 \times V^m$.
By Lemma \ref{prop4.7}, there is some $b^* \in C(a)$ such that $U(b^*/ac)=\textrm{dim}(C(a))$ with $ab^*c \to a'c'b'$.
In particular, $U(b^*/a)=U(b/a)$, and $b^*$ is independent from $c$ over $a$.
Since $C$ is the locus of $(a,b)$, both $(a,b)$ and $(a,b^*)$ are generic on $C$.
By Lemma \ref{samegaloistype}, $t^g(ab/\emptyset)=t^g(ab^*/\emptyset)$, so in particular $\textrm{t}^g(b^*/a)=\textrm{t}^g(b/a)$.
 
Suppose now $a=(a_1, a_2,a_3)$, $a'=(a_1',a_2',a_3')$, $a \to a'$ are as in Definition \ref{good}, and the lemma holds for the specialization $(a_1, a_2) \to (a_1', a_2')$.
Amalgamating over $(a_1, a_2) \to (a_1', a_2')$, we see that there exist $b^*$, $a_3^*$ such that $\textrm{t}^g(b^* a_3^*/a_1 a_2)=\textrm{t}^g(ba_3/a_1 a_2)$, $b^*a_3^*$ is independent from $c$ over $a_1 a_2$, and $a_1 a_2 a_3^*b^*a_3 c \to a_1' a_2' a_3' b' a_3' c'$.
Now, by Definition \ref{good}, $a_3 \in \textrm{bcl}(a_1)$, and as $\textrm{t}^g(b^* a_3^*/a_1 a_2)=\textrm{t}^g(ba_3/a_1 a_2)$, also $a_3^* \in \textrm{bcl}(a_1)$.
By Definition \ref{good}, $a_1 \to a_1'$ is an isomorphism.
From these facts together it follows that $a_1a_3a_3^* \to a_1' a_3' a_3'$ is of rank $0$ and thus an isomorphism.
Hence, $a_3=a_3^*$, so we get $ab^*c \to a'b'c'$ as wanted.
\end{proof}
 
Axiom (8) is proved as follows:
 
\begin{theorem}\label{infamalg}
Let $(a_i:i\in I)$ be independent and strongly indiscernible over $b$, with $I$ infinite. 
Suppose $(a_i':i\in I)$ is strongly indiscernible over $b'$, and $a_ib \to a_i'b'$ for each $i \in I$.
Further suppose $\textrm{rk}(b \to b') \le 1$.
Then, $(ba_i:i \in I) \to (b'a_i':i\in I)$.
\end{theorem}

\begin{proof} 
Since the sequences are strongly indiscernible, we may without loss assume that $I=\omega_1$.
By Theorem \ref{amalg} and induction, there exist elements $a_i^*$, $i<\omega$, independent over $b$ such that $\textrm{t}^g(a_i^*/b)=\textrm{t}^g(a_i/b)$ and $b(a_i^*)_{i<\omega} \to b'(a_i')_{i<\omega}$.
Then, also $(\textrm{exp}(a_i^{*}))_{i <\omega}$ are independent over $\textrm{exp}(b)$.
Moreover, we have for all $i,j \in \omega$,  that $\textrm{tp}(\textrm{exp}(a_i^*/\textrm{exp}(b))=\textrm{tp}(\textrm{exp}(a_j^*/\textrm{exp}(b))$ (for complete first-order types).
Since algebraically closed fields are $\omega$-stable, there are only finitely many free extensions for each complete type.
Thus, there is an infinite $I_0  \subset \omega$ so that $(\textrm{exp}(a_i^{*}))_{i \in I_0}$ are indiscernible (in the field language) over $\textrm{exp}(b)$.
 
Let $\theta$ be a theory consisting of the first-order formulae that express the following (note that we have added the elements of $\textrm{log }(\overline{\mathbb{Q}})$ to our language):
\begin{itemize}
\item The sequence $(\textrm{exp}(x_i))_{i<\omega_1}$ satisfies the same first-order formulae over $\overline{\mathbb{Q}} \cup \textrm{exp}(b)$ as the sequence $(\textrm{exp}(a_i^{*}))_{i \in I_0}$ (note that this is possible as the latter sequence is indiscernible over $ \textrm{exp}(b)$);
\item For each $i<\omega_1$, $x_i$ has the same complete first-order type over $\textrm{bcl}(\emptyset) \cup \{b\}$ as $a_j^*$ for some (and thus every) $j \in I_0$, i.e. $\textrm{tp}(x_i/\textrm{bcl}(\emptyset) b)=\textrm{tp}(a_j^*/\textrm{bcl}(\emptyset) b)$;
\item For each $n$ and each positive, quantifier free first-order formula $\phi$ such that $\neg \phi(b', a_1', \ldots, a_n')$ holds, the theory $\theta$ contains the formula $\neg \phi(b, x_{j_1}, \ldots, x_{j_n})$ for all $j_1< \ldots <j_n<\omega_1$;
\item For any $n$ and any $i_1<\ldots<i_n<\omega_1$, it holds that if $\textrm{exp}(q_1 x_{i_1}+ \ldots + q_n x_{i_n}+c)=1$, where $q_1, \ldots, q_n \in \mathbb{Q}$ and $c$ is some linear combination of elements in $\textrm{log}(\textrm{exp}(\overline{\mathbb{Q}}))\cup b$, then 
$$q_1 x_{i_1}+ \ldots + q_n x_{i_n}+c=q_1 a_1^*+ \ldots + q_n a_n^*+c.$$ (Note that this can be expressed since $q_1 a_1^*+ \ldots + q_n a_n^*+c \in K$ and $K \subset \textrm{log}(\textrm{exp}(\overline{\mathbb{Q}}))$.)
\end{itemize}
This theory is consistent as every finite fragment is realized by the sequence $(a_i^{*})_{i \in I_0}$.
Thus, in a saturated elementary extension $(\mathfrak{V}, \mathfrak{F})$ of the monster model $(V,F)$, we find a sequence $(c_i)_{i<\omega_1}$ realizing $\theta$
(note that since $(\mathfrak{V}, \mathfrak{F})$ is saturated, it is not a model of  $T+(K=\mathbb{Z})$).
From now on, we denote by $K$ the kernel of exp in $V$ (then, $K \cong \mathbb{Z}$) and by $K^{*}$ the kernel of exp in  $\mathfrak{V}$.

Since the $c_i$ satisfy the theory $\theta$, we have $\textrm{span}(c_i)_{i<\omega_1} \cap K^* \subseteq K$. 
Denote now $X=\textrm{span}((c_i)_{i<\omega} \cup \textrm{log}(\mathbb{Q}))$.
We claim that $\textrm{exp}(X)$ is closed under multiplication and inversion and that if for some $a \in \mathfrak{F}$, $a^{\frac{m}{n}} \in \textrm{exp}(X)$ for some choice of the $n$:th root $a^\frac{1}{n}$, where $\frac{m}{n} \in \mathbb{Q}$, then $a \in \textrm{exp}(X)$.
Let $x,y \in \textrm{exp}(X)$.
Then, there are $u,v \in X$ such that $\textrm{exp}(u)=x$, $\textrm{exp}(v)=y$.
Now, $u+v \in X$ and $\textrm{exp}(u+v)=xy$, so $xy \in \textrm{exp}(X)$.
Also, $x^{-1}=\textrm{exp}(-u) \in \textrm{exp}(X)$.
Let $c \in \textrm{exp}(X)$ be such that $a=c^m$ (i.e. $c$ is a choice of a $m$:th root for $a$).
Let  $u \in X$ be such that $\textrm{exp}(u)=c$.
Then, $\textrm{exp}(mu)=a$.
Suppose now $a \in \textrm{exp}(X)$ and let $c$ be a choice for the $m$:th root of $a$.
Now, let $u \in X$ be such that $\textrm{exp}(u)=a$. 
Then, $\textrm{exp}(\frac{u}{m}+\frac{k}{m})=c$ for some $k \in K$.
As $K \subset X$, we have $c \in \textrm{exp}(X)$.  
 
Let now $A=\textrm{acl}(\textrm{exp}(X))$.
Choose $d_0 \in A \setminus \textrm{exp}(X)$ and $x_0 \in \mathfrak{V}$ so that $\textrm{exp}(x_0)=d_0$.
Denote $X_1=\textrm{span}(X \cup \{x_0\})$.
We claim that $X_1 \cap K^{*} \subset K$.
Suppose not.
Let $X_0 \subset X$ be finite and
$$\textrm{exp}(\sum_{v \in X_0} q_v v+q x_0)=1$$
 for $q_v, q \in \mathbb{Q}$, and suppose $\sum_{v \in X_0} q_v v+q x_0 \notin K$.
Then, we must have $q \neq 0$, as otherwise we would have $\sum_{v \in X_0} q_v v+q x_0=\sum_{v \in X_0} q_v v \in K$ (as $\sum_{v \in X_0} q_v v \in X$ and $X \cap Z^{*}=K$).  
This gives us
$$d_0^q \prod_{v \in X_0} (\textrm{exp}(v))^{q_v} =1,$$ and hence $$d_0^q=\prod_{v \in X_0} (\textrm{exp}(v))^{-q_v}$$ (for some suitable choices of the roots in question).
But now $d_0 \in \textrm{exp}(X)$ which is against our assumptions.

We may now repeat the argument, and eventually we will get a set $X' \subset \mathfrak{V}$ such that $\textrm{exp}(X')=A$ and $X' \cap K^{*}=K$.
We have thus constructed a model $(X', A)$ for the theory $T+(K\cong\mathbb{Z})$.
As this theory is categorical, $(X', A)$ is isomorphic to some elementary submodel of our monster model $(V,F)$, and thus we can find the sequence $(c_i)_{i<\omega_1}$ already in $V$.

The sequence $(\textrm{exp}(a_i^{*}))_{i \in I_0}$ is independent over $\textrm{exp}(b)$.
The sequence $(\textrm{exp}(c_i))_{i<\omega_1}$ satisfies the same first-order formulae over $\textrm{exp}(b)$ as the sequence $(\textrm{exp}(a_i^{*}))_{i \in I_0}$, and independence is a local property.
Thus, the sequence $(\textrm{exp}(c_i))_{i<\omega_1}$ is independent over $\textrm{exp}(b)$.
As all ranks are calculated on the field sort, the sequence $(c_i)_{i<\omega_1}$ is independent over $b$.
 
Since the sequence $(c_i)_{i<\omega_1}$ is uncountable, there is some uncountable $J \subset \omega_1$ such that the sequence $(c_i)_{i \in J}$ is Morley over $b$ and thus strongly indiscernible over $b$.
Since the $(a_i')_{i<\omega_1}$ are strongly indiscernible over $b'$ (note that by strong indiscernibility, we may extend the sequence to be arbitrarily long), we have $b(c_i)_{i\in J} \to b' (a_i')_{i<\omega_1}$.

Relabel the indices so that from now  on $(c_i)_{i<\omega_1}$ stands for $(c_i)_{i \in J}$.
The sequence $(c_i)_{i<\omega_1}$ is independent and strongly indiscernible over $b$, and so is $(b, a_i)_{i<\omega_1}$.
For $(b, a_i)_{i<\omega_1} \to (b', a_i')_{i< \omega_1}$, it suffices to show that $t(a_{i_1}, \ldots, a_{i_n}/b)=t(c_{i_1}, \ldots, c_{i_n}/b)$ for all $i_1, \ldots i_n \in \omega_1$.
As the sequences are strongly indiscernible, we may assume $i_1=1, \ldots, i_n=n$.
We have $t(c_1/b)=t(a_1/b)$.
Let $f \in \textrm{Aut}(\M/b)$ be such that $f(a_1)=c_1$, and let $a_2'=f(a_2)$.
Then, $c_1 \da_b a_2'$, and thus $a_2' \da_b c_1$.
Since also $c_2 \da_b c_1$, we have $t(c_2/b c_1)=t(a_2'/b c_1)$, so
$$t(c_1 c_2 /b)=t(c_1 a_2'/b)=t(a_1 a_2/b).$$
Inductively, one shows that $t(a_1, \ldots, a_n/b)=t(c_1, \ldots, c_n/b)$.
\end{proof}

For (9), we still need the following lemma.
The Zariski geometry version was presented originally in \cite{HrZi}.  
Although the basic idea of the proof is the same, we have to do a bit more work.

\begin{lemma}\label{dimthmspec}
Let $a=(a_1, \ldots, a_n), a''=(a_1'', \ldots, a_n'') \in V^n$, $a \to a''$, and suppose $a_1 \neq a_2$, $a_1''=a_2''$.
Then there exists $a'=(a_1', \ldots, a_n') \in V^n$ such that $a_1'=a_2'$, $a \to a' \to a''$, and $U(a)-U(a')=1$.
\end{lemma}

\begin{proof}
Denote $\Delta_{12}^n=\{(v_1, \ldots, v_n) \in V^n \, | \, v_1=v_2\}$.
Let $C$ be the locus of $a$.
Then, $a'' \in C \cap \Delta_{12}^n$.
Hence, $a''$ must lie on some irreducible component $D$ of $C \cap \Delta^n_{12}$.
 By Lemma \ref{dimensio},
$\textrm{dim}( \Delta^n_{12})=n-1$.
Thus, Theorem \ref{dimthm} yields 
$$\textrm{dim}(D) \ge \textrm{dim}(C)+\textrm{dim}( \Delta^n_{12})-n=\textrm{dim}(C)-1.$$ 
As $a_1 \neq a_2$, we have $C \cap \Delta_{12}^n \subsetneq C$, and thus $\textrm{dim}(D) < \textrm{dim}(C)$.
Hence,  $\textrm{dim}(D)=\textrm{dim}(C)-1$.

Next, we will show that $D$ can be defined by a positive, quantifier-free formula over the empty set, in other words that it is in our collection of irreducible sets.
We have $\textrm{exp}(a'') \in \textrm{exp}(C) \cap \Delta_{12}^n$, where by $\Delta_{12}^n$ we denote, abusing the notation, the corresponding diagonal on the field sort.
Then, $\textrm{exp}(a'')$ lies on some irreducible component $X$ of $\textrm{exp}(C) \cap \Delta_{12}^n$.
By the Dimension Theorem for Zariski geometries,
$$\textrm{dim}(X) \le \textrm{dim}(\textrm{exp}(C))-1=\textrm{dim}(C)-1,$$
where the last equality holds by Lemma \ref{dimensio}.
 
Denote $b=\textrm{exp}(a)$.
We claim that $b_1 \neq b_2$.
If not, then we would have $a_1=a_2+k$ for some $0 \neq k \in K$.
Then, there is an element $0 \neq h=(h_1, \ldots, h_n) \in K^n$ such that  $h_1=k$ and $a \in \Delta_{12}^n+h$.
Since we have added the elements of $K$ into the language, $\Delta_{12}^n+h$ is PQF-closed over the empty set.
As $a$ is generic on $C$ (over the empty set), this means that $C \subset \Delta_{12}^n+h$.
But this is impossible, as $a'' \in C \setminus (\Delta_{12}^n+h).$
Thus, $b_1 \neq b_2$.

Hence, we have $\textrm{exp}(a) \notin \textrm{exp}(C) \cap \Delta_{12}^n$, so $X \subsetneq \textrm{exp}(C)$, and thus $\textrm{dim}(X)=\textrm{dim}(C)-1$.
We have $\textrm{exp}(D) \subseteq X$, and thus $\textrm{exp}(D)=X$. 
By Corollary \ref{olennainen}, $D$ is an irreducible component of $\textrm{log}(\textrm{exp}(C) \cap \Delta_{12}^n)$.
The variety $\textrm{exp}(C) \cap \Delta_{12}^n$ is definable over the empty set, and as in the proof of Lemma \ref{generic}, we see that also $D$ is definable over the empty set.
 
Choose $a'$ to be a generic point of $D$.
Then $a'$ is as wanted. 
\end{proof} 

And finally, we prove (9):
 
\begin{theorem}\label{dimthmspec2}
 Let  $a_i, b_i \in V$ with $i < \kappa$, such that $a_0 \neq a_1$ and $b_0=b_1$.
Denote by $K$ the kernel of exp in $V$.
Suppose $(a_i)_{i<\kappa} \to (b_i)_{i <\kappa}$ is a specialization.
Assume there is some unbounded and directed $S \subset \mathcal{P}_{<\omega}(\kappa)$ satisfying the following conditions:
\begin{enumerate}[(i)]
\item  $0,1 \in X$ for all $X \in S$;
\item For all $X,Y \in S$ such that $X \subseteq Y$,  and for all sequences $(c_i)_{i \in Y}$ from $V$, the following holds: 
If $c_0=c_1$, $ (a_i)_{i \in Y} \to (c_i)_{i \in Y} \to (b_i)_{i \in Y},$
and $\textrm{rk}((a_i)_{i \in Y} \to (c_i)_{i \in Y}) \le 1$, then $\textrm{rk}((a_i)_{i \in X} \to (c_i)_{i \in X}) \le 1$.
\end{enumerate}

\vspace{3mm}  
Then, there are $(c_i)_{i<\kappa}$ such that   
$$(a_i)_{i \in \kappa} \to (c_i)_{i \in \kappa} \to (b_i)_{i \in \kappa},$$
$c_0=c_1$ and $\textrm{rk}((a_i)_{i \in X} \to (c_i)_{i \in X}) \le 1$ for all $X \in S$.
\end{theorem}
 
\begin{proof}
Let  $a_i, b_i \in V$ with $i < \kappa$, such that $a_0 \neq a_1$ and $b_0=b_1$.
Suppose $(a_i)_{i<\kappa} \to (b_i)_{i <\kappa}$ is a specialization.
If $S \subset \mathcal{P}_{<\omega}(\kappa)$ is such that $0,1 \in X$ for all $X \in S$, then, by Lemma \ref{dimthmspec} and Remark \ref{ranksama}, for each $X \in S$, there is some sequence $(c_i)_{i \in X} \in V$ so that 
$$(a_i)_{i \in X} \to (c_i)_{i \in X} \to (b_i)_{i \in X},$$
$c_0=c_1$ and $\textrm{rk}((a_i)_{i \in X} \to (c_i)_{i \in X}) \le 1$.

Suppose now that $S$ is unbounded and directed and satisfies condition (ii) from the theorem.
By Compactness, there is a saturated elementary extension $(\mathfrak{V}, \mathfrak{F})$ of $(V,F)$ and elements $c_i \in \mathfrak{V}$ for $i < \kappa$ such that $(a_i)_{i<\kappa} \to (c_i)_{i<\kappa} \to (b_i)_{i <\kappa}$, $c_0=c_1$ and $\textrm{rk}((a_i)_{i \in X} \to (c_i)_{i \in X}) \le 1$ for all $X \in S$.
From now on, we denote by $K$ the kernel of exp in $V$ (then, $K \cong \mathbb{Z}$) and by $K^{*}$ the kernel of exp in  $\mathfrak{V}$.

We now show, using Compactness, that we may choose the $c_i$ so that $\textrm{span}(c_i)_{i < \kappa} \cap K^{*} \subset K$.  
Let $J \subset \kappa$ be finite.  
Then, by Lemma \ref{dimthmspec}, there are $c_i' \in V$ such that $(a_i)_{i \in J} \to (c_i')_{i \in J} \to (b_i)_{i \in J}$
and $\textrm{rk}((a_i)_{i \in X} \to \textrm{rk}(c_i')_{i \in X}) \le 1$ for all $X \in S \cap \mathcal{P}(J)$.
If for some $I_0 \subseteq J$, we have $\textrm{exp}(\sum_{i \in I_0} q_i c_i')=1$, where $q_i \in \mathbb{Q}$, 
then, since the $c_i'$ are in $V$, we have $\sum_{i \in I_0} q_i c_i'=k$ for some $k \in K$.
As there is a specialization from the $c_i'$ to the $b_i$, we must have $\sum_{i \in I_0} q_i b_i=k$.
Thus, for any finite $J \subset \kappa$, we may choose the sequence $(c_i)_{i \in J}$ so that whenever $I_0 \subseteq J$ and  $\textrm{exp}(\sum_{i \in I_0} q_i c_i)=1$, then $\sum_{i \in I_0} q_i c_i=\sum_{i \in I_0} q_i b_i$.
Hence, by Compactness, we may choose the $c_i$ for $i<\kappa$ so that whenever $I_0 \subset \kappa$ is finite and $\textrm{exp}(\sum_{i \in I_0} q_i c_i)=1$, then $\sum_{i \in I_0} q_i c_i=\sum_{i \in J} q_i b_i$.
In particular, $\textrm{span}(c_i)_{i \in I} \cap K^{*} \subset K$.

Now we can show that the sequence $(c_i)_{i<\kappa}$ is in $V$ using the same argument as in the proof of Theorem \ref{infamalg}.
\end{proof}
 
\subsection{Curves on the cover}

\begin{definition}
We say that an irreducible, one-dimensional PQF-closed set  $D$ on $V^n$ is a \emph{curve} on $V^n$.
Note that if $D$ is a curve on $V^n$, then $\textrm{exp}(D)$ is an algebraic curve on $F^n$.

For each $m$, define the \emph{closed sets} on  $D^m$  to be the restrictions of the PQF-closed sets on $V^{mn}$.
Again, we say that a closed set is \emph{irreducible} if it cannot be written as a union of two proper closed subsets.
\end{definition}

We first note that if $D \subset V^n$ is a curve, then each point $x \in D^m$ is also a point of $V^{nm}$ and the locus of $x$ on $D^m$ coincides with the locus of $x$ on $V^{nm}$.
It follows that also ranks coincide and that a map on $D$ is a specialization if and only if the corresponding map on $V$ is one.
From these observations it follow that the axioms (1)-(8) hold for the irreducible closed sets on each $D$.
 
Axiom (9) is more complicated.
Although it holds on the cover, it does not necessary imply that this holds on an arbitrary curve.
If we have a specialization $(a_i)_{i<\kappa} \to (b_i)_{i <\kappa}$ where for each $i$, $a_i, b_i \in D \subset V^n$ and $b_0=b_1$, then we can convert this into a specialization on $V$ by setting $a_i=(a_{i1}, \ldots, a_{in})$ and $b_i=(b_{i1}, \ldots, b_{in})$ and viewing the sequences of the tuples $a_i$ and $b_i$ as sequences of their elements.
Suppose, moreover, that the conditions of the axiom hold.
After reorganizing the indices so that $b_{00}$ is labeled the $0$:th element and $b_{10}$ the $1$:th, we may also think of the set $S$ as a subset of the new index set satisfying the conditions in the axiom in the context of the cover.
Then, Theorem \ref{dimthmspec2} (that is Axiom (9) on the cover) gives us a sequence $(c_{i1}, \ldots, c_{in})_{i<\kappa}$ with $c_{01}=c_{11}$ where
$$(a_{i1}, \ldots, a_{in})_{i<\kappa} \to (c_{i1}, \ldots, c_{in})_{i<\kappa} \to (b_{i1}, \ldots, b_{in})_{i<\kappa}$$
and $\textrm{rk}((a_{i1}, \ldots, a_{in})_{i \in X} \to (c_{i1}, \ldots, c_{in})_{i \in X}) \le 1$ for all $X \in S$.
However, for the statement to hold on $D$, we would need $c_{0j}=c_{1j}$ for each $1 \le j \le n$.
We can of course obtain such a sequence of suitable elements $c_i$ but then the rank  might drop too much for some $X \in S$.
The problem is that the Dimension Theorem does not necessarily hold on arbitrarily curves.

We may, however, remedy the situation by making the extra assumption of $\textrm{exp}(D)$ being a smooth algebraic curve.
Indeed, if  $\textrm{exp}(D)$ is a smooth curve, then the Dimension Theorem  holds on $\textrm{exp}(D)$  and  $\textrm{exp}(D)$ is a Zariski geometry (see \cite{Marker}).
As any closed, irreducible sets on $D^m$ are closed, irreducible on $V^{mn}$, we may prove Theorem \ref{dimthm}, Lemma \ref{dimthmspec} and Theorem \ref{dimthmspec2} using the same arguments as in the case of the cover.


\begin{thebibliography}{9}


\bibitem{monet} Bays, M., Hart, B., Hyttinen, T., Kes\"al\"a, M. and Kirby, J.: \emph{Quasiminimal structures and excellence}. Bull. London Math. Soc. (2014) 46 (1): 155-163. doi: 10.1112/blms/bdt076 First published online: October 24, 2013

\bibitem{bays} Bays, M. and Zilber, B.: \emph{Covers of multiplicative groups of algebraically closed fields of arbitrary characteristic}. Bull. London Math. Soc. 43 (2011): 689-702.

\bibitem{bous1} Bouscaren, E. (ed.): \emph{Model Theory and Algebraic Geometry. An introduction to E. Hrushovski's proof of the geometric Mordell-Lang conjecture}, Springer-Verlag 1999.

\bibitem{Lucy} Burton, L.: \emph{Toric Varieties as Analytic Zariski Structures}, Ph.D.Thesis, University of Oxford, 2007.

\bibitem{Cox} Cox, D., Little, J. and O'Shea, D.: \emph{Ideals, Varieties and Algorithms. An Introduction to Computational Algebraic Geometry and Commutative Algebra}, Third Edition, Springer 2007.

\bibitem{Hr} Hrushovski, E.: \emph{The Mordell-Lang conjecture for function fields}, Journal AMS 9 (1996), 667-690.

\bibitem{HrZi93} Hrushovski, E., Zilber, B.: \emph{Zariski Geometries}, Bulletin AMS 28 (1993), 315-323.  

\bibitem{HrZi} Hrushovski, E., Zilber, B.: \emph{Zariski Geometries}, Journal AMS 9 (Jan 1996),1-56.

\bibitem{Hu} Hulek, K.: \emph{Elementary Algebraic Geometry}, American Mathematical Society 2003.

\bibitem{CatTran} Hyttinen, T. and Kes\"al\"a M.: \emph{Categoricity transfer in simple finitary AECs}. J. Symbolic Logic (2011), 76(3): 759-806.

\bibitem{Hy} Hyttinen, T.: \emph{Finitely generated substructures of a homogeneous structure}. Math. Logic Quarterly, 50 (2004), 77-98. 

\bibitem{Hy2} Hyttinen, T.: \emph{Finiteness of U-rank implies simplicity in homogeneous structures}. Math. Logic Quarterly, 49 (2003), 576-578.

\bibitem{HLS} Hyttinen, T., Lessman, O. and Shelah, S.: \emph{Interpreting groups and fields in some nonelementary classes}. Journal of Mathematical Logic, 5 (2005), 2005, 1-47.

\bibitem{Hy3} Hyttinen, T.: \emph{Locally modular geometries in homogeneous structures}. Math. Logic Quarterly, 51 (2005), 291-298.

\bibitem{HyLess} Hyttinen, T. and Lessman, O.: \emph{Simplicity and uncountable categoricity in excellent classes}. Annals of Pure and Applied Logic, 139 (2006), 110-137.

\bibitem{kir} Kirby, J.: \emph{On quasiminimal excellent classes.} J. Symbolic Logic, 75(2): 551-564, 2010.

\bibitem{MarIntro} Marker, D.: \emph{Model Theory: An Introduction}, Springer-Verlag 2000.

\bibitem{Marker} Marker, D.: \emph{Zariski Geometries}, in Bouscaren, E. (ed.): \emph{Model Theory and Algebraic Geometry. An introduction to E. Hrushovski's proof of the geometric Mordell-Lang conjecture}, Springer-Verlag 1999.

\bibitem{pi} Pillay, A.: \emph{Geometric Stability Theory}. Clarendon Press, Oxford, 1996. 

\bibitem{shafa} Shafarevich, I.R.: \emph{Basic Algebraic Geometry}, Springer-Verlag 1974.

\bibitem{Zilber} Zilber, B.: \emph{Covers of the multiplicative group of an algebraically closed field of characteristic zero}. J. London Math. Soc. 74 (2006), 41-58.

\bibitem{ZilUusi} Zilber, B.: \emph{Model theory, geometry and arithmetic of the cover of a semi-abelian variety}, in \emph{Model Theory and Applications, Proceedings of the Conference in Ravello}, 2000.

\bibitem{Zi} Zilber, B.: \emph{Zariski Geometries. Geometry from the Logician's Point of View}, Cambridge University Press, 2010.

 
 
\end{thebibliography}
\end{document}